\documentclass[11pt,a4paper]{article}
\usepackage[utf8]{inputenc}

 \usepackage[hyphens]{xurl}
 \usepackage{breakurl}

\usepackage{mathrsfs,amsthm,xcolor,verbatim,bbm,amsmath,amsfonts,
amssymb,nicefrac,enumitem,hyperref,bm,mathtools,xparse,etoolbox}
\usepackage[margin=1in]{geometry}
\usepackage[capitalise,sort]{cleveref} % clever reference

\usepackage{textcomp}
\usepackage{sidecap}

\DeclareMathOperator*{\argmin}{arg\,min}

\crefname{enumi}{item}{items}

\crefname{figure}{Figure}{Figures}
\crefname{equation}{}{}
\crefname{subsection}{Subsection}{Subsections}

\hypersetup{
    colorlinks,
    linkcolor={red!50!black},
    citecolor={green},
    urlcolor={blue!80!black}
}

\theoremstyle{plain}
\newtheorem{theorem}{Theorem} [section]
\newtheorem{lemma}[theorem]{Lemma}
\newtheorem{prop}[theorem]{Proposition}
\newtheorem{cor}[theorem]{Corollary}
\newtheorem{setting}[theorem]{Setting}

\theoremstyle{definition}
\newtheorem{definition}[theorem]{Definition}

\crefname{case}{Case}{Cases}
\crefname{cor}{Corollary}{Corollaries}

\numberwithin{equation}{section}

\DeclareFontEncoding{LS1}{}{}
\DeclareFontSubstitution{LS1}{stix}{m}{n}
\DeclareMathAlphabet{\mathscr}{LS1}{stixscr}{m}{n}

\newcommand{\E}{\mathbb{E}}
\renewcommand{\P}{\mathbb{P}}
\newcommand{\Q}{\mathbb{Q}}
\newcommand{\R}{\mathbb{R}}

\newcommand{\bD}{\mathbb{D}}
\newcommand{\N}{\mathbb{N}}
\newcommand{\bbD}{\mathbb{D}}

\newcommand{\bbA}{\mathbb{A}}

\newcommand{\dens}{\mathfrak{p}}

\newcommand{\const}{\fC}
\newcommand{\constt}{\fD}
\newcommand{\consttt}{\fc}
\newcommand{\bL}{\mathbb{L}}
\newcommand{\bB}{\mathbb{B}}

\newcommand{\width}{\mathfrak{h}}

\newcommand{\w}[1]{\mathfrak{w}^{#1}}
\renewcommand{\b}[1]{\mathfrak{b}^{#1}}
\renewcommand{\v}[1]{\mathfrak{v}^{#1}}
\newcommand{\q}[1]{\mathfrak{q}^{#1}}
\renewcommand{\c}[1]{\mathfrak{c}^{#1}}

\newcommand{\smallsum}{\textstyle\sum}

\newcommand{\analytic}{\mathscr{A}}
\newcommand{\polyn}{\scrP}
\newcommand{\ratio}{\scrR}

\newcommand{\with}{\curvearrowleft}

\newcommand{\cA}{\mathcal{A}}
\newcommand{\cB}{\mathcal{B}}

\newcommand{\cD}{\mathcal{D}}
\newcommand{\cE}{\mathcal{E}}
\newcommand{\cF}{\mathcal{F}}
\newcommand{\cG}{\mathcal{G}}

\newcommand{\cL}{\mathcal{L}}

\newcommand{\cN}{\mathcal{N}}

\newcommand{\cX}{\mathcal{X}}

\newcommand{\bfa}{\mathbf{a}}
\newcommand{\bfb}{\mathbf{b}}

\newcommand{\bfd}{\mathbf{d}}

\newcommand{\bfh}{\mathbf{h}}

\newcommand{\bfk}{\mathbf{k}}

\newcommand{\bfm}{\mathbf{m}}
\newcommand{\bfn}{\mathbf{n}}

\newcommand{\bfv}{\mathbf{v}}
\newcommand{\bfw}{\mathbf{w}}

\newcommand{\bfA}{\mathbf{A}}

\newcommand{\bfD}{\mathbf{D}}

\newcommand{\bfH}{\mathbf{H}}

\newcommand{\bfM}{\mathbf{M}}

\newcommand{\bfR}{\mathbf{R}}
\newcommand{\bfS}{\mathbf{S}}

\newcommand{\scrA}{{\scra}}
\newcommand{\scrB}{{\scrb}}
\newcommand{\scrC}{\mathscr{C}}

\newcommand{\scrF}{\mathscr{F}}

\newcommand{\scrJ}{\mathscr{J}}

\newcommand{\scrL}{\mathscr{L}}

\newcommand{\scrN}{\mathscr{N}}

\newcommand{\scrP}{\mathscr{P}}

\newcommand{\scrR}{\mathscr{R}}

\newcommand{\scrV}{\mathscr{V}}

\newcommand{\scra}{\mathscr{a}}
\newcommand{\scrb}{\mathscr{b}}

\newcommand{\fA}{\mathfrak{A}}
\newcommand{\fB}{\mathfrak{B}}
\newcommand{\fC}{\mathfrak{C}}
\newcommand{\fD}{\mathfrak{D}}

\newcommand{\fG}{\mathfrak{G}}

\newcommand{\fI}{\mathfrak{I}}
\newcommand{\fJ}{\mathfrak{J}}

\newcommand{\fL}{\mathfrak{L}}
\newcommand{\fM}{\mathfrak{M}}

\newcommand{\fP}{\mathfrak{P}}

\newcommand{\fS}{\mathfrak{S}}

\newcommand{\fa}{\mathfrak{a}}
\newcommand{\fb}{\mathfrak{b}}
\newcommand{\fc}{\mathfrak{c}}
\newcommand{\fd}{\mathfrak{d}}

\newcommand{\ff}{\mathfrak{f}}
\newcommand{\fg}{\mathfrak{g}}

\newcommand{\fl}{\mathfrak{l}}

\newcommand{\fp}{\mathfrak{p}}
\newcommand{\fq}{\mathfrak{q}}

\newcommand{\fv}{\mathfrak{v}}
\newcommand{\fw}{\mathfrak{w}}
\newcommand{\fx}{\mathscr{x}}

\renewcommand{\emptyset}{\varnothing}

\DeclarePairedDelimiter{\norm}{\lVert}{\rVert}
\DeclarePairedDelimiter{\abs}{\lvert}{\rvert}
\DeclarePairedDelimiter{\rbr}{(}{)}
\DeclarePairedDelimiter{\br}{[}{]}
\DeclarePairedDelimiter{\cu}{\{}{\}}
\DeclarePairedDelimiter{\spro}{\langle}{\rangle}

\newcommand{\Rect}{\mathfrak{R}}
\renewcommand{\d}{ \mathrm{d}}

\newcommand{\qandq}{\quad\text{and}\quad}
\newcommand{\qqandqq}{\qquad\text{and}\qquad}

\newcommand{\indicator}[1]{\mathbbm{1}_{\smash{#1}}}
\newcommand{\realization}[1] {\mathscr{N} ^{ #1  }}

\usepackage{tikz}
\usetikzlibrary{matrix,chains,positioning,decorations.pathreplacing,arrows}
\usetikzlibrary{shapes,arrows}
\tikzset{
	font={\fontsize{9pt}{12}\selectfont}}
\usepackage{adjustbox}

\ExplSyntaxOn

\bool_new:N \g_noteobserve

\NewDocumentCommand{\setnote}{}{
  \bool_gset_true:N \g_noteobserve
}

\NewDocumentCommand{\setobserve}{}{
  \bool_gset_false:N \g_noteobserve
}

\NewDocumentCommand{\nobs}{ o }{
  \IfValueT{#1}{
    \str_if_eq:noTF {note} {#1} {
      \bool_gset_true:N \g_noteobserve
    } {
      \str_if_eq:noTF {Note} {#1} {
        \bool_gset_true:N \g_noteobserve
      } {
        \bool_gset_false:N \g_noteobserve
      }
    }
  }
  \bool_if:nTF { \g_noteobserve } {
    \bool_gset_false:N \g_noteobserve
    note
  } {
    \bool_gset_true:N \g_noteobserve
    observe
  }
  \IfValueF{#1}{~}
}

\NewDocumentCommand{\Nobs}{ o }{
  \IfValueT{#1}{
    \str_if_eq:noTF {note} {#1} {
      \bool_gset_true:N \g_noteobserve
    } {
      \str_if_eq:noTF {Note} {#1} {
        \bool_gset_true:N \g_noteobserve
      } {
        \bool_gset_false:N \g_noteobserve
      }
    }
  }
  \bool_if:nTF { \g_noteobserve } {
    \bool_gset_false:N \g_noteobserve
    Note
  } {
    \bool_gset_true:N \g_noteobserve
    Observe
  }
  \IfValueF{#1}{~}
}

\int_new:N \g_furthermore

\NewDocumentCommand{\Moreover}{ o o }{
  \IfValueT{#1}{
    \str_case:nn {#1} {
      {Furthermore} {\int_set:Nn {\g_furthermore} {0}}
      {Moreover} {\int_set:Nn {\g_furthermore} {1}}
      {In~addition} {\int_set:Nn {\g_furthermore} {2}}
      {note} {\bool_gset_true:N \g_noteobserve}
      {observe} {\bool_gset_false:N \g_noteobserve}
    }
    \IfValueT{#2}{
      \str_case:nn {#2} {
        {Furthermore} {\int_set:Nn {\g_furthermore} {0}}
        {Moreover} {\int_set:Nn {\g_furthermore} {1}}
        {In~addition} {\int_set:Nn {\g_furthermore} {2}}
        {note} {\bool_gset_true:N \g_noteobserve}
        {observe} {\bool_gset_false:N \g_noteobserve}
      }
    }
  }
  \int_case:nn { \int_mod:nn {\g_furthermore} {3} } {
    { 0 } { Furthermore,~\nobs that}
    { 1 } { Moreover,~\nobs that}
    { 2 } { In~addition,~\nobs that}
  }
  \int_incr:N \g_furthermore
  \IfValueF{#1}{~}
}

\bool_new:N \g_hencetherefore

\NewDocumentCommand{\hence}{}{
  \bool_if:nTF { \g_hencetherefore } {
    \bool_gset_false:N \g_hencetherefore
    hence~
  } {
    \bool_gset_true:N \g_hencetherefore
    therefore~
  }
}

\NewDocumentCommand{\Hence}{}{
  \bool_if:nTF { \g_hencetherefore } {
    \bool_gset_false:N \g_hencetherefore
    Hence,~we~obtain~
  } {
    \bool_gset_true:N \g_hencetherefore
    Therefore,~we~obtain~
  }
}

\seq_new:N \g_cflist_loaded
\seq_new:N \g_cflist_pending

\NewDocumentCommand{\cfadd}{ m }
{
	\seq_if_in:NnF \g_cflist_loaded { #1 } {
		\seq_if_in:NnF \g_cflist_pending { #1 } {
			\seq_gput_right:Nn \g_cflist_pending { #1 }
		}
	}
}

\NewDocumentCommand{\cfconsiderloaded}{ m }{
	\seq_gput_right:Nn \g_cflist_loaded {#1}
}

\NewDocumentCommand{\cfremove}{ m }
{
	\seq_gremove_all:Nn \g_cflist_pending { #1 }
}

\NewDocumentCommand{\cfload}{ o }
{
	\seq_if_empty:NTF \g_cflist_pending {\unskip} {
		(cf.\ \cref{\seq_use:Nn \g_cflist_pending {,}})\IfValueTF{#1}{#1~}{\unskip}
		\seq_gconcat:NNN \g_cflist_loaded \g_cflist_loaded \g_cflist_pending
		\seq_gclear:N \g_cflist_pending
	}
}

\NewDocumentCommand{\cfclear} {} {
	\seq_gclear:N \g_cflist_loaded
	\seq_gclear:N \g_cflist_pending
}

\NewDocumentCommand{\cfout}{ o }
{
	\seq_if_empty:NTF \g_cflist_pending {\unskip} {
		(cf.\ \cref{\seq_use:Nn \g_cflist_pending {,}})\IfValueTF{#1}{#1~}{\unskip}
		\seq_gclear:N \g_cflist_pending
	}
}

\NewDocumentCommand{\ifnocf} { m } {
	\seq_if_empty:NT \g_cflist_pending { #1 }
}

\ExplSyntaxOff

\NewDocumentEnvironment{cproof}{m}
{\begin{proof}[Proof of \cref{#1}]}%
	{\noindent The proof of \cref{#1} is thus complete.
\end{proof}}

\NewDocumentEnvironment{cproof2}{m}
{\begin{proof}[Proof of \cref{#1}]}%
	{\noindent This completes the proof of \cref{#1}.
\end{proof}}

\title{On the existence of global minima and 
convergence analyses for gradient descent 
methods in the training of deep neural networks}

\author{
Arnulf Jentzen$^{1,2}$ 
and 
Adrian Riekert$^{3}$
\bigskip
\\
\small{$^1$ School of Data Science and Shenzhen Research Institute of Big Data,}
\vspace{-0.1cm}\\
\small{The Chinese University of Hong Kong, Shenzhen, China, e-mail: \texttt{ajentzen@cuhk.edu.cn}}
\smallskip
\\
\small{$^2$ Applied Mathematics: Institute for Analysis and Numerics,}
\vspace{-0.1cm}\\
\small{University of M{\"u}nster, Germany, e-mail: \texttt{ajentzen@uni-muenster.de}}
\smallskip
\\
\small{$^3$ Applied Mathematics: Institute for Analysis and Numerics,}
\vspace{-0.1cm}\\
\small{University of M{\"u}nster, Germany, e-mail: \texttt{ariekert@uni-muenster.de}}
}

\date{\today}

\begin{document}

\maketitle

\begin{abstract}
Although gradient descent (GD) optimization methods in combination with 
rectified linear unit (ReLU) artificial neural networks (ANNs) 
often supply an impressive performance 
in real world learning problems, till this day it remains -- in 
all practically relevant scenarios -- an open problem of research to 
rigorously prove (or disprove) the conjecture that such 
GD optimization methods do converge in the training of ANNs 
with ReLU activation.

In this article we study fully-connected feedforward deep ReLU ANNs 
with an arbitrarily large number of hidden layers 
and we prove convergence of the risk of the GD optimization method with random 
initializations in the training of such ANNs under the assumption 
that the unnormalized probability density function 
$ \dens \colon [a,b]^d \to [0,\infty) $ of the probability 
distribution of the input data of the considered supervised learning 
problem is piecewise polynomial, under the assumption that 
the target function $ f \colon [a,b]^d \to \R^{ \delta } $ (describing the 
relationship between input data and the output data) is piecewise 
polynomial, and under the assumption that the risk function 
of the considered supervised learning problem admits 
at least one regular global minimum. In addition, in the special 
situation of shallow ANNs with just one hidden layer and one-dimensional 
input we also verify this assumption by proving in 
the training of such shallow ANNs that for every Lipschitz continuous 
target function there exists a global minimum in the risk landscape. 
Finally, in the training of deep ANNs with ReLU activation we also 
study solutions of gradient flow (GF) differential equations and we prove 
% by proving 
that every non-divergent GF trajectory converges with a polynomial rate 
of convergence to a critical point (in the sense of limiting 
Fr\'{e}chet subdifferentiability).

Our mathematical convergence analysis builds up on 
ideas from our previous article
[S.~Eberle, A.~Jentzen, A.~Riekert, \& G.~Weiss, Existence, uniqueness, and convergence rates for gradient flows in the training of artificial neural networks with ReLU activation. \emph{arXiv:2108.08106} (2021)],
on tools from real 
algebraic geometry such as the concept of semi-algebraic functions 
and generalized Kurdyka-\L ojasiewicz inequalities, on tools from functional analysis 
such as the Arzelà--Ascoli theorem on the relative compactness of 
uniformly bounded and equicontinuous sequences of continuous functions, on tools from nonsmooth 
analysis such as the concept of limiting Fr\'{e}chet subgradients, as well as 
on the fact that the set of realization functions of shallow ReLU ANNs 
with fixed architecture forms a closed subset of the set of continuous 
functions revealed in [P.\ Petersen, M.\ Raslan, \& F.\ Voigtlaender, 
Topological properties of the set of functions generated by 
neural networks of fixed size. \emph{Found. Comput. Math.}\ {\bf 21} (2021), no.\ 2, 
375--444]. 
\end{abstract}

\pagebreak 

\tableofcontents

\pagebreak

\section{Introduction and main results}

Even though gradient descent (GD) type optimization methods in combination with 
artificial neural networks (ANNs) often supply an impressive performance 
in real world learning problems, till this day it remains -- in 
all practically relevant scenarios -- an open problem of research to 
rigorously prove (or disprove) the conjecture that such 
GD optimization methods do converge in the training of ANNs. 
Moreover, in the case of ANNs with the widely-used rectified 
linear unit (ReLU) activation function, 
this problem of research receives additional difficulty 
due to the lack of differentiability of the rectifier function 
$ \R \ni x \mapsto \max\{ x, 0 \} \in \R $.

Although the convergence analysis for GD type optimization methods 
in the training of ANNs remains a fundamental open problem of research, 
there are several auspicious approaches in the 
scientific literature which provide interesting 
first steps in this area of research. 
To briefly introduce the reader to this topic of research, 
we now highlight/mention some of those findings 
in a short way and refer to the below mentioned references 
for further details.

In particular,
we refer, for example, to~\cite{AroraDuHuLiWang2019, BachChizatOyallon2019,DuLeeLiWangZhai2019,DuZhaiPoczosSingh2019,EMaWu2020,JacotGabrielHongler2018,JentzenKroeger2021,ZhangMartensGrosse2019,ZouCaoZhouGu2019} for convergence results for gradient flow (GF) and GD processes in the training of ANNs in the so-called overparametrized regime, where the number of ANN parameters highly exceeds the number of considered input-output training data pairs.
As the number of neurons goes to infinity, the corresponding GF processes converge to a measure-valued process called 
Wasserstein gradient flow; 
cf., for instance,~\cite{ZhengdaoRotskoff2020,Chizat2021,ChizatBach2018}, 
\cite[Section 5.1]{EMaWojtowytschWu2020}, and the references mentioned therein.

Regarding abstract results on the convergence of GF and GD processes we refer, 
for example, 
to~\cite{BachMoulines2013, JentzenKuckuckNeufeldVonWurstemberger2021, BachMoulines2011, Nesterov2004, Rakhlin2012} for the case of convex objective functions,
we refer, for instance, to~\cite{AbsilMahonyAndrews2005,AttouchBolte2009,AttouchBolteSvaiter2013,BolteDaniilidis2006,DereichKassing2021,Karimi2020linear,LeePanageasRecht2019,LeeJordanRecht2016,LiMilzarekQiu2021,Lojasiewicz1984,Ochs2019} for convergence results for GF and GD processes under \L ojasiewicz type conditions,
and
 we refer, for instance, to~\cite{BertsekasTsitsiklis2000, FehrmanGessJentzen2020,LeiHuLiTang2020,Patel2021stopping} and the references mentioned therein for further results
 without convexity conditions.
 In general, without global assumptions on the objective function such as convexity, 
 gradient-based methods may converge to non-global local minima or saddle points. 
 It therefore becomes important to analyze critical points of the objective function 
 in the training of ANNs 
 and we refer, for example, 
 to \cite{Cheridito2022Landscape,Safran2018spurious,Swirszcz2017local,ZhangLi2021embedding,Zhang2021embedding} for articles which study the appearance 
 of critical points in the risk landscape in the training of 
 ANNs. 
 The question under which conditions gradient-based optimization algorithms cannot 
 converge to saddle points was investigated, for example, 
 in \cite{Ge2015,LeePanageasRecht2019, LeeJordanRecht2016, PanageasPiliouras2017, Panageas2019firstorder}.
For more detailed overviews and additional references on GD optimization schemes 
we mention, for instance, Bottou et al.~\cite{Bottou2018optimization}, Fehrman et al.~\cite[Section 1.1]{FehrmanGessJentzen2020}, \cite[Section 1]{JentzenKuckuckNeufeldVonWurstemberger2021}, and Ruder~\cite{Ruder2017overview}.

In this article we study the training of fully-connected feedforward 
ANNs with ReLU activation by means of GD type optimization methods 
(we also refer to \cref{figure_shallow} and \cref{figure_deep} 
in this introductory section below for graphical illustrations 
of two example architectures for the ANNs investigated in this work). 
In particular, one of the key contributions of this work is rigorously verify, 
under the assumption that the unnormalized probability density function 
$
  \dens \colon [a,b]^d \to [0,\infty) 
$
of the probability distribution of the 
input data of the considered supervised learning problem 
is piecewise polynomial 
(see \cref{def:multidim:piece:polyn} in \cref{sec:piecewise_polynomial} 
for our precise meaning of a piecewise polynomial function), 
and the assumption that the target function 
$
  f = ( f_1, \dots, f_{ \delta } ) \colon [a,b]^d \to \R^{ \delta }
$
(the function describing the relationship 
between the input data and the output data which one intends to learn approximately)
is piecewise polynomial, 
it holds in the training of \emph{deep ReLU ANNs with an arbitrarily large 
number of hidden layers} 
that the risk function (the function which is to be minimized) 
and its associated generalized gradient function 
satisfy at \emph{every point} of the ANN parameter space a 
\emph{generalized Kurdyka-\L ojasiewicz inequality}; 
see \cref{prop:loss:lojasiewicz} in \cref{subsection:loja} for the precise statement. 
In the previous sentence the quantity $ d \in \N = \{ 1, 2, 3, \dots \} $ 
is an arbitrarily large natural number which describes the dimension 
of the input data, 
the quantity $ \delta \in \N $ is a natural number which describes the dimension 
of the output data, 
and the quantities $ a, b \in \R $ with $ a < b $ are real numbers 
which border the region $ [a,b]^d $ in which the input data 
takes values in. 
\cref{prop:loss:lojasiewicz} in \cref{subsection:loja} in this work generalizes 
Proposition~5.1 in our previous article Eberle et al.~\cite{EberleJentzenRiekertWeiss2021}
where such generalized Kurdyka-\L ojasiewicz inequalities have been 
established in training of ReLU ANNs with one hidden layer. 
The proof of \cref{prop:loss:lojasiewicz} relies on the fact that the considered risk function is semi-algebraic, which we establish in \cref{cor:loss:semialgebraic} below, and the abstract Kurdyka-\L ojasiewicz inequality in Bolte et al.~\cite[Theorem 3.1]{BolteDaniilidis2006}.
A similar result regarding
semi-algebraicity of the empirical risk, which is measured with respect to a finite set of
input-output data pairs, is already known, cf.~Davis et al.~\cite[Corollary 5.11]{Davis2018stochastic}.

We then use the established generalized Kurdyka-\L ojasiewicz inequalities 
in \cref{prop:loss:lojasiewicz}
to prove convergence of GD type optimization methods in the training of deep ReLU ANNs 
where we first focus on \emph{time-continuous} GD optimization methods 
(see \cref{section:gf:loja}) 
and, thereafter, investigate 
\emph{time-discrete} GD optimization methods 
(see \cref{sec:convergence_GD}).

Specifically, in the time-continuous situation 
(see \cref{section:gf:loja} and \cref{sec:intro_GF} in this introductory section) 
we establish in the training of deep ReLU ANNs, 
under the assumption that the unnormalized probability density function 
$
  \dens \colon [a,b]^d \to [0,\infty)
$ 
and the target function 
$
  f \colon [a,b]^d \to \R
$
are both piecewise polynomial,  
that every non-divergent solution of the associated 
\emph{gradient flow (GD) differential equation} 
converges with a strictly positive rate of convergence 
to a \emph{generalized critical point} of the risk function 
(in the sense of the limiting Fr\'{e}chet subdifferential; 
see \cref{def:limit:subdiff} in \cref{sec:subdifferential}) 
and also that the risk of the GF solution converges with rate $ 1 $ 
to the risk of the generalized critical point 
(see \cref{theo:gf:conv:simple} in \cref{sec:convergence_GF_processes} below 
and \cref{thm:main_thm} in \cref{sec:intro_GF} in this introductory section below, respectively, 
for the precise statements). 
This generalizes the approach in 
Eberle et al.~\cite[Subsection~5.2]{EberleJentzenRiekertWeiss2021} 
from shallow ReLU ANNs to deep ReLU ANNs.

Moreover, in the time-discrete situation 
(see \cref{sec:convergence_GD} and 
\cref{subsec:existence_global_minima_intro,susbsec:GD_deep_intro} in this introductory section)
we establish in the training of deep ReLU DNNs, 
under the assumption that 
$
  \dens 
$ 
and 
$
  f 
$
are piecewise polynomial 
and that the risk function of the considered deep supervised learning problem 
admits at least one regular global minimum point, 
that the risk of the plain vanilla GD optimization method with 
random initializations converges in the training of deep ReLU ANNs 
to $ 0 $ as the number of GD steps increases to $ \infty $, 
as the number of random initializations increases to $ \infty $, 
as the step size of the GD method (the learning rate of the GD method) decreases to $ 0 $, 
and as the width of the ANNs increases to $ \infty $; 
see \cref{theo:gd:random:dnn} 
in \cref{sec:convergence_GD} below 
and \cref{theo:intro:convergence} in \cref{susbsec:GD_deep_intro} 
in this introductory section below, respectively, 
for the precise statement.

Another key contribution of this work 
(see \cref{sec:existence_global_minima}) 
is to prove in the special situation 
of shallow ReLU ANNs with just one hidden layer 
and one-dimensional input and output (corresponding to the case $ d = \delta = 1 $) 
that for every Lipschitz continuous target function 
$ 
  f \colon [a,b] \to \R
$
we have that there exist global minimum points of the 
risk function; 
see \cref{theo:existence2} in \cref{subsec:existence_global_minima} 
and \cref{theo:existence_global_minima} in \cref{subsec:existence_global_minima_intro} 
in this introductory section below, respectively, 
for the precise statement. 
In the case of shallow ANNs we thereby verify the above mentioned assumption 
that the risk function of the considered supervised learning problem 
admits at least one regular global minimum point; 
cf.\ \cref{cor:existence:regular} in \cref{subsec:existence_regular_minima} below 
and 
\cref{theo:intro:random:init} in \cref{susbsec:GD_shallow_intro}
in this introductory section below, respectively.

To elucidate the findings of this work more clearly, 
we now present $ 4 $ selected specific results 
(which have already been briefly outlined in the above introductory paragraphs)
regarding the training of ReLU ANNs, 
\cref{theo:existence_global_minima} in \cref{subsec:existence_global_minima_intro}, 
\cref{theo:intro:random:init} in \cref{susbsec:GD_shallow_intro}, 
\cref{thm:main_thm} in \cref{sec:intro_GF}, 
and 
\cref{theo:intro:convergence} in \cref{susbsec:GD_deep_intro}, 
with all details in a self-contained fashion. 
\cref{theo:existence_global_minima}
and 
\cref{theo:intro:random:init}
deal with shallow ReLU ANNs with just one hidden layer 
and one-dimensional input and output ($ d = \delta = 1 $)
and 
\cref{thm:main_thm} 
and 
\cref{theo:intro:convergence} 
treat the situation of deep ReLU ANNs with an arbitrarily large 
number of hidden layers 
and multi-dimensional input and output ($ d, \delta \in \N $).

\subsection{Existence of global minima for shallow artificial neural networks (ANNs)}
\label{subsec:existence_global_minima_intro}

Maybe the most basic question that one can ask regarding the training of ANNs 
is the existence of global minimum points in the risk landscape. 
In particular, without the existence of a global minimum point, 
one can not hope for a GD type optimization method to converge to 
a global minimum point. 
Surprisingly, there is almost no result in the scientific literature 
which actually establishes the existence of global minimum points 
of risk functions in the training of ANNs 
and in our perspective this subject 
is a very important direct of future research.

\cref{theo:existence_global_minima} below proves in the training of 
shallow ANNs with ReLU activation that for every distribution 
$ \mu \colon \cB( [a,b] ) \to [0,\infty] $ 
of the input data of the considered supervised learning problem 
and every Lipschitz continuous target function $ f \colon [a,b] \to \R $
it holds that there exists a global minimum point 
of the risk function. 
The natural number $ \width \in \N $ in \cref{theo:existence_global_minima} 
specifies the number of neurons on the hidden layer of the ANN 
(the dimensionality of the hidden layer of the ANN), 
the natural number $ \fd \in \N $ in \cref{theo:existence_global_minima} 
specifies the overall number of real parameters used to described 
to the considered ANNs, 
and the set $ \cB( [a,b] ) $ is the Borel sigma-algebra 
on the real interval $ [a,b] \subseteq \R $.

\begin{figure}
	\centering
	\begin{adjustbox}{width=\textwidth}
		\begin{tikzpicture}[shorten >=1pt,->,draw=black!50, node distance=\layersep]
		\tikzstyle{every pin edge}=[<-,shorten <=1pt]
		\tikzstyle{input neuron}=[very thick, circle,draw=red, fill=red!30, minimum size=15pt,inner sep=0pt] 
		%\tikzstyle{input neuron}=[neuron];
		\tikzstyle{output neuron}=[very thick, circle,draw=green, fill=green!30, minimum size=20pt,inner sep=0pt]
		\tikzstyle{hidden neuron}=[very thick, circle, draw=blue, fill=blue!30, minimum size=15pt,inner sep=0pt]
		\tikzstyle{annot} = [text width=9em, text centered]
		\tikzstyle{annot2} = [text width=4em, text centered]
		\node[input neuron] (I) at (0,-2.5 cm) {$x$};
		\foreach \name / \y in {1,...,5}
		\path[yshift=0.5cm]
		node[hidden neuron] (H-\name) at (5 cm,-\y cm) {};      
		\path[yshift=-1cm]
		node[output neuron](O) at (10 cm,-1.5 cm) {$\realization{\theta}(x)$};   
		\foreach \dest in {1,...,5}
		\path[line width = 0.8] (I) edge (H-\dest);     
		\foreach \source in {1,...,5}
		\path[line width = 0.8] (H-\source) edge (O);     
		\node[annot,above of=H-1, node distance=1cm, align=center] (hl) {Hidden layer\\(2nd layer)};
		\node[annot, above of=I, node distance = 3cm, align=center] {Input layer\\ (1st layer)};
		\node[annot, above of=O, node distance=3cm, align=center] {Output layer\\(3rd layer)};		
		\node[annot2,below of=H-5, node distance=1cm, align=center] (sl) { $\width=5$};
		\end{tikzpicture}
	\end{adjustbox}
	\caption{Graphical illustration of the considered shallow ANN architecture in \cref{theo:existence_global_minima,theo:intro:random:init} in the special case of an ANN with $\width = 5$ neurons on the hidden layer. In this situation we have for every ANN parameter vector $\theta \in \R^{ \mathfrak{d} } = \R^{16 }$ that the realization function 
		$\R \ni x \mapsto \realization{  \theta }(x) \in \R$ 
		of the considered ANN maps the scalar input $ x \in [a,b] $ 
		to the scalar output 
		$
		  \realization{\theta} ( x ) = \theta_{ \fd } 
  		+ \sum_{ j = 1 }^{ \width } 
		\theta_{ 2 \width + j } 
		\max\{ \theta_j x + \theta_{ \width + j } , 0 \} \in \R $.}
	\label{figure_shallow}
\end{figure}
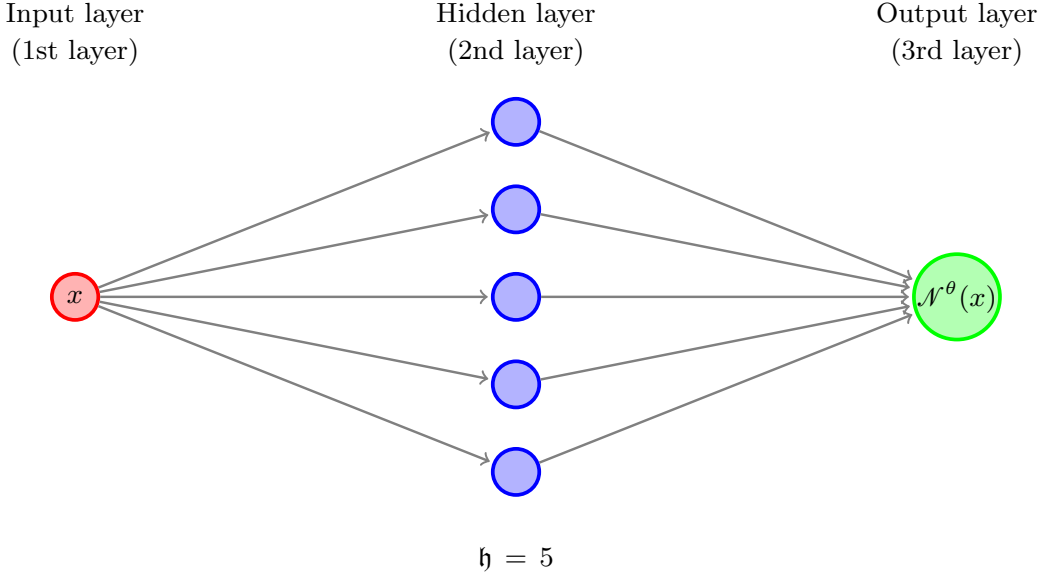

In \cref{theo:existence_global_minima} we thus consider 
ANNs with $ 1 $ neuron on the input layer (with a $ 1 $-dimensional input layer), 
with $ \width $ neurons on the hidden layer (with $ \width $-dimensional hidden layer), 
and with $ 1 $ neuron on the output layer (with a $ 1 $-dimensional output layer). 
There are hence $ \width $ real weight parameters and $ \width $ real bias parameters 
to describe the affine linear transformation between the $ 1 $-dimensional input layer 
and the $ \width $-dimensional hidden layer and $ \width $ real weight parameters 
and $ 1 $ real bias parameter to describe the affine linear transformation between 
the $ \width $-dimensional hidden layer and 
the $ 1 $-dimensional output layer. 
The overall number $ \fd \in \N $ 
of real ANN parameters 
in \cref{theo:existence_global_minima} therefore satisfies 
$
  \fd = ( \width + \width ) + ( \width + 1 ) = 3 \width + 1 
$. 
We also refer to \cref{figure_shallow} for a graphical illustration 
of an example architecture for the ANNs considered in 
\cref{theo:existence_global_minima}.

The function $ \cL \colon \R^{ \fd } \to \R $ 
in \cref{eq:risk_theorem_intro_shallow} 
in \cref{theo:existence_global_minima} 
is the risk function in the considered supervised learning problem 
and the finite measure 
$ \mu \colon \cB( [a,b] ) \to [0, \infty] $ 
is the unnormalized probability distribution 
of the input data of the considered supervised learning problem. 
In \cref{theo:existence_global_minima} 
we considered 
ReLU ANNs and 
the ReLU activation function 
$ \R \ni x \mapsto \max\{ x, 0 \} \in \R $
appears
on the right hand side of \cref{eq:risk_theorem_intro_shallow}. 
In this set-up of shallow ReLU ANNs 
\cref{theo:existence_global_minima} reveals the existence 
of a global minimum point 
$ \theta = ( \theta_1, \dots, \theta_{ \fd } ) \in \R^{ \fd } $ 
of the risk function 
$
  \cL \colon \R^{ \fd } \to \R 
$.

\begin{samepage}
\begin{theorem}
\label{theo:existence_global_minima}
Let $ \width, \fd \in \N $, $ a \in \R $, $ b \in [a, \infty) $
satisfy $ \fd = 3 \width + 1 $, 
let $ f \colon [a,b] \to \R $ be Lipschitz continuous, 
let $ \mu \colon \cB( [a,b] ) \to [0, \infty] $ be a finite measure, 
and let $ \cL \colon \R^{ \fd } \to \R $ 
satisfy for all 
$ \theta = ( \theta_1, \ldots, \theta_{ \fd } ) \in \R^{ \fd } $ 
that
\begin{equation} 
\label{eq:risk_theorem_intro_shallow}
  \cL( \theta ) = 
  \textstyle
  \int_a^b 
  \displaystyle
    ( 
      f( x ) - \theta_{ \fd } 
      - \smallsum_{ j = 1 }^{ \width } 
      \theta_{ 2 \width + j } 
      \max\{ \theta_{ \width + j } + \theta_j x, 0 \} 
    )^2 
  \, \mu( \d x ) .
\end{equation}
Then there exists $ \theta \in \R^{ \fd } $ 
such that 
$
  \cL( \theta ) = 
  \inf_{ \vartheta \in \R^{ \fd } } 
  \cL( \vartheta )
$.
\end{theorem}
\end{samepage}

\cref{theo:existence_global_minima} is an immediate consequence 
of \cref{theo:existence2} in \cref{subsec:existence_global_minima} below. 
\cref{theo:existence_global_minima} proves that there exists 
an ANN parameter vector 
$ \theta = ( \theta_1, \dots, \theta_{ \fd } ) \in \R^{ \fd } $ 
which satisfies that the risk 
$ \cL( \theta ) $
of $ \theta $ 
coincides with the infimum over all risk values
$
  \inf_{ \vartheta \in \R^{ \fd } } 
  \cL( \vartheta )
$. 

As observed in Petersen et al.~\cite{PetersenRaslanVoigtlaender2020}, the existence of global minima has direct implications for the training of ANNs. In particular, if there is no global minimum then any sequence $(\theta_n)_{n  \in \N} \subseteq \R^\fd$ with $\lim_{n \to \infty} \cL ( \theta_n ) = \inf_{ \vartheta \in \R^{ \fd } } 
\cL( \vartheta )$ necessarily diverges to infinity. This behavior is highly undesirable in numerical computations.
If the target function $f$ is not continuous, this divergence phenomenon can indeed be observed in practice, as the results and numerical examples in \cite{GallonJentzen2022} show.
On the other hand,
using our existence result for global minima we are able to establish convergence of GD with random initializations in the training of shallow ANNs if the assumptions of \cref{theo:existence_global_minima} are satisfied, see the next subsection for details.

In the scientific literature 
a similar existence result for ANNs with the Heaviside activation function 
$
  \R \ni x \mapsto \mathbbm{1}_{ [0,\infty) }( x ) \in \R 
$
was established in Kainen et al.~\cite{KAINEN2000695}. 
Moreover, we would like to point out that \cref{theo:existence_global_minima} 
does in general not hold without the Lipschitz continuity assumption on $ f $. 
Indeed, Petersen et al.~\cite[Theorem 3.1]{PetersenRaslanVoigtlaender2020} 
implies in the case where $ \width \geq 2 $ 
and where the measure 
$ \mu $ is non-atomic in the sense that its support is uncountable
that 
the set of realization functions 
\begin{multline}
\textstyle 
  \bigl\{ 
    v \in C( [0,1], \R ) \colon
\\
\textstyle 
    \bigl(
      \exists \, 
      \theta = ( \theta_1, \dots, \theta_{ \fd } ) \in \R^{ \fd } 
      \colon 
      \forall \, x \in [0,1] \colon 
      v(x) = 
      \theta_{ \fd } + \sum_{ j = 1 }^{ \width } 
      \theta_{ 2 \width + j } 
      \max\{ \theta_{ \width + j } + \theta_j x, 0 \}
    \bigr)
  \bigr\} 
\end{multline}
is not closed in the $ L^2 $-space 
$
  L^2( [0,1], \mu ) 
$. 
Specifically, 
Petersen et al.~\cite[Theorem 3.1]{PetersenRaslanVoigtlaender2020} 
shows that there exists $ f \in L^{ \infty }( [0,1], \mu ) $ such that 
$
  \inf_{ \theta \in \R^{ \fd } } \cL( \theta ) = 0
$
and 
$
  \{ \theta \in \R^{ \fd } \colon \cL( \theta ) = 0 \} = \emptyset 
$. 
The function $ f $ constructed in 
Petersen et al.~\cite{PetersenRaslanVoigtlaender2020} is a step function of the form $f (x ) = \indicator{(x^* , 1 ]} ( x )$ for some suitable $x^* \in (0 , 1 )$ depending on the measure $\mu$ and, thus,
does not have a continuous representative. 
A similar non-closedness statement for the logistic activation function 
was proved earlier in Girosi \& Poggio~\cite{Girosi1990bio}. 
In \cref{theo:existence_global_minima} we assume that $ f \colon [a,b] \to \R $ 
is Lipschitz continuous. We guess that the statement remains true 
if one only assumes that $ f $ is continuous.

\subsection{Gradient descent (GD) with random initializations 
in the training of shallow ANNs}
\label{susbsec:GD_shallow_intro}

In \cref{theo:intro:random:init} below we employ 
\cref{theo:existence_global_minima} above 
to establish in the training of shallow ReLU ANNs 
(with $ 1 $ neuron on the input layer, 
$ \width \in \N $ neurons on the hidden layer,
and $ 1 $ neuron on the output layer) 
that the risk of the plain vanilla GD optimization method 
with random initializations 
\emph{converges in probability} to $ 0 $ 
as the number $ n \in \N $ of GD steps increases to $ \infty $, 
as the number $ K \in \N $ of random initializations increases to $ \infty $, 
as the step size $ \gamma \in (0,\infty) $ of the GD optimization method 
(the learning rate of the GD optimization method) decreases to $ 0 $, 
and as the width $ \width \in \N $ of the considered ANNs increases to $ \infty $; 
see \cref{theo:intro:random:init:v4:eq:qresult_shallow} 
in \cref{theo:intro:random:init} below for the precise statement.

In \cref{theo:intro:random:init} we consider the GD optimization method 
with random initializations and the triple $ ( \Omega, \cF, \P ) $ 
in \cref{theo:intro:random:init} serves as the underlying probability space 
for the random initializations. 
Note that the function which maps random variables 
$
  X \colon \Omega \to \R
$
and 
$
  Y \colon \Omega \to \R
$
to the real number 
\begin{equation}
\label{eq:convergence_in_probability}
  \E\bigl[ \min\{ | X - Y |, 1 \} \bigr] 
\end{equation}
is nothing else but one commonly used metric 
which characterizes \emph{convergence in probability} 
(cf., e.g., Klenke~\cite[Theorem 6.7 in Chapter~6]{Klenke2014}) 
and \cref{theo:intro:random:init:v4:eq:qresult_shallow} 
in \cref{theo:intro:random:init} thus indeed 
establishes convergence in probability of the risk 
of the GD optimization method to $ 0 $.

\begin{samepage}
\begin{theorem} 
\label{theo:intro:random:init}
	Let $ N \in \N $, 
	$ \fx_0, \fx_1, \ldots, \fx_N, a, b \in \R $, 
	satisfy 
	$ a = \fx_0 < \fx_1 < \cdots < \fx_N = b $, 
	let 
	$ 
	f \in C( [a,b], \R )
	$, 
	let 
	$ \dens \colon [a,b] \to [0, \infty) $ be a function, 
	assume for all $ n \in \{ 1, \ldots, N \} $ that 
	$ f|_{ ( \fx_{n-1}, \fx_n ) } $ 
	and 
	$ \dens|_{ ( \fx_{n-1}, \fx_n ) } $ 
	are polynomials, 
	for every $ \width \in \N $ 
	let 
	$ \cL_\width \colon \R^{3 \width + 1} \to \R $
	satisfy for all 
	$ \theta = ( \theta_1, \ldots, \theta_{ 3 \width + 1 } ) \in \R^{ 3 \width + 1 } $ 
	that
	\begin{equation} 
	\label{theo:intro:random:init:v4:eq2}
	\cL_\width( \theta ) = 
	\textstyle
	\int_a^b \rbr[\big]{ 
		f( x ) 
		- 
		\theta_{\fd} - 
		\smallsum_{j=1}^\width 
		\theta_{2 \width + j} 
		\max\{ \theta_{j} x + \theta_{\width + j } , 0 \} }^2 
	\dens ( x ) \, \d x,
	\end{equation} 
	for every $ \width \in \N $
	let 
	$ \cG_\width \colon \R^{ 3 \width + 1 } \to \R^{ 3 \width + 1 } $ 
	satisfy for all 
	$ 
	\theta \in 
	\{ 
	\vartheta \in \R^{3 \width + 1} 
	\colon 
	\cL_{ \width }
	\ \text{is}\ \allowbreak \text{differentiable}\ \allowbreak \text{at}\ \vartheta
	\}
	$
	that 
	$ 
	\cG_{ \width }( \theta ) = ( \nabla \cL_{ \width } )( \theta ) 
	$,
	let $ (\Omega, \cF, \P) $ be a probability space, 
	for every 
	$ n, \width, K \in \N_0 $, $ \gamma \in \R $
	let 
	$ \Theta^{ K, \gamma }_{ \width, n } \colon \Omega \to \R^{ 3 \width + 1 } $
	and 
	$ \bfk^{ K, \gamma }_{ \width, n } \colon \Omega \to \N $ 
	be random variables,
	assume for all $ \width \in \N $, $ \gamma \in \R $ that 
	$ \Theta_{ \width, 0 }^{ K, \gamma } $, $ K \in \N $, are i.i.d., 
	assume for all $ \width \in \N $, $ \gamma, r \in (0,1) $, 
	$ \theta \in \R^{ 3 \width + 1 } $ that
	$
	\P( 
	\norm{ \Theta^{ 1 , \gamma }_{ \width, 0 } - \theta } < r
	)
	> 0
	$, 
	and assume for all 
	$ n, \width \in \N_0 $, $ K \in \N $, $ \gamma \in \R $, $ \omega \in \Omega $ 
	that
	\begin{equation} 
	\label{theo:intro:random:init:v4:eq3}
	\Theta_{ \width, n+1 }^{ K, \gamma }( \omega ) 
	= \Theta_{ \width, n }^{ K, \gamma }( \omega ) 
	- \gamma \cG_\width( \Theta_{ \width, n }^{ K, \gamma }( \omega ) ) 
	\quad\text{and}\quad 
	\bfk^{ K, \gamma }_{ \width, n }( \omega) 
	\in 
	\argmin\nolimits_{ \kappa \in \{ 1, \ldots, K \} } 
	\cL_{ \width }( \Theta_{ \width, n }^{ \kappa, \gamma }( \omega ) ).
	\end{equation}
Then
\begin{equation}  
\label{theo:intro:random:init:v4:eq:qresult_shallow}
  \limsup\nolimits_{ \width \to \infty } 
  \limsup\nolimits_{ \gamma \searrow 0 }
  \limsup\nolimits_{ K \to \infty } 
  \limsup\nolimits_{ n \to \infty }
  \mathbb{E}\bigl[ 
    \min\bigl\{
      \cL_{ \width }(
        \Theta^{ \width, \bfk^{ \width, K, \gamma }_n, \gamma }_n 
      )
      , 1 
    \bigr\}
  \bigr] 
  = 0 .
\end{equation}
\end{theorem}
\end{samepage}

\cref{theo:intro:random:init} is a direct consequence of 
\cref{theo:gd:random:init:item2} in \cref{cor:gd:random:shallow} 
in \cref{subsec:convergence_GD_shallow_ANNs} below 
and the reversed version of Fatou's lemma. 
Note that in \cref{theo:intro:random:init:v4:eq3} above the random index $	\bfk^{ K, \gamma }_{ \width, n } ( \omega ) \in \N$ selects the trajectory with the minimal risk after $n \in \N$ gradient steps among the first $K \in \N$ random initializations.
Observe that \cref{theo:intro:random:init:v4:eq:qresult_shallow} 
demonstrates that 
the risk 
$
  \cL_{ \width }(  
    \Theta^{ \width, \bfk^{ \width, K, \gamma }_n, \gamma }_n 
  )
$
of the GD optimization method with random initializations
% $
%   \Theta^{ \width, \bfk^{ \width, K, \gamma }_n, \gamma }_n 
% $
converges in probability 
(see \cref{eq:convergence_in_probability} above) 
to $ 0 $ 
as the number $ n $ of GD steps increases to $ \infty $, 
as the number $ K $ of random initializations increases to $ \infty $, 
as the learning rate $ \gamma $ decreases to $ 0 $,
and 
as the number $ \width $ of neurons on the hidden layer 
increases to $ \infty $. 

Roughly speaking, the proof of \cref{theo:intro:random:init} consists of the following steps.
\begin{enumerate} [label = (\Roman*)]
	\item We strengthen the existence result for global minima from \cref{theo:existence2} by proving in \cref{cor:existence:regular} that each of the risk functions $\cL_\width$, $\width \in \N$, admits a global minimum around which suitable regularity conditions are satisfied.
	\item We establish in \cref{cor:loss:semialgebraic} that the considered risk functions are semi-algebraic. As a consequence, we show in \cref{prop:loss:lojasiewicz} a generalized Kurdyka-\L ojasiewicz inequality for the risk functions.
	\item In \cref{cor:gd:local:simple} below we show an abstract local convergence result to local minima for GD under a Kurdyka-\L ojasiewicz type assumption and a suitable regularity assumption. Specifically, we assume that the considered
	local minimum admits a neighborhood on which the objective function is differentiable with a Lipschitz continuous gradient.
	\item As a consequence, we obtain in \cref{cor:gd:multi:random:init} an abstract convergence result for GD processes with random initializations. Due to the first two steps, \cref{cor:gd:multi:random:init} is applicable to each risk function $\cL_\width$, $\width \in \N$, under the assumptions of \cref{theo:intro:random:init}.
\end{enumerate}
%
%The proof of \cref{theo:intro:random:init} relies on the generalized Kurdyka-\L ojasiewicz inequality for the risk function which we establish in \cref{prop:loss:lojasiewicz} below. This inequality allows us to employ the abstract convergence result for GD processes with random initializations in \cref{cor:gd:multi:random:init} below.
%\cref{cor:gd:multi:random:init}, in turn, relies on the local convergence result to local minima in \cref{cor:gd:local:simple} below. In \cref{cor:gd:local:simple} we require certain regularity conditions for the objective function on a neighborhood of the considered critical point. 
%As a consequence of the results in \cref{sec:existence_global_minima} we have in the situation of \cref{theo:intro:random:init} that the risk function admits a global minimum which satisfies these regularity conditions (see \cref{cor:existence:regular} below for details). Hence, \cref{cor:gd:multi:random:init} is applicable to each risk function $\cL_\width$, $\width \in \N$, under the assumptions of \cref{theo:intro:random:init}.
In \cite[Theorem~1.1]{JentzenRiekert2021rates} a GD convergence result related 
to \cref{theo:intro:random:init} above has been established. 
Roughly speaking, in \cite[Theorem~1.1]{JentzenRiekert2021rates} 
a convergence result similar to \cref{theo:intro:random:init:v4:eq:qresult_shallow} 
has been obtained in the situation 
where the learning rate $ \gamma \in (0,\infty) $ 
must be sufficiently small but may be chosen to be independent 
of the number $ \width \in \N $ of neurons on the hidden layer, 
where the target function
$ f \colon [a,b] \to \R $ must not only be piecewise polynomial 
but even piecewise affine linear, 
and where the unnormalized probability density 
function $ \dens \colon [a,b] \to [0,\infty) $ 
does not necessarily have to be piecewise polynomial 
but instead must be strictly positive and Lipschitz continuous.

The convergence analysis in \cite{JentzenRiekert2021rates} 
follows a completely different strategy than the 
convergence analysis in this work. 
In particular, in contrast to the proof of \cref{theo:intro:random:init} above,  
the proof in \cite[Theorem~1.1]{JentzenRiekert2021rates} 
does not at all use generalized Kurdyka-\L ojasiewicz inequalities 
but instead is based on differential geometric arguments 
and analyses of the Hessian matrices of the risk function 
(cf.\ \cite{FehrmanGessJentzen2020}).

\subsection{Gradient flows (GFs) in the training of deep ANNs}
\label{sec:intro_GF}

In \cref{thm:main_thm} below we demonstrate 
in the training of deep ReLU ANNs 
with an arbitrarily large number of hidden layers, 
under the assumption that the 
unnormalized probability density function 
$
  \dens \colon [a,b]^{ \ell_0 } \to [0,\infty)
$
and the target function 
$ f \colon [a,b]^{ \ell_0 } \to \R^{ \ell_L } $
are piecewise polynomial 
(see \cref{theo:intro:conv:eq0_GF} below for details), 
that every non-divergent solution 
$ \Theta_t $, $ t \in [0,\infty) $, 
of the associated GF differential equation 
converges with a strictly positive rate of convergence 
to a \emph{generalized critical point} $ \vartheta $
(in the sense of the limiting Fr\'{e}chet subdifferential; 
see \cref{def:limit:subdiff} in \cref{sec:subdifferential}) 
and also that the risk 
$ \cL_{ \infty }( \Theta_t ) $, $ t \in [0,\infty) $, 
of the GF solution 
converges with rate $ 1 $ 
to the risk 
$ \cL_{ \infty }( \vartheta ) $
of the generalized critical point; 
see \cref{theo:intro:conv:eq3} in \cref{thm:main_thm} below 
for the precise statement.

The natural number $ L \in \N $ 
in \cref{thm:main_thm} 
specifies the number of affine linear transformations 
in the considered deep ANNs 
(the considered deep ANNs in \cref{thm:main_thm} 
thus consist of $ L - 1 $ hidden layer 
and, including input and output layers, $ L + 1 $ layers overall) 
and the natural numbers 
$ \ell_0, \ell_1, \ell_2, \ldots \in \N $ 
in \cref{thm:main_thm} 
specify the number of neurons of the layers 
in the sense 
that there are $ \ell_0 $ neurons on the input layer 
(the input layer is $ \ell_0 $-dimensional), 
that for every $i \in \cu{1, \ldots, L - 1 }$ there are $ \ell_i $ neurons on the $ i $th hidden layer 
(the $ i $th hidden layer is $ \ell_i $-dimensional), 
and 
that there are $ \ell_L $ neurons on the output layer 
(the output layer is $ \ell_L $-dimensional). 
In the deep ANNs considered in \cref{thm:main_thm}, 
we thus have for every $k \in \cu{1, \ldots, L }$ that there are $ \ell_k \ell_{k-1} $ real weight parameters 
and $ \ell_k $ real bias parameters 
to describe the affine linear transformation 
between the $(k - 1 )$st and the $k$-th layer.
The overall number $ \fd \in \N $ of real ANN parameters 
in \cref{thm:main_thm} thus satisfies 
\begin{equation} 
  \fd = \textstyle \sum_{ k = 1 }^L ( \ell_k \ell_{ k - 1 } + \ell_k ) 
  = \sum_{ k = 1 }^L \ell_k ( \ell_{ k - 1 } + 1 ) 
.
\end{equation} 
We also refer to \cref{figure_deep} 
for a graphical illustration of an example 
architecture for the ANNs considered 
in \cref{thm:main_thm}. 

Because of the lack of differentiability of the ReLU activation function, the risk function $\cL_\infty \colon \R^\fd \to \R$ in \cref{thm:main_thm} is in general not continuously differentiable. 
In order to define an appropriately generalized gradient we approximate, as in \cite{CheriditoJentzenRiekert2022,DNNReLUarXiv,JentzenRiekertFlow,JentzenRiekert2021}, the ReLU function through continuously differentiable functions $\Rect_r \in C^1 ( \R , \R )$, $r \in [1 , \infty ]$ (see \cref{lim_R_thm1} below for details).
 For every $\theta \in \R^\fd$, $r \in [1 , \infty ]$ we define the approximate realization function $\cN^{ L , \theta}_r \colon \R^{ \ell_0 } \to \R^{ \ell_L }$ and the corresponding risk function $\cL_r \colon \R^{ \fd } \to \R$.
  For every parameter vector $\theta \in \R^\fd$ which satisfies that the approximate gradients $(\nabla \cL_r ) ( \theta ) $, $r \in [1 , \infty )$, are convergent as $r \to \infty $ we define the generalized gradient $\cG ( \theta) \in \R^\fd$ as the limit $\lim_{r \to \infty} ( \nabla \cL_r ) ( \theta )$.
In \cref{prop:G} below we verify that this limit, in fact, exists for every $\theta \in \R^\fd$,
and thus the generalized gradient $\cG ( \theta )$ is uniquely defined for every $\theta \in \R^\fd$.
Furthermore, we derive in \cref{prop:G} an explicit formula for the generalized gradient, which agrees with the standard implementation of the gradient via backpropagation.

\def\layersep{2.5cm}
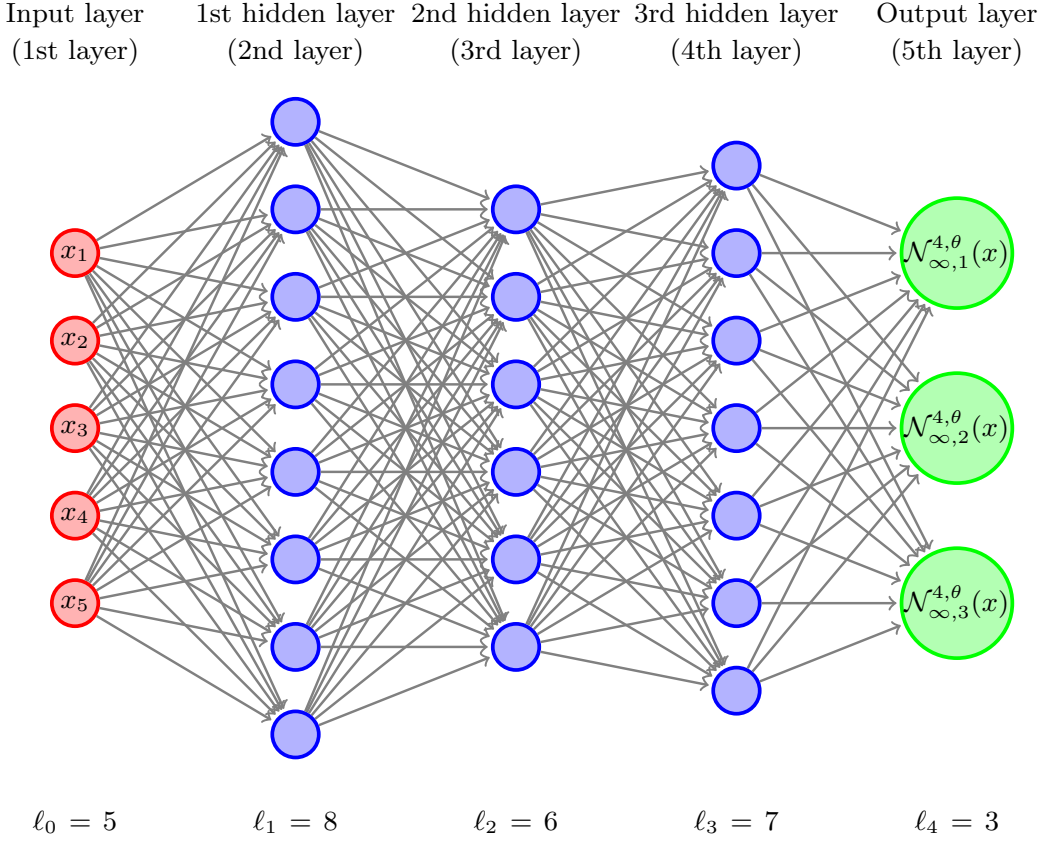
\begin{figure}
	\centering
	\begin{adjustbox}{width=\textwidth}
		\begin{tikzpicture}[shorten >=1pt,->,draw=black!50, node distance=\layersep]
		\tikzstyle{every pin edge}=[<-,shorten <=1pt]
		\tikzstyle{input neuron}=[very thick, circle,draw=red, fill=red!30, minimum size=15pt,inner sep=0pt]
		\tikzstyle{output neuron}=[very thick, circle, draw=green,fill=green!30,minimum size=20pt,inner sep=0pt]
		\tikzstyle{hidden neuron}=[very thick, circle,draw=blue,fill=blue!30,minimum size=15pt,inner sep=0pt]
		\tikzstyle{annot} = [text width=9em, text centered]
		\tikzstyle{annot2} = [text width=4em, text centered]
		
		\foreach \name / \y in {1,...,5}
		\node[input neuron] (I-\name) at (0,-\y) {$x_{\y}$};
		\foreach \name / \y in {1, ..., 8}
		\path[yshift = 1.5cm]
		node[hidden neuron] (H0-\name) at (\layersep, -\y cm) {};
		\foreach \name / \y in {1,...,6}
		\path[yshift=0.5cm]
		node[hidden neuron] (H1-\name) at (2*\layersep,-\y cm) {};      
		\foreach \name / \y in {1,...,7}
		\path[yshift=1cm]
		node[hidden neuron] (H2-\name) at (3*\layersep,-\y cm) {};      
		\foreach \name / \y in {1,2,3}
		\path[yshift=1cm]
		node[output neuron](O-\name) at (4*\layersep,-2*\y cm) {$\cN^{4,\theta}_{\infty,\y}(x)$};   

		\foreach \source in {1,...,5}
		\foreach \dest in {1,...,8}
		\path[line width = 0.8 ] (I-\source) edge (H0-\dest);     
		\foreach \source in {1, ..., 8}
		\foreach \dest in {1, ..., 6 }
		\path[line width = 0.8 ] (H0-\source) edge (H1-\dest);
		\foreach \source in {1,...,6}
		\foreach \dest in {1,...,7}
		\path[line width = 0.8 ] (H1-\source) edge (H2-\dest);     
		\foreach \source in {1,...,7}
		\foreach \dest in {1,...,3}
		\path[line width = 0.8] (H2-\source) edge (O-\dest);
		
		% Annotate the layers
		\node[annot,above of=H0-1, node distance=1cm, align=center] (hl) {1st hidden layer\\(2nd layer)};
		\node[annot,above of=H1-1, node distance=2cm, align=center] (hl2) {2nd hidden layer\\(3rd layer)};
		\node[annot,above of=H2-1, node distance=1.5cm, align=center] (hl3) {3rd hidden layer\\(4th layer)};
		\node[annot,left of=hl, align=center] {Input layer\\ (1st layer)};
		\node[annot,right of=hl3, align=center] {Output layer\\(5th layer)};
		
		\node[annot2,below of=H0-8, node distance=1cm, align=center] (sl) {$\ell_1=8$};
		\node[annot2,below of=H1-6, node distance=2cm, align=center] (sl2) {$\ell_2=6$};
		\node[annot2,below of=H2-7, node distance=1.5cm, align=center] (sl3) {$\ell_3=7$};
		\node[annot2,left of=sl, align=center] {$\ell_0=5$};
		\node[annot2,right of=sl3, align=center] {$\ell_4=3$};
		\end{tikzpicture}
	\end{adjustbox}
	\caption{Graphical illustration of the considered deep ANN architecture in \cref{thm:main_thm} in the special case of a deep ANN with 3 hidden layers (corresponding to $L = 4$ affine linear transformations), with 5 neurons on the input layer (corresponding to $\ell_0 = 5$),
		8 neurons on the 1st hidden layer
		(corresponding to $\ell_1 = 8 $),
		 6 neurons on the 2nd hidden layer (corresponding to $\ell_2 = 6$), 7 neurons on the 3rd hidden layer (corresponding to $\ell_3 = 7$), and 3 neurons on the output layer (corresponding to $\ell_4 = 3$). 
	In this situation the dimension $ \fd $ of the ANN parameter space satisfies 
	$ \fd = \sum_{ i = 1 }^4 \ell_i ( \ell_{i-1}+1 ) = 6 \cdot 8 + 9 \cdot 6 + 7 \cdot 7 + 3 \cdot 8 = 176$.
		Note that we have for every ANN parameter vector $\theta \in \R^{ \mathfrak{d} } = \R^{176 }$ that the realization function 
		$\R^5 \ni x \mapsto \cN^{ 4, \theta }_{ \infty }(x) \in \R^3$ 
		of the considered deep ANN maps the 5-dimensional input vector 
		$ x = ( x_1, x_2, x_3, x_4, x_5 ) \in [a,b]^5$ to the 3-dimensional 
		output vector $\cN^{ 4, \theta }_{ \infty }( x ) = 
		( \cN^{4, \theta}_{\infty, 1}, \cN^{4, \theta}_{\infty, 2}, \cN^{4, \theta}_{\infty, 3} ) \in \R^3$.}
	\label{figure_deep}
\end{figure}

\begin{samepage}
\begin{theorem}
\label{thm:main_thm}
Let $ L, \fd, \fq \in \N $, $ ( \ell_k )_{ k \in \N_0 } \subseteq \N $, 
$ a \in \R $, $ b \in [a,\infty) $ 
satisfy $ \mathfrak{d} = \sum_{ k = 1 }^L \ell_k (\ell_{ k - 1 } + 1 ) $, 
for every $ i \in \cu{ 1, \ldots, \fq } $
let $ \alpha_i \in \R^{ \fq \times d } $,
let $ \beta_i \in \R^{ \fq } $, 
and
let $ P_i \colon \R^{ \ell_0 } \to \R^{ \ell_L + 1 } $ be a polynomial, 
let 
$ f = ( f_1, \dots, f_{ \ell_L } ) \colon [a,b]^{ \ell_0 } \to \R^{ \ell_L } $ 
and 
$ \dens \colon [ a, b ]^{ \ell_0 } \to [0, \infty) $
satisfy for all 
$ x \in [ a, b ]^{ \ell_0 } $ that
\begin{equation} 
\label{theo:intro:conv:eq0_GF}
  (
    f_1(x), f_2(x), \dots, f_{ \ell_L }( x ), \dens(x)
  )
  = 
  \smallsum_{ i = 1 }^{ \fq } 
%   \br*{ 
    P_i( x ) \indicator{ [0, \infty)^{ \fq } } 
    \rbr{ \alpha_i x + \beta_i } 
%   } 
  ,
\end{equation} 
for every 
$ \theta = ( \theta_1, \dots, \theta_{ \fd } ) \in \R^{ \fd } $ 
let 
$ 
  \fw^{ k, \theta } = 
  ( \fw^{ k, \theta }_{ i, j } )_{ 
    (i,j) \in \{ 1, \ldots, \ell_k \} \times \{ 1, \ldots, \ell_{ k - 1 } \}
  }
  \in \R^{ \ell_k \times \ell_{ k - 1 } } 
$, 
$ k \in \N $, 
and 
$ 
  \fb^{ k, \theta } = ( \fb^{ k, \theta }_1, \dots, \fb^{ k, \theta }_{ \ell_k } )
  \in \R^{ \ell_k }
$,
$ k \in \N $, 
satisfy for all  
$ k \in \{ 1, \ldots, L \} $, $ i \in \{1, \ldots, \ell_k \} $, 
$ j \in \{ 1, \ldots, \ell_{ k - 1 } \} $ 
that
\begin{equation}
\label{wb_thm1_GF}
  \fw^{ k, \theta }_{ i, j } 
  =
  \theta_{ 
    (i-1) \ell_{k-1} + j 
    + \sum_{h=1}^{k-1} \ell_h (\ell_{h-1} + 1)
  }
  \qquad \text{and} \qquad
  \fb^{ k, \theta }_i 
  =
  \theta_{ 
    \ell_k \ell_{ k - 1 } + i + \sum_{ h = 1 }^{ k - 1 } \ell_h (\ell_{ h - 1 } + 1 ) 
  } ,
\end{equation}
let 
$ \Rect_r \colon \R \to \R $, 
$ r \in [1, \infty] $,
satisfy for all 
$ r \in [1, \infty) $, 
$ x \in ( - \infty, 2^{ - 1 } r^{ - 1 } ] $, 
$ 
  y \in \R
$, 
$
  z \in [ r^{ - 1 }, \infty ) 
$
that 
\begin{equation}
\label{lim_R_thm1}
  \Rect_r \in C^1( \R, \R ) ,
\quad
  \Rect_r(x) = 0,
\quad 
  0 \leq \Rect_r(y) \leq \Rect_{ \infty }( y ) 
  = 
  \max\{ y, 0 \}
  ,
\qandq
  \Rect_r(z) = z
  ,
\end{equation}
assume 
$
  \sup_{ r \in [1, \infty) }
  \sup_{ x \in \R } 
  | ( \Rect_r )'( x ) | < \infty 
$, 
for every $ r \in [1, \infty] $ 
let 
$ 
  \mathfrak{M}_r \colon ( \cup_{ n \in \N } \R^n )
  \to ( \cup_{ n \in \N } \R^n ) 
$ 
satisfy for all 
$ n \in \N $,
$ x = (x_1, \ldots, x_n ) \in \R^n $
that 
$ 
  \mathfrak{M}_r( x ) = ( \Rect_r( x_1 ), \ldots, \Rect_r( x_n ) ) 
$, 
for every $ r \in [1,\infty] $, 
$ \theta \in \R^{ \fd } $
let 
$ 
  \mathcal{N}^{ k, \theta }_r \colon \R^{ \ell_0 } \to \R^{ \ell_k } 
$, 
$ k \in \N $, 
satisfy for all $ k \in \N $,
$ x \in \R^{ \ell_0 } $
that 
\begin{equation}
\label{def_N_thm1_GF}
  \mathcal{N}^{1,\theta}_r(x) 
  =
  \fb^{ 1, \theta } 
  +
  \fw^{ 1, \theta } x 
  \qquad
  \text{and}
  \qquad
  \mathcal{N}^{ k+1, \theta }_r( x ) 
  =
  \fb^{ k + 1, \theta } 
  +
  \fw^{ k + 1, \theta }(
    \mathfrak{M}_{ r^{ 1 / k } }( \mathcal{N}^{ k, \theta }_r( x ) ) 
  )
  ,
\end{equation}
for every $ r \in [1, \infty] $ 
let $ \cL_r \colon \R^\fd \to \R $
satisfy for all 
$ 
  \theta \in \R^{ \fd } 
$ 
that 
$ 
  \cL_r( \theta ) 
  =
  \int_{ [a,b]^{ \ell_0 } } 
    \norm{ \mathcal{N}^{ L, \theta }_r( x ) - f(x) }^2 
    \,
    \dens(x)
    \, 
  \d x
$,
let $ \cG \colon \R^\fd \to \R^\fd$ satisfy for all
$
  \theta \in 
  \cu{ 
    \vartheta \in \R^\fd \colon ( ( \nabla \cL_r ) ( \vartheta ) )_{ r \in [1,\infty) } 
    \ \allowbreak\text{is}\ \allowbreak\text{convergent} 
  }
$
that 
$
  \cG( \theta ) = 
  \lim_{ r \to \infty } \allowbreak ( \nabla \cL_r )( \theta ) 
$,
and\footnote{Throughout this article we denote by 
$ \norm{ \cdot } \colon ( \cup_{ n \in \N } \R^n ) \to \R $ 
and 
$ 
  \spro{ \cdot , \cdot } \colon( \cup_{ n \in \N }( \R^n \times \R^n ) ) \to \R 
$ 
the functions which
satisfy for all $ n \in \N $,
$ x = ( x_1, \dots, x_n ) $, $ y = ( y_1, \dots, y_n ) \in \R^n $
that $ \norm{x} = ( \sum_{ i = 1 }^n | x_i |^2 )^{ 1 / 2 } $
and $ \spro{ x, y } = \sum_{ i = 1 }^n x_i y_i $.} 
let 
$ \Theta \in C( [0 , \infty) , \R^\fd ) $ satisfy 
$ \liminf_{ t \to \infty } \norm{ \Theta_t } < \infty $ 
and 
$ 
  \forall \, t \in [0, \infty ) \colon 
  \Theta_t = \Theta_0 - \int_0^t \cG( \Theta_s ) \, \d s 
$. 
Then\footnote{In the conclusion of \cref{thm:main_thm} 
we denote by 
$
  ( \bD \cL_{ \infty } )( \vartheta )
$
the limiting Fr\'{e}chet subdifferential of $ \cL_{ \infty } \colon \R^{ \fd } \to \R $ 
at $ \vartheta \in \R^{ \fd } $; see \cref{def:limit:subdiff} below for details.} there exist 
$ \vartheta \in \R^{ \fd } $, 
$ \fC, \beta \in (0, \infty) $ 
with 
$
  0 \in ( \bD \cL_{ \infty } )( \vartheta )
$
such that for all $ t \in [0, \infty) $ 
it holds that 
\begin{equation}
\label{theo:intro:conv:eq3}
  \norm{ \Theta_t - \vartheta } 
  \leq \fC ( 1 + t )^{ - \beta } 
\qqandqq
  \abs{ \cL_\infty ( \Theta_t ) - \cL_\infty ( \vartheta ) } \leq \fC ( 1 + t ) ^{-1} 
  .
\end{equation}
\end{theorem}
\end{samepage}

\cref{thm:main_thm} is an immediate 
consequence of \cref{theo:gf:conv:simple} 
in \cref{sec:convergence_GF_processes} below. 
Note that the first inequality in 
\cref{theo:intro:conv:eq3} in \cref{thm:main_thm} above 
assures that the standard norm 
$
  \| \Theta_t - \vartheta \|
$
of the difference of the GF solution at time $ t $ 
and the generalized critical point 
$ \vartheta $ 
converges with rate $ \beta \in (0,\infty) $ to $ 0 $ 
and note that 
the second inequality in 
\cref{theo:intro:conv:eq3} in \cref{thm:main_thm} above 
assures that the absolute value  
$
  | \cL_{ \infty }( \Theta_t ) - \cL_{ \infty }( \vartheta ) |
$
of the difference of the risks 
of the GF solution at time $ t $ 
and the generalized critical point 
$ \vartheta $ 
converges with rate $ 1 $ to $ 0 $. 
 
In our proof of \cref{thm:main_thm} we combine the generalized KL-inequality for the risk function in \cref{prop:loss:lojasiewicz} with the abstract convergence results for GF processes in \cref{section:gf:loja}. The main regularity condition we need is the chain rule for the risk function $\cL_\infty $, which was established in \cite{DNNReLUarXiv}.
The fact that the limit
$ \vartheta \in \R^{ \fd } $ is a generalized critical point in the sense that 0 is an element of the limiting Fr\'{e}chet subdifferential of $ \cL_{ \infty } \colon \R^{ \fd } \to \R $ 
at $ \vartheta  $ is a consequence of the fact that the generalized gradient we define is an element of the limiting Fr\'{e}chet subdifferential, which we show in \cref{prop:loss:gradient:subdiff} below.

The assumption that the trajectory $(\Theta_t)_{t \in [0 , \infty ) }$ is bounded is necessary and is not implied by the other conditions. Indeed, in \cite{GallonJentzen2022} we show that there are piecewise polynomial target functions for which GF trajectories with certain initialization do diverge to infinity.

In \cite[Theorem 1.2]{EberleJentzenRiekertWeiss2021} 
a GF convergence result related to \cref{thm:main_thm} above 
has been obtained in the case of shallow ANNs with just one hidden layer. 
More specifically, 
in \cite[Theorem 1.2]{EberleJentzenRiekertWeiss2021} 
a GF convergence result similar to \cref{theo:intro:conv:eq3} 
has been established in the situation 
where the target function is additionally continuous 
and where the considered ANNs are not deep but shallow and just consist 
of 3 layers (input layer, output layer, and one hidden layer).

\subsection{Gradient descent (GD) with random initializations 
in the training of deep ANNs}
\label{susbsec:GD_deep_intro}

In \cref{theo:intro:convergence} below 
we establish in the training of deep ReLU ANNs with an arbitrarily large 
number of hidden layers, 
under the assumption that 
the unnormalized probability density function 
$
  \dens \colon [a,b]^d \to [0,\infty)
$
and the target function 
$ f \colon [a,b]^d \to \R^{ \delta } $
are piecewise polynomial 
(see \cref{theo:intro:conv:eq0_GD} below for details)
and that the risk function of the 
considered deep supervised learning problem 
admits at least one regular global minimum point, 
that the risk of the plain vanilla GD optimization method with 
random initializations converges in probability 
to $ 0 $ as the number of GD steps increases to $ \infty $, 
as the number of random initializations increases to $ \infty $, 
as the step size of the GD method (the learning rate of the GD method) decreases to $ 0 $, 
and as the width of the ANNs increases to $ \infty $ 
(see \cref{eq:width_diverges_to_infinity} and \cref{theo:intro:random:init:eqresult_deep} 
below for details).

\begin{samepage}
\begin{theorem} 
\label{theo:intro:convergence}
Let 
$ d, \delta, \fq \in \N $, 
$ a \in \R $, $ b \in [a,\infty) $, 
$ ( \rho_{ \fa } )_{ \fa \in \N } \subseteq ( \N \cap (1,\infty) ) $,
let 
$ 
  \ell^{ \fa } = ( \ell^{ \fa }_0, \ell^{ \fa }_1, \dots, \ell^{ \fa }_{ \rho_{ \fa } } ) \in 
  \{ d \} \times \N^{ \rho_{ \fa } - 1 } \times \{ \delta \}
$, 
$ \fa \in \N $, 
satisfy 
\begin{equation}
\label{eq:width_diverges_to_infinity}
\textstyle
  \liminf_{ \fa \to \infty }
  \min \cu{\ell^\fa_1, \ell^\fa_2, \ldots , \ell^{\fa}_{\rho_\fa - 1 } }
  = \infty
  ,
\end{equation}
for every $ \fa \in \N $
let 
$
  \fd_{ \fa } = 
  \sum_{ k = 1 }^{ \rho_{ \fa } } 
  \ell^{ \fa }_k ( \ell^{ \fa }_{ k - 1 } + 1 )
$, 
for every $ i \in \cu{ 1, \ldots, \fq } $
let $ \alpha_i \in \R^{ \fq \times d } $,
let $ \beta_i \in \R^{ \fq } $, 
and
let $ P_i \colon \R^d \to \R^{ \delta + 1 } $ be a polynomial, 
let 
$ f \colon [a,b]^d \to \R^{ \delta } $ 
and 
$ \dens \colon [ a, b ]^d \to [0, \infty) $
satisfy for all 
$ x \in [ a, b ]^d $ that
\begin{equation}
\label{theo:intro:conv:eq0_GD}
  ( f_1(x), f_2(x), \dots, f_{ \delta }(x), \dens(x) ) 
%   k f( x ) + ( 1 - k ) \dens( x ) 
  = 
  \smallsum_{ i = 1 }^{ \fq } 
%   \br*{ 
    P_i( x ) 
    \indicator{ [0, \infty)^{ \fq } } 
    \rbr{ \alpha_i x + \beta_i } 
%   } 
  ,
\end{equation}
for every $ \fa \in \N $, 
$ k \in \{ 1, \dots, \rho_{ \fa } \} $, 
$ 
  \theta = ( \theta_1, \dots, \theta_{ \fd_{ \fa } } ) \in \R^{ \fd_{ \fa } } 
$
let 
$ 
  \fw^{ k, \theta }_{ \fa }
  = 
  ( 
    \fw^{ k, \theta }_{ \fa, i, j } 
  )_{ 
    (i,j) \in 
    \{ 1, \ldots, \ell_k^{ \fa } \} \times 
    \{ 1, \ldots, \ell_{ k - 1 }^{ \fa } \}
  }
  \allowbreak 
  \in \R^{ \ell_k^{ \fa } \times \ell_{ k - 1 }^{ \fa } } 
$
and 
$ 
  \fb^{ k, \theta }_{ \fa } 
  = 
  ( 
    \fb^{ k, \theta }_{ \fa, 1 }, \dots, \fb^{ k, \theta }_{ \fa, \ell_k^{ \fa } } 
  )
  \in \R^{ \ell_k^{ \fa } }
$
satisfy for all 
$ i \in \{ 1, \ldots, \ell_k^{ \fa } \} $, 
$ j \in \{ 1, \ldots, \ell_{ k - 1 }^{ \fa } \} $ 
that
\begin{equation}
\label{wb_thm1_GD}
  \fw^{ k, \theta }_{ \fa, i, j } 
  =
  \theta_{ 
    (i-1) \ell^\fa_{k-1} + j 
    + \sum_{h=1}^{k-1} \ell^\fa_h (\ell^\fa_{h-1} + 1)
  }
\qqandqq
  \fb^{ k, \theta }_{\fa , i} 
  =
  \theta_{ 
    \ell^\fa_k \ell^\fa_{ k - 1 } + i + \sum_{ h = 1 }^{ k - 1 } \ell^\fa_h (\ell^\fa_{ h - 1 } + 1 ) 
  } ,
\end{equation}
let 
$ 
  \mathfrak{M} \colon ( \cup_{ n \in \N } \R^n )
  \to ( \cup_{ n \in \N } \R^n ) 
$
satisfy for all 
$ n \in \N $,
$ x = (x_1, \ldots, x_n ) \in \R^n $
that 
$ 
  \mathfrak{M}( x ) = ( \max\{ x_1, 0 \}, \ldots, \max\{ x_n, 0 \} ) 
$, 
for every 
$ \fa \in \N $, 
$ \theta \in \R^{ \fd_{ \fa } } $ 
let 
$ 
  \mathcal{N}^{ k, \theta }_{ \fa } 
  \colon \R^d \to \R^{ \ell^{ \fa }_k } 
$, 
$ k \in \N \cap [1, \rho_{ \fa } ]
$, 
satisfy for all 
$ k \in \N \cap [1, \rho_{ \fa } ) $,
$x \in \R^d$
that 
\begin{equation}
\label{def_N_thm1_GD}
  \mathcal{N}^{ 1, \theta }_{ \fa }( x ) 
  =
  \fb^{ 1, \theta }_{ \fa } 
  +
  \fw^{ 1, \theta }_{ \fa } x 
  \qquad
  \text{and}
  \qquad
  \mathcal{N}^{ k+1, \theta }_{ \fa }( x ) 
  =
  \fb^{ k + 1, \theta }_{ \fa } 
  +
  \fw^{ k + 1, \theta }_{ \fa }\big(
    \mathfrak{M}( \mathcal{N}^{ k, \theta }_{ \fa }( x ) ) 
  \big)
  ,
\end{equation} 
for every $ \fa \in \N $
let 
$ \cL_{ \fa } \colon \R^{ \fd_{ \fa } } \to \R $
satisfy for all
$ 
  \theta \in \R^{ \fd_{ \fa } } 
$ 
that 
$ 
  \cL_{ \fa }( \theta ) 
  =
  \int_{ [a,b]^{ d } } 
    \norm{ \mathcal{N}_{ \fa }^{ \rho_\fa, \theta }( x ) - f(x) }^2 
    \, 
    \dens( x )
    \,
  \d x
$, 
for every $ \fa \in \N $ 
let 
$ 
  \vartheta_{ \fa } \in 
  ( \cL_{ \fa } )^{ - 1 }( 
    \{ 
      \inf_{ \theta \in \R^{ \fd_{ \fa } } } 
      \cL_{ \fa }( \theta )
    \} 
  )
$, 
$ \varepsilon_{ \fa } \in (0,1) $ 
satisfy that
$
  \cL_{ \fa }|_{ 
    \{ 
      \theta \in \R^{ \fd_{ \fa } } \colon 
      \norm{ \theta - \vartheta_{ \fa } } 
      < \varepsilon_{ \fa } 
    \}
  } 
$
has a Lipschitz continuous derivative, 
for every $ \fa \in \N $ 
let 
$ \cG_{ \fa } \colon \R^{ \fd_{ \fa } } \to \R^{ \fd_{ \fa } } $
satisfy for all 
$ 
  \theta \in \cu{
\vartheta \in \R^{\fd_\fa} 
\colon 
\cL_{ \fa }
\ \text{is}\ \allowbreak \text{differentiable}\ \allowbreak \text{at}\ \vartheta
}
$
that 
$ 
  \cG_{ \fa }( \theta ) = ( \nabla \cL_{ \fa } )( \theta ) 
$,
let $ ( \Omega, \cF, \P ) $ be a probability space, 
for every 
$ n, \fa, K \in \N_0 $, $ \gamma \in \R $
let 
$ 
  \Theta^{ K, \gamma }_{ \fa, n } \colon \Omega \to  \R^{ \fd_{ \fa } } 
$
and 
$ 
  \bfk^{ K, \gamma }_{ \fa, n } \colon \Omega \to \N 
$ 
be random variables,
assume for all $ \fa \in \N $, $ \gamma \in \R $ 
that $ \Theta_{ \fa, 0 }^{ K, \gamma } $, $ K \in \N $, are i.i.d., 
assume for all $ \fa \in \N $, $ \gamma, r \in (0,1) $, 
$ \theta \in \R^{ \fd_{ \fa } } $ that
$
  \P( 
    \norm{ \Theta^{ 1, \gamma }_{ \fa, 0 } - \theta } < r
  )
  > 0
$, 
and assume for all 
$ n \in \N_0 $, $ \fa, K \in \N $, $ \gamma \in \R $, $ \omega \in \Omega $ 
that
\begin{equation} 
\label{theo:intro:deepGD:random:init:eq3}
    \Theta_{ \fa, n + 1 }^{ K, \gamma }( \omega ) 
    = \Theta_{ \fa, n }^{ K, \gamma }( \omega ) 
    - \gamma \cG_{ \fa }( \Theta_{ \fa, n }^{ K, \gamma }( \omega ) ) 
\quad\text{and}\quad
    \bfk^{ K, \gamma }_{ \fa, n }( \omega) 
    \in 
    \argmin\nolimits_{ \kappa \in \{ 1, \ldots, K \} } 
    \cL_{ \fa }( \Theta_{ \fa, n }^{ \kappa, \gamma }( \omega ) ) .
\end{equation}
Then 
\begin{equation}  
\label{theo:intro:random:init:eqresult_deep}
  \limsup\nolimits_{ \fa \to \infty } 
  \limsup\nolimits_{\gamma \searrow 0 }
  \limsup\nolimits_{ K \to \infty } 
  \limsup\nolimits_{n \to \infty}
  \E\bigl[
    \min\bigl\{
    \cL_{ \fa }\bigl( 
      \Theta^{ 
        \bfk^{ K, \gamma 
        }_{ \fa, n }, 
        \gamma 
      }_{ \fa, n } 
    \bigr)
    , 1  
  \bigr]
  = 0 .
\end{equation}
\end{theorem}
\end{samepage}

\cref{theo:intro:convergence} follows immediately from 
\cref{cor:gd:random:dnn:item2} in \cref{theo:gd:random:dnn} 
in \cref{sec:GD_random_init} below 
and the reversed version of Fatou's lemma. 
Observe that 
\cref{theo:intro:random:init:eqresult_deep} above 
shows that 
the risk 
$
  \cL_{ \fa }(  
    \Theta^{ 
      \bfk^{ K, \gamma 
      }_{ \fa, n }, 
      \gamma 
    }_{ \fa, n } 
  )
$
of the GD optimization method with random initializations
converges in probability (see \cref{eq:convergence_in_probability} above) to $ 0 $ 
as the number $ n $ of GD steps increases to $ \infty $, 
as the number $ K $ of random initializations increases to $ \infty $, 
as the learning rate $ \gamma $ decreases to $ 0 $,
and 
as the width of the ANN 
increases to $ \infty $ 
in the sense of \cref{eq:width_diverges_to_infinity} above. 

The proof of \cref{theo:intro:convergence}
is mostly analogous to the proof of \cref{theo:intro:random:init}.
The main difference is that in the general setting of deep ANNs the existence of global minima is not known.
This is the reason why we assume in \cref{theo:intro:convergence} for every $\fa \in \N$ that the parameter vector $\vartheta_\fa \in \R^{ \fd_\fa}$ is a global minimum of the risk function $\cL_\fa \colon \R^{ \fd_\fa } \to \R$. Additionally, we assume that for every $\fa \in \N$ there exists $\varepsilon_\fa \in (0 , 1 )$ which satisfies that the restriction of $\cL_\fa$ to the neighborhood $\cu{ \theta \in \R^{ \fd_\fa } \colon \norm{ \theta - \vartheta _ \fa } < \varepsilon_\fa } $ is differentiable with a Lipschitz continuous derivative.

\section{Existence of global minima for shallow ANNs}
\label{sec:existence_global_minima}

In this section we establish in 
\cref{theo:existence2} 
in \cref{subsec:existence_global_minima}
below 
in the case where the target function 
$ f \colon [a,b] \to \R $ is Lipschitz continuous 
and where the considered ReLU ANNs 
just consist of a one-dimensional input layer, 
a multi-dimensional hidden layer, and a one-dimensional output layer
that there exists a global minimum point of the risk function. 
\cref{theo:existence_global_minima} in the introduction 
is a direct consequence of \cref{theo:existence2}.

In \cref{cor:existence:regular}
in \cref{subsec:existence_regular_minima} 
we strengthen  
\cref{theo:existence2} 
by showing that there also exists a global minimum point of the risk function such that 
the risk function is continuously differentiable on a neighborhood around 
the global minimum point. Our proof of \cref{cor:existence:regular} 
is based on an application of \cref{theo:existence2} 
as well as on applications of some basic regularity results from our earlier 
article Eberle et al.~\cite[Proposition~2.3 and Corollary~2.7]{EberleJentzenRiekertWeiss2021}.

Our proof of \cref{theo:existence2} can, roughly speaking, be divided 
into three parts. 
\begin{enumerate} [label = (\Roman*)]
	\item In \cref{cor:characterization2} in \cref{ssec:charac_shallow_ANNs} below 
	we establish an 
	explicit characterization for the functions $ f \colon [0,1] \to \R $ 
	which can be exactly represented by a shallow ReLU ANN 
	with $ \width \in \N $ neurons on the hidden layer. 
	\item Thereafter, we employ \cref{cor:characterization2} 
	to prove in \cref{cor:better:approx} 
	in \cref{ssec:structure_preserving_ANNs} below 
	in the case where the target function 
	$ f \colon [a,b] \to \R $ is Lipschitz continuous 
	with Lipschitz constant $ L \in \R $
	and where the considered ReLU ANNs 
	consist of a one-dimensional input layer, 
	an $ \width $-dimensional hidden layer, and 
	a one-dimensional output layer 
	that, roughly speaking, for every ANN parameter vector 
	$ \theta \in \R^{ 3 \width + 1 } $ 
	there exists an ANN parameter vector 
	$ \vartheta \in \R^{ 3 \width + 1 } $
	whose realization function approximates $ f $
	at least as well as the realization function of $ \theta $ 
	but is additionally also Lipschitz continuous 
	with Lipschitz constant at most $ \width L $. 
	\item Finally, 
	we combine \cref{cor:better:approx} 
	with the Arzel\`{a}--Ascoli theorem 
	and the fact that the set of realization functions 
	of shallow ReLU ANNs with fixed architecture forms a closed 
	subset of the set of continuous functions revealed in 
	Petersen et al.~\cite[Theorem 3.8]{PetersenRaslanVoigtlaender2020} 
	to prove \cref{theo:existence2}. 
\end{enumerate}

The question which functions $ f \colon \R^d \to \R $ can be represented by 
a shallow ReLU ANN with a fixed number of neurons on the hidden layer 
has also been investigated in the article Dereich \& Kassing~\cite{DereichKassing2021minimal}
and in Theorem~3.2 in \cite{DereichKassing2021minimal} 
a similar result as \cref{cor:characterization2} 
has been established.

Our proofs of \cref{cor:characterization2} 
and \cref{cor:better:approx} 
also use the elementary 
results and notions regarding piecewise linear functions in 
\cref{ssec:prop_breakpoint} 
%\cref{def:number:of:kinks}, 
%% in \cref{ssec:prop_breakpoint}, 
%\cref{def:piece:linear}, 
%% in \cref{ssec:prop_breakpoint}, 
%\cref{def:slopes}, 
%% in \cref{ssec:prop_breakpoint}, 
%\cref{prop:piece:linear:relations},  
%% in \cref{ssec:prop_breakpoint}, 
%\cref{lem:kinks:triangle:ineq},
%% in \cref{ssec:prop_breakpoint}, 
%\cref{lem:sum:piece:linear},
%% in \cref{ssec:prop_breakpoint}, 
%\cref{def:lip:const},
%% in \cref{ssec:prop_breakpoint}, 
%\cref{prop:piece:linear:lip},
%% in \cref{ssec:prop_breakpoint}, 
as well as the elementary
\cref{lem:realization:l}
% in \cref{ssec:charac_shallow_ANNs}, 
and 
\cref{lem:h-1:kinks},
% in \cref{ssec:charac_shallow_ANNs} 
and only for completeness we include in this section also 
detailed proofs for these results.
%\cref{prop:piece:linear:relations}, 
%\cref{lem:kinks:triangle:ineq}, 
%\cref{lem:sum:piece:linear}, 
%\cref{prop:piece:linear:lip}, 
%and 
%\cref{lem:realization:l}.

\subsection{Mathematical framework for shallow ANNs with ReLU activation}

In \cref{setting:snn} we present our framework for shallow ANNs with ReLU activation which will be employed during the remainder of this section. 

\begin{setting} 
\label{setting:snn}
Let $\width, \fd \in \N$, $L \in \R $, $ f \in C([0,1] , \R)$ satisfy for all $x,y \in [0,1]$ that $\abs{ f(x) - f(y) } \leq L \abs{x-y}$ and $\fd = 3 \width + 1$,
let $\fw  = (( \w{\theta} _ {j}  )_{j \in \{1, \ldots, \width \} })_{ \theta \in \R^{\fd}} \colon \R^{\fd} \to \R^{\width}$,
$\fb =  (( \b{\theta} _ j )_{j \in \{1, \ldots, \width \} })_{ \theta \in \R^{\fd}} \colon \R^{\fd} \to \R^{\width}$,
$\fv = (( \v{\theta} _ j )_{j \in \{1, \ldots, \width \} })_{ \theta \in \R^{\fd}} \colon \R^{\fd} \to \R^{\width}$, 
$\fc = (\c{\theta})_{\theta \in \R^{\fd }} \colon \R^{\fd} \to \R$,
and
$\fq = ( ( \q{\theta}_j)_{j \in \{1, \ldots, \width \} } ) \colon \R^\fd \to (- \infty , \infty] ^\width$
satisfy for all $\theta  = ( \theta_1 ,  \ldots, \theta_{\fd}) \in \R^{\fd}$, $j \in \{1, 2, \ldots, \width \}$ that
	 $\w{\theta}_{ j} = \theta_{ j}$,
	 $\b{\theta}_j = \theta_{\width + j}$, 
	$\v{\theta}_j = \theta_{2 \width + j}$, 
	$\c{\theta} = \theta_{\fd}$,
	and
	\begin{equation}
	\q{\theta}_j = \begin{cases}
	- \nicefrac{\b{\theta}_j}{\w{\theta}_j}, \qquad & \w{\theta}_j \not= 0 \\
	\infty, \qquad & \w{\theta}_j = 0,
	\end{cases}
	\end{equation}
	let $\mu \colon \cB ( [0,1] ) \to [0, \infty]$ be a finite measure,
	and let $\scrN = (\realization{\theta})_{\theta \in \R^{\fd } } \colon \R^{\fd } \to C([0,1] , \R)$ and $ \cL \colon \R^{\fd  } \to \R$
	satisfy for all $\theta \in \R^{\fd}$, $x \in [0,1]$ that
\begin{equation}
\label{eq:def_realization_function}
  \realization{ \theta }(x) 
  = \c{\theta} 
  + \smallsum_{ j = 1 }^{ \width } 
  \v{\theta}_j \max\cu[\big]{ \b{\theta}_j + \w{\theta}_j x , 0 }
\end{equation}
and 
$ 
  \cL( \theta ) = \int_0^1 ( \realization{\theta}( y ) - f( y ) )^2 \, \mu( \d y ) 
$.
\end{setting}

\subsection{Properties of the breakpoint function}
\label{ssec:prop_breakpoint}

\cfclear
\begin{definition}[Breakpoint function]
\label{def:number:of:kinks}
We denote by $ Q \colon C( [0,1], \R ) \to [0, \infty] $ 
the function which satisfies for all $ f \in C( [0,1], \R )$ that
\begin{multline}
\label{eq:def_Q_function}
  Q(f) = 
  \min\bigl( 
    \cu{\infty} 
    \cup 
    \bigl\{ 
      n \in \N_0 \colon \bigl[ \exists \, \fA_1, \fA_2, \ldots, \fA_{n+1}, \fB_1, \fB_2, \ldots, \fB_{n+1},  
      \fq_0, \fq_1, \ldots, \fq_{n+1} \in \R \colon \\
      ( 
        [ 0 = \fq_0 < \fq_1 < \cdots < \fq_{n+1} = 1], \, 
        [ 
          \forall \, j \in \N \cap [1,n+1], \, x \in [\fq_{j-1}, \fq_j] \colon 
%           \\
          f(x) = \fA_j x + \fB_j 
        ] 
      ) 
    \bigr] 
    \bigr\} 
  \bigr)
  .
\end{multline}
\end{definition}

\cfclear
\begin{definition}[Piecewise affine linear functions] 
\label{def:piece:linear} 
\cfadd{def:number:of:kinks}
We denote by $ \scrL \subseteq C([0,1], \R) $ the set given by 
\begin{equation}
  \scrL = \{ f \in C( [0,1], \R ) \colon Q(f) < \infty \} 
\end{equation}
\cfload.
\end{definition}

\cfclear
\begin{definition}[Slopes and axis intercepts for piecewise affine linear functions] 
\label{def:slopes} \cfadd{def:piece:linear}
Let $ f \in \scrL $ \cfload. 
Then we denote by 
$ 
  A_1(f), A_2(f), \ldots, A_{Q(f)+1}(f), 
  \allowbreak B_1 (f), B_2 (f), \ldots, B_{Q(f) + 1}(f), 
  \allowbreak 
  q_0(f), q_1(f), \ldots, q_{ Q(f) + 1 }(f) \in \R 
$ 
the real numbers which satisfy 
$ 0 = q_0(f) < q_1(f) < \cdots < q_{ Q(f) + 1 }(f) = 1 $ 
and 
\begin{equation} \label{eq:def:slopes} \cfadd{def:number:of:kinks} 
  \forall \, j \in \N \cap [1,Q(f)+1] 
%   \{ 1, 2, \ldots, Q(f) + 1 \}
  , \, x \in[q_{j-1} ( f ) , q_j ( f ) ] \colon 
  f(x) = A_j (f) x + B_j ( f ) 
\end{equation}
\cfout.
\end{definition}

\cfclear
\begin{prop} 
\label{prop:piece:linear:relations} 
\cfadd{def:piece:linear,def:number:of:kinks}
Let $ f \in \scrL $, $ i \in \{ 1, 2, \dots, Q(f) \} $ \cfload.\cfadd{def:slopes} 
Then 
\begin{enumerate}[label = (\roman*)]
\item 
it holds that 
$ A_{ i + 1 }( f ) \not= A_i( f ) $, 
\item 
it holds that 
$ B_{ i + 1 }( f ) = B_i( f ) - ( A_{ i + 1 }( f ) - A_i( f ) ) q_i( f ) $, 
and 
\item 
it holds that 
$
  B_{ i + 1 }( f ) = B_1( f ) - \sum_{ j = 1 }^i ( A_{ j + 1 }( f ) - A_j( f ) ) q_j( f )
$. 
\end{enumerate}
\cfadd{def:slopes}\cfadd{def:lip:const}\cfout.
\end{prop}
\begin{cproof}{prop:piece:linear:relations}
\Nobs that \cref{eq:def_Q_function} ensures 
that $ A_i( f ) \not= A_{ i + 1 }( f ) $. 
Next \nobs that the fact that 
for all $ j \in \{ 1, 2, \ldots, Q(f) + 1 \} $, 
$ x \in [ q_{ j - 1 }( f ), q_j( f ) ] $ 
it holds that $ f(x) = A_j( f ) x + B_j( f ) $ 
proves that 
for all $ j \in \{ 1, 2, \ldots, Q(f) \} $ we have that
\begin{equation}
  A_j( f ) q_j( f ) + B_j( f ) = A_{ j + 1 }( f ) q_j( f ) + B_{ j + 1 }( f ).
\end{equation}
Hence, we obtain for all $ j \in \{ 1, 2, \ldots, Q(f) \} $ that 
$ B_{ j + 1 }( f ) = B_j( f ) - ( A_{ j + 1 }( f ) - A_j( f ) ) q_j( f ) $. 
Induction hence establishes that for all 
$ j \in \{ 1, 2, \ldots, Q(f) \} $ it holds that 
$ 
  B_{ j + 1 }( f ) = 
  B_1( f ) - \sum_{ k = 1 }^j ( A_{ k + 1 }( f ) - A_k( f ) ) q_k( f )
$.
\end{cproof}

\cfclear 
\begin{lemma}[Subadditivity of the breakpoint function] 
\label{lem:kinks:triangle:ineq}
Let $ f, g \in C( [0, 1] , \R ) $. 
Then 
\begin{equation}
  Q(f + g) \le Q(f) + Q(g) 
\end{equation} \cfadd{def:number:of:kinks}\cfload.
\end{lemma}
\begin{cproof}{lem:kinks:triangle:ineq}
Throughout this proof assume without loss of generality 
that 
\begin{equation}
\label{eq:obdA_breakpoint_function}
  Q(f) + Q(g) < \infty 
  .
\end{equation}
\Nobs that \cref{eq:obdA_breakpoint_function} 
implies that there exist 
$ N \in \N_0 \cap [ 0, Q(f) + Q(g) ] $, 
$
  \fq_0, \fq_1, \ldots, \fq_{ N + 1 } \in \R
$
which satisfy 
\begin{equation}
\label{eq:q_strictly_increasing}
  0 = \fq_0 < \fq_1 < \cdots < \fq_{ N + 1 } = 1 
\end{equation}
and 
\begin{equation}
\label{eq:union_of_f_and_g}
  \{ 
    \fq_0, \fq_1, \dots, \fq_{ N + 1 }
  \}
  =
  \{
    q_0(f), q_1(f), \dots, q_{ Q(f)+1 }(f)
  \} 
  \cup 
  \{
    q_0(g), q_1(g), \dots, q_{ Q(g)+1 }(g)
  \}
  .
\end{equation}
\cfadd{def:slopes}\cfload. 
\Nobs that 
\cref{eq:q_strictly_increasing}
and 
\cref{eq:union_of_f_and_g} 
ensure that 
for all $ i \in \{ 0, 1, \ldots, N \} $ 
it holds that 
$ 
  (f + g)|_{ [ \fq_i, \fq_{ i + 1 } ] }
$ 
is affine linear. 
Hence, we obtain that 
$
  Q( f + g ) \leq N \leq Q(f) + Q(g)
$.
\end{cproof}

\cfclear
\begin{cor} 
\label{lem:sum:piece:linear}
\cfadd{def:piece:linear}
Let $ f, g \in \scrL $ \cfload. Then
\begin{enumerate}[label = (\roman*)]
\item it holds that $ Q(f + g) \leq Q(f) + Q(g) $ and
\item it holds that $ f + g \in \scrL $
\cfadd{def:number:of:kinks}
\end{enumerate}
\cfout.
\end{cor}
\begin{cproof2}{lem:sum:piece:linear}
\Nobs that 
\cref{lem:kinks:triangle:ineq} and 
the assumption that $ f, g \in \scrL $ 
assure that $ Q(f+g) \le Q(f) + Q( g ) < \infty $.
\end{cproof2}

\begin{definition}[Lipschitz constant] 
\label{def:lip:const} 
We denote by $ \operatorname{Lip} \colon C( [0,1], \R) \to [0, \infty] $ 
the function which satisfies for all $ f \in C([0,1], \R) $ 
that
\begin{equation}
  \operatorname{Lip} ( f ) = 
  \sup_{ 
    \substack{
      x, y \in [0,1], 
    \\
      x \neq y
    }
  } 
  \rbr*{ 
    \frac{ 
      \abs{ f(x) - f(y) } 
    }{ 
      \abs{ x - y }
    }
  }
  .
\end{equation}
\end{definition}

\cfclear 
\begin{lemma}
\label{prop:piece:linear:lip}
\cfadd{def:piece:linear}
Let $ f \in \scrL $ \cfload. 
Then 
\begin{equation}
\label{eq:Lipschitz}
\textstyle
  \operatorname{Lip}( f ) = 
  \max_{ i \in \{ 1, 2, \ldots, Q( f ) + 1 \} } 
  \abs{ A_i( f ) }
\end{equation}
\cfadd{def:number:of:kinks}\cfadd{def:slopes}\cfadd{def:lip:const}\cfout.
\end{lemma}
\begin{cproof}{prop:piece:linear:lip}
\Nobs that the fact that 
$
  0 = q_0(f) < q_1(f) < \dots < q_{ Q(f) + 1 }(f) = 1
$
and the fact that for all 
$ j \in \{ 1, 2, \dots, Q(f) + 1 \} $, 
$ x \in [ q_{ j - 1 }(f), q_j(f) ] $
it holds that 
$
  f(x) = A_j(f) x + B_j(f)
$
ensure that
\begin{equation}
\textstyle
  \sup\nolimits_{ 
    x \in [a,b] \backslash \{ q_0(f), q_1(f), \dots, q_{ Q(f) + 1 }(f) \} 
  }
  \abs{ f'(x) }
  =  
  \max_{ i \in \{ 1, 2, \ldots, Q( f ) + 1 \} } 
  \abs{ A_i( f ) }
  .
\end{equation}
This and the fundamental theorem of calculus establish \cref{eq:Lipschitz}. 
\end{cproof}

\subsection{Characterization results for realization functions of shallow ANNs}
\label{ssec:charac_shallow_ANNs}

The objective of this subsection is to establish \cref{cor:characterization2}, which provides a complete characterization of all functions in $C([0,1] , \R)$
 that can be represented by a shallow ANN with ReLU activation and $\width \in \N$ hidden neurons.
We first prove in \cref{lem:realization:l} a simple necessary condition: All representable functions are piecewise linear with at most $\width \in \N$ breakpoints.

\cfclear
\begin{lemma}
\label{lem:realization:l}
Assume \cref{setting:snn} and let $\theta \in \R^\fd$. 
Then
\begin{enumerate}[label = (\roman*)] 
\cfadd{def:piece:linear}\cfadd{def:number:of:kinks}
  \item \label{lem:realization:l:item1} 
  it holds that $ \realization{ \theta } \in \scrL $ and
  \item \label{lem:realization:l:item2} 
  it holds that $ Q( \realization{\theta} ) \leq \width $
\end{enumerate}
\cfload.
\end{lemma}
\begin{cproof}{lem:realization:l}
Throughout this proof 
let $ g_j \in C( [0,1] , \R ) $, $ j \in \{ 0, 1, \ldots, \width \} $, 
satisfy for all $ j \in \{ 1, 2, \ldots, \width\} $, $ x \in [0,1] $ 
that 
\begin{equation}
\label{eq:def_gj_charac}
  g_j( x ) 
  = \v{ \theta }_j 
  \max \cu[\big]{ \w{\theta}_j x + \b{\theta}_j , 0 }
  \qqandqq
  g_0 ( x ) = \c{\theta} 
  .
\end{equation}
\Nobs that \cref{eq:def_gj_charac} ensures
for all $ j \in \{ 1, 2, \ldots, \width \} $ 
that 
\begin{equation}
  g_j \in \scrL 
  \qqandqq 
  Q(g_j) \in \{ 0, 1 \} .
\end{equation}
Furthermore, \nobs that \cref{eq:def_gj_charac} 
demonstrates that $ g_0 \in \scrL $ and $ Q(g_0) = 0 $. 
Combining this, the fact that 
for all $ x \in [0,1] $ it holds that 
$ \realization{\theta}( x ) = \sum_{ j = 0 }^{ \width } g_j(x) $, 
\cref{lem:sum:piece:linear}, 
and induction establishes \cref{lem:realization:l:item1,lem:realization:l:item2}. 
\end{cproof}

Moreover, every piecewise linear function with at most $\width - 1 $ breakpoints is representable, as we show in \cref{lem:h-1:kinks}.

\cfclear
\begin{lemma} 
\label{lem:h-1:kinks}
\cfadd{def:piece:linear} \cfadd{def:number:of:kinks}
Assume \cref{setting:snn} and let $ g \in \scrL $ satisfy $ Q(g) \leq \width - 1 $ \cfload. 
Then there exists $ \theta \in \R^\fd $ such that $ \realization{\theta} = g $.
\end{lemma}
\begin{cproof}{lem:h-1:kinks}
Throughout this proof let $ \theta \in \R^\fd $ satisfy 
for all $ j \in \{ 1, 2, \ldots, \width \} $ that 
\begin{equation}
\label{eq:def_w}
  \w{\theta}_j = 
  \begin{cases}
    1 
  &
    \colon 
    j \leq Q(g) + 1
  \\
    0
  &
    \colon 
    j > Q(g) + 1
    ,
  \end{cases}
  \qquad 
  \b{\theta}_j 
  =
  \begin{cases}
    - q_{ j - 1 }(g) 
  &
    \colon 
    j \leq Q(g) + 1
  \\
    0
  &
    \colon 
    j > Q(g) + 1
    ,
  \end{cases}
\end{equation}
\begin{equation}
\label{eq:def_v}
  \v{\theta}_j 
  =
  \begin{cases}
    A_1(g)
  &
    \colon
    j = 1
  \\ 
    A_j(g) - A_{ j - 1 }(g) 
  &
    \colon 
    1 < j \leq Q(g) + 1
  \\
    0
  &
    \colon 
    j > Q(g) + 1 ,
  \end{cases}
\end{equation}
and 
$ \c{\theta} = B_1(g) $. 
\Nobs that \cref{eq:def_realization_function},  
\cref{eq:def_w}, and \cref{eq:def_v} ensure 
for all 
$ x \in [0,1] $ that 
\begin{equation}
\begin{split}
\textstyle
  \realization{\theta}(x)
& =
\textstyle
  \c{\theta} 
  +
  \sum\limits_{ j = 1 }^{ \width }
  \v{\theta}_j
  \max\{ 
    \w{\theta}_j x
    +
    \b{\theta}_j
  , 0 \}
  =
  B_1(g)
  +
  \sum\limits_{ j = 1 }^{ Q(g) + 1 }
  \v{\theta}_j
  \max\{ 
    \w{\theta}_j x
    +
    \b{\theta}_j
  , 0 \}
\\ & =
\textstyle
  B_1(g)
  +
  \sum\limits_{ j = 1 }^{ Q(g) + 1 }
  \v{\theta}_j
  \max\{ 
    x
    -
    q_{ j - 1 }(g)
  , 0 \}
\\ & =
\textstyle
  B_1(g)
  +
  A_1(g) x
  +
  \sum\limits_{ j = 2 }^{ Q(g) + 1 }
  (
    A_j(g) - A_{j-1}(g)
  )
  \max\{ 
    x
    -
    q_{ j - 1 }(g)
  , 0 \}
\\ & =
\textstyle
  B_1(g)
  +
  A_1(g) x
  +
  \sum\limits_{ j = 1 }^{ Q(g) }
  (
    A_{ j + 1 }(g) - A_j(g)
  )
  \max\{ 
    x
    -
    q_j(g)
  , 0 \}
  .
\end{split}
\end{equation}
Combining this with \cref{prop:piece:linear:relations} establishes 
for all $ i \in \{ 0, 1, \ldots, Q(g) \}$, $ x \in [ q_i(g) , q_{i+1}(g) ] $ 
that
\begin{equation}
\begin{split}
  \realization{\theta}( x ) 
& 
  = 
  B_1(g) + A_1(g) x 
  + \sum_{ j = 1 }^i ( A_{ j + 1 }(g) - A_j(g) )( x - q_j(g) ) 
\\
&
  = A_{i+1}( g ) x + B_1 (g) - \sum_{ j = 1 }^i (A_{j+1}(g) - A_j( g ) ) q_j(g) 
  = A_{i+1}(g) x + B_{i+1}(g) = g(x)
  .
\end{split}
\end{equation}
\end{cproof}

For piecewise linear functions with exactly $\width$ breakpoints, the situation is more involved: They are only representable by a shallow ANN with $\width$ hidden neurons
if the slopes fulfill a certain linear relation;
 see \cref{lem:necessary_condition:eq:defs} below for details.
  In \cref{lem:necessary_condition} we establish that this condition is necessary for a function to be representable, and afterwards we show in \cref{lem:sufficient_condition} that it is also sufficient. Both proofs proceed by induction on the number of breakpoints.

\cfclear
\begin{lemma}
\label{lem:necessary_condition}
For every $ \width \in \N_0 $, 
$ \theta = ( \theta_1, \dots, \theta_{ 3 \width + 1 } ) \in \R^{ 3 \width + 1 } $ 
let 
$ 
  \realization{\theta} 
  \colon [0,1] \to \R
$
satisfy for all 
$ x \in [0,1] $
that 
\begin{equation}
  \realization{ \theta }( x ) 
= 
  \theta_{ 3 \width + 1 }
  + 
  \sum_{ j = 1 }^{ \width }
  \theta_{ 2 \width + j }
  \max\{ 
    \theta_{ \width + j }
    +
    \theta_j x
    , 0 
  \}
  ,
\end{equation}
for every $ \width \in \N_0 $ let 
$ \bfR_{ \width } \subseteq C( [0,1], \R ) $ 
satisfy 
$
  \bfR_{ \width }
  =
  \{ 
    f \in 
    Q^{ - 1 }( \{ \width \} )
    \colon 
    [
      \exists \, \theta \in \R^{ 3 \width + 1 } \colon 
      f = \realization{\theta}
    ]
  \}
$,
and 
for every $ \width \in \N_0 $ 
let 
$
\cfadd{def:number:of:kinks,def:slopes}
  \bfS_{ \width } \subseteq C( [0,1], \R ) 
$
satisfy 
\begin{multline} \label{lem:necessary_condition:eq:defs}
  \bfS_{ \width } 
  = 
  \Big\{ 
    f \in Q^{ - 1 }( \{ \width \} )
    \colon 
    \Big(
      \exists \, k \in \N, i_1, i_2, \dots, i_k \in \N \colon
\\
\textstyle
      \big[ 
        \big( 
          \tfrac{ k }{ 2 } \notin \N 
        \big) ,
        \big(
          i_1 < i_2 < \dots < i_k \leq \width + 1 
        \big) ,
        \big(
          \sum_{ j = 1 }^k 
          ( - 1 )^j 
          A_{ i_j }( f )
          = 0
        \big)
      \big]
    \Big)
  \Big\}
\end{multline}
\cfout.
Then it holds for all $ \width \in \N_0 $
that 
\begin{equation}
\label{eq:necessary_condition}
  \bfR_{ \width } \subseteq \bfS_{ \width } 
  .
\end{equation}
\end{lemma}
\begin{cproof}{lem:necessary_condition}
Throughout this proof 
let 
$ 
  \operatorname{sgn} \colon \R \to \R 
$ 
satisfy for all $ x \in (0,\infty) $, $ k \in \{ - 1, 0, 1 \} $ that 
$
  \operatorname{sgn}( k x ) = k
$
and 
for every $ \width \in \N $, 
$ \theta = ( \theta_1, \dots, \theta_{ 3 \width + 1 } ) \in \R^{ 3 \width + 1 } $, 
$ j \in \{ 1, 2, \dots, \width \} $ 
let 
$ \q{ \theta }_j \in (-\infty, \infty] $
satisfy
\begin{equation}
\label{eq:def_q_in_necessary_condition}
  \q{ \theta }_j 
  =
  \begin{cases}
    -
    \frac{
      \theta_{ \width + j }
    }{
      \theta_j
    }
  &
    \colon \theta_j \neq 0
  \\
    \infty
  &
    \colon \theta_j = 0 .
  \end{cases}
\end{equation}
Observe that 
\cref{eq:def_q_in_necessary_condition} 
ensures that
for all 
$ \width \in \N $,
$ \theta \in \R^{ 3 \width + 1 } $
it holds that 
\begin{equation}
\label{eq:Q_upper_estimate}
  Q( \realization{\theta} ) 
  \leq 
  \big| 
    \big(
      \{ 
        \q{\theta}_1 ,
        \q{\theta}_2 ,
        \dots ,
        \q{\theta}_{ \width }
      \}
      \cap \R
    \big)
  \big|
  \leq 
  \width .
\end{equation}
We prove \cref{eq:necessary_condition} by induction 
on $ \width \in \N_0 $. 
For the base case 
$ \width = 0 $ observe that for all $ \theta \in \R $, $ x \in [0,1] $ 
it holds that 
$
  \realization{\theta}(x) = \theta
$. 
Hence, we obtain that for all $ \theta \in \R $ that 
$
  Q( \realization{\theta} ) = 0
$, 
$ 
  A_1( \realization{\theta} ) = 0
$,
and 
$
  B_1( \realization{\theta} ) = \theta
$. 
Therefore, we obtain that 
$
  \bfR_0 = ( \cup_{ \theta \in \R } \{ \realization{\theta} \} ) \subseteq \bfS_0
$.
This establishes \cref{eq:necessary_condition} 
in the base case $ \width = 0 $. 
For the induction step let $ \width \in \N_0 $ 
satisfy $ \bfR_{ \width } \subseteq \bfS_{ \width } $
and let $ F \in \bfR_{ \width + 1 } $. 
We intend to prove that $ F \in \bfS_{ \width + 1 } $.
Observe that the fact that $ F \in \bfR_{ \width + 1 } $ 
ensures that 
there exists 
$ \Xi \in \R^{ 3 ( \width + 1 ) + 1 } = \R^{ 3 \width + 4 } $ 
which satisfies 
$
  \realization{\Xi} = F 
$.
\Nobs that \cref{eq:Q_upper_estimate} 
and the fact that $ Q( F ) = \width + 1 $ 
demonstrate that 
$
  \q{\Xi}_1, \q{\Xi}_2, \dots, \q{\Xi}_{ \width + 1 } 
  \in \R
$
and 
$
  | \{ \q{\Xi}_1, \q{\Xi}_2, \dots, \q{\Xi}_{ \width + 1 } \} | 
  = \width + 1
$. 
This shows that there exists a bijective 
$ p \colon \{ 1, 2, \dots, \width + 1 \} \to \{ 1, 2, \dots, \width + 1 \} $ 
which satisfies 
\begin{equation}
\label{eq:properties_of_q}
  - \infty < 
  \q{ \Xi }_{ p(1) } < \q{ \Xi }_{ p(2) } 
  < \dots < \q{ \Xi }_{ p(\width+1) } 
  < \infty
  .
\end{equation}
In the following let 
$ 
  \Theta = ( \Theta_1, \dots, \Theta_{ 3 \width + 4 } ) 
$ 
satisfy 
for all 
$ j \in \{ 1, 2, \dots, \width + 1 \} $ 
that 
\begin{equation}
\label{eq:def_vartheta_necessary_condition}
  \Theta_j 
  = \Xi_{ p(j) }
  ,
  \quad 
  \Theta_{ \width + 1 + j }
  = \Xi_{ \width + 1 + p(j) }
  ,
  \quad 
  \Theta_{ 2 \width + 2 + j }
  = \Xi_{ 2 \width + 2 + p(j) } 
  ,
  \qandq
  \Theta_{ 3 \width + 4 }
  = \Xi_{ 3 \width + 4 }
  .
\end{equation}
\Nobs that 
\cref{eq:properties_of_q}, 
\cref{eq:def_vartheta_necessary_condition}, 
and the fact that 
$
  F = \realization{\Xi}
$
ensure that 
\begin{equation}
\label{eq:new_q_condition:necessary_condition}
  \realization{\Theta} = \realization{\Xi} = F
\qqandqq
  - \infty 
  < 
  \q{ \Theta }_1
  < 
  \q{ \Theta }_2 
  < 
  \dots 
  <
  \q{ \Theta }_{ \width + 1 }
  < 
  \infty 
  .
\end{equation}
In the following let 
$ 
  \theta = ( \theta_1, \dots, \theta_{ 3 \width + 1 } ) 
  \in \R^{ 3 \width + 1 }
$
satisfy 
\begin{equation}
\label{eq:def:varTheta:necessary_condition}
  \theta = 
  ( 
    \Theta_1, \dots, \Theta_{ \width }, 
    \Theta_{ \width + 2 }, \dots, \Theta_{ 2 \width + 1 }, 
    \Theta_{ 2 \width + 3 }, \dots, \Theta_{ 3 \width + 2 }, 
    \Theta_{ 3 \width + 4 } 
  )
\end{equation}
and let 
$ f \in C( [0,1], \R ) $ 
satisfy 
$ f = \realization{\theta} $. 
\Nobs that \cref{eq:def:varTheta:necessary_condition} 
ensures
for all $ x \in [0,1] $ 
that 
\begin{equation}
\label{eq:relation_f_and_F:necessary_condition}
  f(x)
  =
  \theta_{ 3 \width + 1 }
  +
  \sum_{ j = 1 }^{ \width }
  \theta_{ 2 \width + j }
  \max\{ 
    \theta_{ \width + j }
    +
    \theta_j x
    , 0 
  \}
  = 
  F(x)
  -
  \Theta_{ 3 \width + 3 }
  \max\{ 
    \Theta_{ 2 \width + 2 }
    +
    \Theta_{ \width + 1 } x
    , 0 
  \} 
  .
\end{equation} 
Next \nobs that \cref{eq:new_q_condition:necessary_condition} 
assures that 
$
  Q( f ) = \width 
$
and 
$
  - \infty 
  < \q{ \theta }_1 = \q{ \Theta }_1 
  < \q{ \theta }_2 = \q{ \Theta }_2
  < \dots 
  < \q{ \theta }_{ \width } = \q{ \Theta }_{ \width } 
  < \q{ \Theta }_{ \width + 1 }
  < \infty 
$. 
Combining this with the fact that 
$ 
  f = \realization{\theta} 
$ 
demonstrates that $ f \in \bfR_{ \width } $. The induction hypothesis 
that $ \bfR_{ \width } \subseteq \bfS_{ \width } $ 
\hence assures
that $ f \in \bfS_{ \width } $. This proves that there exist 
$ k \in \N $, $ i_1, i_2, \dots, i_k \in \N $ 
which satisfy 
\begin{equation}
\label{eq:necessary_condition:ik_properties}
\textstyle
  \frac{ k }{ 2 } \notin \N,
  \qquad
  i_1 < i_2 < \dots < i_k \leq \width + 1 ,
  \qqandqq 
  \sum_{ j = 1 }^k ( - 1 )^j A_{ i_j }( f ) = 0
  .
\end{equation}
Next let $ K \in \N $, 
$ I_1, I_2, \dots, I_K \in \N $ 
satisfy     
\begin{equation}
\label{eq:necessary_condition:def_IK}
  \{ I_1, I_2, \dots, I_K \}
  =
  \begin{cases}
    \{ i_1, i_2, \dots, i_k \}
  &
    \colon 
    \Theta_{ \width + 1 } > 0
  \\
    ( \cup_{ l = 1 }^{ k - 1 } \{ i_l \} )
    \cup 
    \{ \width + 2 \} 
  &
    \colon 
    \Theta_{ \width + 1 } < 0 = \width + 1 - i_k
  \\
    \{ i_1, \dots, i_k, \width + 1, \width + 2 \}
  &
    \colon 
    \Theta_{ \width + 1 } < 0 < \width + 1 - i_k .
  \end{cases}
\end{equation}
\Nobs that 
\cref{eq:necessary_condition:ik_properties} 
and
\cref{eq:necessary_condition:def_IK} ensure that 
$ \frac{ K }{ 2 } \notin \N $ 
and 
$ I_1 < I_2 < \dots < I_K \leq \width + 2 $. 
In order to prove that $ F \in \bfS_{ \width + 1 } $, 
it is thus sufficient to verify that 
\begin{equation}
\label{eq:necessary_condition:remains_to_prove}
\textstyle
  \sum_{ j = 1 }^K ( - 1 )^j A_{ I_j }( F ) = 0
  .
\end{equation}
For this \nobs that 
\cref{eq:relation_f_and_F:necessary_condition} 
assures for all $ x \in [0,1] $ 
that
\begin{equation}
\label{eq:relation_f_and_F_2:necessary_condition}
\begin{split}
  F( x )
& =
  f(x)
  +
  |
    \Theta_{ \width + 1 }
  |
  \Theta_{ 3 \width + 3 }
  \max\{ 
    |
      \Theta_{ \width + 1 }
    |^{ - 1 }
    \Theta_{ 2 \width + 2 }
    +
    |
      \Theta_{ \width + 1 }
    |^{ - 1 }
    \Theta_{ \width + 1 } x
    , 0 
  \} 
\\ & = 
  f(x)
  +
  \Theta_{ \width + 1 }
  \Theta_{ 3 \width + 3 }
  \operatorname{sgn}( \Theta_{ \width + 1 } ) 
  \max\{ 
    ( 
      x 
      - \q{ \Theta }_{ \width + 1 }
    )
    \operatorname{sgn}( \Theta_{ \width + 1 } ) 
    , 0 
  \} 
  .
\end{split}
\end{equation}
In the following we distinguish between the case 
$
  \Theta_{ \width + 1 } > 0
$,
the case 
$
  \Theta_{ \width + 1 } < 0 = \width + 1 - i_k
$, 
and the case 
$
  \Theta_{ \width + 1 } < 0 < \width + 1 - i_k 
$. We first prove \cref{eq:necessary_condition:remains_to_prove} 
in the case 
\begin{equation}
\label{eq:necessary_condition:case1}
  \Theta_{ \width + 1 } > 0
  .
\end{equation}
\Nobs that 
\cref{eq:relation_f_and_F_2:necessary_condition}
and \cref{eq:necessary_condition:case1} 
demonstrate for all $ x \in [0,1] $ 
that 
$
  F( x )
=
  f(x)
  +
  \Theta_{ \width + 1 }
  \Theta_{ 3 \width + 3 }
  \max\{ 
      x 
      - \q{ \Theta }_{ \width + 1 }
    , 0 
  \} 
$.
Hence, we obtain 
for all $ j \in \{ 1, 2, \ldots, \width + 1 \} $ 
that $ A_j( F ) = A_j( f ) $. 
Combining this with \cref{eq:necessary_condition:def_IK} implies that 
$
  \sum_{ j = 1 }^K (-1)^j A_{ I_j }( F ) 
= 
  \sum_{ j = 1 }^k (-1)^j A_{ i_j }( f ) 
= 0
$. 
This establishes \cref{eq:necessary_condition:remains_to_prove} 
in the case $ \Theta_{ \width + 1 } > 0 $. 
In the next step we prove 
\cref{eq:necessary_condition:remains_to_prove} 
in the case 
\begin{equation}
\label{eq:necessary_condition:case2}
  \Theta_{ \width + 1 } < 0 = \width + 1 - i_k
  .
\end{equation}
\Nobs that 
\cref{eq:relation_f_and_F_2:necessary_condition}
and \cref{eq:necessary_condition:case2} 
show for all $ x \in [0,1] $ 
that 
$
  F( x )
=
  f(x)
  +
  \Theta_{ \width + 1 }
  \Theta_{ 3 \width + 3 }
  \min\{ 
      x 
      - \q{ \Theta }_{ \width + 1 }
    , 0 
  \} 
$. 
Therefore, we obtain 
for all $ j \in \{ 1, 2, \ldots, \width + 1 \} $ 
that 
$ A_j(F) = A_j(f) + \Theta_{ \width + 1 } \Theta_{ 3 \width + 3 } $ 
and $ A_{ \width + 2 }( F ) = A_{ \width + 1 }(f) $. 
Combining this with 
\cref{eq:necessary_condition:def_IK}, 
\cref{eq:necessary_condition:case2}, 
and the fact that $ \frac{ k }{ 2 } \notin \N $
shows that 
\begin{equation}
\begin{split}
\textstyle
  \sum_{ j = 1 }^K (-1)^j A_{ I_j }( F ) 
& =
\textstyle
  \left[ 
    \sum_{ j = 1 }^{ K - 1 } 
    (-1)^j A_{ I_j }( F ) 
  \right]
  +
  (-1)^K  A_{ I_K }( F ) 
\\ & =
\textstyle
  \left[ 
    \sum_{ j = 1 }^{ k - 1 } 
    (-1)^j 
    \left( 
      A_{ i_j }( f ) 
      +
      \Theta_{ \width + 1 } 
      \Theta_{ 3 \width + 3 }
    \right)
  \right]
  -
  A_{ \width + 2 }( F ) 
\\ & =
\textstyle
  \left[ 
    \sum_{ j = 1 }^{ k - 1 } 
    (-1)^j 
    A_{ i_j }( f ) 
  \right]
  +
  \left[
    \sum_{ j = 1 }^{ k - 1 } 
    (-1)^j 
      \Theta_{ \width + 1 } 
      \Theta_{ 3 \width + 3 }
  \right]
  -
  A_{ \width + 1 }( f ) 
\\ & =
\textstyle
  \left[ 
    \sum_{ j = 1 }^{ k - 1 } 
    (-1)^j 
    A_{ i_j }( f ) 
  \right]
  +
  \left[
    \sum_{ j = 1 }^{ k - 1 } 
    (-1)^j 
  \right]
      \Theta_{ \width + 1 } 
      \Theta_{ 3 \width + 3 }
  +
  ( - 1 )^k
  A_{ i_k }( f ) 
\\ & =
\textstyle
    \sum_{ j = 1 }^k 
    (-1)^j 
    A_{ i_j }( f ) 
= 0
  .
\end{split}
\end{equation}
This establishes 
\cref{eq:necessary_condition:remains_to_prove} 
in the case 
$
  \Theta_{ \width + 1 } < 0 = \width + 1 - i_k
$. 
Next we prove \cref{eq:necessary_condition:remains_to_prove}
in the case 
\begin{equation}
\label{eq:necessary_condition:case3}
  \Theta_{ \width + 1 } < 0 < \width + 1 - i_k
  .
\end{equation}
\Nobs that 
\cref{eq:relation_f_and_F_2:necessary_condition}
and \cref{eq:necessary_condition:case3} 
demonstrate for all $ x \in [0,1] $ 
that 
$
  F( x )
=
  f(x)
  +
  \Theta_{ \width + 1 }
  \Theta_{ 3 \width + 3 }
  \min\{ 
      x 
      - \q{ \Theta }_{ \width + 1 }
    , 0 
  \} 
$. 
Hence, we obtain 
for all $ j \in \{ 1, 2, \ldots, \width + 1 \} $ 
that 
$ A_j(F) = A_j(f) + \Theta_{ \width + 1 } \Theta_{ 3 \width + 3 } $ 
and $ A_{ \width + 2 }( F ) = A_{ \width + 1 }(f) $. 
Combining this with 
\cref{eq:necessary_condition:def_IK}, 
\cref{eq:necessary_condition:case3}, 
and the fact that $ \frac{ k }{ 2 } \notin \N $
shows that 
\begin{equation}
\begin{split}
&
\textstyle
  \sum_{ j = 1 }^K (-1)^j A_{ I_j }( F ) 
\\
& =
\textstyle
  \left[ 
    \sum_{ j = 1 }^{ k } 
    (-1)^j A_{ i_j }( F ) 
  \right]
  +
  ( - 1 )^{ k + 1 }
  A_{ I_{ k + 1 } }( F ) 
  +
  ( - 1 )^{ k + 2 }
  A_{ I_{ k + 2 } }( F ) 
\\ & =
\textstyle
  \left[ 
    \sum_{ j = 1 }^{ k } 
    (-1)^j 
    \left( 
      A_{ i_j }( f ) 
      +
      \Theta_{ \width + 1 }
      \Theta_{ 3 \width + 3 }
    \right)
  \right]
  +
  A_{ I_{ k + 1 } }( F ) 
  -
  A_{ I_{ k + 2 } }( F ) 
\\ & =
\textstyle
  \left[ 
    \sum_{ j = 1 }^{ k } 
    (-1)^j 
    A_{ i_j }( f ) 
  \right]
  +
  \left[ 
    \sum_{ j = 1 }^{ k } 
    (-1)^j 
    \Theta_{ \width + 1 }
    \Theta_{ 3 \width + 3 }
  \right]
  +
  A_{ \width + 1 }( F ) 
  -
  A_{ \width + 2 }( F ) 
\\ & =
\textstyle
  \left[ 
    \sum_{ j = 1 }^{ k } 
    (-1)^j 
  \right]
  \Theta_{ \width + 1 }
  \Theta_{ 3 \width + 3 }
  +
  A_{ \width + 1 }( f ) 
  +
  \Theta_{ \width + 1 }
  \Theta_{ 3 \width + 3 }
  -
  A_{ \width + 1 }( f ) 
  =
  0
  .
\end{split}
\end{equation}
This establishes 
\cref{eq:necessary_condition:remains_to_prove} 
in the case 
$
  \Theta_{ \width + 1 } < 0 < \width + 1 - i_k
$. 
\Nobs that 
\cref{eq:necessary_condition:remains_to_prove}, 
the fact that $ \frac{ K }{ 2 } \notin \N $, 
and the fact that 
$
  I_1 < I_2 < \dots < I_K \leq \width + 2 
$
prove that 
$
  F \in \bfS_{ \width + 1 }
$. 
Induction thus establishes \cref{eq:necessary_condition}. 
\end{cproof}

\begin{lemma}
\label{lem:sufficient_condition}
For every $ \width \in \N_0 $, 
$ \theta = ( \theta_1, \dots, \theta_{ 3 \width + 1 } ) \in \R^{ 3 \width + 1 } $ 
let 
$ 
  \realization{\theta} 
  \colon [0,1] \to \R
$
satisfy for all 
$ x \in [0,1] $
that 
\begin{equation}
  \realization{ \theta }( x ) 
= 
  \theta_{ 3 \width + 1 }
  + 
  \sum_{ j = 1 }^{ \width }
  \theta_{ 2 \width + j }
  \max\{ 
    \theta_{ \width + j }
    +
    \theta_j x
    , 0 
  \}
  ,
\end{equation}
for every $ \width \in \N_0 $ let 
$ \bfR_{ \width } \subseteq C( [0,1], \R ) $ 
satisfy 
$
  \bfR_{ \width }
  =
  \{ 
    f \in 
    Q^{ - 1 }( \{ \width \} )
    \colon 
    [
      \exists \, \theta \in \R^{ 3 \width + 1 } \colon 
      f = \realization{\theta}
    ]
  \}
$,
and 
for every $ \width \in \N_0 $ 
let 
$
\cfadd{def:number:of:kinks,def:slopes}
  \bfS_{ \width } \subseteq C( [0,1], \R ) 
$
satisfy 
\begin{multline}
\label{eq:sufficient_condition:defS}
  \bfS_{ \width } 
  = 
  \Big\{ 
    f \in Q^{ - 1 }( \{ \width \} ) 
    \colon 
    \Big(
      \exists \, k \in \N, i_1, i_2, \dots, i_k \in \N \colon
\\
\textstyle
      \big[ 
        \big( 
          \tfrac{ k }{ 2 } \notin \N 
        \big) ,
        \big(
          i_1 < i_2 < \dots < i_k \leq \width + 1 
        \big) ,
        \big(
          \sum_{ j = 1 }^k 
          ( - 1 )^j 
          A_{ i_j }( f )
          = 0
        \big)
      \big]
    \Big)
  \Big\}
\end{multline}
\cfout.
Then it holds for all $ \width \in \N_0 $
that 
\begin{equation}
\label{eq:sufficient_condition}
  \bfS_{ \width } 
  \subseteq 
  \bfR_{ \width } 
  .
\end{equation}
\end{lemma}
\begin{cproof}{lem:sufficient_condition}
We prove \cref{eq:sufficient_condition} by induction on $ \width \in \N_0 $. 
For the base case $ \width = 0 $ 
\nobs that \cref{eq:sufficient_condition:defS} 
ensures that
\begin{equation}
\begin{split}
  \bfS_0 
& = 
  \left\{ 
    f \in Q^{ - 1 }( \{ 0 \} ) \colon A_1(f) = 0
  \right\}
  =
  \left\{ 
    f \in C( [0,1], \R ) 
    \colon 
    ( \forall \, x \in [0,1] \colon f(x) = f(0) )
  \right\}
\\ &
=
  \cup_{ \theta \in \R }
  \{ \realization{\theta} \} 
  .
\end{split}
\end{equation}
This establishes \cref{eq:sufficient_condition:defS} in the base case $ \width = 0 $. 
For the induction step let 
$ \width \in \N_0 $ 
satisfy 
$ \bfS_{ \width } \subseteq \bfR_{ \width } $ 
and let $ F \in \bfS_{ \width + 1 } $. 
We intend to prove that $ F \in \bfR_{ \width + 1 } $. 
Note that \cref{eq:sufficient_condition:defS} ensures that 
there exist $ K \in \N $, 
$ I_1, I_2, \dots, I_K \in \N $ 
which satisfy 
\begin{equation}
\label{eq:properties_F:sufficient_condition}
\textstyle
  Q( F ) = \width + 1 , 
\quad
  \frac{ K }{ 2 } \notin \N ,
\quad
  I_1 < I_2 < \dots < I_K \leq \width + 2 ,
\qandq
  \sum_{ j = 1 }^K
  ( - 1 )^j
  A_{ I_j }( F )
  = 0
  .
\end{equation}
Next let 
$ f \colon [0,1] \to \R $ satisfy for all $ x \in [0,1] $ that 
\begin{equation} 
\label{eq:sufficient_condition:def_induction_function}
  f(x) 
  = 
  \begin{cases}
    F(x) - ( A_{ \width + 2 }( F ) - A_{ \width + 1 }( F ) ) 
    \max\cu{ x - q_{\width + 1 }( F ) , 0 } 
  &
    \colon 
    I_K - 2 \neq \width
\\
    F(x) - ( A_{ \width + 2 }( F ) - A_{ \width + 1 }( F ) ) 
    \max\cu{ q_{\width + 1 }( F ) - x , 0 } 
  &
    \colon 
    I_K - 2 = \width
  \end{cases}
\end{equation}
and let 
$ k \in \N $, $ i_1, i_2, \dots, i_k \in \N $
satisfy 
\begin{equation}
\label{eq:sufficient_condition:def_ik}
  \{ i_1, i_2, \dots, i_k \} 
  =
  \begin{cases}
    \{ I_1, I_2, \dots, I_K \}
  &
    \colon 
    I_K - 2 \neq \width 
  \\
    \{ I_1, I_2, \dots, I_{ K - 2 } \}
  &
    \colon
    I_K - 2 = \width 
    < 
    \min\{ 
      I_{ \max\{ K - 1 , 1 \} } ,
      K + \width - 1
    \} 
  \\
    ( \cup_{ l = 1 }^{ K - 1 } \{ I_l \} ) \cup \{ \width + 1 \}
  &
    \colon 
    I_K - 2 = \width 
    \geq 
    \min\{ 
      I_{ \max\{ K - 1 , 1 \} } ,
      K + \width - 1
    \} 
    .
  \end{cases}
\end{equation}
\Nobs that 
\cref{eq:properties_F:sufficient_condition} and 
\cref{eq:sufficient_condition:def_ik} 
assure that 
$ 
  \frac{ k }{ 2 } \notin \N 
$
and 
$
  i_1 < i_2 < \dots < i_k \leq \width + 1
$. 
Moreover, \nobs that 
\cref{eq:properties_F:sufficient_condition} and 
\cref{eq:sufficient_condition:def_induction_function}
ensure that 
$
  f \in C( [0,1], \R ) 
$,
$
  Q(f) = \width 
$, 
and 
\begin{equation}
\label{eq:relation_f_and_F:sufficient_condition}
  \left(
    \forall \, i \in \{ 1, 2, \dots, \width + 1 \} \colon 
    A_i( f ) 
    = 
    \begin{cases}
      A_i( F )
    &
      \colon 
      I_K - 2 \neq \width 
    \\
      A_i( F ) 
      +
      A_{ \width + 2 }( F ) - A_{ \width + 1 }( F )
    & 
      \colon
      I_K - 2 = \width 
    \end{cases}
  \right)
  .
\end{equation}
In the next step we prove that 
\begin{equation}
\label{eq:sufficient_condition:to_prove}
\textstyle
  \sum_{ j = 1 }^k (-1)^j A_{ i_j }( f ) = 0
  .
\end{equation}
In the following we distinguish between 
the case 
$
  I_K - 2 \neq \width  
$, 
the case 
$
  I_K - 2 = \width 
  < 
  \min\{ 
    I_{ \max\{ K - 1 , 1 \} } ,
\allowbreak 
    K + \width - 1
  \} 
$, 
and the case 
$
  I_K - 2 = \width 
  \geq 
  \min\{ 
    I_{ \max\{ K - 1 , 1 \} } ,
    K + \width - 1
  \} 
$. 
We first prove 
\cref{eq:sufficient_condition:to_prove} 
in the case 
\begin{equation}
\label{eq:sufficient_condition:case1}
  I_K - 2 \neq \width  
  .
\end{equation}
\Nobs that 
\cref{eq:sufficient_condition:def_ik}
and 
\cref{eq:sufficient_condition:case1} 
ensure that 
\begin{equation}
\textstyle
  \sum_{ j = 1 }^k (-1)^j A_{ i_j }( f ) 
=
  \sum_{ j = 1 }^K (-1)^j A_{ I_j }( F ) 
= 
  0
  .
\end{equation}
This establishes 
\cref{eq:sufficient_condition:to_prove} in the case 
$
  I_K - 2 \neq \width  
$. 
In the next step we prove 
\cref{eq:sufficient_condition:to_prove} 
in the case 
\begin{equation}
\label{eq:sufficient_condition:case2}
  I_K - 2 = \width 
  < 
  \min\{ 
    I_{ \max\{ K - 1 , 1 \} } ,
    K + \width - 1
  \} 
  .
\end{equation}
\Nobs that 
\cref{eq:properties_F:sufficient_condition}, 
\cref{eq:sufficient_condition:def_ik}, 
\cref{eq:relation_f_and_F:sufficient_condition}, 
and 
\cref{eq:sufficient_condition:case2} 
assure that 
\begin{equation}
\begin{split}
&
\textstyle
  \sum_{ j = 1 }^k (-1)^j A_{ i_j }( f ) 
=
  \sum_{ j = 1 }^{ K - 2 } (-1)^j A_{ I_j }( f ) 
\\ & 
\textstyle
=
  \sum_{ j = 1 }^{ K - 2 } 
  (-1)^j 
  \left( 
    A_{ I_j }( F ) 
    +
    A_{ \width + 2 }( F )
    -
    A_{ \width + 1 }( F )
  \right)
\\ & 
\textstyle
= 
  \left[
    \sum_{ j = 1 }^{ K - 2 } 
    (-1)^j 
    A_{ I_j }( F ) 
  \right]
  +
  \left[
    \sum_{ j = 1 }^{ K - 2 } 
    (-1)^j 
  \right]
  \left[
    A_{ \width + 2 }( F )
    -
    A_{ \width + 1 }( F )
  \right]
\\ &
\textstyle
=
  \left[
    \sum_{ j = 1 }^K 
    (-1)^j 
    A_{ I_j }( F ) 
  \right]
  -
  \left[ 
    \sum_{ j = K - 1 }^K
    (-1)^j 
    A_{ I_j }( F ) 
  \right]
  -
  \left[
    A_{ \width + 2 }( F )
    -
    A_{ \width + 1 }( F )
  \right]
\\ &
\textstyle
=
  -
  \left[ 
    A_{ I_{ K - 1 } }( F ) 
    -
    A_{ I_K }( F ) 
  \right]
  -
  \left[
    A_{ \width + 2 }( F )
    -
    A_{ \width + 1 }( F )
  \right]
\\ &
\textstyle
=
  -
  \left[ 
    A_{ \width + 1 }( F ) 
    -
    A_{ \width + 2 }( F ) 
  \right]
  -
  \left[
    A_{ \width + 2 }( F )
    -
    A_{ \width + 1 }( F )
  \right]
=
  0
  .
\end{split}
\end{equation}
This establishes 
\cref{eq:sufficient_condition:to_prove} in the case 
$
  I_K - 2 = \width 
  < 
  \min\{ 
    I_{ \max\{ K - 1 , 1 \} } ,
    K + \width - 1
  \} 
$. 
In the next step we prove 
\cref{eq:sufficient_condition:to_prove} 
in the case 
\begin{equation}
\label{eq:sufficient_condition:case3}
  I_K - 2 = \width 
  \geq 
  \min\{ 
    I_{ \max\{ K - 1 , 1 \} } ,
    K + \width - 1
  \} 
  .
\end{equation}
\Nobs that 
\cref{eq:properties_F:sufficient_condition}, 
\cref{eq:sufficient_condition:def_ik}, 
\cref{eq:relation_f_and_F:sufficient_condition}, 
and 
\cref{eq:sufficient_condition:case3} 
assure that 
\begin{equation}
\begin{split}
&
\textstyle
  \sum_{ j = 1 }^k (-1)^j A_{ i_j }( f ) 
=
  \left[ 
    \sum_{ j = 1 }^{ K - 1 } (-1)^j A_{ I_j }( f ) 
  \right]
  +
  ( - 1 )^k
  A_{ \width + 1 }( f )
\\ & 
\textstyle
  =
  \left[ 
    \sum_{ j = 1 }^{ K - 1 } (-1)^j 
    \big( 
      A_{ I_j }( F ) + A_{ \width + 2 }( F ) - A_{ \width + 1 }( F ) 
    \big)
  \right]
  -
  \left[
    A_{ \width + 1 }( F )
    +
    A_{ \width + 2 }( F )
    -
    A_{ \width + 1 }( F )
  \right]
\\ & 
\textstyle
=
  \left[ 
    \sum_{ j = 1 }^{ K - 1 } 
    (-1)^j 
    A_{ I_j }( F ) 
  \right]
  +
  \left[
    \sum_{ j = 1 }^{ K - 1 } 
    (-1)^j 
  \right]
  \left( 
    A_{ \width + 2 }( F ) - A_{ \width + 1 }( F ) 
  \right)
  -
  A_{ \width + 2 }( F )
\\ & 
\textstyle
=
  \left[ 
    \sum_{ j = 1 }^K 
    (-1)^j 
    A_{ I_j }( F ) 
  \right]
  +
  A_{ I_K }( F ) 
  -
  A_{ \width + 2 }( F )
=
  A_{ \width + 2 }( F ) 
  -
  A_{ \width + 2 }( F )
=
  0
  .
\end{split}
\end{equation}
This establishes 
\cref{eq:sufficient_condition:to_prove} in the case 
$
  I_K - 2 = \width 
  \geq 
  \min\{ 
    I_{ \max\{ K - 1 , 1 \} } ,
    K + \width - 1
  \} 
$. 
\Nobs that \cref{eq:sufficient_condition:to_prove} 
implies that 
$ f \in \bfS_{ \width } $. 
The induction hypothesis that 
$
  \bfS_{ \width } \subseteq \bfR_{ \width }
$
hence assures that 
$
  f \in \bfR_{ \width }
$. 
Combining this with 
\cref{eq:sufficient_condition:def_induction_function} 
and the fact that $ Q( F ) = \width + 1 $ 
shows that 
$
  F \in \bfR_{ \width + 1 }
$. 
Induction thus establishes \cref{eq:sufficient_condition}. 
\end{cproof}

Finally, in \cref{cor:characterization2} we combine the previous results to obtain the promised characterization.

\cfclear
\begin{cor}
\label{cor:characterization2}
Let 
$ \width \in \N_0 $, 
for every 
$ 
  \theta = ( \theta_1, \dots, \theta_{ 3 \width + 1 } ) \in \R^{ 3 \width + 1 } 
$
let 
$ 
  \realization{\theta} \colon [0,1] \to \R 
$
satisfy 
$ x \in [0,1] $
that 
$
  \realization{ \theta }( x ) 
= 
  \theta_{ 3 \width + 1 }
  + 
  \sum_{ j = 1 }^{ \width }
  \theta_{ 2 \width + j }
  \max\{ 
    \theta_{ \width + j }
    +
    \theta_j x
    , 0 
  \}
$, 
and let 
$ f \in C( [0,1], \R ) $.
Then the following two statements are equivalent: 
\begin{enumerate}[label = (\roman*)]
\item 
\label{charac2:item_i}
Then exists $ \theta \in \R^{ 3 \width + 1 } $ such that $ \realization{\theta} = f $. 
\item 
\label{charac2:item_ii}
There exist $ k \in \N $, $ i_1, i_2, \ldots, i_k \in \N $ such that
$ \frac{ k }{ 2 } \notin \N$,
$ i_1 < i_2 < \cdots < i_k \leq \width + 1 $,
$ Q(f) \leq \width $, 
and 
\begin{equation}
\cfadd{def:number:of:kinks,def:slopes}
\textstyle
  ( \width - Q(f) - 1 ) \bigl| \sum_{ j = 1 }^k (-1)^j A_{ \min\{ i_j, Q(f) + 1 \} }(f) \bigr| \geq 0 
\end{equation} \cfout.
\end{enumerate}
\end{cor}
\begin{cproof}{cor:characterization2} 
\Nobs that 
\cref{lem:realization:l} and 
\cref{lem:necessary_condition} 
prove that 
(\cref{charac2:item_i}\allowbreak$ \rightarrow $\allowbreak\cref{charac2:item_ii}). 
\Moreover 
\cref{lem:h-1:kinks} and \cref{lem:sufficient_condition} 
establish that 
(\cref{charac2:item_ii}\allowbreak$ \rightarrow $\allowbreak\cref{charac2:item_i}). 
\end{cproof}

\subsection{Structure preserving approximations 
for piecewise affine linear functions}

The next elementary lemma is an immediate consequence of the definitions in \cref{ssec:prop_breakpoint}. It will be employed in the sequel to switch the endpoints of the domain $[0,1]$ and thereby make some simplifying assumptions.

\cfclear
\begin{lemma} 
\label{lem:breakpoint:transformation}
Let $ L \in \R $, $ f \in C( [0,1], \R ) $ 
satisfy for all $ x, y \in [0,1] $ that 
$
  \abs{ f(x) - f(y) } \leq L \abs{ x - y } 
$, 
let $ g \in \scrL $, $ i \in \{ 1, 2, \ldots, Q(g) + 1 \} $, $ \bfa \in \R $ 
satisfy $ L \le \abs{\bfa} \le \abs{ A_i(g) } $ and $ \bfa A_i( g ) > 0 $,
let $F \colon [0,1] \to \R$ and $G \colon [0,1] \to \R$ 
satisfy for all $x \in [0,1]$ that $F ( x ) = - f ( 1 - x ) $
 and $G ( x ) = - g ( 1 - x ) $,
and let $I \in \N$ satisfy $I = Q(g) + 2 - i $
\cfadd{def:piece:linear}\cfadd{def:number:of:kinks}\cfadd{def:slopes}\cfload. 
Then
\begin{enumerate} [label = (\roman*)]
\item \label{lem:breakpoint:transformation:item1}
it holds that $F \in C([0,1], \R)$,
\item \label{lem:breakpoint:transformation:item2} 
it holds for all $x , y \in [0,1]$ that 	$
\abs{ F ( x ) - F ( y ) } \leq L \abs{ x - y } 
$, 
\item \label{lem:breakpoint:transformation:item3}
it holds that $G \in \scrL$,
\item \label{lem:breakpoint:transformation:item4}
it holds that $Q(G) = Q(g)$,
\item \label{lem:breakpoint:transformation:item4b}
it holds that $ I \in \{ 1, 2, \dots, Q(G) + 1 \} $,
\item \label{lem:breakpoint:transformation:item5}
it holds for all $j \in \cu{0, 1, \ldots, Q(g) + 1 }$ that $ q_j( G ) = 1 - q_{Q(g) + 1 - j } ( g )$,
\item \label{lem:breakpoint:transformation:item8}
it holds for all $j \in \cu{1, 2, \ldots, Q(g) + 1 }$ that $ A_j( G ) = A_{Q ( g ) + 2 - j } ( g )$,
\item \label{lem:breakpoint:transformation:item6}
it holds that $L \le \abs{\bfa} \le \abs{A_I(G)} = \abs{A_i ( g ) }$,
and
\item \label{lem:breakpoint:transformation:item7}
it holds that $ \bfa A_I( G ) = \bfa A_i( g ) > 0$.
\end{enumerate}
\end{lemma}
\begin{cproof}{lem:breakpoint:transformation}
\Nobs that \cref{eq:def_Q_function,eq:def:slopes} establish 
\cref{lem:breakpoint:transformation:item1,lem:breakpoint:transformation:item2,lem:breakpoint:transformation:item3,lem:breakpoint:transformation:item4,lem:breakpoint:transformation:item4b,lem:breakpoint:transformation:item5,lem:breakpoint:transformation:item6,lem:breakpoint:transformation:item7,lem:breakpoint:transformation:item8}.
\end{cproof}

Our next goal is to prove in \cref{lem:piecewise:approx:1} below that for any piecewise linear function $g \in \scrL$ there exists a piecewise linear $h \in \scrL$ which has at most as many breakpoints as $g$, approximates a given Lipschitz continuous target function $f \in C([0,1], \R)$ as least as well as $g$, and has a Lipschitz constant bounded by the Lipschitz constant of $f$.
To show this we will, roughly speaking, adjust the slopes of the piecewise linear function $g$ one by one and apply induction. Loosely speaking, the following three results, \cref{lem:piecewise:approx:case1,lem:piecewise:approx:case2,lem:piecewise:approx:case3} all consider different cases depending on the slope to be adjusted in each step. The various cases are illustrated in the figures.

\cfclear
\begin{lemma}
\label{lem:piecewise:approx:case1}
Let $ L \in (0, \infty) $, $ f \in C( [0,1], \R ) $ 
satisfy for all $ x, y \in [0,1] $ that 
$
  \abs{ f(x) - f(y) } \leq L \abs{ x - y } 
$, 
let $ g \in \scrL $, $ i \in \{ 1, 2, \ldots, Q(g) + 1 \} $, $ \bfa \in \R $ 
satisfy $ L \le \bfa \le  A_i(g)  $,
assume for all $x \in (q_{i-1} ( g ) , q_i ( g )) $
that
$ g(x) \not= f(x) $,
and let $ \mu \colon \cB( [0,1] ) \to [0, \infty] $ 
be a finite measure \cfadd{def:piece:linear}\cfadd{def:number:of:kinks}\cfadd{def:slopes}\cfload. 
Then there exists $ h \in \scrL $ such that 
$
  \int_0^1 (h(y) - f(y) ) ^2 \, \mu ( \d y ) \leq \int_0^1 (g(y) - f(y) )^2 \, \mu ( \d y ) 
$, 
$
  Q( h ) \leq Q( g )
$,
and
\begin{equation}
\label{eq:lem:piecewise:approx:case1_to_prove}
\textstyle
  \left( 
    Q(g) - Q(h) - 1 
  \right) 
  \big( 
    \sum_{ j = 1 }^{ Q(h) + 1 }
    | 
      A_j( h )
      - A_j(g) \indicator{ \N \backslash \{ i \} }( j ) 
      - \bfa \indicator{ \{ i \} }( j )
    |
  \big)
  \geq 0
  .
\end{equation}
\end{lemma}
\begin{cproof}{lem:piecewise:approx:case1}
Throughout this proof assume without loss of generality that
$ \bfa < A_i( g ) $, 
let $ \fq_0, \fq_1, \ldots, \fq_{ Q(g) + 1 } \in \R $ 
satisfy for all 
$ j \in \cu{ 0, 1, \ldots, Q(g) + 1 } $ 
that $ \fq_j = q_j( g ) $, 
and assume without loss of 
generality\footnote{Otherwise 
the fact that $ f $ and $ g $ 
are continuous ensures that 
$
  \forall \, x \in ( \fq_{ i - 1 }, \fq_i ) \colon f(x) > g(x) 
$
and we can consider 
$
  f \with ( [0,1] \ni x \mapsto -f( 1 - x ) \in \R ) 
$, 
$  
  g \with ( [0,1] \ni x \mapsto - g ( 1 - x ) \in \R ) 
$, and 
$
  i \with Q(g) + 2 - i 
$ 
(cf.~\cref{lem:breakpoint:transformation}).} 
that 
$
  \forall \, x \in ( \fq_{ i - 1 }, \fq_i ) 
  \colon f(x) < g(x) 
$.
\Nobs that 
the fact that $ f $ and $ g $ are continuous 
proves that $ f( \fq_{ i - 1 } ) \leq g( \fq_{ i - 1 } ) $ 
and $ f( \fq_i ) \leq g( \fq_i ) $.
In the following we distinguish between several cases:
\begin{enumerate}[label = (\Roman*)]
\item 
\begin{SCfigure}
	\begin{tikzpicture}[
	thick,
	>=stealth',
	dot/.style = {
		draw,
		fill = white,
		circle,
		inner sep = 0pt,
		minimum size = 4pt
	}
	]
	\coordinate (O) at (0,0);
	\draw[->] (-0.3,0) -- (5,0) coordinate[label = {below:$x$}] (xmax);
	\draw[->] (0,-0.3) -- (0,6) coordinate[label = {right:$f(x)$}] (ymax);
	\draw (0.3, 1.5) -- (2, 2) node[midway, above] {$g$};
	\draw (2,2) -- (4, 5);
	\draw[gray, thin] (2,2) -- (2,0) node[below] {$ \fq_{ Q ( g ) }$};
	\draw[gray, thin] (4,5) -- (4,0) node[below] {$ \fq_{Q ( g ) + 1 } $};
	\draw[blue] (2,2) -- (4,4.5) node[midway, below]  {$h$};
	\draw[red] plot[smooth] coordinates {(0.3,1) (2,0.5) (4,1.8)} node[right] {$f$};
	\end{tikzpicture}
	\caption{Case (I) in \cref{lem:piecewise:approx:case1}. Note that $\fq_i = \fq_{Q ( g ) + 1 } = 1$. The new function $h \in \scrL$ is linear on $[\fq_{i-1} , \fq_{i}] = [ \fq_{Q(g)} , \fq_{Q ( g ) + 1 }]$ with slope $\mathbf{a}$ and agrees with $g$ on $[ 0 , \fq_{Q(g)} ]$.}
	\label{fig:piecewise:approx:1.1}
\end{SCfigure}
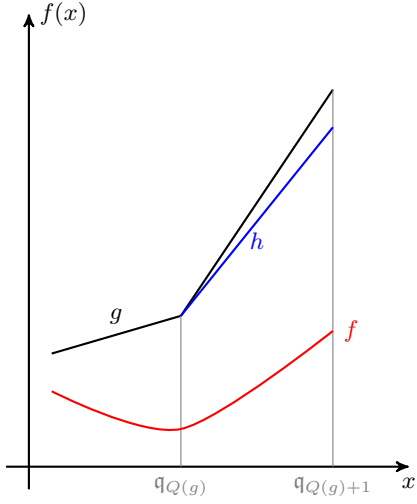
We first prove \cref{eq:lem:piecewise:approx:case1_to_prove} in the case
\begin{equation}
  i = Q(g) + 1 
\end{equation}
(cf.~\cref{fig:piecewise:approx:1.1}).
Let $ h \in \scrL $ satisfy 
for all $ x \in [0, \fq_{ Q(g) } ] $, 
$ y \in [ \fq_{ Q(g) }, 1 ] $ 
that 
$ h(x) = g(x) $ 
and 
$ 
  h(y) = g( \fq_{ Q(g) } ) + \bfa ( y - \fq_{ Q(g) } ) 
$.
\Nobs that for all 
$ j \in \{ 1, 2, \ldots, Q(g) \} $ 
that 
$ 
  h|_{ [ \fq_{ j - 1 }, \fq_j ] } 
$ 
is affine-linear with slope $ A_j( g ) $. 
\Moreover $ h|_{ [ \fq_{ Q(g) }, 1 ] } $ 
is affine-linear with slope $ \bfa $. 
\Moreover 
$ 
  \bigl( 
    ( A_{ Q( g ) }( g ) = \bfa ) 
    \rightarrow 
    ( Q(h) = Q( g ) - 1 < Q( g ) )
  \bigr)
$
and 
\begin{multline}
  \bigl( 
    ( A_{ Q( g ) }( g ) \neq \bfa ) 
    \rightarrow 
    [
      ( Q(h) = Q(g) )
      \wedge 
      ( A_{ Q(g) + 1 }( h ) = \bfa )
\\
      \wedge 
      ( 
        \forall \, j \in \N \cap [1, Q(h)] \colon A_j ( h ) = A_j ( g )
      )
    ]
  \bigr)
  .
\end{multline}
\Moreover the fact that 
$ \forall \, x, y \in [0,1] \colon \abs{ f(x) - f(y) } \leq L \abs{ x - y } $,
the fact that 
$ 
  \forall \, y \in [ \fq_{ Q(g) }, 1 ] \colon 
  g(y) = g( \fq_{ Q( g ) } ) + A_i ( g ) ( y - \fq_{ Q( g ) } ) 
$, 
and the fact that 
$
  \bfa \in [ L, A_{ i }( g ) ) 
$ 
ensure that for all $ y \in [ \fq_{ Q(g) }, 1 ] $ 
we have that
\begin{equation}
\begin{split}
  f(y) \leq 
  f( \fq_{ Q( g ) } ) + L ( y - \fq _{Q ( g ) } ) 
  \leq 
  g( \fq_{ Q( g ) } ) + \bfa ( y - \fq_{ Q( g ) } ) 
  = h( y ) \leq g( y ) .
\end{split}
\end{equation}
This implies that 
$
  \int_0^1 ( h(y) - f(y) )^2 \, \mu( \d y ) 
  \leq 
  \int_0^1 ( g(y) - f(y) )^2 \, \mu( \d y ) 
$,
 which establishes \cref{eq:lem:piecewise:approx:case1_to_prove} in the case 
$ i = Q(g) + 1 $.
\item 
Next we prove \cref{eq:lem:piecewise:approx:case1_to_prove} 
in the case 
\begin{equation}
  ( i < Q(g) + 1 ) \wedge  
  ( A_{ i + 1 }( g ) > A_i( g ) ) \wedge 
  ( 
    g( \fq_{ i + 1 } ) - g( \fq_{ i - 1 } ) 
    \geq \bfa ( \fq_{ i + 1 } - \fq_{ i - 1 } )
  )
\end{equation}
\begin{SCfigure}
	\begin{tikzpicture}[
	thick,
	>=stealth',
	dot/.style = {
		draw,
		fill = white,
		circle,
		inner sep = 0pt,
		minimum size = 4pt
	}
	]
	\coordinate (O) at (0,0);
	\draw[->] (-0.3,0) -- (6,0) coordinate[label = {below:$x$}] (xmax);
	\draw[->] (0,-0.3) -- (0,6.5) coordinate[label = {right:$f(x)$}] (ymax);
	\draw (0.3, 1.5) -- (2, 2);
	\draw (2,2) -- (4, 4) node[midway, above] {$g$};
	\draw (4,4) -- (5,6);
	\draw[gray, thin] (2,2) -- (2,0) node[below] {$\fq_{i - 1 }$};
	\draw[gray, thin] (4,4) -- (4,0) node[below] {$ \fq_i$};
	\draw[gray, thin] (5,6) -- (5,0) node[below] {$ \fq_{i + 1 }$};
	\draw[blue] (3.5,3) -- (4,4) ;
	\draw[blue] (2,2) -- (3.5,3) node[midway, below] {$h$};
	\draw[gray, thin] (3.5,3) -- (3.5,0) node[below] {$u$};
	\draw[red] plot[smooth] coordinates {(0.3,1) (2,0.5) (4,1.8) (5,3)} node[right] {$f$};
	\end{tikzpicture}
	\caption{Case (II) in \cref{lem:piecewise:approx:case1}. The new function $h \in \scrL$ is linear on $[u, \fq_i]$ with slope $A_{i+1} ( g )$,
		linear on $[ \fq_{i-1}, u ]$ with slope $\mathbf{a}$,
		 and agrees with $g$ outside of $[ \fq_{i-1}, \fq_i]$.}
	 \label{fig:piecewise:approx:1.2}
\end{SCfigure}
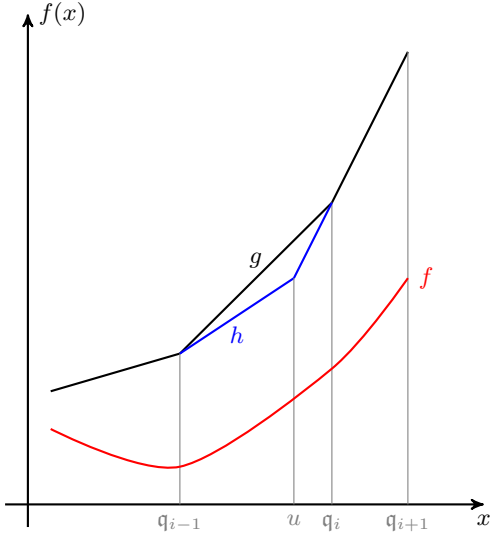
(cf.~\cref{fig:piecewise:approx:1.2}).
\Nobs that the fact that 
$ \bfa \in [ L, A_i( g ) ) $ 
shows that
\begin{equation}
\label{eq:intermediate_value_theorem_assumption_border_B}
\begin{split}
&
  g( \fq_{ i + 1 } ) 
  + A_{ i + 1 }( g ) ( \fq_i - \fq_{ i + 1 } ) 
  - g( \fq_{ i - 1 } ) - \bfa ( \fq_i - \fq_{ i - 1 } ) 
\\ &
  = g( \fq_i ) - g( \fq_{ i - 1 } ) - \bfa ( \fq_i - \fq_{ i - 1 } ) 
  > 0 .
\end{split}
\end{equation} 
\Moreover the fact that $ A_{ i + 1 }( g ) > A_i( g ) $ 
shows that
\begin{equation}
  A_{ i + 1 }( g ) 
  > 
  \bigl[
    \tfrac{ 
      \fq_{ i + 1 } - \fq_i 
    }{ 
      \fq_{ i + 1 } - \fq_{ i - 1 } 
    } 
  \bigr]
  A_{ i + 1 }( g ) 
  + 
  \bigl[ 
    \tfrac{ 
      \fq_i - \fq_{ i - 1 } 
    }{
      \fq_{ i + 1 } - \fq_{ i - 1 } 
    } 
  \bigr] 
  A_i( g ) 
  = 
  \tfrac{ 
    g( \fq_{ i + 1 } ) - g( \fq_{ i - 1 } ) 
  }{ 
    \fq_{ i + 1 } - \fq_{ i - 1 } 
  } 
  .
\end{equation}
\Hence 
$ 
  g( \fq_{ i + 1 } ) 
  + A_{ i + 1 }( g ) ( \fq_{ i - 1 } - \fq_{ i + 1 } ) 
  - g( \fq_{ i - 1 } ) < 0
$.
The intermediate value theorem 
and \cref{eq:intermediate_value_theorem_assumption_border_B}
\hence 
assure that there exists $ u \in ( \fq_{ i - 1 }, \fq_i ) $ 
which satisfies
\begin{equation}
  g( \fq_{ i - 1 } ) 
  + \bfa ( u - \fq_{ i - 1 } ) 
  = 
  g( \fq_{ i + 1 } ) + A_{ i + 1 }( g ) ( u - \fq_{ i + 1 } ) 
  .
\end{equation}
Let $ h \in \scrL $ satisfy 
for all $ x \in [ 0, \fq_{ i - 1 } ] \cup [ \fq_{ i + 1 }, 1 ] $, 
$ y \in [ \fq_{ i - 1 }, u ] $, $ z \in [ u, \fq_{ i + 1 } ] $ 
that $ h(x) = g(x) $, 
$ h(y) = g( \fq_{ i - 1 } ) + \bfa ( y - \fq_{ i - 1 } ) $, 
and $ h(z) = g( \fq_{ i + 1 } ) + A_{ i + 1 }( g ) ( z - \fq_{ i + 1 } ) $. 
\Nobs that
\begin{multline}
  \bigl(
    [ 
      ( i = 1 ) \vee ( A_{ \max\{ i - 1 , 1 \} }( g ) \not= \bfa ) 
    ] 
    \rightarrow 
    [
      ( Q(h) = Q(g) )
\\
      \wedge
      ( 
        \forall \, j \in ( \N \cap [1,Q(h)+1] ) \backslash \{ i \} \colon 
        A_j( h ) = A_j( g )
      )
    ]
  \bigr)
\end{multline}
and 
$
  \bigl(
    [
      ( i > 1 )
      \wedge 
      ( A_{ \max\{ i - 1 , 1 \} }( g ) = \bfa )
    ]
    \rightarrow 
    (
      Q(h) < Q(g)
    )
  \bigr)
$.
\Moreover the fact that $ \bfa \in [ L, A_i( g ) ] $ implies 
for all $ y \in [ \fq_{ i - 1 }, \fq_{ i + 1 } ] $ 
that $ h(y) \leq g(y) $. 
\Moreover 
the fact that 
$
  f( \fq_{ i - 1 } ) \leq g( \fg_{ i - 1 } )
$, 
the fact that 
$ L \leq \bfa \leq A_{ i+1 }( g ) $, 
and 
the fact that 
$ 
  \forall \, x, y \in [0,1] \colon \abs{ f(x) - f(y) } \leq L \abs{ x - y } 
$ 
prove for all $ y \in [ \fq_{ i - 1 }, \fq_{ i + 1 } ] $ that $ f(y) \leq h(y) $. 
\Hence 
$
  \int_0^1 ( h(y) - f(y) )^2 \, \mu( \d y ) 
  \leq \int_0^1 ( g(y) - f(y) )^2 \, \mu( \d y )
$. 
This establishes \cref{eq:lem:piecewise:approx:case1_to_prove} in the case 
$
  ( i < Q(g) + 1 ) \wedge  
  ( A_{ i + 1 }( g ) > A_i( g ) ) \wedge 
  ( 
    g( \fq_{ i + 1 } ) - g( \fq_{ i - 1 } ) 
    \geq \bfa ( \fq_{ i + 1 } - \fq_{ i - 1 } )
  )
$.
\item 
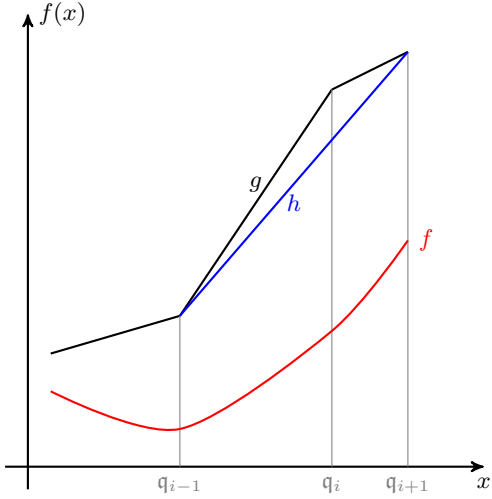
\begin{SCfigure}
	\begin{tikzpicture}[
	thick,
	>=stealth',
	dot/.style = {
		draw,
		fill = white,
		circle,
		inner sep = 0pt,
		minimum size = 4pt
	}
	]
	\coordinate (O) at (0,0);
	\draw[->] (-0.3,0) -- (6,0) coordinate[label = {below:$x$}] (xmax);
	\draw[->] (0,-0.3) -- (0,6) coordinate[label = {right:$f(x)$}] (ymax);
	\draw (0.3, 1.5) -- (2, 2);
	\draw (2,2) -- (4, 5) node[midway, above] {$g$};
	\draw (4,5) -- (5,5.5);
	\draw[gray, thin] (2,2) -- (2,0) node[below] {$ \fq_{i - 1 }$};
	\draw[gray, thin] (4,5) -- (4,0) node[below] {$ \fq_i$};
	\draw[gray, thin] (5, 5.5) -- (5,0) node[below] {$ \fq_{i + 1 }$};
	\draw[blue] (2,2) -- (5,5.5) node[midway, below]  {$h$};
	\draw[red] plot[smooth] coordinates {(0.3,1) (2,0.5) (4,1.8) (5,3)} node[right] {$f$};
	\end{tikzpicture}
	\caption{Case (III) in \cref{lem:piecewise:approx:case1}. The new function $h \in \scrL$ is linear on $[ \fq_{i-1} , \fq_{i+1}]$ with slope $ \frac{g ( \fq_{i+1} ) - g ( \fq_{i-1} ) }{ \fq_{i+1} - \fq_{i-1} } \ge \mathbf{a}$ and agrees with $g$ outside of $[ \fq_{i-1}, \fq_{i + 1 }]$. It thus satisfies $Q(h) < Q(g)$.}
	\label{fig:piecewise:approx:1.3}
\end{SCfigure}
Next we prove \cref{eq:lem:piecewise:approx:case1_to_prove} in the case 
\begin{equation}
  ( i < Q(g) + 1 )
  \wedge
  ( A_{ i + 1 }( g ) < A_i( g ) )
  \wedge 
  ( 
    g( \fq_{ i + 1 } ) - g( \fq_{ i - 1 } ) \geq \bfa ( \fq_{ i + 1 } - \fq_{ i - 1 } ) 
  )
\end{equation}
(cf.~\cref{fig:piecewise:approx:1.3}).
Let $ h \in \scrL $ satisfy 
for all $ x \in [ 0, \fq_{ i - 1 } ] \cup [ \fq_{ i + 1 }, 1 ] $, 
$ y \in [ \fq_{ i - 1 }, \fq_{ i + 1 } ] $ that 
$
  h(x) = g(x)
$ 
and 
$
  h(y) = 
  g( \fq_{ i - 1 } ) 
  + 
  \bigl[
    \tfrac{ 
      g( \fq_{ i + 1 } ) - g( \fq_{ i - 1 } ) 
    }{
      \fq_{ i + 1 } - \fq_{ i - 1 } 
    } 
  \bigr]
  ( y - \fq_{ i - 1 } ) 
$. 
Clearly, we have that $ h \in \scrL $ and $ Q(h) < Q(g) $. 
\Moreover the fact that 
$
  A_{ i + 1 }( g ) < A_i( g )
$ 
shows that
\begin{equation}
  A_{ i + 1 }( g )
  \leq 
  \bigl[
    \tfrac{ \fq_{ i + 1 } - \fq_i }{ 
      \fq_{ i + 1 } - \fq_{ i - 1 } 
    } 
  \bigr]
  A_{ i + 1 }( g ) 
  + 
  \bigl[ 
    \tfrac{ \fq_i - \fq_{ i - 1 } 
    }{ 
      \fq_{ i + 1 } - \fq_{ i - 1 } 
    } 
  \bigr]
  A_i( g ) 
  = 
  \tfrac{
    g( \fq_{ i + 1 } ) - g( \fq_{ i - 1 } ) 
  }{ 
    \fq_{ i + 1 } - \fq_{ i - 1 } 
  } 
  \leq A_i ( g ) .
\end{equation}
\Hence for all $ y \in [ \fq_{ i - 1 }, \fq_{ i + 1 } ] $ 
that $ h(y) \leq g(y) $. 
\Moreover 
the fact that 
$
  f( \fq_{ i - 1 } ) \leq g( \fg_{ i - 1 } )
$, 
the fact that 
$ 
  L 
  \leq 
  \frac{
    g( \fq_{ i + 1 } ) - g( \fg_{ i - 1 } )
  }{
    \fq_{ i + 1 } - \fg_{ i - 1 }
  }
%   \bfa 
$, 
and 
the fact that 
$ 
  \forall \, x, y \in [0,1] \colon \abs{ f(x) - f(y) } \leq L \abs{ x - y } 
$ 
prove for all 
$ y \in [ \fq_{ i - 1 }, \fq_{ i + 1 } ] $ 
that $ f(y) \leq h(y) $.
\Hence
$ 
  \int_0^1 ( h(y) - f(y) )^2 \, \mu( \d y ) 
  \leq \int_0^1 (g(y) - f(y) ) ^2 \, \mu( \d y )
$. 
This establishes \cref{eq:lem:piecewise:approx:case1_to_prove} 
in the case 
$
  ( i < Q(g) + 1 ) \wedge  
  ( A_{ i + 1 }( g ) > A_i( g ) ) \wedge 
  ( 
    g( q_{ i + 1 } ) - g( \fq_{ i - 1 } ) \geq \bfa ( \fq_{ i + 1 } - \fq_{ i - 1 } )
  )
$.
% \end{case}
% 
% 
% 
% \begin{case}
\item
\begin{SCfigure} 
	\begin{tikzpicture}[
	thick,
	>=stealth',
	dot/.style = {
		draw,
		fill = white,
		circle,
		inner sep = 0pt,
		minimum size = 4pt
	}
	]
	\coordinate (O) at (0,0);
	\draw[->] (-0.3,0) -- (6,0) coordinate[label = {below:$x$}] (xmax);
	\draw[->] (0,-0.3) -- (0,6) coordinate[label = {right:$f(x)$}] (ymax);
	\draw (0.3, 1.5) -- (2, 2);
	\draw (2,2) -- (4, 5) node[midway, above] {$g$};
	\draw (4,5) -- (5,4.5);
	\draw[gray, thin] (2,2) -- (2,0) node[below] {$ \fq_{i - 1 }$};
	\draw[gray, thin] (4,5) -- (4,0) node[below] {$ \fq_i$};
	\draw[gray, thin] (5, 4.5) -- (5,0) node[below] {$ \fq_{i + 1 }$};
	\draw[blue] (2,2) -- (4.5,4.75) node[midway, below]  {$h$};
	\draw[gray, thin] (4.5,4.75) -- (4.5, 0) node[below]{$z$};
	\draw[red] plot[smooth] coordinates {(0.3,1) (2,0.5) (4,1.8) (5,3)} node[right] {$f$};
	\end{tikzpicture}
	\caption{Case (IV) in \cref{lem:piecewise:approx:case1}. The new function $h \in \scrL$ is linear on $[ \fq_{i-1} , z]$ with slope $\mathbf{a}$ and agrees with $g$ outside of $[ \fq_{i-1}, z]$.}
	\label{fig:piecewise:approx:1.4}
\end{SCfigure}
Finally, we prove 
\cref{eq:lem:piecewise:approx:case1_to_prove} 
in the case 
\begin{equation}
  ( i < Q(g) + 1 ) 
  \wedge 
  ( 
    g( \fq_{ i + 1 } ) - g( \fq_{ i - 1 } ) 
    < \bfa ( \fq_{ i + 1 } - \fq_{ i - 1 } ) 
  )
\end{equation}
(cf.~\cref{fig:piecewise:approx:1.4}).
\Nobs that the fact that 
$ 
  g( \fq_i ) - g( \fq_{ i - 1 } ) 
  = A_i(g) ( \fq_i - \fq_{ i - 1 } ) 
  > \bfa ( \fq_i - \fq_{ i - 1 } ) 
$ 
and the intermediate value theorem demonstrate 
that there exists 
$ z \in ( \fq_i, \fq_{ i + 1 } ) $ 
which satisfies 
$
  g(z) = g( \fq_{ i - 1 } ) + \bfa ( z - \fq_{ i - 1 } ) 
$. 
Let $ h \in \scrL $ satisfy 
for all $ x \in [ 0, \fq_{ i - 1 } ] \cup [z,1] $, 
$ y \in [ \fq_{ i - 1 }, z] $ that 
$ h(x) = g(x) $ and 
$ 
  h(y) = g( \fq_{ i - 1 } ) + \bfa ( y - \fq_{ i - 1 } ) 
$.
\Nobs that
$ Q(h) = Q(g) $, 
$ A_i( h ) = \bfa $, and 
$ 
  \forall \, j \in \{ 1, 2, \ldots, Q( g ) + 1 \} \backslash \{ i \} \colon 
  A_j(h) = A_j(g)
$. 
\Moreover 
the fact that 
$
  f( \fq_{ i - 1 } ) \leq g( \fg_{ i - 1 } )
$, 
the fact that 
$ L \leq \bfa $, 
and 
the fact that 
$
  \forall \, x, y \in [0,1] \colon | f(x) - f(y) | \leq L | x - y |
$
prove 
for all $ y \in [ \fq_{ i - 1 }, z ] $ 
that $ f(y) \leq h(y) \leq g(y) $. 
\Hence 
$ 
  \int_0^1 ( h(y) - f(y) )^2 \, \mu( \d y ) 
  \leq \int_0^1 ( g(y) - f(y) )^2 \, \mu( \d y )
$. 
This establishes 
\cref{eq:lem:piecewise:approx:case1_to_prove} 
in the case 
$
  ( i < Q(g) + 1 ) 
  \wedge 
  ( 
    g( \fq_{ i + 1 } ) - g( \fq_{ i - 1 } ) 
    < \bfa ( \fq_{ i + 1 } - \fq_{ i - 1 } ) 
  )
$. 
\end{enumerate}
% \end{case}
\end{cproof}

\setcounter{case}{0}
\begin{lemma}
\label{lem:piecewise:approx:case2}
Let $ L \in (0, \infty) $, $ f \in C( [0,1], \R ) $ 
satisfy for all $ x, y \in [0,1] $ that 
$
  \abs{ f(x) - f(y) } \leq L \abs{ x - y } 
$, 
let $ g \in \scrL $, $ i \in \{ 1,  Q(g) + 1 \} $, $ \bfa \in \R $ 
satisfy $ L \le \bfa \le A_i(g) $,
let $ z \in ( q_{ i - 1 }( g ), q_i( g ) ) $
satisfy $ g(z) = f(z) $,
and let 
$ 
  \mu \colon \cB( [0,1] ) \to [0,\infty] 
$ 
be a finite measure 
\cfadd{def:piece:linear}\cfadd{def:number:of:kinks}\cfadd{def:slopes}\cfload. 
Then there exists $ h \in \scrL $ such that 
$
  \int_0^1 ( h(y) - f(y) )^2 \, \mu( \d y ) 
  \leq \int_0^1 ( g(y) - f(y) )^2 \, \mu( \d y ) 
$, 
$
  Q( h ) \leq Q( g )
$,
and
\begin{equation}
\label{eq:piecewise:approx:case2:to_prove}
\textstyle
  \left( 
    Q(g) - Q(h) - 1 
  \right) 
  \big( 
    \sum_{ j = 1 }^{ Q(h) + 1 }
    | 
      A_j( h )
      - A_j(g) \indicator{ \N \backslash \{ i \} }( j ) 
      - \bfa \indicator{ \{ i \} }( j )
    |
  \big)
  \geq 0
  .
\end{equation}
\end{lemma}
\begin{cproof}{lem:piecewise:approx:case2}
Throughout this proof assume without loss of generality that
$ \bfa < A_i( g ) $, 
assume without loss of generality that 
$ i = 1 $ (cf.~\cref{lem:breakpoint:transformation}), 
and let $ \fq_0, \fq_1, \ldots, \fq_{ Q(g) + 1 } \in \R $ 
satisfy for all 
$ j \in \cu{ 0, 1, \ldots, Q(g) + 1 } $ that 
$ \fq_j = q_j( g ) $. 
In the following we distinguish between several cases:  
\begin{enumerate}[label = (\Roman*)]
\item 
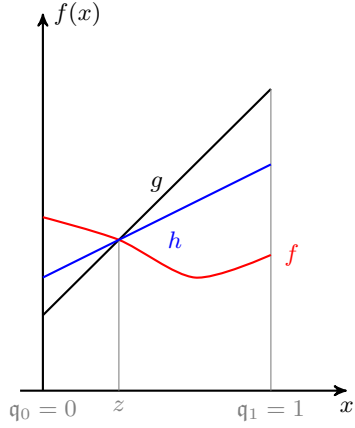
\begin{SCfigure}
	\begin{tikzpicture}[
	thick,
	>=stealth',
	dot/.style = {
		draw,
		fill = white,
		circle,
		inner sep = 0pt,
		minimum size = 4pt
	}
	]
	\coordinate (O) at (0,0);
	\draw[->] (-0.3,0) -- (4,0) coordinate[label = {below:$x$}] (xmax);
	\draw[->] (0,0)  node[below, gray, thin] {$ \fq_{0} = 0$} -- (0,5) coordinate[label = {right:$f(x)$}] (ymax);
	\draw (0, 1) -- (3, 4) node[midway, above] {$g$};
	\draw[gray, thin] (3,4) -- (3,0) node[below] {$ \fq_{1  } = 1$};
	\draw[red] plot[smooth] coordinates {(0,2.3) (1,2) (2,1.5) (3,1.8)} node[right = 1pt] {$f$};
	\draw[gray, thin] (1,2) -- (1,0) node[below] {$z$};
	\draw[blue] (0,3/2) -- (3, 3)  node[midway, below right] {$h$} ;
	\end{tikzpicture}
	\caption{Case (I) in \cref{lem:piecewise:approx:case2}. 
		Here $Q(g)=0$, so $g$ is linear. The new function $h \in \scrL$ is linear on $[\fq_0, \fq_1] = [0,1]$ with slope $\mathbf{a}$ and satisfies $h(z) = g(z) = f(z)$.}
	\label{fig:piecewise:approx:2.1}
\end{SCfigure}
We first prove \cref{eq:piecewise:approx:case2:to_prove} 
in the case 
\begin{equation}
  1 = i = Q ( g ) + 1 
\end{equation}
(cf.~\cref{fig:piecewise:approx:2.1}).
Let $ h \in \scrL $ satisfy for all $ x \in [0,1] $ 
that 
$
  h(x) = g(z) + \bfa ( x - z ) 
$.
\Nobs that 
$ h \in \scrL $, $ Q(h) = 0 = Q(g) $, and $ A_1( h ) = \bfa $.
\Moreover the assumption that $ A_1(g) > \bfa \ge L $ 
and the fact that 
$ \forall \, x, y \in [0,1] \colon \abs{ f(x) - f(y) } \leq L \abs{ x - y } $ 
prove for all 
$ x \in [0,z] $, $ y \in [z,1] $ 
that $ g(x) \leq h(x) \leq f(x) $ and $ f(y) \leq h(y) \leq g(y) $.
This implies  
$ 
  \int_0^1 ( h(y) - f(y) )^2 \, \mu( \d y ) \leq \int_0^1 ( g(y) - f(y) )^2 \, \mu( \d y ) 
$. 
This establishes \cref{eq:piecewise:approx:case2:to_prove} in the case 
$
  1 = i = Q ( g ) + 1 
$. 
\item 
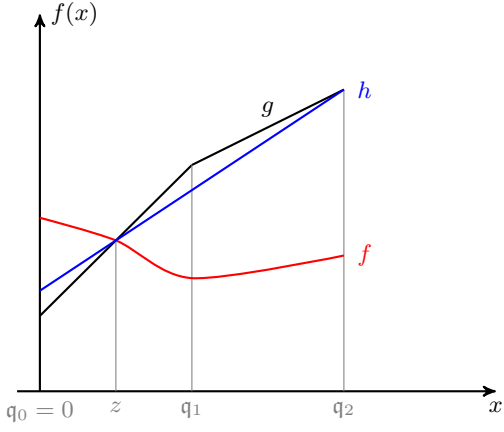
\begin{SCfigure} 
	\begin{tikzpicture}[
	thick,
	>=stealth',
	dot/.style = {
		draw,
		fill = white,
		circle,
		inner sep = 0pt,
		minimum size = 4pt
	}
	]
	\coordinate (O) at (0,0);
	\draw[->] (-0.3,0) -- (6,0) coordinate[label = {below:$x$}] (xmax);
	\draw[->] (0,0)  node[below, gray, thin] {$ \fq_{0} = 0$} -- (0,5) coordinate[label = {right:$f(x)$}] (ymax);
	\draw (0, 1) -- (2, 3);
	\draw (2,3) -- (4, 4) node[midway, above] {$g$};
	\draw[gray, thin] (2,3) -- (2,0) node[below] {$ \fq_{1  }$};
	\draw[gray, thin] (4,4) -- (4,0) node[below] {$ \fq_{2 }$};
	\draw[red] plot[smooth] coordinates {(0,2.3) (1,2) (2,1.5) (4,1.8)} node[right = 1pt] {$f$};
	\draw[gray, thin] (1,2) -- (1,0) node[below] {$z$};
	\draw[blue] (0, 4/3) -- (4,4) node[right=1pt]  {$h$};
	\end{tikzpicture}
	\caption{Case (II) in \cref{lem:piecewise:approx:case2}. The new function $h \in \scrL$ satisfies $h(z) = f(z) = g(z)$, is linear on $[ \fq_0, \fq_2]$ with slope $ \frac{g ( \fq_2) - g(z) }{ \fq_2 - z } \ge \mathbf{a}$, and agrees with $g$ on $[ \fq_{2} , 1]$.}
	\label{fig:piecewise:approx:2.2}
\end{SCfigure}
Next we prove \cref{eq:piecewise:approx:case2:to_prove} 
in the case 
\begin{equation} 
  ( 1 = i < Q( g ) + 1 ) \wedge 
  ( A_2( g ) < A_1( g ) ) \wedge  
  ( g( \fq_2 ) - g( z ) \geq \bfa ( \fq_2 - z ) )
\end{equation}
(cf.~\cref{fig:piecewise:approx:2.2}).
Let $ h \in \scrL $ satisfy for all 
$ x \in [0, \fq_2] $, $ y \in [\fq_2, 1] $ that 
$
  h(x) = g(z) + \bigl[ \frac{ g( \fq_2 ) - g( z ) }{ \fq_2 - z } \bigr] ( x - z) 
$ 
and 
$ 
  h(y) = g(y)
$.
Clearly, we have that $ h \in \scrL $ and $ Q( h ) < Q( g ) $. 
\Moreover the fact that
\begin{equation}
  A_1( g ) 
  > 
  \bigl[
    \tfrac{ \fq_2 - \fq_1 }{ \fq_2 - z } 
  \bigr]
  A_2( g ) 
  + 
  \bigl[ 
    \tfrac{ \fq_1 - z }{ \fq_2 - z } 
  \bigr] 
  A_1( g ) 
  = 
  \tfrac{ g( \fq_2 ) - g( z ) }{ \fq_2 - z } 
  \ge \max\cu{ A_2( g ) , \bfa }
  \ge \max\cu{ A_2( g ) , L }
\end{equation}
and the fact $ \forall \, x, y \in [0,1] \colon \abs{ f(x) - f(y) } \leq L \abs{x-y} $ 
prove that for all $ x \in [0, z] $, $ y \in [ z, \fq_2 ] $ 
we have that 
$
  g(x) \leq h(x) \leq f(x)
$ 
and 
$
  f(y) \leq h(y) \leq g(y)
$. 
\Hence 
$ 
  \int_0^1 ( h(y) - f(y) )^2 \, \mu( \d y ) \leq \int_0^1 ( g(y) - f(y) )^2 \, \mu( \d y ) 
$. 
This establishes \cref{eq:piecewise:approx:case2:to_prove} 
in the case 
$
  ( 1 = i < Q( g ) + 1 ) \wedge 
  ( A_2( g ) < A_1( g ) ) \wedge  
  ( g( \fq_2 ) - g( z ) \geq \bfa ( \fq_2 - z ) )
$.
\item 
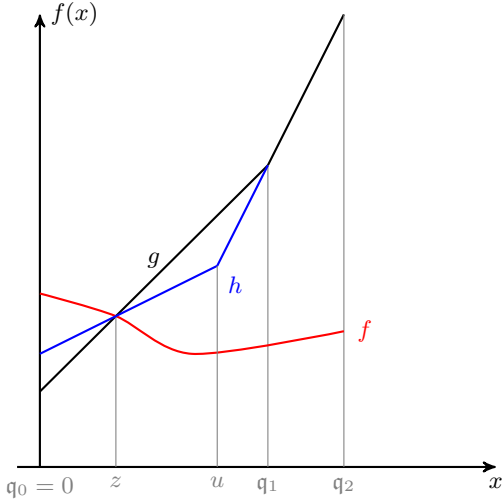
\begin{SCfigure}
	\begin{tikzpicture}[
	thick,
	>=stealth',
	dot/.style = {
		draw,
		fill = white,
		circle,
		inner sep = 0pt,
		minimum size = 4pt
	}
	]
	\coordinate (O) at (0,0);
	\draw[->] (-0.3,0) -- (6,0) coordinate[label = {below:$x$}] (xmax);
	\draw[->] (0,0)  node[below, gray, thin] {$ \fq_{0} = 0$} -- (0,6) coordinate[label = {right:$f(x)$}] (ymax);
	\draw (0, 1) -- (3, 4) node[midway, above] {$g$};
	\draw (3,4) -- (4, 6) ;
	\draw[gray, thin] (3,4) -- (3,0) node[below] {$ \fq_{1  }$};
	\draw[gray, thin] (4,6) -- (4,0) node[below] {$ \fq_{2 }$};
	\draw[red] plot[smooth] coordinates {(0,2.3) (1,2) (2,1.5) (4,1.8)} node[right = 1pt] {$f$};
	\draw[gray, thin] (1,2) -- (1,0) node[below] {$z$};
	\draw[blue] (0,3/2) -- (7/3, 8/3)  node[ below right] {$h$} ;
	\draw[blue] (7/3, 8/3) -- (3,4) ;
	\draw[gray, thin] (7/3, 8/3) -- (7/3, 0) node[below]{$u$} ;
	\end{tikzpicture}
	\caption{Case (III) in \cref{lem:piecewise:approx:case2}. The new function $h \in \scrL$ is linear on $[0 , u]$ with slope $\mathbf{a}$, linear on $[u, \fq_1]$ with slope $A_2(g)$, and agrees with $g$ on $[ \fq_1 , 1]$.}
	\label{fig:piecewise:approx:2.3}
\end{SCfigure}
Next we prove \cref{eq:piecewise:approx:case2:to_prove} 
in the case 
\begin{equation}
  ( 1 = i < Q(g) + 1 ) \wedge 
  ( A_2( g ) > A_1( g ) ) \wedge
  ( g( \fq_2 ) - g( z ) \geq \bfa ( \fq_2 - z ) )
\end{equation}
(cf.~\cref{fig:piecewise:approx:2.3}).
\Nobs that the intermediate value theorem ensures 
that there exists $ u \in [z, \fq_1] $ 
which satisfies 
$ 
  g(z) + \bfa ( u - z ) = g( \fq_2 ) + A_2( g ) ( u - \fq_2 ) 
$. 
Let $ h \in \scrL $ satisfy 
for all $ x \in [0, u] $, $ y \in [u, \fq_1] $,
$z \in [\fq_1 , 1]$ 
that 
$ 
  h(x) = g(z) + \bfa ( x - z ) 
$,
$h(y) = h(u ) + A_2(g) ( y - u )$,
and 
$
  h( z ) = g ( z ) 
$.
\Nobs that $ Q(h) = Q(g) $, $ A_1(h) = \bfa $, and 
$ 
  \forall \, j \in \cu{ 2, 3, \ldots, Q(g) + 1 } \colon A_j( h ) = A_j( g ) 
$. 
\Moreover for all $ x \in [0, z] $, $ y \in [z, \fq_1 ] $ 
it holds that 
$ 
  g(x) \le h(x) \le f(x)
$ 
and 
$
  f(y) \le h(y) \le g(y) 
$. 
\Hence  
$
  \int_0^1 ( h(y) - f(y) )^2 \, \mu( \d y ) 
  \leq \int_0^1 (g(y) - f(y) ) ^2 \, \mu( \d y ) 
$. 
This establishes \cref{eq:piecewise:approx:case2:to_prove} in the case 
$
  ( 1 = i < Q(g) + 1 ) \wedge 
  ( A_2( g ) > A_1( g ) ) \wedge
  ( g( \fq_2 ) - g( z ) \geq \bfa ( \fq_2 - z ) )
$. 
\item 
% \label{lem:piecewise:approx:case:2:case3} 
\begin{SCfigure}  
	\begin{tikzpicture}[
	thick,
	>=stealth',
	dot/.style = {
		draw,
		fill = white,
		circle,
		inner sep = 0pt,
		minimum size = 4pt
	}
	]
	\coordinate (O) at (0,0);
	\draw[->] (-0.3,0) -- (6,0) coordinate[label = {below:$x$}] (xmax);
	\draw[->] (0,0)  node[below, gray, thin] {$ \fq_{0} = 0$} -- (0,5) coordinate[label = {right:$f(x)$}] (ymax);
	\draw (0, 1) -- (3, 4) node[midway, above] {$g$};
	\draw (3,4) -- (4, 2);
	\draw[gray, thin] (3,4) -- (3,0) node[below] {$ \fq_{1  }$};
	\draw[gray, thin] (4, 2) -- (4,0) node[below] {$ \fq_{2 }$};
	\draw[red] plot[smooth] coordinates {(0,2.3) (1,2) (2,1.5) (4,1.8)} node[right = 1pt] {$f$};
	\draw[gray, thin] (1,2) -- (1,0) node[below] {$z$};
	\draw[blue] (0,3/2) -- (17/5, 16/5)  node[midway, below right] {$h$} ;
	\draw[gray, thin] (17/5, 16/5) -- (17/5, 0) node[below]{$u$} ;
	\end{tikzpicture}
	\caption{Case (IV) in \cref{lem:piecewise:approx:case2}. The new function $h \in \scrL$ is linear on $[0 , u]$ with slope $\mathbf{a}$ and agrees with $g$ on $[u , 1]$.}
	\label{fig:piecewise:approx:2.4}
\end{SCfigure}
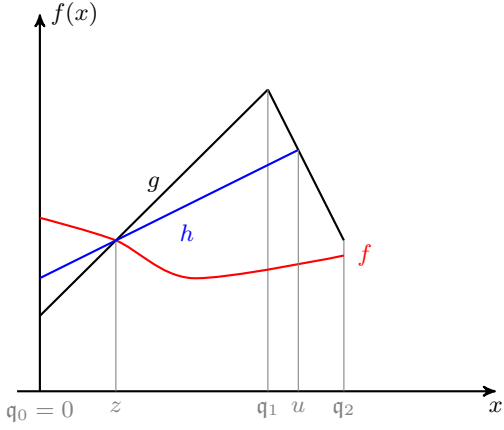
Finally, we prove \cref{eq:piecewise:approx:case2:to_prove} 
in the case 
\begin{equation}
  (
    1 = i < Q( g ) + 1
  )
  \wedge 
  (
    g( \fq_2 ) - g( z ) < \bfa ( \fq_2  - z )
  )
\end{equation}
(cf.~\cref{fig:piecewise:approx:2.4}).
\Nobs that the fact that 
$
  g( \fq_1 ) - g( z ) 
  = A_1( g ) ( \fq_1 - z ) > \bfa ( \fq_1 - z ) 
$ 
and the intermediate value theorem demonstrate 
that there exists $ u \in ( \fq_1, \fq_2 ) $ 
which satisfies 
$ 
  g(u) = g(z) + \bfa ( u - z )
$.
Let $ h \in \scrL $ satisfy for all 
$ x \in [0, u] $, $ y \in [u,1] $ that 
$
  h(x) = g(z) + \bfa ( x - z)
$ 
and 
$
  h(y) = g(y) 
$.
\Nobs that 
$ \bigl( ( A_2(g) = \bfa ) \rightarrow ( Q(h) < Q(g) ) \bigr) $
and 
\begin{equation}
  \bigl( 
    ( A_2(g) \neq \bfa ) 
    \rightarrow 
    [ 
      ( Q(h) = Q(g) ) \wedge ( A_1(h) = \bfa ) \wedge 
      (
        \forall \, j \in \N \cap (1,Q(g) + 1] \colon A_j(h) = A_j(g)
      ) 
    ] 
  \bigr)
  .
\end{equation}
\Moreover the assumption that 
$ A_1(g) > \bfa $ and 
the fact that 
$ 
  \forall \, x, y \in [0,1] \colon \abs{ f(x) - f(y) } 
  \leq L \abs{ x - y } 
$ 
prove that for all $ x \in [0, z] $, $ y \in [z, u] $ 
we have that 
$ g(x) \leq h(x) \leq f(x) $ and $ f(y) \leq h(y) \leq g(y) $.
This implies 
$
  \int_0^1 ( h(y) - f(y) )^2 \, \mu( \d y ) 
  \leq 
  \int_0^1 ( g(y) - f(y) )^2 \, \mu( \d y )
$. 
This establishes \cref{eq:piecewise:approx:case2:to_prove} 
in the case 
$ 
  (
    1 = i < Q( g ) + 1
  )
  \wedge 
  (
    g( \fq_2 ) - g( z ) < \bfa ( \fq_2  - z )
  )
$.
\end{enumerate}
\end{cproof}

\setcounter{case}{0}
\begin{lemma}
\label{lem:piecewise:approx:case3}
Let $ L \in (0, \infty) $, $ f \in C( [0,1], \R ) $ 
satisfy for all $ x, y \in [0,1] $ that 
$
  \abs{ f(x) - f(y) } \leq L \abs{ x - y } 
$, 
let $ g \in \scrL $, $ i \in \N \cap (1, Q(g)] $, $ \bfa \in \R $ 
satisfy $ L \le \bfa \le  A_i(g) $,
let $ z \in ( q_{ i - 1 }( g ), q_i( g ) ) $
satisfy
$ g(z) = f(z) $,
and let $ \mu \colon \cB( [0,1] ) \to [0, \infty] $ 
be a finite measure 
\cfadd{def:piece:linear}\cfadd{def:number:of:kinks}\cfadd{def:slopes}\cfload. 
Then there exists $ h \in \scrL $ such that 
$
  \int_0^1 ( h(y) - f(y) )^2 \, \mu( \d y ) 
  \leq \int_0^1 (g(y) - f(y) )^2 \, \mu( \d y ) 
$, 
$
  Q( h ) \leq Q( g )
$,
and
\begin{equation}
\label{eq:lem:piecewise:approx:case3:to_prove}
\textstyle
  \left( 
    Q(g) - Q(h) - 1 
  \right) 
  \big( 
    \sum_{ j = 1 }^{ Q(h) + 1 }
    | 
      A_j( h )
      - A_j(g) \indicator{ \N \backslash \{ i \} }( j ) 
      - \bfa \indicator{ \{ i \} }( j )
    |
  \big)
  \geq 0
  .
\end{equation}
\end{lemma}
\begin{cproof}{lem:piecewise:approx:case3}
Throughout this proof assume 
without loss of generality that $ \bfa < A_i( g ) $
and let 
$ \fq_0, \fq_1, \ldots, \fq_{ Q(g) + 1 } \in \R $ 
satisfy 
for all $ j \in \cu{ 0, 1, \ldots, Q(g) + 1 } $ 
that $ \fq_j = q_j( g ) $. 
In the following we distinguish between several 
cases: 
\begin{enumerate}[label = (\Roman*)]
\item 
\begin{SCfigure}
	\begin{tikzpicture}[
	thick,
	>=stealth',
	dot/.style = {
		draw,
		fill = white,
		circle,
		inner sep = 0pt,
		minimum size = 4pt
	}
	]
	\coordinate (O) at (0,0);
	\draw[->] (-0.3,0) -- (6,0) coordinate[label = {below:$x$}] (xmax);
	\draw[->] (0,-0.3) -- (0,6.5) coordinate[label = {right:$f(x)$}] (ymax);
	\draw (1, 0) -- (2, 2)  node[midway, above] {$g$};
	\draw (2,2) -- (4, 4);
	\draw (4,4) -- (5,6);
	\draw[gray, thin] (2,2) -- (2,0) node[below] {$ \fq_{i - 1 }$};
	\draw[gray, thin] (4,4) -- (4,0) node[below] {$ \fq_i$};
	\draw[gray, thin] (5,6) -- (5,0) node[below] {$ \fq_{i + 1 }$};
	\draw[gray, thin] (1, 0) -- (1,0) node[below] { $ \fq_{i-2}$};
	\draw[gray, thin] (3,3) -- (3,0) node[below] {$z$};
	\draw[blue] (2,2) -- (7/3, 8/3) node[left] {$h$};
	\draw[blue] (7/3, 8/3) -- (11/3, 10/3) ;
	\draw[blue] (11/3, 10/3) -- (4,4) ;
	\draw[gray , thin] (7/3, 8/3) -- (7/3, 0) node[above right]{$u$};
	\draw[gray , thin] (11/3, 10/3) -- (11/3, 0) node[above left]{$v$};
	\draw[red] plot[smooth] coordinates {(1, 3.3) (2,3.2) (3,3) (4,2.2) (5,2)} node[right] {$f$};
	\end{tikzpicture}
	\caption{Case (I) in \cref{lem:piecewise:approx:case3}. The new function $h \in \scrL$ is linear on $[ \fq_{i-1} , u]$ with slope $A_{i-1} ( g )$,
		linear on $[u,v]$ with slope $\mathbf{a}$,
		linear on $[v, \fq_i]$ with slope $A_{i+1} ( g )$,
		and agrees with $g$ outside of $[ \fq_{i-1}, \fq_i]$.}
	\label{fig:piecewise:approx:3.1}
\end{SCfigure}
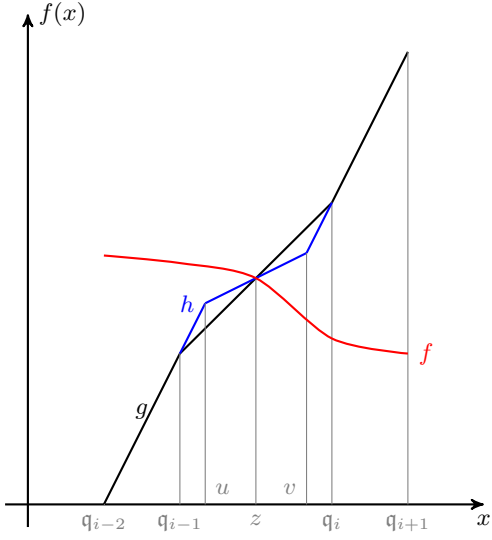
We first prove \cref{eq:lem:piecewise:approx:case3:to_prove} in the case
\begin{equation}
  A_i ( g ) 
  < \min\{ A_{ i - 1 }( g ), A_{ i + 1 }( g ) \}
\end{equation}
(cf.~\cref{fig:piecewise:approx:3.1}).
\Nobs that the fact that 
$ \bfa < A_i( g ) $ and 
the intermediate value theorem assure 
that there exist 
$ u \in ( \fq_{ i - 1 }, z ) $, $ v \in ( z, \fq_i ) $ 
which satisfy 
$
  g( \fq_{ i - 2 } ) + A_{ i - 1 }( g ) ( u - \fq_{ i - 2 } ) 
  = g(z) + \bfa ( u - z)
$ 
and 
$
  g( \fq_{ i + 1 } ) + A_{ i + 1 }( g ) ( v - \fq_{ i + 1 } ) 
  = g(z) + \bfa ( v - z )
$. 
Let $ h \in \scrL $ satisfy 
for all $ x \in [ 0, \fq_{ i - 1 } ] \cup [ \fq_{ i + 1 }, 1] $, 
$ y_1 \in [ \fq_{ i - 1 }, u ] $, 
$ y_2 \in [u , v ] $,  
$ y_3 \in [v, \fq_i] $ 
that 
$ h(x) = g(x) $, 
$ 
  h(y_1) = g( \fq_{i-1} ) + A_{ i - 1 }( g ) ( y_1 - \fq_{ i - 1 } ) 
$, 
$ 
  h( y_2 ) = g(z) + \bfa ( y_2 - z ) 
$, 
and 
$ 
  h( y_3 ) = g( \fq_i ) + A_{ i + 1 }( g ) ( y_3 - \fq_i ) 
$. 
\Nobs that $ Q(h) = Q(g) $, 
$ A_i( h ) = \bfa $, 
and 
$
  \forall \, j \in \{ 1, 2, \ldots, Q( g ) + 1 \} \backslash \{ i \} \colon 
  A_j( h ) = A_j( g ) 
$.  
\Moreover the assumption that $ A_i(g) > \bfa $ 
and 
the fact that 
$ 
  \forall \, x, y \in [0,1] \colon \abs{ f(x) - f(y) } \leq L \abs{ x - y } 
$ 
demonstrate for all $ x \in [ \fq_{ i - 1 } , z ] $, 
$ y \in [ z, \fq_i ] $ that 
$ g(x) \leq h(x) \leq f(x) $ 
and 
$ f(y) \leq h(y) \leq g(y) $. 
\Hence 
$ 
  \int_0^1 ( h(y) - f(y) )^2 \, \mu( \d y ) 
  \leq \int_0^1 (g(y) - f(y) ) ^2 \, \mu( \d y )
$. 
This establishes \cref{eq:lem:piecewise:approx:case3:to_prove} in the case 
$
  A_i ( g ) 
  < \min\{ A_{ i - 1 }( g ), A_{ i + 1 }( g ) \}
$. 
\item 
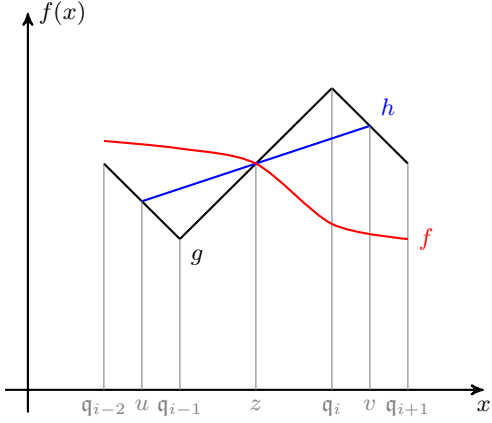
\begin{SCfigure} 
	\begin{tikzpicture}[
	thick,
	>=stealth',
	dot/.style = {
		draw,
		fill = white,
		circle,
		inner sep = 0pt,
		minimum size = 4pt
	}
	]
	\coordinate (O) at (0,0);
	\draw[->] (-0.3,0) -- (6,0) coordinate[label = {below:$x$}] (xmax);
	\draw[->] (0,-0.3) -- (0,5) coordinate[label = {right:$f(x)$}] (ymax);
	\draw (1, 3) -- (2, 2)  node[below right] {$g$};
	\draw (2,2) -- (4, 4);
	\draw (4,4) -- (5,3);
	\draw[gray, thin] (2,2) -- (2,0) node[below] {$ \fq_{i - 1 }$};
	\draw[gray, thin] (4,4) -- (4,0) node[below] {$ \fq_i$};
	\draw[gray, thin] (5,3) -- (5,0) node[below] {$ \fq_{i + 1 }$};
	\draw[gray, thin] (1, 3) -- (1,0) node[below] { $ \fq_{i-2}$};
	\draw[gray, thin] (3,3) -- (3,0) node[below] {$z$};
	\draw[blue] (1.5, 2.5) -- (4.5, 3.5) node[above right] {$h$};
	\draw[gray , thin] (1.5, 2.5) -- (1.5, 0) node[below]{$u$};
	\draw[gray , thin] (4.5, 3.5) -- (4.5, 0) node[below]{$v$};
	\draw[red] plot[smooth] coordinates {(1, 3.3) (2,3.2) (3,3) (4,2.2) (5,2)} node[right] {$f$};
	\end{tikzpicture}
	\caption{Case (II) in \cref{lem:piecewise:approx:case3}. The new function $h \in \scrL$ is linear on $[u,v]$ with slope $\mathbf{a}$
		and agrees with $g$ outside of $[u, v]$.}
	\label{fig:piecewise:approx:3.2}
\end{SCfigure}
Next we prove \cref{eq:lem:piecewise:approx:case3:to_prove} in the case
\begin{equation}
  \max\bigl\{ 
    \tfrac{ 
      g( \fq_{ i + 1 } ) - g( z ) 
    }{ 
      \fq_{ i + 1 } - z 
    } , 
    \tfrac{ 
      g( \fq_{ i - 2 } ) - g( z ) 
    }{
      \fq_{ i - 2 } - z 
    }
  \bigr\} < \bfa 
\end{equation}
(cf.~\cref{fig:piecewise:approx:3.2}).
\Nobs that the fact that $ A_i( g ) > \bfa $ 
proves that 
$
  \max\{ A_{ i - 1 }( g ), A_{ i + 1 }( g ) \} < \bfa 
$. 
\Moreover the fact that $ A_i( g ) > \bfa $ 
and the intermediate value theorem assure 
that there exist 
$ u \in ( \fq_{ i - 2 }, \fq_{ i - 1 } ) $, 
$ v \in ( \fq_i, \fq_{ i + 1 } ) $ 
which satisfy
$
  \frac{ g(u) - g(z) }{ u - z } 
  = \frac{ g(v) - g(z) }{ v - z } 
  = \bfa 
$. 
Let $ h \in \scrL $ satisfy 
for all $ x \in [0, u] \cup [v, 1] $, $ y \in [u, v] $ 
that $ h(x) = g(x) $ and $ h(y) = g(z) + \bfa ( y - z ) $. 
\Nobs that the fact that 
$ A_{ i - 1 }( g ) \not= \bfa $ 
and the fact that $ A_{ i + 1 }( g ) \not= \bfa $ 
show that 
$
  Q(h) = Q(g)
$, 
$ A_i( h ) = \bfa $, 
and 
$ 
  \forall \, j \in \{ 1, 2, \ldots, Q( g ) + 1 \} \backslash \{ i \} \colon A_j(h) = A_j(g) 
$. 
\Moreover the assumption that $ A_i(g) > \bfa $ and 
the fact that 
$ 
  \forall \, x,y \in [0,1] \colon \abs{ f(x) - f(y) } \leq L \abs{ x - y } 
$ 
demonstrate for all $ x \in [u, z] $, $ y \in [z, v] $ that 
$ 
  g(x) \leq h(x) \leq f(x) 
$ 
and 
$
  f(y) \leq h(y) \leq g(y)
$. 
\Hence 
$
  \int_0^1 ( h(y) - f(y) )^2 \, \mu( \d y ) 
  \leq \int_0^1 ( g(y) - f(y) )^2 \, \mu( \d y ) 
$. 
This establishes \cref{eq:lem:piecewise:approx:case3:to_prove}
in the case 
$
  \max\bigl\{ 
    \frac{ 
      g( \fq_{ i + 1 } ) - g( z ) 
    }{ 
      \fq_{ i + 1 } - z 
    } , 
    \frac{ 
      g( \fq_{ i - 2 } ) - g( z ) 
    }{
      \fq_{ i - 2 } - z 
    }
  \bigr\} < \bfa 
$.
\item 
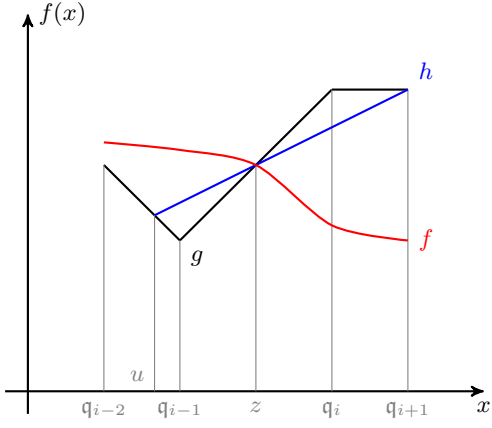
\begin{SCfigure}
	\begin{tikzpicture}[
	thick,
	>=stealth',
	dot/.style = {
		draw,
		fill = white,
		circle,
		inner sep = 0pt,
		minimum size = 4pt
	}
	]
	\coordinate (O) at (0,0);
	\draw[->] (-0.3,0) -- (6,0) coordinate[label = {below:$x$}] (xmax);
	\draw[->] (0,-0.3) -- (0,5) coordinate[label = {right:$f(x)$}] (ymax);
	\draw (1, 3) -- (2, 2)  node[below right] {$g$};
	\draw (2,2) -- (4, 4);
	\draw (4,4) -- (5,4);
	\draw[gray, thin] (2,2) -- (2,0) node[below] {$ \fq_{i - 1 }$};
	\draw[gray, thin] (4,4) -- (4,0) node[below] {$ \fq_i$};
	\draw[gray, thin] (5,4) -- (5,0) node[below] {$ \fq_{i + 1 }$};
	\draw[gray, thin] (1, 3) -- (1,0) node[below] { $ \fq_{i-2}$};
	\draw[gray, thin] (3,3) -- (3,0) node[below] {$z$};
	\draw[blue] (5/3 , 7/3) -- (5,4) node[above right] {$h$};
	\draw[gray , thin] (5/3, 7/3) -- (5/3 , 0) node[above left]{$u$};
	\draw[red] plot[smooth] coordinates {(1, 3.3) (2,3.2) (3,3) (4,2.2) (5,2)} node[right] {$f$};
	\end{tikzpicture}
	\caption{Case (III) in \cref{lem:piecewise:approx:case3}. The new function $h \in \scrL$ is linear on $[u , \fq_{i+1}]$ with slope $\frac{g ( \fq_{i+1} ) - g(z)}{ \fq_{i+1} - z } \ge \mathbf{a}$
		and agrees with $g$ outside of $[u, \fq_{i+1}]$.}
	\label{fig:piecewise:approx:3.3}
\end{SCfigure}
Next we prove \cref{eq:lem:piecewise:approx:case3:to_prove} 
in the case 
\begin{equation}
  \bigl( 
    A_i( g ) > \max\{ A_{ i - 1 }( g ), A_{ i + 1 }( g ) \} 
  \bigr)
  \wedge 
  \bigr(
    \max\bigl\{ 
      \tfrac{ 
        g( \fq_{ i + 1 } ) - g( z ) 
      }{
        \fq_{ i + 1 } - z 
      } , 
      \tfrac{ 
        g( \fq_{ i - 2 } ) - g( z ) 
      }{
        \fq_{ i - 2 } - z 
      } 
    \bigr\} 
    \geq \bfa
  \bigl)
\end{equation}
(cf.~\cref{fig:piecewise:approx:3.3}).
In the following 
we assume without loss of generality that 
$
  \frac{ g( \fq_{ i + 1 } ) - g( z ) }{ \fq_{ i + 1 } - z } 
  \geq 
  \frac{ 
    g( \fq_{ i - 2 } ) - g( z )
  }{
    \fq_{ i - 2 } - z
  }
$ (cf.~\cref{lem:breakpoint:transformation}). 
\Nobs that the fact that 
$ 
  A_i( g ) > \max\{ A_{ i - 1 }( g ), A_{ i + 1 }( g ) \} 
$ 
shows that
\begin{equation}
  A_i( g ) 
  = 
  \frac{ g( \fq_{ i - 1 } ) - g( z ) }{ \fq_{ i - 1 } - z } 
  >
%   \geq 
  \frac{ g( \fq_{ i + 1 } ) - g( z ) }{ \fq_{ i + 1 } - z } 
  \geq 
  \frac{ g( \fq_{ i - 2 } ) - g( z ) }{ \fq_{ i - 2 } - z }
  .
\end{equation}
The intermediate value theorem \hence proves 
that there exists 
$
  u \in [ \fq_{ i - 2 }, \fq_{ i - 1 } ) 
$ 
which satisfies 
$
  g(u) = 
  g(z) 
  + 
  \bigl[ 
    \frac{ g( \fq_{ i + 1 } ) - g( z ) }{ \fq_{ i + 1 } - z } 
  \bigr] ( u - z )
$.  
Let $ h \in \scrL $ satisfy 
for all $ x \in [0, u] \cup [ \fq_{ i + 1 }, 1] $, $ y \in [u, \fq_{ i + 1 } ]$ 
that $ h(x) = g(x) $ and 
$
  h(y) = 
  g(z) + \bigl[ \frac{ g( \fq_{ i + 1 } ) - g( z ) }{ \fq_{ i + 1 } - z } \bigr] (y-z) 
$. 
\Nobs that $ Q(h) < Q(g) $.
\Moreover the assumption that $ A_i( g ) > \bfa $ 
and the fact that 
$ 
  \forall \, x, y \in [0,1] \colon \abs{ f(x) - f(y) } \leq L \abs{ x - y } 
$ 
demonstrate for all 
$ x \in [u, z] $, $ y \in [z, \fq_{ i + 1 } ] $ 
that 
$
  g(x) \leq h(x) \leq f(x)
$ 
and 
$
  f(y) \leq h(y) \leq g(y) 
$. 
\Hence 
$
  \int_0^1 ( h(y) - f(y) )^2 \, \mu( \d y ) 
  \leq \int_0^1 ( g(y) - f(y) )^2 \, \mu( \d y ) 
$. 
This establishes \cref{eq:lem:piecewise:approx:case3:to_prove} in the case 
$
  \bigl( 
    A_i( g ) > \max\{ A_{ i - 1 }( g ), A_{ i + 1 }( g ) \} 
  \bigr)
  \wedge 
  \bigr(
    \max\bigl\{ 
      \tfrac{ 
        g( \fq_{ i + 1 } ) - g( z ) 
      }{
        \fq_{ i + 1 } - z 
      } , 
      \tfrac{ 
        g( \fq_{ i - 2 } ) - g( z ) 
      }{
        \fq_{ i - 2 } - z 
      } 
    \bigr\} 
    \geq \bfa
  \bigl)
$.
\item 
\begin{SCfigure} 
	\begin{tikzpicture}[
	thick,
	>=stealth',
	dot/.style = {
		draw,
		fill = white,
		circle,
		inner sep = 0pt,
		minimum size = 4pt
	}
	]
	\coordinate (O) at (0,0);
	\draw[->] (-0.3,0) -- (6,0) coordinate[label = {below:$x$}] (xmax);
	\draw[->] (0,-0.3) -- (0,6.5) coordinate[label = {right:$f(x)$}] (ymax);
	\draw (1, 2) -- (2, 2)  node[below right] {$g$};
	\draw (2,2) -- (4, 4);
	\draw (4,4) -- (5,6);
	\draw[gray, thin] (2,2) -- (2,0) node[below] {$ \fq_{i - 1 }$};
	\draw[gray, thin] (4,4) -- (4,0) node[below] {$ \fq_i$};
	\draw[gray, thin] (5,6) -- (5,0) node[below] {$ \fq_{i + 1 }$};
	\draw[gray, thin] (1, 2) -- (1,0) node[below] { $ \fq_{i-2}$};
	\draw[gray, thin] (3,3) -- (3,0) node[below] {$z$};
	\draw[blue] (1,2) -- (11/3, 10/3) node[midway, above] {$h$};
	\draw[blue] (11/3, 10/3) -- (4,4);
	\draw[gray , thin] (11/3, 10/3) -- (11/3 , 0) node[above left]{$u$};
	\draw[red] plot[smooth] coordinates {(1, 3.3) (2,3.2) (3,3) (4,2.2) (5,2)} node[right] {$f$};
	\end{tikzpicture}
	\caption{Case (IV) in \cref{lem:piecewise:approx:case3}. The new function $h \in \scrL$ is linear on $[ \fq_{i-2} , u ]$ with slope $\frac{g ( z ) - g ( \fq_{i-2} ) }{z - \fq_{i-2}} \ge \mathbf{a}$,
		linear on $[u, \fq_i]$ with slope $A_{i+1} ( g )$,
		and agrees with $g$ outside of $[ \fq_{i-2} , \fq_i]$.}
	\label{fig:piecewise:approx:3.4}
\end{SCfigure}
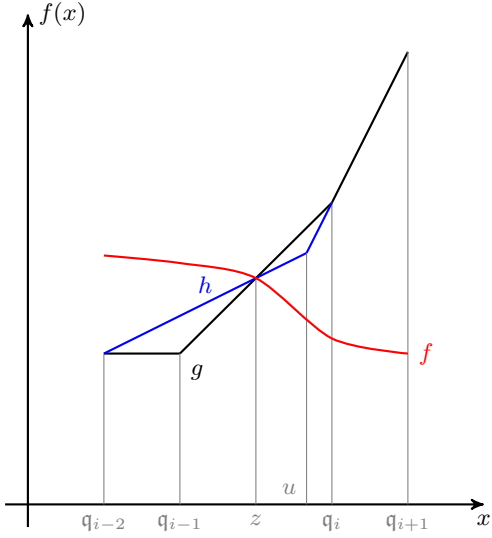
Next we prove \cref{eq:lem:piecewise:approx:case3:to_prove} in the case
\begin{multline}
  \bigl(
    \max\{ 
      A_{ i - 1 }( g ), A_{ i + 1 }( g ) 
    \} 
    > A_i( g ) 
    > \min\{ A_{ i - 1 }( g ), A_{ i + 1 }( g ) \} 
  \bigr)
\\
  \wedge
  \bigl(
    \min\bigl\{ 
      \tfrac{ g( \fq_{ i + 1 } ) - g( z ) }{ \fq_{ i + 1 } - z } , 
      \tfrac{ g( \fq_{ i - 2 } ) - g( z ) }{ \fq_{ i - 2 } - z } 
    \bigr\} \geq \bfa
  \bigr)
\end{multline}
(cf.~\cref{fig:piecewise:approx:3.4}).
In the following we assume without loss of generality that 
$ A_{ i - 1 }( g ) < A_i( g ) < A_{ i + 1 }( g ) $ 
(cf.~\cref{lem:breakpoint:transformation}). 
\Nobs that the fact 
$
  \frac{ g( \fq_{ i - 2 } ) - g( z ) }{ \fq_{ i - 2 } - z } < A_i( g ) 
$ 
and the intermediate value theorem imply that there exists 
$
  u \in (z, \fq_i)
$ 
which satisfies 
$
  g(z) 
  + 
  \bigl[ 
    \frac{ g( \fq_{ i - 2 } ) - g( z ) }{ \fq_{ i - 2 } - z } 
  \bigr] ( u - z ) 
  = g(q_i ) + A_{ i + 1 }( g ) ( u - \fq_i )
$. 
Let $ h \in \scrL $ satisfy 
for all $ x \in [ 0, \fq_{ i - 2 } ] \cup [ \fq_i, 1 ] $, 
$ y_1 \in [ \fq_{ i - 2 }, u ] $, 
$ y_2 \in [ u, \fq_i ] $ 
that 
$ h(x) = g(x) $, 
$ 
  h( y_1 ) = 
  g( \fq_{ i - 2 } ) 
  + 
  \bigl[ \frac{ g( \fq_{ i - 2 } ) - g( z ) }{ \fq_{ i - 2 } - z } 
  \bigr] ( y_1 - \fq_{ i - 2 } ) 
$, 
and 
$
  h( y_2 ) = g( \fq_i ) + A_{ i + 1 }( g ) ( y_2 - \fq_i ) 
$. 
\Nobs that $ Q(h) < Q(g) $ 
and 
$ 
  \int_0^1 ( h(y) - f(y) )^2 \, \mu( \d y ) \leq \int_0^1 ( g(y) - f(y) )^2 \, \mu( \d y ) 
$. 
This establishes \cref{eq:lem:piecewise:approx:case3:to_prove} in the case 
$
  \bigl(
    \max\{ 
      A_{ i - 1 }( g ), A_{ i + 1 }( g ) 
    \} 
    > A_i( g ) 
    > \min\{ A_{ i - 1 }( g ), A_{ i + 1 }( g ) \} 
  \bigr)
  \wedge
  \bigl(
    \min\bigl\{ 
      \frac{ g( \fq_{ i + 1 } ) - g( z ) }{ \fq_{ i + 1 } - z } , 
      \frac{ g( \fq_{ i - 2 } ) - g( z ) }{ \fq_{ i - 2 } - z } 
    \bigr\} \geq \bfa
  \bigr)
$.
\item 
\begin{SCfigure} 
	\begin{tikzpicture}[
	thick,
	>=stealth',
	dot/.style = {
		draw,
		fill = white,
		circle,
		inner sep = 0pt,
		minimum size = 4pt
	}
	]
	\coordinate (O) at (0,0);
	\draw[->] (-0.3,0) -- (6,0) coordinate[label = {below:$x$}] (xmax);
	\draw[->] (0,-0.3) -- (0,6.5) coordinate[label = {right:$f(x)$}] (ymax);
	\draw (1, 3) -- (2, 2)  node[below right] {$g$};
	\draw (2,2) -- (4, 4);
	\draw (4,4) -- (5,6);
	\draw[gray, thin] (2,2) -- (2,0) node[below] {$ \fq_{i - 1 }$};
	\draw[gray, thin] (4,4) -- (4,0) node[below] {$ \fq_i$};
	\draw[gray, thin] (5,6) -- (5,0) node[below] {$ \fq_{i + 1 }$};
	\draw[gray, thin] (1, 3) -- (1,0) node[below] { $ \fq_{i-2}$};
	\draw[gray, thin] (3,3) -- (3,0) node[below] {$z$};
	\draw[blue] (5/3, 7/3) node[above right =2pt]{$h$} -- (11/3, 10/3) ;
	\draw[blue] (11/3, 10/3) -- (4,4);
	\draw[gray , thin] (5/3, 7/3) -- (5/3 , 0) node[above left]{$u$};
	\draw[gray , thin] (11/3, 10/3) -- (11/3 , 0) node[above left]{$v$};
	\draw[red] plot[smooth] coordinates {(1, 3.3) (2,3.2) (3,3) (4,2.2) (5,2)} node[right] {$f$};
	\end{tikzpicture}
	\caption{Case (V) in \cref{lem:piecewise:approx:case3}. The new function $h \in \scrL$ is linear on $[u,v]$ with slope $\mathbf{a}$,
		linear on $[v, \fq_i ]$ with slope $A_{i+1} ( g )$,
		and agrees with $g$ outside of $[v , \fq_i ]$.}
	\label{fig:piecewise:approx:3.5}
\end{SCfigure}
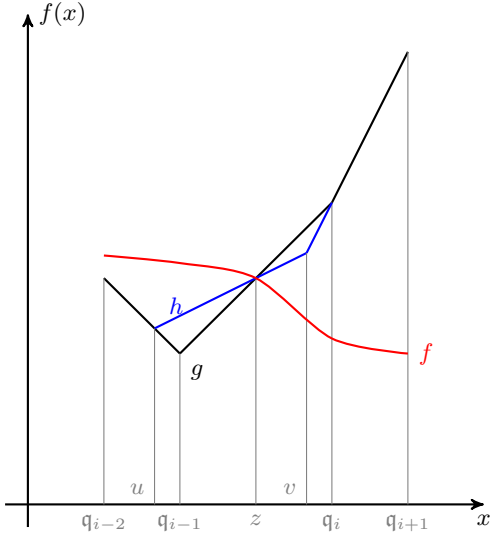
Finally, we prove \cref{eq:lem:piecewise:approx:case3:to_prove} in the case
\begin{multline}
\label{eq:last_case}
  \bigl(
    \max\{ 
      A_{ i - 1 }( g ), 
      A_{ i + 1 }( g ) 
    \} 
    > A_i( g ) 
    > 
    \min\{ 
      A_{ i - 1 }( g ), A_{ i + 1 }( g ) 
    \}
  \bigr)
\\
  \wedge 
  \bigl(
    \min\bigl\{ 
      \tfrac{ 
        g( \fq_{ i + 1 } ) - g( z ) 
      }{
        \fq_{ i + 1 } - z 
      } , 
      \tfrac{ 
        g( \fq_{ i - 2 } ) - g( z ) 
      }{
        \fq_{ i - 2 } - z 
      } 
    \bigr\} < \bfa
  \bigr)
\end{multline}
(cf.~\cref{fig:piecewise:approx:3.5}).
In the following we assume without loss of generality that 
\begin{equation}
\label{eq:obdA_in_proof}
  A_{ i - 1 }( g ) < A_i( g ) < A_{ i + 1 }( g ) 
\end{equation}
(cf.~\cref{lem:breakpoint:transformation}). 
\Nobs that \cref{eq:last_case,eq:obdA_in_proof}
show that 
$
  \frac{ 
    g( \fq_{ i - 2 } ) - g( z ) 
  }{
    \fq_{ i - 2 } - z 
  } 
  = 
  \min\bigl\{ 
    \frac{ 
      g( \fq_{ i + 1 } ) - g( z ) 
    }{
      \fq_{ i + 1 } - z 
    } , \allowbreak 
    \frac{ 
      g( \fq_{ i - 2 } ) - g( z ) 
    }{
      \fq_{ i - 2 } - z 
    } 
  \bigr\} 
  < \bfa 
$.
The intermediate value theorem therefore implies that 
there exist 
$
  u \in ( \fq_{ i - 2 }, \fq_{ i - 1 } ) 
$, 
$
  v \in ( z, \fq_i ) 
$ 
which satisfy 
$
  \frac{ g( u ) - g( z ) 
  }{
    u - z 
  } 
  = \bfa 
$ and 
$
  g( z ) + \bfa ( v - z ) = g( \fq_i ) + A_{ i + 1 }( g ) ( v - \fq_i ) 
$. 
Let $ h \in \scrL $ satisfy 
for all $ x \in [ 0, u ] \cup [ \fq_i, 1 ] $, 
$ y_1 \in [ u , v ] $, $ y_2 \in [ v , \fq_i ] $ 
that 
$
  h( x ) = g( x ) 
$, 
$
  h( y_1 ) = g(z) + \bfa ( y_1 - z ) 
$, 
and 
$
  h( y_2 ) = g( \fq_i ) + A_{ i + 1 }( g ) ( y_2 - \fq_i ) 
$. 
\Nobs that 
$ Q(h) = Q(g) $, 
$ A_i( h ) = \bfa $, 
and 
$
  \forall \, j \in \{ 1, 2, \ldots, Q( g ) + 1 \} \backslash \{ i \} \colon 
  A_j( h ) = A_j( g ) 
$.  
\Moreover the fact that 
$
  \forall \, x, y \in [0,1] \colon \abs{ f(x) - f(y) } \leq L \abs{ x - y }
$
and the fact that 
$
  \bfa < A_i( g )
$
demonstrate that
$ 
  \int_0^1 ( h(y) - f(y) )^2 \, \mu( \d y ) 
  \leq 
  \int_0^1 ( g(y) - f(y) )^2 \, \mu( \d y ) 
$. 
This establishes \cref{eq:lem:piecewise:approx:case3:to_prove} in the case 
$
  \bigl(
    \max\{ 
      A_{ i - 1 }( g ), 
      A_{ i + 1 }( g ) 
    \} 
    > A_i( g ) 
    > 
    \min\{ 
      A_{ i - 1 }( g ), A_{ i + 1 }( g ) 
    \}
  \bigr)
  \wedge 
  \bigl(
    \min\bigl\{ 
      \tfrac{ 
        g( \fq_{ i + 1 } ) - g( z ) 
      }{
        \fq_{ i + 1 } - z 
      } , 
      \tfrac{ 
        g( \fq_{ i - 2 } ) - g( z ) 
      }{
        \fq_{ i - 2 } - z 
      } 
    \bigr\} < \bfa
  \bigr)
$.
\end{enumerate}
\end{cproof}

Next, we summarize \cref{lem:piecewise:approx:case1,lem:piecewise:approx:case2,lem:piecewise:approx:case3} in the following corollary.

\cfclear
\begin{cor}
\label{lem:piecewise:approx:0}
Let $ L \in (0, \infty) $, $ f \in C( [0,1], \R ) $ 
satisfy for all $ x, y \in [0,1] $ that 
$
  \abs{ f(x) - f(y) } \leq L \abs{ x - y } 
$, 
let $ g \in \scrL $, $ i \in \{ 1, 2, \ldots, Q(g) + 1 \} $, $ \bfa \in \R $ 
satisfy $ L \le \abs{\bfa} \le \abs{ A_i(g) } $ and $ \bfa A_i( g ) > 0 $,
and let $ \mu \colon \cB( [0,1] ) \to [0, \infty] $ 
be a finite measure \cfadd{def:piece:linear}\cfadd{def:number:of:kinks}\cfadd{def:slopes}\cfload. 
Then there exists $ h \in \scrL $ such that 
$
  \int_0^1 (h(y) - f(y) ) ^2 \, \mu ( \d y ) \leq \int_0^1 (g(y) - f(y) )^2 \, \mu ( \d y ) 
$, 
$
  Q( h ) \leq Q( g )
$,
and
\begin{equation} \label{lem:piecewise:approx:0:eqclaim}
\textstyle
  \left( 
    Q(g) - Q(h) - 1 
  \right) 
  \big( 
    \sum_{ j = 1 }^{ Q(h) + 1 }
    | 
      A_j( h )
      - A_j(g) \indicator{ \N \backslash \{ i \} }( j ) 
      - \bfa \indicator{ \{ i \} }( j )
    |
  \big)
  \geq 0
  .
\end{equation}
\end{cor}
\begin{cproof} {lem:piecewise:approx:0}
Throughout this proof 
assume without loss 
of generality\footnote{Otherwise 
we consider $ f \with - f $, $ g \with -g $, $ \bfa \with - \bfa $.} 
that $ L \le \bfa \le A_i( g ) $. 
\Nobs that \cref{lem:piecewise:approx:case1} 
establishes \cref{lem:piecewise:approx:0:eqclaim} 
in the case 
\begin{equation}
  \bigl[
    \forall \, x \in ( q_{ i - 1 }( g ) , q_i( g ) ) \colon f(x) \not= g(x) 
  \bigr]
  .
\end{equation}
\Moreover \cref{lem:piecewise:approx:case2} 
establishes \cref{lem:piecewise:approx:0:eqclaim} 
in the case 
\begin{equation}
  \bigl[
    ( 
      \exists \, z \in ( q_{ i - 1 }( g ), q_i( g ) ) \colon f(z) = g(z) 
    ) 
    \wedge 
    ( 
      i \in \cu{ 1, Q(g) + 1 } 
    )
  \bigr]
  .
\end{equation}
\Moreover \cref{lem:piecewise:approx:case3} establishes \cref{lem:piecewise:approx:0:eqclaim} 
in the case 
\begin{equation}
  \bigl[
    ( 
      \exists \, z \in ( q_{ i - 1 }( g ), q_i( g ) ) \colon f(z) = g(z) 
    ) 
    \wedge 
    ( 
      i \notin \cu{ 1, Q(g) + 1 } 
    )
  \bigr]
  .
\end{equation}
\end{cproof}

The following two results are a consequence of \cref{cor:piecewise:approx:0} and induction.

\cfclear
\begin{cor}
\label{cor:piecewise:approx:0}
Let $ L \in (0, \infty) $, $ f \in C( [0,1], \R ) $ 
satisfy for all $ x, y \in [0,1] $ that 
$
  \abs{ f(x) - f(y) } \leq L \abs{ x - y } 
$, 
let $ g \in \scrL $, 
let 
$ \bbA \subseteq \{ 1, 2, \dots, Q(g) + 1 \} $ 
be a set, 
let 
$ \bfa = ( \bfa_j )_{ j \in \bbA } \colon \bbA \to \R $ 
satisfy for all $ j \in \bbA $ that 
$ L \le \abs{\bfa_j} \le \abs{ A_j(g) } $ and $ \bfa_j A_j( g ) > 0 $,
and let $ \mu \colon \cB( [0,1] ) \to [0, \infty] $ 
be a finite measure \cfadd{def:piece:linear}\cfadd{def:number:of:kinks}\cfadd{def:slopes}\cfload. 
Then there exists $ h \in \scrL $ such that 
$
  \int_0^1 (h(y) - f(y) ) ^2 \, \mu ( \d y ) \leq \int_0^1 (g(y) - f(y) ) ^2 \, \mu ( \d y ) 
$, 
$
  Q( h ) \leq Q( g )
$,
and
\begin{equation}
\label{eq:cor:piecewise:approx:0}
\textstyle
  \left( 
    Q(g) - Q(h) - 1 
  \right) 
  \big( 
    \sum_{ j = 1 }^{ Q(h) + 1 }
    | 
      A_j( h )
      - A_j(g) \indicator{ \N \backslash \bbA }( j ) 
      - \bfa_j \indicator{ \bbA }( j )
    |
  \big)
  \geq 0
  .
\end{equation}
\end{cor}
\begin{cproof}{cor:piecewise:approx:0}
\Nobs that induction and \cref{lem:piecewise:approx:0} 
establish \cref{eq:cor:piecewise:approx:0}. 
\end{cproof}

\cfclear
\begin{lemma}
\label{lem:piecewise:approx:1}
Let $ f \in C( [0,1], \R ) $
and 
let $ \mu \colon \cB( [0,1] ) \to [0, \infty] $ be a finite measure. 
Then for all $ g \in \scrL $ 
there exists $ h \in \scrL $ such that 
\begin{equation}
\label{eq:lem:piecewise:approx:1a_to_prove}
\textstyle
  Q(h) \leq Q(g) ,
  \quad
  \operatorname{Lip}( h ) 
  \leq
  \operatorname{Lip}( f ) ,
  \qandq
  \int_0^1 ( h(y) - f(y) )^2 \, \mu ( \d y ) 
  \leq 
  \int_0^1 ( g(y) - f(y) )^2 \, \mu ( \d y ) 
\end{equation}
\cfadd{def:piece:linear}\cfadd{def:lip:const}
\cfadd{def:number:of:kinks}\cfadd{def:slopes}\cfout.
\end{lemma}
\begin{cproof}{lem:piecewise:approx:1}
Throughout this proof 
assume without loss of generality that 
$
  0 < \operatorname{Lip}( f ) < \infty 
$
and let 
$ 
  \Q \colon \scrL \to \N_0
$ 
satisfy for all 
$
  g \in \scrL 
$
that
\begin{equation}
\label{eq:def_bbQ_functional}
  \Q(g) 
  =
  ( Q(g) + 1 )^2 
  +
  \#\!\left(
    \left\{ 
      i \in \{ 1, 2, \dots, Q(g) + 1 \} 
      \colon
      | A_i(g) | > \operatorname{Lip}( f )
    \right\}
  \right)
  .
\end{equation}
\Nobs that \cref{eq:def_bbQ_functional} 
assures for all $ g_1, g_2 \in \scrL $ 
with $ Q( g_1 ) < Q( g_2 ) $ 
that 
\begin{equation}
\label{eq:g1g2_growth_condition}
\begin{split}
  \Q( g_1 ) 
& \leq 
  ( Q( g_1 ) + 1 )^2 + Q(g_1) + 1 
  < 
  ( Q( g_1 ) + 1 )^2 + 2 ( Q(g_1) + 1 ) + 1
\\ &
  = ( Q( g_1 ) + 2 )^2 \leq ( Q( g_2 ) + 1 )^2 
  \leq \Q( g_2 ) 
  .
\end{split}
\end{equation}
Next we claim that for all $ k \in \N_0 $, 
$
  g \in \Q^{ - 1 }( \{ k \} )
$ 
there exists $ h \in \scrL $ such that 
\begin{equation}
\label{eq:induction_preserve_Lipschitz}
\textstyle
  Q(h) \leq Q(g) , 
  \quad
  \operatorname{Lip}( h ) \leq \operatorname{Lip}( f ) ,
  \qandq
  \int_0^1 ( h(y) - f(y) )^2 \, \mu ( \d y ) 
  \leq 
  \int_0^1 ( g(y) - f(y) )^2 \, \mu ( \d y ) 
  .
\end{equation}
We now prove \cref{eq:induction_preserve_Lipschitz} by 
induction on $ k \in \N_0 $. 
For the base case $ k = 0 $ we \nobs that 
$
  \Q^{ - 1 }( \{ 0 \} ) = \emptyset
$. 
This establishes \cref{eq:induction_preserve_Lipschitz} 
in the base case $ k = 0 $. 
For the induction step let $ k \in \N_0 $ satisfy 
for all $ g \in \Q^{ - 1 }( \{ 0, 1, \dots, k \} ) $ that there exists 
$ h \in \scrL $ such that 
\begin{equation}
\label{eq:induction_hypothesis:approx1a}
\textstyle
  Q(h) \leq Q(g) , 
  \quad
  \operatorname{Lip}( h ) \leq \operatorname{Lip}( f ) ,
  \qandq
  \int_0^1 ( h(y) - f(y) )^2 \, \mu ( \d y ) 
  \leq 
  \int_0^1 ( g(y) - f(y) )^2 \, \mu ( \d y ) 
\end{equation}
and let $ g \in \Q^{ - 1 }( \{ k + 1 \} ) $ 
satisfy 
\begin{equation}
\label{eq:Lip_large}
  \operatorname{Lip}( g ) > \operatorname{Lip}( f )
  .
\end{equation}
\Nobs that 
\cref{prop:piece:linear:lip}
and 
\cref{eq:Lip_large} ensure that 
there exists $ i \in \{ 1, 2, \dots, Q(g) + 1 \} $ 
which satisfies 
\begin{equation}
\label{eq:A_i_g_Lip_f}
  | A_i( g ) | > \operatorname{Lip}( f )
  .
\end{equation}
\Nobs that \cref{eq:A_i_g_Lip_f} shows that there exists 
$ \bfa \in \R $ 
which satisfies 
\begin{equation}
\label{eq:preparation_for_applying_piecewise:approx:0}
  \operatorname{Lip}( f ) = | \bfa | \leq | A_i( g ) |
  \qqandqq
  \bfa A_i( g ) > 0
  .
\end{equation}
\Nobs that 
\cref{eq:preparation_for_applying_piecewise:approx:0}, 
the fact that $ \operatorname{Lip}( f ) \in (0,\infty) $, 
and 
\cref{lem:piecewise:approx:0}
demonstrate that there exists 
$ \fg \in \scrL $ 
which satisfies 
$
  \int_0^1 ( \fg(y) - f(y) ) ^2 \, \mu ( \d y ) 
  \leq \int_0^1 (g(y) - f(y) ) ^2 \, \mu ( \d y ) 
$, 
$
  Q( \fg) \leq Q( g )
$,
and
\begin{equation}
\label{eq:frak_g_condition}
\textstyle
  \left( 
    Q(g) - Q(\fg) - 1 
  \right) 
  \big( 
    \sum_{ j = 1 }^{ Q(\fg) + 1 }
    | 
      A_j( \fg)
      - A_j(g) \indicator{ \N \backslash \{ i \} }( j ) 
      - \bfa \indicator{ \{ i \} }( j )
    |
  \big)
  \geq 0
  .
\end{equation}
\Nobs that \cref{eq:g1g2_growth_condition} 
and \cref{eq:frak_g_condition} assure that 
$
  \Q( \fg ) < \Q( g ) = k + 1
$. 
\Hence  
$ 
  \fg \in \Q^{ - 1 }( \{ 0, 1, \dots, k \} )
$. 
Combining this with \cref{eq:induction_hypothesis:approx1a} 
and \cref{eq:frak_g_condition} demonstrates that 
there exists $ h \in \scrL $ such that
$ Q(h) \leq Q( \fg ) \leq Q ( g ) $, 
$ \operatorname{Lip}( h ) \leq \operatorname{Lip}( f ) $, 
and
\begin{equation}
\textstyle
  \int_0^1 ( h(y) - f(y) ) ^2 \, \mu ( \d y ) 
  \leq \int_0^1 ( \fg(y) - f(y) ) ^2 \, \mu ( \d y ) 
  \leq \int_0^1 (g(y) - f(y) ) ^2 \, \mu ( \d y ) 
  .
\end{equation}
Induction thus establishes \cref{eq:induction_preserve_Lipschitz}. 
\Nobs that \cref{eq:induction_preserve_Lipschitz} 
implies \cref{eq:lem:piecewise:approx:1a_to_prove}.
\end{cproof}

\cref{lem:piecewise:approx:1} is not yet sufficient to establish \cref{prop:better:approx} since, as mentioned before, not every piecewise linear function with at most $\width \in \N$ breakpoints is representable by an ANN with $\width$ hidden neurons. Thus we need to ensure that the linear relation for the slopes (cf.~\cref{cor:characterization2}) is also preserved by our inductive construction. This is the content of \cref{lem:piecewise:approx:2}, which is again a consequence of \cref{lem:piecewise:approx:0} and induction.

\cfclear
\begin{lemma} 
\label{lem:piecewise:approx:2}
\cfadd{def:piece:linear}\cfadd{def:number:of:kinks}\cfadd{def:slopes}
Let 
$ f \colon [0,1] \to \R $
be Lipschitz continuous 
and 
let $ \mu \colon \cB( [0,1] ) \to [0, \infty] $ 
be a finite measure. 
Then for all 
$ g \in \scrL $, $ k \in \N $, $ i_1, i_2, \ldots, i_k \in \N $ 
with 
$ \frac{ k }{ 2 } \notin \N $, 
$ i_1 < i_2 < \cdots < i_k \leq Q ( g ) + 1 $, 
and 
$ \sum_{ j = 1 }^k ( - 1 )^j A_{ i_j }( g ) = 0 $ 
there exists $ h \in \scrL $ such that 
$ 
  \int_0^1 ( h(y) - f(y) )^2 \, \mu ( \d y ) 
  \leq \int_0^1 ( g(y) - f(y) )^2 \, \mu ( \d y ) 
$, 
$
  Q( h ) \leq Q(g) 
$,
and
\begin{equation} 
\label{eq:lem:piecewise:approx:2_A0}
\textstyle 
  (
    Q ( g ) - Q(h) - 1
  )
  \big(
    \bigl| 
      \sum_{ j = 1 }^k (-1)^j A_{ \min\{ i_j, Q(h)+1 \} }( h ) 
    \bigr|  
    +  
    \max\{ 
      \operatorname{Lip}( h )
      - Q( g ) \operatorname{Lip}(f)
      , 0
    \}
  \big)
  \geq 0
\end{equation}
\cfload. 
\end{lemma}
\begin{cproof}{lem:piecewise:approx:2}
\cfclear
Throughout this proof assume without loss of generality 
that $ \operatorname{Lip}( f ) > 0 $, 
let 
$
  \operatorname{sgn} \colon \R \to \R
$
satisfy for all $ x \in (0,\infty) $ that
$
  \operatorname{sgn}( x ) = 1 
$,
$
  \operatorname{sgn}( - x ) = -1
$, 
and 
$
  \operatorname{sgn}( 0 ) = 0
$,
and let 
$
  \Q \colon \scrL \to \N_0
$
satisfy for all 
$
  g \in \scrL 
$
that
\cfadd{def:piece:linear}\cfadd{def:number:of:kinks}\cfadd{def:slopes}
\begin{equation}
\label{eq:def_bbQ_functional_0}
  \Q(g)
  =
  \left( Q(g) + 1 \right)^2 
  +
  \#\bigl(
    \{ 
      i \in \{ 1, 2, \dots, Q(g) + 1 \} 
      \colon
      | A_i(g) | > Q(g) \operatorname{Lip}(f)
    \}
  \bigr)
\end{equation}
\cfload. 
\Nobs that \cref{eq:def_bbQ_functional_0} 
assures for all $ g_1, g_2 \in \scrL $ 
with $ Q( g_1 ) < Q( g_2 ) $ 
that 
\begin{equation}
\label{eq:g1g2_growth_condition_0}
\begin{split}
  \Q( g_1 ) 
& \leq 
  ( Q( g_1 ) + 1 )^2 + Q(g_1) + 1 
  < 
  ( Q( g_1 ) + 1 )^2 + 2 ( Q(g_1) + 1 ) + 1
\\ &
  = ( Q( g_1 ) + 2 )^2 \leq ( Q( g_2 ) + 1 )^2 
  \leq \Q( g_2 ) 
  .
\end{split}
\end{equation}
We claim that for all 
$ n \in \N_0 $, 
$ g \in \Q^{ - 1 }( \{ n \} ) $, 
$ k \in \N $,
$ i_1, i_2, \dots, i_k \in \N $
with 
$ \frac{ k }{ 2 } \notin \N $, 
$ i_1 < i_2 < \cdots < i_k \leq Q ( g ) + 1 $, 
and 
$ \sum_{ j = 1 }^k ( - 1 )^j A_{ i_j }( g ) = 0 $ 
there exists $ h \in \scrL $ such that 
$ 
  \int_0^1 ( h(y) - f(y) )^2 \, \mu ( \d y ) 
  \leq \int_0^1 ( g(y) - f(y) )^2 \, \mu ( \d y ) 
$, 
$
  Q( h ) \leq Q(g) 
$,
and
\begin{equation} 
\label{eq:lem:piecewise:approx:2_A00}
\textstyle 
  (
    Q ( g ) - Q(h) - 1
  )
  \big(
    \bigl| 
      \sum_{ j = 1 }^k (-1)^j A_{ \min\{ i_j, Q(h) + 1 \} }( h ) 
    \bigr|  
    +  
    \max\{ 
      \operatorname{Lip}( h )
      - Q( g ) \operatorname{Lip}(f)
      , 0
    \}
  \big)
  \geq 0 
  .
\end{equation}
We now prove \cref{eq:lem:piecewise:approx:2_A00} by induction 
on $ n \in \N_0 $. 
For the base case $ n = 0 $ \nobs that 
$ \Q^{ - 1 }( \{ 0 \} ) = \emptyset $. 
This establishes \cref{eq:lem:piecewise:approx:2_A00} in the base case $ n = 0 $. 
For the induction step let $ n \in \N_0 $ satisfy for all 
$ g \in \Q^{ - 1 }( \{ 0, 1, \dots, n \} ) $, 
$ k \in \N $,
$ i_1, i_2, \dots, i_k \in \N $
with 
$ \frac{ k }{ 2 } \notin \N $, 
$ i_1 < i_2 < \cdots < i_k \leq Q ( g ) + 1 $, 
and 
$ \sum_{ j = 1 }^k ( - 1 )^j A_{ i_j }( g ) = 0 $ 
that there exists 
$ h \in \scrL $ such that 
$ 
  \int_0^1 ( h(y) - f(y) )^2 \, \mu ( \d y ) 
  \leq \int_0^1 ( g(y) - f(y) )^2 \, \mu ( \d y ) 
$, 
$
  Q( h ) \leq Q(g) 
$,
and
\begin{equation} 
\label{eq:lem:piecewise:approx:2_A0_induction_hypothesis}
\textstyle 
  (
    Q( g ) - Q(h) - 1
  )
  \big(
    \bigl| 
      \sum_{ j = 1 }^k (-1)^j A_{ \min\{ i_j, Q(h) + 1 \} }( h ) 
    \bigr|  
    +  
    \max\{ 
      \operatorname{Lip}( h )
      - Q( g ) \operatorname{Lip}(f)
      , 0
    \}
  \big)
  \geq 0 
  ,
\end{equation}
and\footnote{\Nobs that we could choose 
$ h = g $ in \cref{eq:lem:piecewise:approx:2_A0_to_prove} 
if we would have  
$
  \operatorname{Lip}( g ) \leq Q(g) \operatorname{Lip}(f)
$.}    
let $ g \in \Q^{ - 1 }( \{ n + 1 \} ) $, 
$ k \in \N $,
$ i_1, i_2, \dots, i_k \in \N $
satisfy
\begin{equation}
\label{eq:g_condition_for_induction}
\textstyle
  \frac{ k }{ 2 } \notin \N ,
\qquad
  i_1 < i_2 < \cdots < i_k \leq Q ( g ) + 1 , 
\qquad
  \sum_{ j = 1 }^k ( - 1 )^j A_{ i_j }( g ) = 0 ,
\end{equation}
and 
$
  \operatorname{Lip}( g ) > Q(g) \operatorname{Lip}(f)
$. 
We now prove that there exists 
$ h \in \scrL $ such that 
$ 
  \int_0^1 ( h(y) - f(y) )^2 \, \mu ( \d y ) 
  \leq \int_0^1 ( g(y) - f(y) )^2 \, \mu ( \d y ) 
$, 
$
  Q( h ) \leq Q(g) 
$,
and
\begin{equation} 
\label{eq:lem:piecewise:approx:2_A0_to_prove}
\textstyle 
  (
    Q( g ) - Q(h) - 1
  )
  \big(
    \bigl| 
      \sum_{ j = 1 }^k (-1)^j A_{ \min\{ i_j, Q(h) + 1 \} }( h ) 
    \bigr|  
    +  
    \max\{ 
      \operatorname{Lip}( h )
      - Q( g ) \operatorname{Lip}(f)
      , 0
    \}
  \big)
  \geq 0 
  .
\end{equation}
\Nobs that 
\cref{prop:piece:linear:lip}
% and 
% \cref{eq:def_bbQ_functional_0}
and 
the fact that 
$
  \operatorname{Lip}( g ) > Q( g ) \operatorname{Lip}( f )
$
ensure that there exist
$ \fI \in \{ 1, 2, \dots, Q(g) + 1 \} $, $ s \in \{ - 1, 1 \} $ 
which satisfy  
\begin{equation}
\label{eq:bigJ}
  s A_{ \fI }( g ) = 
  | A_{ \fI }( g ) | 
  > Q(g) \operatorname{Lip}( f )
  .
\end{equation}
In the following we distinguish between the case 
$ \fI \notin \{ i_1, i_2, \dots, i_k \} $
and the case 
$ 
  \fI \in \{ i_1, i_2, \dots, i_k \} 
$.
We first prove \cref{eq:lem:piecewise:approx:2_A0_to_prove} 
in the case 
\begin{equation}
\label{eq:induction_case1}
  \fI \notin \{ i_1, i_2, \dots, i_k \} 
  .
\end{equation}
\Nobs that 
\cref{eq:bigJ} and 
\cref{lem:piecewise:approx:0} 
assure that there exists 
$ \fg \in \scrL $ 
which satisfies 
$ 
  \int_0^1 ( \fg(y) - f(y) )^2 \, \mu ( \d y ) 
  \leq \int_0^1 ( g(y) - f(y) )^2 \, \mu ( \d y ) 
$, 
$
  Q( \fg ) \leq Q(g) 
$,
and
\begin{equation}
\label{eq:application_Lipschitz_move_lemma_1}
\textstyle
  \left( 
    Q(g) - Q(\fg) - 1 
  \right) 
  \big( 
    \sum_{ j = 1 }^{ Q( \fg ) + 1 }
    | 
      A_j( \fg )
      - A_j(g) \indicator{ \N \backslash \{ \fI \} }( j ) 
      - 
      s Q(g) \! \operatorname{Lip}( f )
      \indicator{ \{ \fI \} }( j )
    |
  \big)
  \geq 0
  .
\end{equation}
Moreover, \nobs that 
\cref{eq:g_condition_for_induction} 
ensures that 
$
  \Q( g ) = n + 1 
$. 
Combining this with 
\cref{eq:def_bbQ_functional_0},
\cref{eq:g1g2_growth_condition_0}, 
and 
\cref{eq:bigJ}
demonstrates that 
$
  \Q( \fg ) < \Q( g ) = n + 1
$. 
\Hence  
\begin{equation}
\label{eq:fg_ready_for_induction_hypothesis}
  \fg \in \Q^{ - 1 }( \{ 0, 1, \dots, n \} ) 
  .
\end{equation}
In addition, \nobs that 
\cref{eq:g_condition_for_induction}, 
% \cref{eq:bigJ}, 
\cref{eq:induction_case1}, 
and 
\cref{eq:application_Lipschitz_move_lemma_1} 
show that 
$ Q( \fg ) \leq Q( g ) $,
$ 
  \int_0^1 ( \fg(y) - f(y) )^2 \, \mu ( \d y ) 
  \leq \int_0^1 ( g(y) - f(y) )^2 \, \mu ( \d y ) 
$, 
and 
\begin{equation}
\label{eq:fg_satisfy_g_condition}
\textstyle
  (
    Q( g ) - Q( \fg ) - 1
  )
    \bigl| 
      \sum_{ j = 1 }^k (-1)^j A_{ \min\{ i_j, Q(\fg) + 1 \} }( \fg ) 
    \bigr|  
  \geq 0 
  .
\end{equation}
Combining 
\cref{eq:lem:piecewise:approx:2_A0_induction_hypothesis}, 
\cref{eq:g_condition_for_induction}, 
and 
\cref{eq:fg_ready_for_induction_hypothesis} 
hence\footnote{\Nobs that we 
can choose $ h = \fg $ in \cref{eq:lem:piecewise:approx:2_A0_to_prove} 
in the case where $ Q( \fg ) < Q( g ) $.} 
establishes 
\cref{eq:lem:piecewise:approx:2_A0_to_prove} 
in the case $ \fI \notin \{ i_1, i_2, \dots, i_k \} $. 
In the next step we prove 
\cref{eq:lem:piecewise:approx:2_A0_to_prove} 
in the case 
\begin{equation} 
\label{eq:induction_case2}
  \fI \in \{ i_1, i_2, \dots, i_k \} 
  .
\end{equation}
\Nobs that \cref{eq:induction_case2} demonstrates that there exist 
$
  \scrJ \in \{ 1, 2, \dots, k \}
$, 
$
  \fS \in \{ - 1, 1 \}
$
which satisfy 
\begin{equation}
\label{eq:def_big_S}
  i_{ \scrJ } = \fI
\qqandqq
  \fS 
  = s (-1)^{ \scrJ } 
  = \operatorname{sgn}( ( - 1 )^{ \scrJ } A_{ i_{ \scrJ } } ) 
  = \operatorname{sgn}( ( - 1 )^{ \scrJ } A_{ \fI } ) 
  . 
\end{equation}
In the following let $ \alpha_v \in \R $, $ v \in \{ - 1, 1 \} $, 
satisfy for all $ v \in \{ - 1, 1 \} $ that
\begin{equation}
\label{eq:def:alpha_v}
\textstyle
  \alpha_v 
  =
%   \left[
    \sum_{ 
      \substack{
        j \in \N \cap [1,k] , 
        \,
%       \\
        \operatorname{sgn}(
          ( - 1 )^j A_{ i_j }( g ) 
        )
        = v \fS 
      }
    } 
    | A_{ i_j }( g ) | 
%   \right]
  .
\end{equation}
\Nobs that \cref{eq:g_condition_for_induction} 
and \cref{eq:def:alpha_v}
ensure that
\begin{equation}
\begin{split}
\textstyle
&
  \fS
  \left( 
    \alpha_1
    -
    \alpha_{ - 1 }
  \right)
=
\textstyle
  \fS \alpha_1 
  -
  \fS \alpha_{ - 1 }
  =
  \sum_{ v \in \{ - 1, 1 \} }
  \left[ 
    v \fS \alpha_v
  \right]
\\ & =
\textstyle 
  \sum_{ v \in \{ - 1, 1 \} }
    \sum_{ 
      \substack{
        j \in \N \cap [1,k] , 
        \,
        \operatorname{sgn}(
          ( - 1 )^j A_{ i_j }( g ) 
        )
        = v \fS 
      }
    } 
  \left[ 
    v \fS 
    | A_{ i_j }( g ) | 
  \right]
\\ & =
\textstyle 
  \sum_{ v \in \{ - 1, 1 \} }
    \sum_{ 
      \substack{
        j \in \N \cap [1,k] , 
        \,
        \operatorname{sgn}(
          ( - 1 )^j A_{ i_j }( g ) 
        )
        = v \fS 
      }
    } 
  \left[ 
    v \fS 
    | ( - 1 )^j A_{ i_j }( g ) | 
  \right]
\\ & =
\textstyle 
  \sum_{ v \in \{ - 1, 1 \} }
    \sum_{ 
      \substack{
        j \in \N \cap [1,k] , 
        \,
        \operatorname{sgn}(
          ( - 1 )^j A_{ i_j }( g ) 
        )
        = v \fS 
      }
    } 
  \left[ 
    ( - 1 )^j A_{ i_j }( g ) 
  \right]
% \\ &
% \textstyle
  =
  \sum_{ j \in \N \cap [1,k] }
  ( - 1 )^j 
  A_{ i_j }( g )
  = 
  0 .
\end{split}
\end{equation}
\Hence  
$
  \alpha_1 = \alpha_{ - 1 } 
$. 
Next \nobs that 
% \cref{eq:induction_case2}, 
\cref{eq:def_big_S} and \cref{eq:def:alpha_v} assure that 
$ \alpha_1 \geq | A_{ \fI }(g) | $. 
Combining this with 
\cref{eq:bigJ} 
and the fact that 
$
  \alpha_{ - 1 } = \alpha_1
$
demonstrates that 
\begin{equation}
  \alpha_{ - 1 } 
  \geq 
  \big( 
    | A_{ \fI }( g ) | - Q(g) \! \operatorname{Lip}( f ) 
  \big) 
  +
  Q(g) \! \operatorname{Lip}( f )
  >
  Q(g) \! \operatorname{Lip}( f )
  .
\end{equation}
\Hence that
there exist 
$ l \in \N $, 
$ J_1, J_2, \dots, J_l \in \N $, 
$ r_1, r_2, \dots, r_l \in [0,\infty) $ 
which satisfy for all $ v \in \{ 1, 2, \dots, l \} $ that 
\begin{equation}
\label{eq:introduction_of_Jv}
  J_1 < J_2 < \dots < J_l \leq k,
\qquad
  \operatorname{sgn}( 
    ( - 1 )^{ J_v }
    A_{ i_{ J_v } }( g )
  )
  =
  - \fS
  ,
\qquad
  | A_{ i_{ J_v } }( g ) | - r_v \geq \operatorname{Lip}( f )
  ,
\end{equation}
and 
$
  r_1 + r_2 + \ldots + r_l
  =
  | A_{ \fI }( g ) | - Q(g) \! \operatorname{Lip}( f ) 
$. 
In the following let 
$ 
  \bbA \subseteq \{ 1, 2, \dots, Q(g)+1 \} 
$
satisfy
\begin{equation}
  \bbA = 
  \{ i_{ J_1 }, i_{ J_2 }, \dots, i_{ J_l } \} \cup \{ \fI \}
\end{equation}
and let 
$ 
  \bfa = ( \bfa_j )_{ j \in \bbA } \colon \bbA \to \R 
$
satisfy for all 
$ v \in \{ 1, 2, \dots, l \} $ 
that 
\begin{equation}
\label{eq:introduction_of_bfa_function}
  \bfa_{ i_{ J_v } }
  =
  \left( 
    | 
      A_{ i_{ J_v } }
    |
    -
    r_v
  \right)
  \operatorname{sgn}( A_{ i_{ J_v } } )
\qqandqq
  \bfa_{ \fI }
  =
  s 
  Q(g) \! \operatorname{Lip}( f )
  .
\end{equation}
Note that
\cref{eq:introduction_of_Jv}
and 
\cref{eq:introduction_of_bfa_function}
ensure
for all $ v \in \{ 1, 2, \dots, l \} $ 
that 
\begin{equation}
\label{eq:strictly_positive_bfa_condition}
\begin{split}
  \bfa_{ i_{ J_v } }
  A_{
    i_{ J_v }
  }( g )
& =
  \left( 
    | 
      A_{ i_{ J_v } }( g )
    |
    -
    r_v
  \right)
  \operatorname{sgn}( A_{ i_{ J_v } }( g ) )
  A_{
    i_{ J_v }
  }( g )
  =
  \left( 
    | 
      A_{ i_{ J_v } }( g )
    |
    -
    r_v
  \right)
  |
    A_{
      i_{ J_v }( g )
    }
  |
\\ & 
\geq 
  \operatorname{Lip}( f )
  |
    A_{
      i_{ J_v }
    }( g )
  | 
\geq
  \operatorname{Lip}( f )
  (
    |
      A_{
        i_{ J_v }
      }( g )
    | 
    - r_v 
  )
  \geq 
  \left[ 
    \operatorname{Lip}( f )
  \right]^2
  > 0
  .
\end{split}
\end{equation}
Next \nobs that \cref{eq:bigJ} implies that 
$ 
  | A_{ \fI }( g ) | > Q( g ) \operatorname{Lip}( f ) \geq 0
$. 
\Hence $ A_{ \fI }( g ) \neq 0 $. 
This and \cref{eq:g_condition_for_induction} prove that 
$  
  Q(g) > 0
$. 
Combining 
\cref{eq:bigJ}, 
\cref{eq:introduction_of_Jv}, 
and 
\cref{eq:introduction_of_bfa_function}
therefore shows that 
\begin{equation}
  \bfa_{ \fI } 
  A_{ \fI }( g ) 
  =
  \big[ 
    s 
    Q(g) \! \operatorname{Lip}( f )
  \big]
  \big[ 
    s^{ - 1 } 
    | A_{ \fI }( g ) |
  \big]
  =
    | A_{ \fI }( g ) |
    Q(g) \!
    \operatorname{Lip}( f )
  > 
  \left[
    Q(g) \!
    \operatorname{Lip}( f )
  \right]^2
  > 0
  .
\end{equation}
This and \cref{eq:strictly_positive_bfa_condition} 
assure for all $ j \in \bbA $ that 
\begin{equation}
\label{eq:strictly_positive_bfa_condition_joined}
  \bfa_j A_j( g ) > 0
  .
\end{equation}
Furthermore, \nobs that 
\cref{eq:bigJ}, 
\cref{eq:introduction_of_Jv}, 
\cref{eq:introduction_of_bfa_function}, 
and the fact that 
$
  Q(g) \geq 1
$
demonstrate for all $ v \in \{ 1, 2, \dots, l \} $ 
that 
\begin{equation}
  \operatorname{Lip}( f )
  \leq 
  | A_{ i_{ J_v } } | - r_v
  =
  | \bfa_{ i_{ J_v } } | 
  \leq 
  | A_{ i_{ J_v } } | 
\qqandqq
  \operatorname{Lip}( f )
  \leq 
  Q( g ) \!
  \operatorname{Lip}( f )
  =
  | \bfa_{ \fJ } |
  < | A_{ \fJ }( g ) |
  .
\end{equation}
\Hence for all $ j \in \bbA $ that
$
  \operatorname{Lip}( f ) \leq | \bfa_j | \leq A_j( g ) 
$.
Combining this with 
\cref{eq:strictly_positive_bfa_condition_joined}
enables us to apply 
\cref{cor:piecewise:approx:0} 
to obtain that there exists 
$ \fg \in \scrL $ which satisfies  
$
  \int_0^1 ( \fg(y) - f(y) )^2 \, \mu( \d y ) 
  \leq \int_0^1 ( g(y) - f(y) )^2 \, \mu( \d y ) 
$, 
$
  Q( \fg ) \leq Q( g )
$,
and
\begin{equation}
\label{eq:corollary_consequence_bfa_substracted}
\textstyle
  \left( 
    Q(g) - Q(\fg) - 1 
  \right) 
  \big( 
    \sum_{ j = 1 }^{ Q(\fg) + 1 }
    | 
      A_j( \fg )
      - A_j( g ) \indicator{ \N \backslash \bbA }( j ) 
      - \bfa_j \indicator{ \bbA }( j )
    |
  \big)
  \geq 0
  .
\end{equation}
\Nobs that 
\cref{eq:def_bbQ_functional_0}, 
\cref{eq:g1g2_growth_condition_0}, 
\cref{eq:bigJ}, 
\cref{eq:introduction_of_bfa_function}, 
and 
\cref{eq:corollary_consequence_bfa_substracted}
show that 
$ \Q( \fg ) < \Q( g ) $.
This and \cref{eq:g_condition_for_induction} 
show that 
$ 
  \Q( \fg ) \leq \Q( g ) - 1 = ( n + 1 ) - 1 = n 
$.
\Hence  
\begin{equation}
\label{eq:fg_ready_for_induction_hypothesis_case2}
  \fg \in \Q^{ - 1 }( \{ 0, 1, \dots, n \} ) 
  .
\end{equation}
Moreover, \nobs that 
\cref{eq:g_condition_for_induction}, 
% \cref{eq:introduction_of_Jv}, 
\cref{eq:introduction_of_bfa_function}, 
% and 
\cref{eq:corollary_consequence_bfa_substracted}, 
and the fact that 
$
  r_1 + r_2 + \ldots + r_l = | A_{ \fJ }( g ) | - Q( g ) \operatorname{Lip}( f )
$
assure that 
$
  \int_0^1 ( \fg(y) - f(y) )^2 \, \mu( \d y ) 
  \leq \int_0^1 ( g(y) - f(y) )^2 \, \mu( \d y ) 
$, 
$
  Q( \fg ) \leq Q( g )
$, 
and 
\begin{equation}
\textstyle
  \left( 
    Q(g) - Q(\fg) - 1 
  \right) 
  |
    \sum_{ j = 1 }^k
    ( - 1 )^j
    A_{ \min\{ i_j, Q( \fg ) + 1 \} }( \fg )
  |
  \geq 0
  .
\end{equation}
Combining 
\cref{eq:lem:piecewise:approx:2_A0_induction_hypothesis}, 
\cref{eq:g_condition_for_induction}, 
and 
\cref{eq:fg_satisfy_g_condition} 
hence\footnote{\Nobs 
that we can choose $ h = \fg $ in \cref{eq:lem:piecewise:approx:2_A0_to_prove} 
in the case where $ Q( \fg ) < Q( g ) $.} establishes 
\cref{eq:lem:piecewise:approx:2_A0_to_prove} 
in the case $ \fI \in \{ i_1, i_2, \dots, i_k \} $. 
Induction thus proves \cref{eq:lem:piecewise:approx:2_A00}. 
\Nobs that \cref{eq:lem:piecewise:approx:2_A00} 
establishes \cref{eq:lem:piecewise:approx:2_A0}. 
\end{cproof}

\subsection{Structure preserving approximations 
for realization functions of shallow ANNs}
\label{ssec:structure_preserving_ANNs}

In this subsection we employ \cref{lem:piecewise:approx:1} and \cref{lem:piecewise:approx:2} above to prove in \cref{prop:better:approx} the announced result about the existence of a better ANN approximation which is additionally Lipschitz with a constant depending only on the width $\width \in \N$ and the target function $f$.

\cfclear
\begin{lemma}
\label{lem:better_approx}
\cfadd{def:piece:linear}\cfadd{def:number:of:kinks}
Assume \cref{setting:snn} 
and let $ \theta \in \R^{ \fd } $ satisfy 
$
  (
    \width - Q( \realization{\theta} ) - 1
  )
    \max\{ 
      \operatorname{Lip}( \realization{\theta} ) 
      - \width L
      , 0
    \}
  \geq 0
$
\cfload. 
\cfadd{def:lip:const}
Then there exists $ \vartheta \in \R^\fd $ such that 
\begin{equation}
\label{eq:to_prove_better_approx_lemma}
  \cL( \vartheta ) \leq 
  \cL( \theta )
  ,
  \qquad 
  Q( \realization{\vartheta} ) \leq Q( \realization{\theta} ) ,
  \qquad
  \text{and}
  \qquad
  \operatorname{Lip}( \realization{\vartheta} ) \leq \width L 
\end{equation}
\cfload. 
\end{lemma}
\begin{cproof}{lem:better_approx} 
In the following we distinguish between the case 
$
  Q( \realization{\theta} ) = \width 
$ 
and the case 
$
  Q( \realization{\theta} ) < \width 
$.
We first prove \cref{eq:to_prove_better_approx_lemma} 
in the case 
\begin{equation}
\label{eq:lem:better_approx_case1}
  Q( \realization{\theta} ) = \width 
  .
\end{equation}
\Nobs that \cref{eq:lem:better_approx_case1} and the assumption that 
$
  (
    \width - Q( \realization{\theta} ) - 1
  )
    \max\{ 
      \operatorname{Lip}( \realization{\theta} ) 
      - \width L
      , 0
    \}
  \geq 0
$
ensure that 
$
  -
    \max\{ 
      \operatorname{Lip}( \realization{\theta} ) 
      - \width L
      , 0
    \}
  \geq 0
$. 
\Hence  
$
  \operatorname{Lip}( \realization{\theta} ) 
  \leq 
  \width L
$. 
This establishes \cref{eq:to_prove_better_approx_lemma} 
in the case $ Q( \realization{\theta} ) = \width $. 
In the next step we prove \cref{eq:to_prove_better_approx_lemma} 
in the case 
\begin{equation}
\label{eq:lem:better_approx_case2}
  Q( \realization{\theta} ) < \width 
  .
\end{equation}
\Nobs that \cref{eq:lem:better_approx_case2} 
and \cref{lem:piecewise:approx:1} 
% and \cref{prop:piece:linear:lip} 
prove that there exists 
$ h \in \scrL $ 
which satisfies $ Q(h) \leq Q( \realization{\theta} ) \leq \width - 1 $, 
$ \operatorname{Lip}(h) \leq L \leq \width L $, 
and
\begin{equation}
\label{eq:better_approx_case1}
  \int_0^1 ( h(y) - f(y) )^2 \, \mu ( \d y ) 
\leq 
  \int_0^1 ( \realization{\theta}( y ) - f( y ) )^2 \, \mu ( \d y ) 
= 
  \cL( \theta ) .
\end{equation}
\Nobs that \cref{lem:h-1:kinks} 
and the fact that $ Q( h ) \leq \width - 1 $ 
show that there exists $ \vartheta \in \R^\fd $ 
which satisfies 
\begin{equation}
\label{eq:better_approx_existence_of_realization_case1}
  \realization{\vartheta} = h 
  .
\end{equation}
\Nobs that 
\cref{eq:better_approx_case1} 
and \cref{eq:better_approx_existence_of_realization_case1} 
ensure that 
$
  \cL( \vartheta ) \leq \cL ( \theta) 
$, 
$
  Q( \realization{\vartheta} ) = Q( h ) \leq Q( \realization{\theta} )
$, 
and 
$ 
  \operatorname{Lip}( \realization{\vartheta} ) = \operatorname{Lip}( h ) \leq \width L 
$.
This establishes \cref{eq:to_prove_better_approx} in the case 
$ 
  Q( \realization{\theta} ) < \width 
$.
\end{cproof}

\cfclear
\begin{prop}
\label{prop:better:approx} 
\cfadd{def:lip:const}
Assume \cref{setting:snn} and let $\theta \in \R^\fd$. 
\cfadd{def:number:of:kinks}
Then there exists $ \vartheta \in \R^\fd $ such that 
\begin{equation}
\label{eq:to_prove_better_approx}
  \cL( \vartheta ) \leq \cL ( \theta ) ,
  \qquad 
  Q( \realization{\vartheta} ) \leq Q( \realization{\theta} ) ,
  \qquad
  \text{and}
  \qquad
  \operatorname{Lip} ( \realization{\vartheta} ) \leq \width L 
\end{equation}
\cfload. 
\end{prop}
\begin{cproof}{prop:better:approx} 
\cfadd{def:piece:linear}\cfadd{def:number:of:kinks}
\Nobs that \cref{lem:realization:l} proves that $\realization{\theta} \in \scrL$ and 
$Q( \realization{\theta} ) \leq \width$ \cfload. 
In the following we distinguish between the case 
$ Q( \realization{\theta} ) < \width $
and the case 
$ Q( \realization{\theta} ) = \width $. 
We first prove \cref{eq:to_prove_better_approx} in the case 
\begin{equation}
\label{eq:case1_better_approx}
  Q( \realization{\theta} ) < \width 
  .
\end{equation}
\Nobs that \cref{eq:case1_better_approx} and \cref{lem:better_approx} show that 
there exists $ \vartheta \in \R^{ \fd } $ such that 
$
  \cL( \vartheta ) \leq \cL( \theta )
$,
$
  Q( \realization{\vartheta} ) \leq Q( \realization{\theta} )
$,
and 
$
  \operatorname{Lip}( \realization{\vartheta} ) \leq \width L
$.
This establishes \cref{eq:to_prove_better_approx} in the case 
$ 
  Q( \realization{\theta} ) < \width 
$.
In the next step we prove \cref{eq:to_prove_better_approx} in the case 
\begin{equation}
\label{eq:better_approx_case2_formulation_of_case}
  Q( \realization{\theta} ) = \width  
  .
\end{equation}  
\Nobs that \cref{eq:better_approx_case2_formulation_of_case} 
and \cref{cor:characterization2} imply that there exists 
$ k \in \N $, $ i_1, i_2, \ldots, i_k \in \N $ which satisfy  
$ \frac{ k }{ 2 } \notin \N $, $ 1 \leq i_1 < i_2 < \cdots < i_k \leq \width + 1 $, 
and 
$
  \sum_{ j = 1 }^k (-1)^j A_{ i_j }( \realization{\theta} ) 
$ 
\cfadd{def:slopes}\cfload. 
Combining this with \cref{lem:piecewise:approx:2} ensures 
that there exists $ h \in \scrL $ which satisfies 
\begin{equation}
\textstyle
  \int_0^1 ( h(y) - f(y) )^2 \, \mu ( \d y ) 
  \leq \int_0^1 ( \realization{\theta}(y) - f(y) )^2 \, \mu( \d y ) 
  = \cL( \theta )
  ,
  \qquad 
  Q( h ) \leq Q( \realization{\theta} ) 
  ,
\end{equation}
and
\begin{equation} 
\textstyle 
  (
    Q( \realization{\theta} ) - Q( h ) - 1
  )
  \big(
    \bigl| 
      \sum_{ j = 1 }^k (-1)^j A_{ i_j }( h ) 
    \bigr|  
    +  
    \max\{ 
      \operatorname{Lip}( h ) 
      - Q( \realization{\theta} ) L
      , 0
    \}
  \big)
  \geq 0
  .
\end{equation}
Hence, we obtain that 
\begin{equation}
\label{eq:lem:piecewise:approx:1a}
\textstyle
  \int_0^1 ( h(y) - f(y) )^2 \, \mu ( \d y ) 
  \leq \cL( \theta )
  ,
  \qquad 
  Q( h ) \leq \width
  ,
\end{equation}
\begin{equation} 
\label{eq:lem:piecewise:approx:2BBB}
\textstyle 
  \text{and}
  \qquad
  (
    \width - Q( h ) - 1
  )
  \big(
    \bigl| 
      \sum_{ j = 1 }^k (-1)^j A_{ i_j }( h ) 
    \bigr|  
    +  
    \max\{ 
      \operatorname{Lip}( h ) 
      - \width L
      , 0
    \}
  \big)
  \geq 0
  .
\end{equation}
Therefore, we get that
$
  (
    \width - Q( h ) - 1
  )
    \bigl| 
      \sum_{ j = 1 }^k (-1)^j A_{ i_j }( h ) 
    \bigr|  
  \geq 0
$. 
Combining \cref{cor:characterization2} and 
\cref{eq:lem:piecewise:approx:1a} hence shows that there exist 
$ \psi \in \R^{ \fd } $ which satisfies 
\begin{equation}
\label{eq:existence_realization}
  \realization{\psi} = h 
  .
\end{equation}
Note that \cref{eq:lem:piecewise:approx:1a}, 
\cref{eq:lem:piecewise:approx:2BBB}, and \cref{eq:existence_realization} demonstrate that 
\begin{equation}
  \cL( \psi ) \leq \cL( \theta ) ,
  \quad 
  Q( \realization{\psi} ) \leq \width ,
  \quad
  \text{and}
  \quad
  (
    \width - Q( \realization{\psi} ) - 1
  )
    \max\{ 
      \operatorname{Lip}( \realization{\psi} ) 
      - \width L
      , 0
    \}
  \geq 0
  .
\end{equation}
\cref{lem:better_approx} hence implies that there exists 
$ \vartheta \in \R^{ \fd } $ which satisfies 
\begin{equation}
  \cL( \vartheta ) \leq \cL( \psi ) \leq \cL( \theta ) ,
  \qquad 
  Q( \realization{\vartheta} ) \leq Q( \realization{\psi} ) \leq 
  \width = Q( \realization{\theta} ) ,
  \qquad
  \text{and}
  \qquad
  \operatorname{Lip}( \realization{\vartheta} ) \leq \width L
  .
\end{equation}
This proves \cref{eq:to_prove_better_approx} in the case $ Q( \realization{\theta} ) = \width $. 
\end{cproof}

As a simple consequence of \cref{prop:better:approx} we obtain in \cref{cor:better:approx} below
that the new network parameter vector $\vartheta \in \R^\fd$ can be chosen in such a way that in addition to the Lipschitz constant also the supremum norm of its realization function $\realization{\vartheta}$ is bounded by a constant depending only on $\width$ and the target function $f$.

\cfclear
\begin{cor} 
\label{cor:better:approx} 
\cfadd{def:lip:const}
Assume \cref{setting:snn} and let $ \theta \in \R^{ \fd } $. 
Then there exists $ \vartheta \in \R^\fd $ such that 
$ \cL( \vartheta ) \leq \cL( \theta ) $, 
$
  Q( \realization{\vartheta} ) \leq Q( \realization{\theta} )
$,
$
  \sup_{ x \in [0,1] } 
  \abs{ 
    \realization{ \vartheta }( x ) 
  } \leq \width L + \sup_{ x \in [0,1] } \abs{ f(x) } 
$, 
and 
$
  \operatorname{Lip}( \realization{\vartheta} ) \leq \width L 
$ \cfload.
\end{cor}
\begin{cproof}{cor:better:approx}
\Nobs that \cref{prop:better:approx} establishes that there exist $\psi \in \R^\fd$,
$r \in [0, \infty)$
which satisfy 
\begin{equation} 
\label{cor:better:approx:eqr}
\textstyle
  \cL( \psi ) \leq \cL( \theta ) ,
  \quad 
  Q( \realization{\psi} ) \leq Q( \realization{\theta} ) ,
  \quad 
  \operatorname{Lip} ( \realization{\psi} ) \leq \width L,
  \qandq 
  r = 
  \inf_{ x \in [0,1] } 
  \abs{ 
    \realization{\psi}( x ) - f( x ) 
  } 
  .
\end{equation}
\Nobs that \cref{cor:better:approx:eqr} assures 
that there exist $ y \in [0,1] $, $ k \in \cu{-1, 1} $ 
which satisfy $ \realization{\psi}( y ) - f( y ) = kr $. 
In the following let $ \vartheta \in \R^\fd $ satisfy 
for all $ i \in \{ 1, 2, \ldots, 3 \width \} $ that 
$ \vartheta_i = \psi_i $ and $ \vartheta_{ \fd } = \psi_{ \fd } - k r $. 
\Nobs that the fact that for all $ x \in [0,1] $ it holds that 
$ \realization{\vartheta}( x ) = \realization{\psi}( x ) - k r $ 
and \cref{cor:better:approx:eqr}
show that 
$ 
  \operatorname{Lip}( \realization{ \vartheta } ) 
  = \operatorname{Lip}( \realization{ \psi } ) 
  \leq \width L 
$
and 
$
  Q( \realization{ \vartheta } ) = Q( \realization{\psi} ) \leq Q( \realization{\theta} )
$. 
The fact that 
$ 
  \realization{\vartheta}( y ) = \realization{\psi}( y ) - k r = f( y ) 
$
and the triangle inequality
therefore imply for all $z \in [0,1]$
that
\begin{equation}
\begin{split}
  \abs{ \realization{\vartheta } ( z ) } 
& 
  \leq \abs{ \realization{\vartheta } ( y ) } + \abs{ \realization{ \vartheta } ( z ) - \realization{\vartheta } ( y ) } = \abs {f( y ) } + \abs{ \realization{\vartheta } ( z ) - \realization{\vartheta } ( y ) } \\
&
  \leq \sup\nolimits_{ x \in [0,1] } \abs{ f(x) } + \width L \abs{z - y } 
  \leq \sup\nolimits_{ x \in [0,1] } \abs{ f(x) } + \width L .
\end{split}
\end{equation}
It remains to prove that 
$ \cL( \vartheta ) \leq \cL( \psi ) $.
For this we assume without loss of generality that
\begin{equation}
\label{eq:obdA_r_stricly_positive}
  r > 0 .
\end{equation}
\Nobs that 
\cref{cor:better:approx:eqr}, 
\cref{eq:obdA_r_stricly_positive},  
and the fact that $ [0, 1] \ni x \mapsto \realization{\psi}( x ) - f( x ) \in \R $
is continuous imply for all $ x \in [0,1] $ 
that $ k ( \realization{\psi}( x ) - f( x ) ) \ge r $.
Hence, we obtain for all $ x \in [0,1] $ that
\begin{equation}
\begin{split}
  \abs{\realization{\vartheta}( x ) - f( x ) } 
& 
  = \abs{ \realization{\psi}( x ) - f( x ) - k r }
  = \abs{ k ( \realization{\psi}( x ) - f( x ) ) - r } 
\\
&
  = k ( \realization{\psi}( x ) - f( x ) ) - r
  \le k ( \realization{\psi}( x ) - f( x ) ) 
  \le \abs{ \realization{\psi}( x ) - f( x ) } .
\end{split}
\end{equation}
This demonstrates that $ \cL( \vartheta ) \leq \cL( \psi ) $.
\end{cproof}

\subsection{Existence of global minima for shallow ANNs}
\label{subsec:existence_global_minima}

In this subsection we establish in \cref{theo:existence} the existence of a global minimizer of the risk function under the assumptions of \cref{setting:snn}. For the proof, we combine \cref{cor:better:approx} with the Arzel\`{a}--Ascoli theorem to extract a convergent subsequence from a minimizing sequence.
Due to the fact that the set of realization functions 
of shallow ReLU ANNs with fixed architecture is closed in the set of continuous functions with respect to the supremum norm (cf.~Petersen et al.~\cite[Theorem 3.8]{PetersenRaslanVoigtlaender2020})
the limit is again equal to the realization function of a suitable ANN.

\cfclear
\begin{prop} 
\label{theo:existence}
\cfadd{def:lip:const}
Assume \cref{setting:snn}. 
Then there exists $ \theta \in \R^\fd $ such that 
$
  \cL( \theta ) = \inf_{ \vartheta \in \R^\fd } \cL( \vartheta ) 
$,
$
  \operatorname{Lip}( \realization{ \theta } ) \leq \width L 
$,
and 
$ 
  \sup_{ x \in [0,1] } \abs{ \realization{ \theta }(x) } 
  \leq \width L + \sup_{ x \in [0,1] } \abs{ f(x) } 
$ \cfload.
\end{prop}
\begin{cproof}{theo:existence}
\Nobs that there exists 
$ \phi = ( \phi_n )_{ n \in \N } \colon \N \to \R^\fd $ 
which satisfies 
\begin{equation}
\label{eq:approximating_sequence}
\textstyle
  \limsup_{ n \to \infty } \cL( \phi_n ) 
  = \inf_{ \vartheta \in \R^\fd } \cL( \vartheta  )
  .
\end{equation}
\Nobs that \cref{cor:better:approx} implies that there exists  
$ \psi = ( \psi_n )_{ n \in \N } \colon \N \to \R^\fd $
which satisfies for all $ n \in \N $ that 
\begin{equation}
\label{eq:improved_approx}
\textstyle
  \cL( \psi_n ) \leq \cL( \phi_n ) ,
  \quad
  \sup_{ x \in [0,1] } \abs{ \realization{ \psi_n }(x) } 
  \leq \width L + \sup_{ x \in [0,1] } \abs{ f(x) } 
  , 
  \quad
  \text{and}
  \quad 
  \operatorname{Lip}( \realization{ \psi_n } ) \leq \width L 
  .
\end{equation}
\Nobs that \cref{eq:approximating_sequence} and \cref{eq:improved_approx} show that
\begin{equation}
\textstyle 
  \inf_{ \vartheta \in \R^{ \fd } } 
  \cL( \vartheta  ) 
  \leq \limsup_{ n \to \infty } \cL( \psi_n ) 
  \leq \limsup_{ n \to \infty } \cL( \phi_n ) 
  = \inf_{ \vartheta \in \R^\fd } \cL( \vartheta ) 
  .
\end{equation}
Hence, we obtain that 
$ 
  \lim_{ n \to \infty } \cL( \psi_n ) = \inf_{ \vartheta \in \R^\fd } \cL( \vartheta ) 
$. 
Furthermore, \nobs that 
\cref{eq:improved_approx} and 
the Arzela--Ascoli theorem demonstrate that 
there exist $ g \in C( [0,1], \R ) $ and 
a strictly increasing $ k \colon \N \to \N $ 
such that 
\begin{equation}
\label{eq:convergence_uniform_sequence}
\textstyle
  \limsup_{ n \to \infty } 
  \sup_{ x \in [0,1] } 
  \abs{ \realization{ \psi_{ k(n) } }(x) - g(x) } = 0
  .
\end{equation} 
Next \nobs that 
Petersen et al.~\cite[Theorem 3.8]{PetersenRaslanVoigtlaender2020} 
assures that 
$
  \{ 
    h \in C( [0,1], \R ) \colon
    (
      \exists \, 
      \vartheta \in \R^{ \fd } \colon 
      \realization{ \vartheta } = h
    )
  \}
$
is a closed subset of $C([0,1], \R)$ with respect to the supremum norm on $ C( [0,1], \R ) $. 
Combining this with \cref{eq:convergence_uniform_sequence} implies  
that there exists $ \theta \in \R^\fd $ 
which satisfies 
\begin{equation}
\label{eq:existence_of_realization}
  \realization{\theta} = g 
  .
\end{equation}
\Nobs that 
\cref{eq:convergence_uniform_sequence}, 
\cref{eq:existence_of_realization}, 
and Lebesgue's theorem of dominated convergence ensure that
\begin{equation}
\begin{split}
  \cL( \theta ) 
& = 
  \int_0^1 ( \realization{\theta }( y ) - f ( y ) )^2 \, \mu ( \d y ) 
= 
  \int_0^1 ( g(y) - f(y) )^2 \, \mu ( \d y ) \\
& = 
  \int_0^1 
  \left[ 
    \lim_{ n \to \infty } 
    ( \realization{ \psi _{ k( n ) } }(y) - f(y) )^2 
  \right]  
  \mu( \d y ) 
  = 
  \lim_{ n \to \infty } 
  \br*{ 
    \int_0^1 ( \realization{\psi _{k ( n ) } }(y) - f(y) )^2 
  \, \mu ( \d y ) } 
\\
& 
  = \lim_{ n \to \infty } \cL( \psi_{ k(n) } ) 
  = \inf_{ \vartheta \in \R^\fd } \cL( \vartheta  ) .
\end{split}
\end{equation}
Furthermore, \nobs that \cref{eq:improved_approx}, \cref{eq:convergence_uniform_sequence}, 
and \cref{eq:existence_of_realization} demonstrate that 
$
  \sup_{ x \in [0,1] } \abs{ \realization{ \theta }(x) } 
  \leq \width L + \sup_{ x \in [0,1] } \abs{ f(x) } 
$
and 
$
  \operatorname{Lip}( \realization{ \theta } ) \leq \width L 
$.
\end{cproof}

\cref{theo:existence} is formulated only for the input domain $[0,1]$. 
In \cref{theo:existence2} we generalize this result to a general input interval $[a,b] \subseteq \R$ by employing a suitable coordinate transformation.

\cfclear
\begin{theorem} 
\label{theo:existence2}
Let $ \width, \fd \in \N $, $ L, a \in \R $, $ b \in (a, \infty) $, 
$ f \in C( [a,b], \R ) $ 
satisfy for all $ x, y \in [a,b] $ that 
$ \fd = 3 \width + 1 $
and 
$
  \abs{ f(x) - f(y) } \leq L \abs{ x - y }
$, 
let $ \mu \colon \cB( [a,b] ) \to [0, \infty] $ be a finite measure, 
for every $ \theta = ( \theta_1, \dots, \theta_{ \fd } ) \in \R^{ \fd } $ let 
$ 
  \realization{ \theta } \colon \R \to \R 
$ 
satisfy for all $ x \in \R $ that
\begin{equation}
\textstyle
  \realization{ \theta }( x ) 
  =
  \theta_{ \fd } 
  +
  \smallsum_{ j = 1 }^{ \width } 
  \theta_{ 2 \width + j } 
  \max\{ \theta_{ \width + j } + \theta_j x , 0 \} 
  ,
\end{equation}
and let $ \cL \colon \R^{ \fd } \to \R $ 
satisfy for all 
$ \theta \in \R^{ \fd } $ 
that
$
  \cL( \theta ) = 
  \textstyle
  \int_a^b 
  \displaystyle
    ( 
      f( x ) 
      - 
      \realization{ \theta }( x )
    )^2 
  \, \mu( \d x ) 
$.
Then there exists 
$ \theta \in \R^{ \fd } $ 
such that 
$
  \cL( \theta ) = 
  \inf_{ \vartheta \in \R^{ \fd } } 
  \cL( \vartheta )
$, 
$ 
  \sup_{ x \in [a,b] } \abs{ \realization{ \theta }(x) } 
  \leq \width L ( b - a ) + \sup_{ x \in [a,b] } \abs{ f(x) } 
$, 
and 
\begin{equation}
\textstyle
  \sup_{ x, y \in [a,b], x \neq y }
  \big[
    \frac{
      \abs{ \realization{\theta} ( x ) - \realization{\theta} ( y ) }
    }{
      \abs{ x - y }
    }
  \big]
  \leq \width L 
  .
\end{equation}
\end{theorem}
\begin{cproof}{theo:existence2}
Throughout this proof let $ {\bf f} \colon [0,1] \to [a,b] $ 
and $ {\bf F} \colon \R^\fd \to \R^\fd $ 
satisfy for all 
$ x \in [0,1] $, 
$ \theta = ( \theta_1, \ldots, \theta_\fd) \in \R^\fd $ 
that $ {\bf f}( x ) = a + (b-a) x $ 
and
\begin{equation}
  {\bf F}( \theta ) = ( (b-a) \theta_1, \ldots, (b-a) \theta_{ \width } ,
  \theta_{ \width + 1 } + a \theta_1, \ldots, \theta_{2 \width } + a \theta_\width ,
  \theta_{ 2 \width + 1 }, \ldots, \theta_{ 3 \width }, \theta_{ 3 \width + 1 } ),
\end{equation}
let $ g \in C( [0,1], \R ) $ satisfy 
for all $ x \in [0,1] $ that 
$ g(x) = f( {\bf f}( x ) ) $,
and let $ \nu \colon \cB( [0,1] ) \to [0, \infty] $ 
satisfy for all $ E \in \cB( [0,1] ) $ that 
$ \nu( E ) = \mu( {\bf f}( E ) )$.
\Nobs that $ {\bf f} $ and $ {\bf F} $ are bijective.
Moreover, \nobs that for all $ x \in [0,1] $, 
$\theta = ( \theta_1, \ldots, \theta_\fd) \in \R^{ \fd } $ 
it holds that
\begin{equation}
\label{theo:existence2:eq:realization}
\begin{split}
  \realization{ 
    {\bf F}( \theta ) 
  }( x ) 
&
= \theta_{ \fd } 
+
\smallsum_{ j = 1 }^{ \width } 
\theta_{ 2 \width + j } 
\max\{ (b-a) \theta_j x + \theta_{ \width + j } + a \theta_j , 0 \} \\
&
= \theta_{ \fd } 
+
\smallsum_{ j = 1 }^{ \width } 
\theta_{ 2 \width + j } 
\max\{ \theta_j {\bf f}( x ) + \theta_{ \width + j } , 0 \} 
  = \realization{\theta}( {\bf f}( x ) )
  .
\end{split}
\end{equation}
In addition, \nobs that for all $ x, y \in [0,1] $ we have that
\begin{equation}
  \abs{ g(x) - g(y) } 
= 
  \abs{ 
    f( {\bf f}( x ) ) 
    - 
    f( {\bf f}( y ) ) 
  } 
\le 
  L \abs{ {\bf f}( x ) - {\bf f}( y ) } 
= 
  L ( b - a ) \abs{ x - y } .
\end{equation}
\cref{theo:existence} \hence demonstrates that 
there exists $ \psi \in \R^\fd $ 
which satisfies 
$
  \operatorname{Lip}( \realization{\psi} ) 
  \le \width L ( b - a ) 
$, 
$
  \sup_{ x \in [0,1] } 
  \abs{ \realization{ \psi }( x ) } 
  \le \width L ( b - a ) + \sup_{ x \in [0,1] } \abs{ g(x) } 
$,
and
\begin{equation} 
\label{eq:global_infimum}
\cfadd{def:lip:const}
  \int_0^1 ( \realization{\psi}( x ) - g( x ) )^2 \, \nu( \d x ) 
  = \inf_{ \vartheta \in \R^\fd } 
  \br*{ 
    \int_0^1 ( \realization{ \vartheta }( x ) - g( x ) )^2 \, \nu( \d x ) 
  }
\end{equation}
\cfload.
In the following let $ \theta \in \R^{ \fd } $ 
satisfy 
$ 
  \theta = {\bf F}^{ - 1 }( \psi ) 
$. 
\Nobs that \cref{theo:existence2:eq:realization} 
and the integral transformation theorem assure 
for all $ \vartheta \in \R^{ \fd } $ that
\begin{equation}
\begin{split}
  \cL( \vartheta ) 
&
= 
  \int_a^b ( \realization{\vartheta}( x ) - f( x ) )^2 \, \mu( \d x ) 
=
  \int_0^1 
    ( \realization{\vartheta} ( {\bf f} ( x ) ) - f( {\bf f} ( x ) ) )^2 \, 
  \nu( \d x ) 
\\
&
  = 
  \int_0^1 ( \realization{ {\bf F } ( \vartheta ) }( x ) - g(x) )^2 \, \nu( \d x ) .
\end{split}
\end{equation}
Combining this with \cref{eq:global_infimum} 
and the fact that $ {\bf F} $ is bijective shows that
\begin{equation}
\begin{split}
  \cL( \theta ) 
&= 
  \int_0^1 ( \realization{\psi}( x ) - g ( x ) )^2 \, \nu( \d x ) 
= 
  \inf_{ \vartheta \in \R^{ \fd } } 
  \br*{ 
    \int_0^1 ( \realization{ \vartheta }( x ) - g( x ) )^2 \, \nu( \d x ) 
  } 
\\
&
= 
  \inf_{ \vartheta \in \R^{ \fd } } 
  \br*{ 
    \int_0^1 ( \realization{ {\bf F}( \vartheta ) }( x ) - g( x ) )^2 \, \nu( \d x ) 
  }
= 
  \inf_{ \vartheta \in \R^{ \fd } } \cL( \vartheta ) .
\end{split}
\end{equation}
In addition, \nobs that \cref{theo:existence2:eq:realization} 
ensures for all $ x \in [a,b] $ that
\begin{equation}
  \abs{ \realization{ \theta }( x ) } 
= 
  \abs{ \realization{ \psi }( {\bf f}^{ - 1 }( x ) ) } 
\le 
  \width L ( b - a ) + \sup\nolimits_{ y \in [0,1] } \abs{ g(y) } 
= 
  \width L ( b - a ) + \sup\nolimits_{ y \in [a,b] } \abs{ f(y) } .
\end{equation}
Finally, \nobs that \cref{theo:existence2:eq:realization} demonstrates
for all $ x, y \in [a,b] $ that
\begin{equation}
\begin{split}
  \abs{ 
    \realization{ \theta }( x ) - \realization{ \theta }( y ) 
  } 
& 
  = 
  \abs{ 
    \realization{\psi}( {\bf f}^{ - 1 }( x ) ) 
    - 
    \realization{\psi}( {\bf f}^{ - 1 }( y ) ) 
  }
\le 
  \width L ( b - a ) 
  \abs{ 
    {\bf f}^{ - 1 }( x ) - {\bf f}^{ - 1 }( y ) 
  } 
\\
& 
  = \width L \abs{ x - y } 
  .
\end{split}
\end{equation}
\end{cproof}

\subsection{Existence of regular global minima for shallow ANNs}
\label{subsec:existence_regular_minima}

In the final result of this section, \cref{cor:existence:regular}, 
we strengthen  
\cref{theo:existence2} 
by showing that there also exists a global minimizer of the risk function which admits a neighborhood on which
the risk function is continuously differentiable.
Furthermore, the gradient on this neighborhood can be obtained from a sequence of approximate realization functions using suitable differentiable approximations of the ReLU function, as outlined in the introduction.
The proof relies on regularity results from our previous
article Eberle et al.~\cite{EberleJentzenRiekertWeiss2021}.

\begin{cor} 
\label{cor:existence:regular}
Let $ \width, \fd \in \N $, $ L, a \in \R $, $ b \in (a, \infty) $, 
$ f \in C( [a,b], \R ) $ 
satisfy for all $ x, y \in [a,b] $ that 
$ \fd = 3 \width + 1 $
and 
$
  \abs{ f(x) - f(y) } \leq L \abs{ x - y }
$, 
let $ \dens \colon  [a,b]  \to [0, \infty)$ be bounded and measurable, 
let 
$
  \Rect_r \colon \R \to \R 
$, 
$ 
  r \in \N \cup \{ \infty \} 
$,
satisfy for all $ x \in \R $ that 
$ 
  ( \cup_{ r \in \N } \{ \Rect_r \} ) \subseteq C^1( \R, \R ) 
$,
$ 
  \Rect_{ \infty }( x ) = \max\{ x, 0 \} 
$,
$
  \sup_{ r \in \N }
  \sup_{ y \in [ - |x|, |x| ] } 
  | (\Rect_r)'( y ) | < \infty 
$, 
and
\begin{equation}
  \limsup\nolimits_{ r \to \infty 
  }( 
    | 
      \Rect_r(x) - \Rect_\infty( x ) 
    |
    +
    | 
      ( \Rect_r )'( x ) - \indicator{ (0,\infty) }( x ) 
    | 
  ) = 0 
  ,
\end{equation}
for every 
$ r \in \N \cup \cu{ \infty } $ 
let $ \cL_r \colon \R^{ \fd } \to \R $ 
satisfy for all 
$ \theta = ( \theta_1, \dots, \theta_{ \fd } ) \in \R^{ \fd } $ 
that 
\begin{equation}
\textstyle
  \cL_r( \theta ) = 
  \int_a^b 
  \displaystyle
    ( 
      f( x ) 
      - 
%       \realization{ \theta }_r( x ) 
      \theta_{ \fd } 
      -
      \smallsum_{ j = 1 }^{ \width } 
      \theta_{ 2 \width + j } 
      [ 
        \Rect_r( \theta_{ \width + j } + \theta_j x )
      ]
    )^2 
    \, \dens ( x ) \, 
  \d x
  ,
\end{equation}
let $ U \subseteq \R^{ \fd } $ satisfy
\begin{equation}
\label{eq:def_set_U_in_proof_regular_minimum}
  U = 
  \bigl\{
    \theta = ( \theta_1, \ldots, \theta_ \fd ) \in \R^{ \fd } \colon
    \bigl(
      \forall \, i \in \cu{ 1, 2, \ldots, \width } \colon 
      \abs{ \theta_i } + \abs{ \theta_{ \width + i } } > 0 
    \bigr)
  \bigr\}
  ,
\end{equation}
and let 
$
  \cG \colon \R^{ \fd } \to \R^{ \fd } 
$ 
satisfy for all 
$
  \theta \in 
  \cu{ 
    \vartheta \in \R^{ \fd } \colon 
    ( ( \nabla \cL_r )( \vartheta ) )_{ r \in \N } 
    \text{ is convergent} 
  } 
$ 
that 
$
  \cG( \theta ) = \lim_{ r \to \infty } ( \nabla \cL_r )( \theta )
$. 
Then 
\begin{enumerate}[label = (\roman*)]
\item 
\label{cor:regular_minimum_item_i}
it holds that
$ U $ 
is open, 
\item 
\label{cor:regular_minimum_item_ii}
it holds that
$
  ( \cL_{ \infty } )|_U 
  \in C^1( U, \R )
$, 
\item 
\label{cor:regular_minimum_item_iii}
it holds that
$
  U \ni \theta \mapsto ( \nabla \cL_{ \infty } )( \theta ) \in \R^{ \fd } 
$
is locally Lipschitz continuous,  
\item 
\label{cor:regular_minimum_item_iv}
it holds for all $ \theta \in U $ that 
$
  ( \nabla \cL_{ \infty } )( \theta ) = \cG( \theta ) 
$,
and
\item 
\label{cor:regular_minimum_item_v}
it holds that there exists $ \theta \in U $ 
such that 
$
  \cL_{ \infty }( \theta ) = 
  \inf_{ \vartheta \in \R^{ \fd } } 
  \cL_{ \infty }( \vartheta )
$, 
$ 
  \sup_{ x \in [a,b] } \abs{ \realization{ \theta }(x) } 
  \leq \width L ( b - a ) + \sup_{ x \in [a,b] } \abs{ f(x) } 
$, 
and 
$
  \sup_{ x, y \in [a,b], x \neq y }
  (
    | x - y |^{ - 1 }
    | \realization{\theta} ( x ) - \realization{\theta} ( y ) |
  )
  \leq \width L 
$.
\end{enumerate}
\end{cor}

\begin{cproof}{cor:existence:regular}
Throughout this proof 
for every $ \theta = ( \theta_1, \dots, \theta_{ \fd } ) \in \R^{ \fd } $,
$ r \in \N \cup \cu{ \infty } $ 
let 
$ \realization{ \theta }_r \colon \R \to \R $ 
satisfy for all $ x \in \R $ that
\begin{equation}
\textstyle
  \realization{ \theta }_r ( x ) 
  =
  \theta_{ \fd } 
  +
  \smallsum_{ j = 1 }^{ \width } 
  \theta_{ 2 \width + j } 
  [ 
    \Rect_r( \theta_{ \width + j } + \theta_j x )
  ]
  .
\end{equation}
\Nobs that \cref{eq:def_set_U_in_proof_regular_minimum} 
proves \cref{cor:regular_minimum_item_i}. 
\Nobs that Eberle et al.~\cite[Proposition 2.3]{EberleJentzenRiekertWeiss2021} 
establishes \cref{cor:regular_minimum_item_ii,cor:regular_minimum_item_iv}.
\Nobs that Eberle et al.~\cite[Corollary 2.7]{EberleJentzenRiekertWeiss2021} 
and \cref{cor:regular_minimum_item_iv}
prove \cref{cor:regular_minimum_item_iii}. 
\Nobs that  
\cref{theo:existence2} 
(applied with\footnote{Here and in the remainder of this article, when applying another theorem/lemma/proposition we use the notation $\with$ to indicate which values are assigned to the variables in the applied result.
	In this particular case,
	 \cref{theo:existence2}, where $\mu$ is an arbitrary finite measure on $[a,b]$, is applied with the specific measure $\cB ( [a,b]) \ni E \mapsto \int_E \dens ( x ) \, \d x \in [0, \infty ]$.}
$ 
  \mu \with (\cB ( [a,b]) \ni E \mapsto \int_E \dens ( x ) \, \d x \in [0, \infty ]) 
$ 
in the notation of \cref{theo:existence2}) 
proves that there exists 
$ \psi \in \R^{ \fd } $ 
which satisfies	
$ 
  \sup_{ x \in [a,b] } \abs{ \realization{ \psi }(x) } 
  \leq \width L ( b - a ) + \sup_{ x \in [a,b] } \abs{ f(x) } 
$, 
$
  \sup_{ x, y \in [a,b], x \neq y }
  (
    | x - y |^{ - 1 }
    | \realization{ \psi } ( x ) - \realization{ \psi } ( y ) |
  )
  \leq \width L 
$,
and 
\begin{equation}
\label{eq:existence_global_minima_in_proof}
\textstyle
  \cL_{ \infty }( \psi ) = 
  \inf_{ \vartheta \in \R^{ \fd } } 
  \cL_{ \infty }( \vartheta )
  .
\end{equation}
In the following let 
$
  \theta = ( \theta_1, \ldots, \theta_{ \fd } ) 
$ 
satisfy for all 
$
  i \in 
  \N \cap ( [1, \width] \cup (2 \width, \fd] )
$, 
$
  j \in \N \cap (\width, 2 \width]
$ 
that
\begin{equation}
\label{eq:proof_construction_of_theta}
  \theta_i = \psi_i 
\qqandqq
  \theta_j =
  \psi_j 
  - 
  \mathbbm{1}_{ \{ 0 \} }( 
    \abs{ \psi_{ j - \width } } + \abs{ \psi_j }
  )
  .
\end{equation}
% \Nobs that \cref{eq:proof_construction_of_theta} 
% establishes \cref{cor:regular_minimum_item_i}. 
\Nobs that \cref{eq:proof_construction_of_theta} 
shows
for all $ i \in \cu{ 1, 2, \ldots, \width } $, $ x \in \R $ 
that 
$
  \max\cu{
    \theta_{ \width + i } + \theta_i x , 0 
  }
  =
  \max\cu{
    \psi_{ \width + i } + \psi_i x , 0 
  } 
$. 
\Hence 
for all $ x \in \R $ 
that
$
  \realization{ \theta }_{ \infty }( x ) = \realization{ \psi }_{ \infty }( x )
$. 
Combining this with \cref{eq:existence_global_minima_in_proof} 
establishes \cref{cor:regular_minimum_item_v}. 
\end{cproof}

\section{Regularity analysis for generalized gradients in the training 
of deep ANNs}
\label{section:risk:diff}

In this section we introduce in \cref{setting:dnn} 
in \cref{subsection:dnn:framework} below 
our mathematical framework for deep ReLU ANNs. 
As in \cite{CheriditoJentzenRiekert2022,DNNReLUarXiv,JentzenRiekert2021} 
we approximate the ReLU activation function 
$ \Rect_{ \infty } \colon \R \to \R $ through 
continuously differentiable functions 
$  
  \Rect_r \colon \R \to \R 
$, 
$
  r \in [1, \infty)
$, 
in order to define an appropriate 
generalized gradient 
$ \cG \colon \R^{ \fd } \to \R^{ \fd } $ 
of the risk function 
$
  \cL_{ \infty } \colon \R^{ \fd } \to \R 
$; 
see \cref{lim_R} and \cref{eq:def_risk_function} 
in \cref{setting:dnn}.

In \cref{prop:G} in \cref{ssec:properties} below (explicit representation 
and pointwise approximations for $ \cG $), 
in \cref{lem:realization:lip} in \cref{subsection:risk:local:lip} below (local Lipschitz 
continuity of $ \cL_{ \infty } $), 
and 
in \cref{lem:approx:gradient:bounded} in \cref{subsection:risk:weak:diff} below 
(uniform local boundedness for $ \nabla \cL_r $, $ r \in [1,\infty) $)
we then state several important regularity properties of 
the risk function 
$
  \cL \colon \R^{ \fd } \to \R
$
and its generalized gradient function 
$ \cG \colon \R^{ \fd } \to \R^{ \fd } $. 
\cref{prop:G} is proved in 
Hutzenthaler et al.~\cite[Theorem~2.9]{DNNReLUarXiv}, 
\cref{lem:realization:lip} follows, e.g., from 
Hutzenthaler et al.~\cite[Lemma~2.10]{DNNReLUarXiv}, 
and \cref{lem:approx:gradient:bounded} is a consequence from 
Hutzenthaler et al.~\cite[Lemma~3.6]{DNNReLUarXiv}.

In \cref{cor:weak_derivative} in \cref{subsection:risk:weak:diff} 
we show that the risk function $ \cL_{ \infty } \colon \R^{ \fd } \to \R $ 
is weakly differentiable with the generalized gradient function 
$ \cG \colon \R^{ \fd } \to \R^{ \fd } $ serving as a weak gradient function.

In \cref{prop:risk:ae:diff} 
in \cref{subsection:risk:strong:diff} below 
we establish that the risk function 
$ \cL_{ \infty } \colon \R^{ \fd } \to \R $ 
is differentiable Lebesgue almost everywhere 
with its gradients agreeing Lebesgue almost everywhere 
with the generalized gradient function $ \cG \colon \R^{ \fd } \to \R^{ \fd } $.
% and that its gradient agrees with the generalized gradient $\cG$ almost everywhere. 
% 
% 
% 
Our proof of \cref{prop:risk:ae:diff} relies on 
\cref{lem:realization:lip}, 
\cref{cor:weak_derivative}, 
and
well-known results on weak derivatives of locally Lipschitz continuous functions (cf.~Evans~\cite{Evans2010}).

In \cref{sec:subdifferential} below we gather several known notions and lemmas 
regarding Fr\'{e}chet subdifferentials.
In particular, in the scientific literature \cref{def:limit:subdiff} 
can be found, e.g., as Rockafellar \& Wets~\cite[Definition 8.3]{RockafellarWets1998} 
and Bolte et al.~\cite[Definition 2.10]{BolteDaniilidis2006},
\cref{item:subdiff_item_iii,item:subdiff_item_iv,item:subdiff_item_v} 
in \cref{lem:subdifferential:c1} are proved, e.g., 
as \cite[Theorem 8.6 and Exercise 8.8]{RockafellarWets1998},
and
\cref{lem:limiting_derivatives} is a reformulation of the well-known fact that the limiting Fr\'{e}chet subdifferential 
of a continuous function has a closed graph 
(see, e.g.,~\cite[Proposition 8.7]{RockafellarWets1998}).

Finally, in \cref{prop:loss:gradient:subdiff} 
in \cref{subsection:gen:grad:frechet} below (the main result of this section) 
we establish that for every ANN parameter vector 
$ \theta \in \R^{ \fd } $ 
it holds that the generalized gradient $ \cG( \theta ) $ 
is a limiting subgradient of the risk function $ \cL_\infty $ at $ \theta $. 
Our proof of \cref{prop:loss:gradient:subdiff} relies 
on \cref{prop:risk:ae:diff},
on \cref{lem:subdifferential:c1}, 
as well as
on the continuity type result 
for the generalized gradient function 
in \cref{lem:gradient:approx:sequence} 
in \cref{subsection:gen:grad:continuous}. 
Our proof of \cref{lem:gradient:approx:sequence}, in turn, 
is based on 
local underestimate type result in 
\cref{lem:gradient:left:approx} in \cref{subsection:local:under}
as well as on 
the conditional continuity result for the generalized gradient function 
in \cref{lem:gradient:convergence} 
in \cref{subsection:gen:grad:continuous} below.

\subsection{Mathematical framework for deep ANNs with ReLU activation}
\label{subsection:dnn:framework}

\begin{setting}
\label{setting:dnn}
Let 
$ a \in \R $, $ b \in [a,\infty) $, 
$ \scrA \in (0,\infty) $, $ \scrB \in ( \scrA, \infty) $, 
$ ( \ell_k )_{ k \in \N_0 } \subseteq \N $, 
$ L, \fd \in \N $
satisfy 
$
  \mathfrak{d} = \sum_{ k = 1 }^L \ell_k ( \ell_{ k - 1 } + 1 ) 
$, 
for every 
$ 
  \theta = ( \theta_1, \dots, \theta_{ \fd } ) \in \R^{ \fd } 
$
let 
$ 
  \fw^{ k, \theta } = 
  ( \fw^{ k, \theta }_{ i, j } )_{ 
    (i,j) \in \{ 1, \ldots, \ell_k \} \times \{ 1, \ldots, \ell_{ k - 1 } \} 
  }
  \in \R^{ \ell_k \times \ell_{ k - 1 } }
$, 
$
  k \in \N 
$, 
and 
$
  \fb^{ k, \theta } 
  = 
  ( \fb^{ k, \theta }_1, \dots, \fb^{ k, \theta }_{ \ell_k} )
  \in \R^{ \ell_k } 
$,
$ k \in \N $, 
satisfy for all 
$ k \in \{ 1, \dots, L \} $, 
$ i \in \{ 1, \ldots, \ell_k \} $,
$ j \in \{ 1, \ldots, \ell_{ k - 1 } \} $ 
that
\begin{equation}
\label{wb}
  \fw^{ k, \theta }_{ i, j }
  = 
  \theta_{ ( i - 1 ) \ell_{ k - 1 } + j 
  + 
  \sum_{ h = 1 }^{ k - 1 } \ell_h ( \ell_{ h - 1 } + 1 ) }
\qqandqq
  \fb^{ k, \theta }_i 
  =
  \theta_{ \ell_k \ell_{ k - 1 } + i 
  + 
  \sum_{ h = 1 }^{ k - 1 } \ell_h ( \ell_{ h - 1 } + 1 ) } 
  ,
\end{equation}
for every 
$ k \in \N $, 
$ \theta \in \R^{ \fd } $ 
let 
$
  \cA_k^{ \theta }
  =
  ( 
    \cA_{ k, 1 }^{ \theta }, \ldots, \cA_{ k, \ell_k }^{ \theta } 
  )
  \colon \R^{ \ell_{ k - 1 } } \to \R^{ \ell_k } 
$
satisfy 
for all 
$ x \in \R^{ \ell_{ k - 1 } } $
that 
\begin{equation}
\label{eq:def_Ak_transformation_deep_ANNs}
  \cA_k^{ \theta }( x ) 
  = 
  \fb^{ k, \theta } + \fw^{ k, \theta } x 
  ,
\end{equation}
let 
$ \Rect_r \colon \R \to \R $, 
$ r \in [1, \infty] $,
satisfy for all 
$ r \in [1, \infty) $, 
$ x \in ( - \infty, \scrA r^{ - 1 } ] $, 
$ 
  y \in \R
$, 
$
  z \in [ \scrB r^{ - 1 }, \infty ) 
$
that 
\begin{equation}
\label{lim_R}
  \Rect_r \in C^1( \R, \R ) ,
\quad
  \Rect_r(x) = 0,
\quad 
  0 \leq \Rect_r(y) \leq \Rect_{ \infty }( y ) 
  = 
  \max\{ y, 0 \}
  ,
\qandq
  \Rect_r(z) = z
  ,
\end{equation}
assume 
$
  \sup_{ r \in [1, \infty) }
  \sup_{ x \in \R } 
  | ( \Rect_r )'( x ) | < \infty 
$, 
for every 
$ r \in [ 1, \infty ] $, 
$ k \in \N $
let 
$ \mathfrak{M}_{ r, k } \colon \R^{ \ell_k } \to \R^{ \ell_k } $ 
satisfy for all
$ x = ( x_1, \ldots, x_{ \ell_k } ) \in \R^{ \ell_k } $ 
that 
\begin{equation}
\textstyle
  \mathfrak{M}_{ r, k }( x ) 
  = ( \Rect_r( x_1 ), \ldots, \Rect_r( x_{ \ell_k } ) )
  ,
\end{equation}
for every 
$ \theta \in \R^{ \fd } $
let 
$
  \mathcal{N}^{ k, \theta }_r 
  = 
  ( 
    \mathcal{N}^{ k, \theta }_{ r, 1 }, \ldots, \mathcal{N}^{ k, \theta }_{ r, \ell_k } 
  )
  \colon \R^{ \ell_0 } \to \R^{ \ell_k } 
$, 
$ k \in \N $, 
$ r \in [1,\infty] $, 
and 
$ 
  \mathcal{X}^{ k, \theta }_i \subseteq \R^{ \ell_0 } 
$, 
$ k, i \in \N $,
satisfy 
for all 
$ k \in \N $, 
$ r \in [1,\infty] $,
$ i \in \{ 1, \ldots, \ell_k \} $
that
\begin{equation}
\label{eq:def_NN_realization}
  \mathcal{N}^{ 1, \theta }_r = \cA^{ \theta }_1, 
  \quad
  \mathcal{N}^{ k + 1, \theta }_r 
  = \cA_{ k + 1 }^{ \theta } \circ 
  \mathfrak{M}_{ r^{ 1 / k}, k } \circ 
  \mathcal{N}^{ k, \theta }_r ,
  \quad 
  \text{and}
  \quad 
  \mathcal{X}^{ k, \theta }_i 
  = 
  \{ 
    x \in[a,b]^{ \ell_0 } \colon
    \mathcal{N}^{ k, \theta }_{ \infty, i }( x ) > 0 
  \}
  ,
\end{equation}
let
$ 
  f = 
  ( f_1, \ldots, f_{ \ell_L } ) 
  \colon [a,b]^{ \ell_0 } \to \R^{ \ell_L } 
$
be measurable,
let $ \mu \colon \cB( [a,b]^{ \ell_0 } ) \to [0, \infty] $ be a finite measure, 
for every 
$ r \in [1, \infty] $ 
let $ \cL_r \colon \R^{ \fd } \to \R $ 
satisfy 
for all $ \theta \in \R^{ \fd } $
that
\begin{equation}
\label{eq:def_risk_function}
\textstyle 
  \cL_r( \theta ) 
  = 
  \int_{ [a,b]^{ \ell_0 } } 
  \norm{ \mathcal{N}_r^{ L, \theta }( x ) - f(x) }^2 \, \mu( \d x )
  ,
\end{equation}
and let 
$
  \cG = ( \cG_1, \ldots, \cG_{ \fd } ) \colon \R^{ \fd } \to \R^{ \fd } 
$ 
satisfy
for all 
$
  \theta \in 
  \{
    \vartheta \in \R^{ \fd } \colon 
    ( ( \nabla\cL_r )( \vartheta ) )_{ r \in [1,\infty) }
    \text{ is convergent} 
  \}
$ 
that
$
  \cG( \theta ) = \lim_{ r \to \infty }( \nabla\cL_r )( \theta ) 
$.
\end{setting}

\subsection{Explicit representations for the generalized gradient function}
\label{ssec:properties}

In \cref{prop:G} we show that the approximating sequence of gradients $(\nabla \cL_r ) ( \theta )$, $r \in [1, \infty )$,
in \cref{setting:dnn} converges for every $\theta \in \R^\fd$.
Furthermore, we derive
in \cref{G4',G4''} explicit formulae for the limit $\cG ( \theta )$.
This explicit representation of $\cG$ agrees with the standard generalized gradient obtained by formally defining the derivative of the ReLU as the left derivative $\indicator{(0 , \infty )}$ and applying the chain rule.

\begin{prop}
\label{prop:G}
Assume \cref{setting:dnn} and let 
$
  \theta = ( \theta_1, \ldots, \theta_{ \fd } ) \in \R^{ \fd }
$. 
Then
\begin{enumerate}[label = (\roman*)]
\item
\label{G1}
it holds for all $ r \in [1,\infty) $ that $ \cL_r \in C^1( \R^{ \fd }, \R ) $,
\item 
\label{G3}
it holds that 
$
  \limsup\nolimits_{ r \to \infty }
  \left(
    \abs{ 
      \cL_r( \theta ) - \cL_{ \infty }( \theta ) 
    } 
    +
    \norm{ 
      ( \nabla\cL_r )( \theta ) - \cG( \theta ) 
    } 
  \right) 
  = 0 
$,
\item  
\label{G4'}
it holds for all 
$ k \in \{ 1, \ldots, L \} $,
$ i \in \{ 1, \ldots, \ell_k \} $, 
$ j \in \{ 1, \ldots, \ell_{ k - 1 } \} $ 
that
\begin{equation}
\begin{split}
&
  \cG_{ 
    ( i - 1 ) \ell_{ k - 1 } + j 
    + 
    \sum_{ h = 1 }^{ k - 1 } \ell_h ( \ell_{ h - 1 } + 1 ) 
  }( \theta ) 
\\
&
  =
      \sum_{ 
        \substack{
          v_k, v_{ k + 1 }, \ldots, v_L \in \N, 
        \\
          \forall \, w \in \N \cap [k,L] \colon 
          v_w \leq \ell_w
        }
      }
  \int_{ [a,b]^{ \ell_0 } }
    2
    \br*{
      \Rect_{ \infty }(
        \cN^{ \max\{ k - 1, 1 \}, \theta }_{ \infty, j }( x ) 
      ) 
      \indicator{ (1, L] }( k )
      + 
      x_j \indicator{ \cu{1} }( k ) 
    }
%     \Biggl[
\\
\label{G_w}
&
        \cdot
%       \biggl[
        \Bigl[
          \indicator{ \{ i \} }( v_k )
        \Bigr]
        \Bigl[
          \cN_{ 
            \infty, v_L
          }^{ L, \theta }( x ) 
          - 
          f_{ v_L }( x ) 
        \Bigr]
        \Bigl[
          \textstyle
          \prod_{ n = k + 1 }^L
          \bigl(
            \fw^{ n, \theta }_{ v_n, v_{ n - 1 } }
            \indicator{
              \cX^{ n - 1, \theta }_{ v_{ n - 1 } } 
            }( x )
%             \indicator{ [0, \ell_n] }( v_n )
          \bigr)
        \Bigr]
%       \biggr]
%     \Biggr] 
  \,
  \mu( \d x ) ,
\end{split}
\end{equation}
and
\item  
\label{G4''}
it holds for all 
$ k \in \{ 1, \ldots, L \} $,
$ i \in \{ 1, \ldots, \ell_k \} $ that
\begin{equation}
\label{G_b}
\begin{split}
&
  \cG_{ 
    \ell_k \ell_{ k - 1 } + i 
    + 
    \sum_{ h = 1 }^{ k - 1 } 
    \ell_h ( \ell_{ h - 1 } + 1 ) 
  }( \theta )
% \\&
  =
      \sum_{ 
        \substack{
          v_k, v_{ k + 1 }, \ldots, v_L \in \N, 
        \\
          \forall \, w \in \N \cap [k,L] \colon 
          v_w \leq \ell_w
        }
      }
  \int_{ [a,b]^{ \ell_0 } } 
    2
        \Bigl[
          \indicator{ \{ i \} }( v_k )
        \Bigr]
\\&
%         \quad 
        \cdot 
        \Bigl[
          \cN^{ L, \theta }_{ \infty, v_L }( x ) 
          - 
          f_{ v_L }( x )
        \Bigr]
        \Bigl[
          \textstyle{ \prod }_{ n = k + 1 }^L 
          \bigl(
            \fw^{ n, \theta }_{ v_n, v_{ n - 1 } }
            \indicator{ 
              \cX^{ n - 1, \theta }_{ v_{ n - 1 } } 
            }( x )
          \bigr)
        \Bigr]
  \, 
  \mu( \d x) .
\end{split}
\end{equation}
\end{enumerate}
\end{prop}
\begin{cproof}{prop:G}
\Nobs that \cite[Items~(i), (iv), (v), and (vi) in Theorem~2.9]{DNNReLUarXiv} 
establishes 
\cref{G1,G3,G4',G4''}
\end{cproof}

\subsection{Local Lipschitz continuity of the risk function}
\label{subsection:risk:local:lip}

\begin{lemma}
\label{lem:realization:lip}
Assume \cref{setting:dnn} and let $ K \subseteq \R^{ \fd } $ be compact.
Then there exists $ \fC \in \R $ such that 
for all $ \theta, \vartheta \in K $ it holds that
\begin{equation}
\label{eq:local_Lip_estimate}
  \abs{
    \cL_{ \infty }( \theta ) - \cL_{ \infty }( \vartheta ) 
  } 
  +
  \bigl( 
    \sup\nolimits_{ x \in [a, b]^{ \ell_0 } } 
    \norm{ 
      \cN_{ \infty }^{ L, \theta }( x ) - \cN_{ \infty }^{ L, \vartheta }( x ) 
    } 
  \bigr) 
  \le 
  \fC \norm{ \theta - \vartheta } .
\end{equation}
\end{lemma}
\begin{cproof}{lem:realization:lip}
\Nobs that, e.g., \cite[Lemma 2.10]{DNNReLUarXiv} 
establishes \cref{eq:local_Lip_estimate}. 
\end{cproof}

\subsection{Weak differentiability properties of the risk function}
\label{subsection:risk:weak:diff}

\begin{lemma} 
\label{lem:approx:gradient:bounded}
Assume \cref{setting:dnn} and let $ K \subseteq \R^\fd $ be non-empty and compact. 
Then
\begin{equation}
\label{eq:a_priori_bound}
  \sup\nolimits_{ \theta \in K 
  } 
  \sup\nolimits_{ r \in [1,\infty)
  } 
  \rbr*{ 
    \abs{ 
      \cL_r( \theta ) 
    } 
    +
    \abs{ 
      \cL_{ \infty }( \theta )
    }
    + 
    \norm{ 
      ( \nabla \cL_r )( \theta ) 
    } 
    +
    \norm{
      \cG( \theta )
    }
  } 
  < \infty .
\end{equation}
\end{lemma}
\begin{cproof}{lem:approx:gradient:bounded}
\Nobs that \cite[Lemma 3.6]{DNNReLUarXiv} and 
\cref{G1} in \cref{prop:G} 
show for all $ s \in (0,\infty) $ that
$
  \sup_{ 
    \theta \in \{ \vartheta \in \R^{ \fd } \colon \| \vartheta \| \leq s \}
  } 
  \sup_{ 
    r \in [1,\infty)
  } 
  \norm{ 
    ( \nabla \cL_r )( \theta ) 
  } 
  < \infty 
$. 
The fundamental theorem of calculus 
and the fact that for all $ r \in [1,\infty) $
it holds that 
$
  \cL_r( 0 ) = \cL_{ \infty }( 0 )
$
\hence 
demonstrate that for all $ s \in (0,\infty) $
we have that 
$
  \sup_{ 
    \theta \in \{ \vartheta \in \R^{ \fd } \colon \| \vartheta \| \leq s \}
  } 
  \sup_{ 
    r \in [1,\infty)
  } 
  \rbr*{ 
    \abs{ 
      \cL_r( \theta ) 
    } 
    +
    \norm{ 
      ( \nabla \cL_r )( \theta ) 
    } 
  } 
  < \infty 
$. 
Combining this with 
\cref{G3} in \cref{prop:G} 
establishes \cref{eq:a_priori_bound}.
\end{cproof}

As a consequence of \cref{prop:G} and the uniform boundedness result in \cref{lem:approx:gradient:bounded} we obtain in \cref{cor:weak_derivative} that the generalized gradient $\cG$ serves as a weak gradient of the risk function $\cL_\infty $.

\begin{cor}[Weak differentiability]
\label{cor:weak_derivative}
Assume \cref{setting:dnn}, 
let 
$ 
  \varphi = 
  ( 
    \varphi( \theta ) 
  )_{ 
    \theta = ( \theta_1, \dots, \theta_{ \fd } ) \in \R^{ \fd }
  } 
  \colon \R^{ \fd } \to \R
%   \in C^{ \infty }( \R^n , \R ) 
$ 
be compactly supported and continuously differentiable, 
and let 
% $ \theta = ( \theta_1, \dots, \theta_{ \fd } ) \in \R^{ \fd } $, 
$ i \in \{ 1, 2, \dots, \fd \} $. 
Then 
$
  \int_{ \R^{ \fd } } 
  |
    \cL_{ \infty }( \theta ) 
    \,
    ( \tfrac{ \partial }{ \partial \theta_i } \varphi )( \theta ) 
  |
  +
  |
    \cG_i( \theta ) 
    \,
    \varphi( \theta )
  |
  \, 
  \d \theta 
  < \infty
$
and 
\begin{equation}
\label{eq:weak_derivative_of_risk_function}
  \int_{ \R^{ \fd } } 
  \cL_{ \infty }( \theta ) 
  \,
  ( \tfrac{ \partial }{ \partial \theta_i } \varphi )( \theta ) \, \d \theta 
  =
  -
  \int_{ \R^{ \fd } } 
  \cG_i( \theta ) 
  \,
  \varphi( \theta )
  \, 
  \d \theta 
  .
\end{equation}
\end{cor}
\begin{cproof}{cor:weak_derivative}
\Nobs that the assumption that $ \varphi $ has a compact support 
ensures that there exists $ R \in (0,\infty) $ which satisfies 
for all 
$ \theta \in \R^{ \fd } \backslash [ - R, R ]^{ \fd } $ 
that
\begin{equation}
\label{eq:K_compact_set_defining_property}
  \varphi( \theta ) = 0 .
\end{equation}
\Nobs that \cref{lem:approx:gradient:bounded}
demonstrates that 
\begin{equation}
\label{eq:cG_cL_local_boundedness}
  \sup\nolimits_{ \theta \in [-R,R]^{ \fd } } 
  \sup\nolimits_{ r \in [1,\infty) }
  \bigl(
    | \cL_r( \theta ) | 
    +
    | \cL_{ \infty }( \theta ) |
    +
    \| ( \nabla \cL_r )( \theta ) \|
    +
    \| \cG( \theta ) \|
  \bigr)
  < \infty .
\end{equation}
This and \cref{eq:K_compact_set_defining_property} 
assure that for all $ r \in [1,\infty) $, 
$ \theta = ( \theta_1, \dots, \theta_{ \fd } ) \in \R^{ \fd } $ it holds that
\begin{equation}
\label{eq:Lebesgues_Majorante}
\begin{split}
&
  |
    \cL_r( \theta ) 
    ( \tfrac{ \partial }{ \partial \theta_i } \varphi )( \theta ) 
  |
  +
  |
    ( \tfrac{ \partial }{ \partial \theta_i } \cL_r )( \theta ) 
    \varphi( \theta ) 
  |
\\ &
  \leq 
  \bigl[ 
    \sup\nolimits_{ \vartheta \in [ - R, R ]^{ \fd } } 
    \sup\nolimits_{ s \in [1,\infty) }
    \bigl(
      | \cL_s( \vartheta ) | 
      +
      \| ( \nabla \cL_s )( \vartheta ) \|
      +    
      | \varphi( \vartheta ) | 
      +    
      \| ( \nabla \varphi )( \vartheta ) \| 
    \bigr)
  \big]
  \mathbbm{1}_{ [ - R, R ]^{ \fd } }( \theta )
  < \infty 
  .
\end{split}
\end{equation}
\cref{prop:G}, Lebesgue's dominated convergence theorem, 
\cref{eq:K_compact_set_defining_property}, and 
integration by parts \hence ensure that  
\begin{equation}
\begin{split}
&
  \int_{ \R^{ \fd } } 
  \cL_{ \infty }( \theta ) 
  \, 
  ( \tfrac{ \partial }{ \partial \theta_i } \varphi )( \theta ) 
  \, \d \theta 
  = 
  \lim_{ r \to \infty } 
  \left[ 
    \int_{ \R^{ \fd } } 
    \cL_r( \theta ) 
    \,
    ( \tfrac{ \partial }{ \partial \theta_i } \varphi )( \theta ) 
    \, \d \theta 
  \right]
\\ &
  = 
  \lim_{ r \to \infty } 
  \left[ 
    \int_{ [ - R, R ]^{ \fd } } 
    \cL_r( \theta ) 
    \,
    ( \tfrac{ \partial }{ \partial \theta_i } \varphi )( \theta ) 
    \, \d \theta 
  \right]
  = 
  - 
  \left(
    \lim_{ r \to \infty } 
    \left[ 
      \int_{ [ - R, R ]^{ \fd } } 
      ( \tfrac{ \partial }{ \partial \theta_i } \cL_r)( \theta ) 
      \, \varphi( \theta ) \, \d \theta 
    \right]
  \right)
  .
\end{split}
\end{equation}
\cref{prop:G}, 
\cref{eq:Lebesgues_Majorante}, 
% \cref{lem:approx:gradient:bounded}, 
and 
Lebesgue's dominated convergence theorem 
\hence show that
\begin{equation}
\begin{split}
  \int_{ \R^{ \fd } } 
  \cL_{ \infty }( \theta ) 
  \, 
  ( \tfrac{ \partial }{ \partial \theta_i } \varphi )( \theta ) 
  \, \d \theta 
  = 
  - 
  \int_{ 
    [ - R, R ]^{ \fd }
  } 
    \Bigl[
      \lim_{ r \to \infty } 
      ( \tfrac{ \partial }{ \partial \theta_i } \cL_r )( \theta ) 
    \Bigr]
    \varphi( \theta ) 
  \, \d \theta  
  = 
  - \int_{ \R^\fd } \cG_i( \theta ) \, \varphi ( \theta ) \, \d \theta .
\end{split}
\end{equation}
This, \cref{eq:K_compact_set_defining_property}, 
and 
\cref{eq:cG_cL_local_boundedness} 
establish \cref{eq:weak_derivative_of_risk_function}. 
\end{cproof}

\subsection{Strong differentiability properties of the risk function}
\label{subsection:risk:strong:diff}

We next establish in \cref{prop:risk:ae:diff} that the risk function is a.e.~strongly differentiable with gradient $\cG$. The proof relies on \cref{cor:weak_derivative},
the local Lipschitz continuity result in \cref{lem:realization:lip},
Rademacher's theorem,
and the fact that locally Lipschitz continuous functions are weakly differentiable with the weak gradient a.e.~equal to the strong gradient (cf.~Evans~\cite{Evans2010}).

\begin{prop} 
\label{prop:risk:ae:diff}
Assume \cref{setting:dnn}. Then there exists 
$ E \in \cB( \R^{ \fd } ) $ such that 
\begin{enumerate}[label = (\roman*)]
\item 
\label{item:lem_strong_diff_1}
it holds that 
$
  \int_{ \R^{ \fd } \backslash E } 1 \, \d \theta = 0
$, 
\item 
\label{item:lem_strong_diff_2}
it holds for all $ \theta \in E $ that 
$
  \cL_{ \infty }
$
is differentiable at $ \theta $, 
and 
\item 
\label{item:lem_strong_diff_3}
it holds for all $ \theta \in E $ that
$
  ( \nabla \cL_{ \infty } )( \theta ) = \cG( \theta )
$. 
\end{enumerate}
\end{prop}
\begin{cproof}{prop:risk:ae:diff}
Throughout this proof let 
$
  G = ( G_1, \dots, G_{ \fd } ) \colon \R^{ \fd } \to \R^{ \fd }
$
satisfy for all $ \theta \in \R^{ \fd } $ that
\begin{equation}
\label{eq:def_G_derivative}
  G( \theta )
  =
  \begin{cases}
    ( \nabla \cL_{ \infty } )( \theta )
  &
    \colon 
    \cL_{ \infty } \text{ is differentiable at } \theta 
  \\
    0
  &
    \colon
    \cL_{ \infty } \text{ is not differentiable at } \theta 
  \end{cases}
  .
\end{equation}
\Nobs that 
\cref{eq:def_G_derivative}, 
the fact that for all measurable 
$ g_n \colon \R^{ \fd } \to \R^{ \fd } $, $ n \in \N $, 
it holds that 
$
  \{ 
    \theta \in \R^{ \fd } \colon 
    ( g_n( \theta ) )_{ n \in \N }
    \text{ is a Cauchy sequence}
  \}
$
is measurable, 
and 
the fact that for all measurable and 
pointwise convergent 
$ g_n \colon \R^{ \fd } \to \R^{ \fd } $, $ n \in \N $, 
it holds that 
$
  \R^{ \fd } \ni \theta \mapsto \lim_{ n \to \infty } g_n( \theta ) \in \R^{ \fd }
$
is measurable establish that 
$ 
  G 
$
is measurable. 
\Moreover \cref{lem:realization:lip} 
ensures that $ \cL_{ \infty } $ is locally Lipschitz continuous.
Rademacher's theorem (cf.~Evans~\cite[Theorem 5.8.6]{Evans2010})
\hence demonstrates that there exists 
$ 
  \cE 
  \in \{ 
    A \in \cB( \R^{ \fd } ) \colon 
    \int_{ \R^{ \fd } \backslash A } 1 \, \d \theta = 0
  \} 
$ 
which satisfies for all $ \theta \in \cE $ that 
$ \cL_{ \infty } $ is differentiable at $ \theta $. 
\cref{lem:realization:lip}, 
Evans~\cite[Theorems 5.8.4 and 5.8.5]{Evans2010},
and 
\cref{eq:def_G_derivative}
\hence show
for all compactly supported $ \varphi \in C^{ \infty }( \R^{ \fd }, \R ) $
and all $ i \in \{ 1, 2, \dots, \fd \} $ 
that 
$
  \int_{ \R^{ \fd } } 
  |
    \cL_{ \infty }( \theta ) 
    \,
    ( \tfrac{ \partial }{ \partial \theta_i } \varphi )( \theta ) 
  |
  +
  |
    G_i( \theta ) 
    \,
    \varphi( \theta )
  |
  \, 
  \d \theta 
  < \infty
$
and
\begin{equation}
\label{eq:Rademacher_consequence}
  \int_{ \R^{ \fd } } 
  \cL_{ \infty }( \theta ) 
  \,
  ( \tfrac{ \partial }{ \partial \theta_i } \varphi )( \theta ) \, \d \theta 
  =
  -
  \int_{ \R^{ \fd } } 
  G_i( \theta ) 
  \,
  \varphi( \theta )
  \, 
  \d \theta 
  .
\end{equation}
Combining this with \cref{cor:weak_derivative} and the fundamental lemma of calculus of variations (cf., e.g., Hörmander~\cite[Theorem 1.2.5]{Hormander2003})
implies that
there exists 
$ 
  {\bf E} 
  \in \{ 
    A \in \cB( \R^{ \fd } ) \colon 
    \int_{ \R^{ \fd } \backslash A } 1 \, \d \theta = 0
  \} 
$ 
which satisfies for all $ \theta \in {\bf E} $ that 
\begin{equation}
\label{eq:G_and_cG_agree}
  G( \theta ) = \cG( \theta )
  .
\end{equation}
\Nobs that 
\cref{eq:def_G_derivative}, 
\cref{eq:G_and_cG_agree}, 
and the fact that for all 
$ \theta \in \cE $ 
it holds that 
$ \cL_{ \infty } $ 
is differentiable at $ \theta $ 
assure that for all 
$ \theta \in ( {\bf E} \cap \cE ) $
it holds that
\begin{equation}
\label{eq:cG_on_E_cap_cE}
  \cG( \theta ) = ( \nabla \cL_{ \infty } )( \theta ) .  
\end{equation}
\Moreover the fact that 
$
  {\bf E} 
  \in \{ 
    A \in \cB( \R^{ \fd } ) \colon 
    \int_{ \R^{ \fd } \backslash A } 1 \, \d \theta = 0
  \} 
$
and the fact that 
$
  \cE
  \in \{ 
    A \in \cB( \R^{ \fd } ) \colon 
    \int_{ \R^{ \fd } \backslash A } 1 \, \d \theta = 0
  \} 
$
ensure that 
$
  ( {\bf E} \cap \cE ) 
  \in \{ 
    A \in \cB( \R^{ \fd } ) \colon 
    \int_{ \R^{ \fd } \backslash A } 1 \, \d \theta = 0
  \} 
$. 
Combining this
and the fact that for all 
$ \theta \in ( {\bf E} \cap \cE ) $
it holds that 
$
  \cL_{ \infty }
$
is differentiable at $ \theta $ 
with \cref{eq:cG_on_E_cap_cE} 
establishes \cref{item:lem_strong_diff_1,item:lem_strong_diff_2,item:lem_strong_diff_3}. 
\end{cproof}

\subsection{Fr\'{e}chet subdifferentials and limiting Fr\'{e}chet subdifferentials}
\label{sec:subdifferential}

\begin{definition}[Fr\'{e}chet subdifferentials and limiting Fr\'{e}chet subdifferentials]
\label{def:limit:subdiff}
Let $ n \in \N $, $ f \in C( \R^n, \R) $, $ x \in \R^n $.
Then we denote by 
$ (\cD f)(x) \subseteq \R^n $ 
the set given by
\begin{equation}
  ( \cD f)( x )
  = 
  \cu*{ 
    y \in \R^n \colon 
    \left[
      \liminf_{\R^n \backslash \cu{  0 } \ni h \to 0 } 
      \rbr*{ \frac{f(x + h ) - f ( x ) - \spro{y , h } }{\norm{h}} } \geq 0  
    \right]
  } 
\end{equation}
and we denote by 
$
  (\bD f)(x) \subseteq \R^n 
$ 
the set given by
\begin{equation}
\label{def:limit:subdiff:eq}
  (\bD f)( x ) =
  \textstyle\bigcap_{ \varepsilon \in (0, \infty) } 
  \overline{
    \br*{
      \cup_{ 
        y \in 
        \cu{ 
          z \in \R^n \colon \norm{ x - z } < \varepsilon 
        }
      } 
      ( \cD f )( y ) 
    } 
  }
  .
\end{equation}
\end{definition}

\cfclear
\begin{lemma}[Properties of Fr\'{e}chet subdifferentials]
\label{lem:subdifferential:c1}
\cfadd{def:limit:subdiff}
Let $ n \in \N $, $ f \in C( \R^n , \R ) $. 
Then 
\begin{enumerate}[label = (\roman*)]
\item 
\label{item:subdiff_item_i}
it holds for all $ x \in \R^n $ that
\begin{multline} 
\label{lem:subdiff:eq} 
\cfadd{def:limit:subdiff}
  ( \bD f )( x ) 
  = 
  \bigl\{ 
    y \in \R^n \colon 
    \bigl[
      \exists \, z = (z_1, z_2) \colon \N \to \R^n \times \R^n \colon 
      \bigl( 
        \br[\big]{
          \forall \, k \in \N \colon z_2( k ) \in ( \cD f )( z_1(k) ) 
        } 
\\ 
        \wedge
          \br[\big]{
            \limsup\nolimits_{ k \to \infty } 
            ( 
              \norm{ z_1(k) - x } + \norm{ z_2(k) - y } 
            ) 
            = 0
          } 
      \bigr) 
    \bigr]
  \bigr\}
  ,
\end{multline}
\item 
\label{item:subdiff_item_ii}
it holds for all $ x \in \R^n $ that
$ 
  ( \cD f)( x )
  \subseteq 
  ( \bD f)( x )
$,
\item 
\label{item:subdiff_item_iii}
it holds for all 
$ x \in \{ y \in \R^n \colon f \text{ is differentiable at } y \} $
that
$
  ( \cD f )( x ) 
  = \cu{ ( \nabla f )( x ) }  
$,
\item 
\label{item:subdiff_item_iv}
it holds for all 
$
  x \in \cup_{ U \subseteq \R^n, \, U \text{ is open}, \, f|_U \in C^1( U , \R) } U $
that
$ ( \bD f )( x ) = \cu{ ( \nabla f ) ( x ) } $,
and
\item 
\label{item:subdiff_item_v}
it holds for all $ x \in \R^n $ that 
$
  ( \bD f )( x ) 
$
is closed. 
\end{enumerate}
\cfload. 
\end{lemma}
\begin{cproof}{lem:subdifferential:c1}
Throughout this proof let 
$
  Z^{ x, y } = ( Z^{ x, y }_1, Z^{ x, y }_2 ) 
  \colon  
  \N \to \R^n \times \R^n
$, 
$ x, y \in \R^n $, 
satisfy for all $ x, y \in \R^n $, $ k \in \N $ 
that 
\begin{equation}
\label{eq:def_Z_process}
  Z^{ x, y }_1( k ) = x
\qqandqq
  Z^{ x, y }_2( k ) = y
  .
\end{equation}
\Nobs that \cref{def:limit:subdiff:eq} establishes \cref{item:subdiff_item_i}. 
\Nobs that \cref{eq:def_Z_process} shows
for all $ x \in \R^n $, $ y \in ( \cD f )( x ) $ 
that
\begin{equation}
% \textstyle
  \biggl[
    \forall \, k \in \N \colon 
    \Bigl(
      Z^{ x, y }_2( k ) 
      \in ( \cD f )( Z^{ x, y }_1( k ) )
    \Bigr)
  \biggr]
  \wedge
  \biggl[
    \limsup_{ k \to \infty }
    \bigl(
      \| Z^{ x, y }_1( k ) - x \|
      +
      \| Z^{ x, y }_2( k ) - y \|
    \bigr)
    = 0
  \biggr]
  .
\end{equation}
This establishes \cref{item:subdiff_item_ii}.
Next \nobs that Rockafellar \& Wets~\cite[Exercise 8.8]{RockafellarWets1998}
establishes \cref{item:subdiff_item_iii,item:subdiff_item_iv}. 
Finally, \cite[Theorem 8.6]{RockafellarWets1998}
 establishes \cref{item:subdiff_item_v}. 
\end{cproof}

\cfclear
\begin{lemma}[Limits of limiting Fr\'{e}chet subgradients]
\label{lem:limiting_derivatives}
\cfadd{def:limit:subdiff}
Let $ n \in \N $, $ f \in C( \R^n , \R ) $, 
let 
$ ( x_k )_{ k \in \N_0 } \subseteq \R^n $
and 
$ ( y_k )_{ k \in \N_0 } \subseteq \R^n $
satisfy 
$
  \limsup_{ k \to \infty }
  ( 
    \| x_k - x_0 \| 
    +
    \| y_k - y_0 \| 
  )
  = 0
$, 
and assume for all $ k \in \N $ that 
$
  y_k \in ( \bbD f)( x_k )
$
\cfload. 
Then 
$
  y_0 \in ( \bbD f)( x_0 )
$.
\end{lemma}

\begin{cproof}{lem:limiting_derivatives}
\Nobs that, e.g., \cite[Proposition 8.7]{RockafellarWets1998} implies that
$y_0 \in ( \bbD f)( x_0 )$.
\end{cproof}

\subsection{Local underestimates for the realization functions of DNNs}
\label{subsection:local:under}

Next we establish in \cref{lem:gradient:left:approx} a technical lemma that will be used in the proof of \cref{lem:gradient:approx:sequence} below. Roughly speaking, since we work with the left derivative of the ReLU function we need to approximate the realization functions from below to obtain convergence of the generalized gradients.

\begin{lemma} 
\label{lem:gradient:left:approx}
Assume \cref{setting:dnn} and let 
$ \theta \in \R^{ \fd } $,
$ \varepsilon \in (0, \infty) $. 
Then there exists a non-empty and open 
$
  U \subseteq \R^{ \fd } 
$ 
such that for all $ \vartheta \in U $, 
$ k \in \cu{ 1, 2, \ldots, L } $, 
$ i \in \cu{ 1, 2, \ldots, \ell_k } $, 
$ x \in [a,b]^{ \ell_0 } $
it holds that 
\begin{equation}
  \norm{ \vartheta - \theta } < \varepsilon 
\qqandqq
  \cN_{ \infty, i }^{ k, \vartheta }( x ) \le \cN_{ \infty, i }^{ k, \theta }( x ) 
  .
\end{equation}
\end{lemma}
\begin{cproof}{lem:gradient:left:approx}
Throughout this proof let 
$ \fC_k \in (0, \infty) $, 
$ k \in \N $, 
satisfy 
for all $ k \in \N $ that
$
  \fC_1 = \max \cu{\ell_0 \abs{a}, \ell_0 \abs{b},  1 }
$ 
and
\begin{equation}
\label{eq:def_Ck_constants}
%    \forall \, k \in \cu{1,  \ldots, L - 1} \colon 
   \fC_{ k + 1 } = 
   2 
   \fC_k 
   ( k + 1 ) 
   \max\cu{ 1, \abs{a}, \abs{b} } 
   ( \max\cu{ 1, \norm{ \theta } + 2 \fC_k } )^k 
   \bigl[ 
     \textstyle\prod_{ j = 0 }^k ( \ell_j + 1 ) 
   \bigr]
   ,
\end{equation}
let $ \delta \in (0, \infty) $ 
satisfy 
$
  \delta = \min\cu{ 1, \varepsilon ( 2 \fC_L \fd )^{ - 1 } } 
$, 
and let 
$ U \subseteq \R^{ \fd } $ 
satisfy
\begin{multline}
\label{eq:def_set_U}
  U = 
  \biggl\{
    \vartheta \in \R^{ \fd } \colon 
    \biggl(
      \Bigl[ 
        \forall \, k \in \cu{ 1, \ldots, L }, 
        i \in \cu{ 1, \ldots, \ell_k }, 
        j \in \cu{ 1, \ldots, \ell_{ k - 1 } } \colon 
        \abs{ \w{ k, \vartheta }_{ i, j } - \w{ k, \theta }_{ i, j } } < \delta  
      \Bigr]
%       \Bigr. 
\\
        \qquad 
        \wedge
%         \Bigl. 
        \br*{
          \forall \, k \in \cu{ 1, \ldots, L }, i \in \cu{ 1, \ldots, \ell_k } 
          \colon 
            \b{ k, \theta }_i - 2 \fC_k \delta 
            <
            \b{ k, \vartheta }_i 
            <
            \b{ \theta, k }_i - \fC_k \delta 
        } 
    \biggr)
  \biggr\} 
  .
\end{multline}
\Nobs that \cref{eq:def_set_U} ensures 
that $ U \subseteq \R^{ \fd } $ is non-empty and open. 
\Moreover \cref{eq:def_Ck_constants} shows
for all 
$ k \in \N $ 
that 
$
  \fC_{ k + 1 } > 2 \fC_k 
$ 
and 
$ \fC_k \geq 1 $. 
Combining this with \cref{eq:def_set_U} 
assures 
for all 
$ 
  \vartheta \in U 
$, 
$ i \in \cu{ 1, 2, \ldots, \fd } $ 
that 
\begin{equation}
  \abs{ \vartheta_i - \theta_i } 
  < \max\{ \delta, 2 \fC_1 \delta, 2 \fC_2 \delta, \dots, 2 \fC_L \delta \} 
  = 2 \fC_L \delta 
  \leq 
  2 \fC_L 
  \bigl( \varepsilon ( 2 \fC_L \fd )^{ - 1 } \bigr)
  = \fd^{ - 1 } \varepsilon 
  .
\end{equation}
\Hence for all $ \vartheta \in U $ that
\begin{equation}
\textstyle
  \norm{\vartheta - \theta } 
  =
  \bigl[
    \sum_{ i = 1 }^{ \fd }
    \abs{ \vartheta_i - \theta_i }^2
  \bigr]^{ 1 / 2 }
  \le 
  \fd 
  \bigl[
    \max\nolimits_{ i \in \cu{ 1, 2, \ldots, L } } 
    \abs{ \vartheta_i - \theta_i } 
  \bigr]
  < \varepsilon .
\end{equation}
\Moreover 
\cref{eq:def_Ak_transformation_deep_ANNs,eq:def_NN_realization,eq:def_set_U} 
ensure
for all 
$ \vartheta \in U $, 
$ i \in \cu{ 1, 2, \dots, \ell_1 } $, 
$ x = ( x_1, \dots, x_{ \ell_0 } ) \in [a,b]^{ \ell_0 } $
that
\begin{equation}
\begin{split}
  \cN_{ \infty, i }^{ 1, \vartheta }( x ) - 
  \cN_{ \infty, i }^{ 1, \theta }( x ) 
&
  =
  ( 
    \b{ 1, \vartheta }_i - \b{ 1, \theta }_i 
  ) 
  + 
  \smallsum_{ j = 1 }^{ \ell_0 } 
  ( 
    \w{ 1, \vartheta }_{ i, j } - \w{ 1, \theta }_{ i, j } 
  ) x_j 
\\
& 
  < 
  - \fC_1 \delta 
  + 
  \smallsum_{ j = 1 }^{ \ell_0 } 
  |
    ( 
      \w{ 1, \vartheta }_{ i, j } - \w{ 1, \theta }_{ i, j } 
    ) 
    x_j 
  |
\\ &
  \leq
    - \fC_1 \delta 
    + \delta \bigl( \smallsum_{ j = 1 }^{ \ell_0 } \abs{ x_j } \bigr)
  \leq 
    - \fC_1 \delta + \ell_0 \delta \max\cu{ \abs{a}, \abs{b} } 
  \le 0 .
\end{split}
\end{equation}
It thus remains to prove an analogous inequality for the subsequent layers. 
For this let
$ 
  \vartheta \in U 
$,
$
  k \in \N \cap (0,L) 
%   \cu{ 1, 2, \ldots, L - 1 } 
$, 
$
  i \in \cu{ 1, 2, \ldots, \ell_{ k + 1 } }
$, 
$ x \in [a,b]^{ \ell_0 } $, 
let 
$ 
  \bfd \in \N 
$
satisfy 
$
  \bfd = \ell_{ k + 1 } \ell_k + 1 + \sum_{ j = 1 }^{ k } \ell_j ( \ell_{ j - 1 } + 1 ) 
$, 
let 
$
  \bfD \in \N
$
satisfy 
$
  \bfD = \sum_{ j = 1 }^{ k + 1 } \ell_j ( \ell_{ j - 1 } + 1 ) 
$, 
and let 
$
  \psi \in \R^{ \fd } 
$
satisfy 
\begin{equation}
\label{eq:def_psi_zwischenvector}
  \psi = 
  ( 
    \vartheta_1, \vartheta_2, \ldots, \vartheta_{ \bfd - 1 }, \theta_{ \bfd }, 
    \theta_{ \bfd + 1 }, \dots, \theta_{ \fd }
  )
  .
\end{equation}
\Nobs that 
\cref{eq:def_Ak_transformation_deep_ANNs,eq:def_NN_realization,eq:def_set_U,eq:def_psi_zwischenvector} 
show that
\begin{equation}
\label{eq:vartheta_estimate_from_above_general_k}
\begin{split}
  \cN_{ \infty, i }^{ k + 1, \vartheta }(x) 
& = 
  \cN_{ \infty, i }^{ k + 1, \psi }( x ) 
  + 
  ( 
    \b{ k + 1, \vartheta }_i - \b{ k + 1, \psi }_i 
  ) 
\\ & =
  \cN_{ \infty, i }^{ k + 1, \psi }( x ) 
  + 
  ( 
    \b{ k + 1, \vartheta }_i - \b{ k + 1, \theta }_i 
  ) 
< 
  \cN_{ \infty, i }^{ k + 1, \psi }( x ) - \fC_{ k + 1 } \delta 
  .
\end{split}
\end{equation}
Next note that, e.g., \cite[Theorem 2.1]{BeckJentzenKuckuck2022}
(applied with 
$
  a \with a
$, 
$ b \with b $, 
$
  d \with \bfD
$,
$
  L \with k + 1 
$, 
$
  \ell \with ( \ell_0, \ell_1, \ldots, \ell_{ k + 1 } ) 
$
in the notation of \cite[Theorem 2.36]{BeckJentzenKuckuck2022})
demonstrates that
\begin{equation}
\label{eq:DNN_lipschitz_estimate}
\begin{split}
&
  \abs{ 
    \cN_{ \infty, i }^{ k + 1, \theta }( x ) 
    - 
    \cN_{ \infty, i }^{ k + 1, \psi }( x ) 
  } 
\\
&
  \leq 
  ( k + 1 ) \max\cu{ 1, \abs{a}, \abs{b} } 
  \br*{ 
    \max\cu*{ 
      1, 
      \max\nolimits_{ 
        i \in \cu{ 1, 2, \ldots, \bfD } 
      } 
      \abs{ \theta_i } ,
      \max\nolimits_{ 
        i \in \cu{ 1, 2, \ldots, \bfD } 
      } 
      \abs{ \psi_i } 
    } 
  }^k 
\\
& 
\quad 
  \cdot 
  \bigl[ 
    \textstyle\prod_{ m = 0 }^k ( \ell_m + 1 ) 
  \bigr] 
  \br*{ 
    \max\nolimits_{ i \in \cu{ 1, 2, \ldots, \bfD } } 
    \abs{ \theta_i - \psi_i } 
  }
  .
\end{split}
\end{equation}
\Moreover \cref{eq:def_set_U} ensures that 
\begin{equation}
  \max\nolimits_{ i \in \cu{ 1, 2, \ldots, \bfD } } 
  \abs{ \theta_i - \psi_i } 
  =
  \max\nolimits_{ i \in \cu{ 1, 2, \ldots, \bfd - 1 } } 
  \abs{ \theta_i - \vartheta_i } 
  \leq
  \max\{ \delta, 2 \fC_1 \delta, 2 \fC_2 \delta, \dots, 2 \fC_k \delta \}
  = 
  2 \fC_k \delta 
  .
\end{equation}
Combining this with \cref{eq:def_Ck_constants,eq:DNN_lipschitz_estimate} proves that 
\begin{equation}
\begin{split}
&
  \abs{ 
    \cN_{ \infty, i }^{ k + 1, \theta }( x ) 
    - 
    \cN_{ \infty, i }^{ k + 1, \psi }( x ) 
  } 
\\
& 
\leq 
  ( k + 1 ) 
  \max\cu{ 1, \abs{a}, \abs{b} } 
  \bigl[ 
    \max\cu*{ 
      1, 
      2 \fC_k \delta 
      + 
      \max\nolimits_{ i \in \cu{ 1, 2, \ldots, \bfD } } \abs{ \theta_i }, 
      \max\nolimits_{ i \in \cu{ 1, 2, \ldots, \bfD } } \abs{ \theta_i } 
    } 
  \bigr]^k 
\\
& \quad 
  \cdot 
  \bigl[ 
    \textstyle\prod_{ m = 0 }^k ( \ell_m + 1 ) 
  \bigr] 
  \bigl[ 
    2 \fC_k \delta 
  \bigr]
\\
&
\leq 
  2 \fC_k \delta ( k + 1 ) 
  \max\cu{ 1, \abs{a}, \abs{b} } 
  \bigl[ 
    \textstyle\prod_{ m = 0 }^k ( \ell_m + 1 ) 
  \bigr] 
  \bigl[ 
    \max\cu{ 1, \norm{ \theta } + 2 \fC_k } 
  \bigr]^k 
= 
  \fC_{ k + 1 } \delta .
\end{split}
\end{equation}
This and \cref{eq:vartheta_estimate_from_above_general_k} 
assure that
\begin{equation}
\begin{split}
  \cN_{ \infty, i }^{ k + 1, \vartheta }(x) 
& 
<  
  \cN_{ \infty, i }^{ k + 1, \psi }( x ) - \fC_{ k + 1 } \delta 
=
  \cN_{ \infty, i }^{ k + 1, \theta }( x ) 
  + 
  \bigl( 
    \cN_{ \infty, i }^{ k + 1, \psi }(x) - \cN_{ \infty, i }^{ k + 1, \theta }( x )  
  \bigr)
  - \fC_{k+1} \delta 
\\
&  
\leq
  \cN_{ \infty, i }^{ k + 1, \theta }( x ) 
  + \abs{ \cN_{ \infty, i }^{ k + 1, \psi }(x) - \cN_{ \infty, i }^{ k + 1, \theta }( x ) } 
  - \fC_{k+1} \delta 
\leq 
  \cN_{ \infty, i }^{ k + 1, \theta }( x ) .
\end{split}
\end{equation}
\end{cproof}

\subsection{Continuity properties for the generalized gradient function}
\label{subsection:gen:grad:continuous}

\begin{lemma}[Continuity points of the generalized gradient function]
\label{lem:gradient:convergence}
Assume \cref{setting:dnn}
and 
let 
$ 
  \theta = ( \theta_n )_{ n \in \N_0 } \colon \N_0 \to \R^{ \fd } 
$ 
satisfy 
for all 
$ k \in \cu{ 1, 2, \ldots, L } $, 
$ i \in \cu{ 1, 2, \ldots, \ell_k } $, 
$ x \in [a,b]^{ \ell_0 } $
that 
\begin{equation}
\label{eq:limit_theta_assumption}
  \limsup\nolimits_{ n \to \infty } 
  \bigl(
    \| \theta_n - \theta_0 \| 
    +
    |
      \indicator{ (0, \infty) }( 
        \cN_{ \infty, i }^{ k, \theta_n }( x ) 
      ) 
      - 
      \indicator{ (0, \infty) 
      }( 
        \cN_{ \infty, i }^{ k, \theta_0 }( x ) 
      ) 
    | 
  \bigr)
  = 0
  .
\end{equation}
Then 
$
  \limsup_{ n \to \infty } 
  \| \cG( \theta_n ) - \cG ( \theta_0 ) \| = 0 
$.
\end{lemma}
\begin{cproof}{lem:gradient:convergence}
\Nobs that \cref{lem:realization:lip} (applied 
for every $ k \in \{ 1, 2, \dots, L \} $ 
with 
$
  L \with k
$
in the notation of \cref{lem:realization:lip}), 
\cref{wb,eq:def_Ak_transformation_deep_ANNs,eq:def_NN_realization}
assure that 
for all $ k \in \{ 1, 2, \dots, L \} $, 
$ j \in \{ 1, 2, \dots, \ell_k \} $
it holds that 
\begin{equation}
\label{eq:continuous_realization_function}
\textstyle
  \limsup_{ n \to \infty }
  \sup_{ x \in [a,b]^{ \ell_0 } }
  |
    \cN_{ \infty, j }^{ k, \theta_n }( x ) 
    -
    \cN_{ \infty, j }^{ k, \theta_0 }( x ) 
  |
  = 0 .
\end{equation}
\Moreover \cref{eq:def_NN_realization,eq:limit_theta_assumption}
ensure for all $ x \in [a,b]^{ \ell_0 } $, 
$ k \in \cu{ 1, 2, \ldots, L } $, 
$ i \in \cu{ 1, 2, \ldots, \ell_k } $ 
 that 
$
  \limsup_{ n \to \infty } 
  |
    \indicator{ \cX_i^{ k, \theta_n } }( x ) 
    -
    \indicator{ \cX_i^{ k, \theta_0 } }( x )
  |
  = 0
$. 
Combining this, 
\cref{eq:limit_theta_assumption},
and 
\cref{eq:continuous_realization_function} 
with \cref{prop:G} and 
Lebesgue's dominated convergence theorem establishes that
$
  \limsup_{ n \to \infty } 
  \| \cG( \theta_n ) - \cG( \theta_0 ) \| = 0
$.
\end{cproof}

As a consequence of \cref{lem:gradient:left:approx,lem:gradient:convergence} we show in \cref{lem:gradient:approx:sequence} that, loosely speaking,
 the generalized gradient $\cG( \theta)$ at an arbitrary point $\theta \in \R^\fd$ can be represented as the limit of generalized gradients of a sequence $\vartheta_n \to \theta$, even after removing an arbitrary set of zero measure.

\begin{lemma} 
\label{lem:gradient:approx:sequence}
Assume \cref{setting:dnn}
and let $ \theta \in \R^{ \fd } $, $ E \in \cB( \R^{ \fd } ) $ satisfy 
$
  \int_{ \R^{ \fd } \backslash E } 1 \, \d \vartheta = 0
$. 
Then there exists 
$
  \vartheta = ( \vartheta_n )_{ n \in \N } \colon \N \to E
$ 
such that
\begin{equation}
\label{eq:vartheta_to_prove_construction_sequence}
  \limsup\nolimits_{
    n \to \infty
  } 
  \bigl( 
    \norm{ \vartheta_n - \theta } 
    + 
    \norm{ \cG( \vartheta_n ) - \cG( \theta ) } 
  \bigr)
  = 0 .
\end{equation}
\end{lemma}
\begin{cproof}{lem:gradient:approx:sequence}
\Nobs that \cref{lem:gradient:left:approx} assures 
that there exist 
non-empty and open $ U_n \subseteq \R^{ \fd } $, $ n \in \N $, 
which satisfy for all 
$ n \in \N $, 
$ \vartheta \in U_n $,
$ k \in \cu{ 1, 2, \ldots, L} $, 
$ i \in \cu{ 1, 2, \ldots, \ell_k } $, 
$ x \in [a,b]^{ \ell_0 } $ 
that 
\begin{equation}
\label{eq:property_sets_Un}
  \norm{ \vartheta - \theta } < \tfrac{ 1 }{ n } 
\qqandqq
  \cN_{ \infty, i }^{ k, \vartheta }(x) 
  \le 
  \cN_{ \infty, i }^{ k, \theta }( x )
\end{equation}
\Nobs that the assumption that 
$
  \int_{ \R^{ \fd } \backslash E } 1 \, \d \vartheta = 0
$ 
implies for all $ n \in \N $ that 
$
  ( U_n \cap E ) \not= \emptyset
$. 
In the following let 
$
  \vartheta = ( \vartheta_n )_{ n \in \N } \colon \N \to E 
$ 
satisfy for all $ n \in \N $ that 
\begin{equation}
\label{eq:construction_of_sequence_vartheta}
  \vartheta_n \in U_n
  .
\end{equation}
\Nobs that \cref{eq:property_sets_Un}
assures
for all $ n \in \N $ 
that 
$
  \norm{ \vartheta_n - \theta } < \frac{ 1 }{ n } 
$. 
\Hence that 
\begin{equation}
\label{eq:vartheta_n_is_convergent}
\textstyle
  \limsup_{ n \to \infty } \| \vartheta_n - \theta \| = 0 
  .
\end{equation}
\cref{lem:realization:lip} (applied for every $ k \in \{ 1, 2, \dots, L \} $ 
with 
$
  L \with k
$
in the notation of \cref{lem:realization:lip}) 
\hence implies that for all 
$
  k \in \cu{ 1, 2, \ldots, L }
$, 
$
  i \in \cu{ 1, 2, \ldots, \ell_k } 
$, 
$
  x \in [a,b]^{ \ell_0 }
$
we have that 
\begin{equation}
\label{eq:realization_functions_sequence_is_convergent}
\textstyle
  \limsup_{ n \to \infty } 
  |
    \cN_{ \infty, i }^{ k, \vartheta_n }( x ) 
    -
    \cN_{ \infty, i }^{ k, \theta }( x ) 
  |
  = 0
  .
\end{equation}
\Moreover 
\cref{eq:property_sets_Un,eq:construction_of_sequence_vartheta} assure
for all 
$
  n \in \N
$, 
$
  k \in \cu{ 1, 2, \ldots, L }
$, 
$
  i \in \cu{ 1, 2, \ldots, \ell_k }
$,  
$
  x \in [a,b]^{ \ell_0 }
$
that 
$
  \cN_{ \infty, i }^{ k, \vartheta_n }( x ) 
  \le 
  \cN_{ \infty, i }^{ k, \theta }( x ) 
$. 
Combining this and 
\cref{eq:realization_functions_sequence_is_convergent} 
with the fact that the function 
$
  \R \ni x \mapsto \indicator{ (0,\infty) }( x ) \in \R
$
is left continuous 
demonstrates 
for all 
$ 
  x \in [a,b]^{ \ell_0 } 
$, 
$
  k \in \cu{ 1, 2, \ldots, L } 
$, 
$
  i \in \cu{ 1, 2, \ldots, \ell_k } 
$ 
that 
$
  \limsup_{ n \to \infty } 
  |
    \indicator{ (0, \infty ) }( 
      \cN_{ \infty, i }^{ k, \vartheta_n }( x ) 
    ) 
    -
    \indicator{ (0, \infty) }( 
      \cN_{ \infty, i }^{ k, \theta }( x ) 
    )
  | = 0
$.
\cref{lem:gradient:convergence} and 
\cref{eq:vartheta_n_is_convergent} \hence show 
that 
$
  \limsup_{ n \to \infty } 
  \| \cG( \vartheta_n ) - \cG( \theta ) \| = 0
$. 
Combining this with \cref{eq:vartheta_n_is_convergent}
establishes \cref{eq:vartheta_to_prove_construction_sequence}. 
\end{cproof}

\subsection{Generalized gradients as limiting Fr\'{e}chet subdifferentials}
\label{subsection:gen:grad:frechet}

We next employ the differentiability result from \cref{prop:risk:ae:diff},
the approximation result for the generalized gradient from \cref{lem:gradient:approx:sequence},
and the definition of the limiting Fr\'{e}chet subdifferential to establish in \cref{prop:loss:gradient:subdiff} the main result of this section: For every $\theta \in \R^\fd$, the generalized gradient $\cG ( \theta )$ is an element of the limiting Fr\'{e}chet subdifferential
$( \bD \cL_{ \infty } )( \theta )$.

\cfclear
\begin{prop}
\label{prop:loss:gradient:subdiff}
\cfadd{def:limit:subdiff}
Assume \cref{setting:dnn} and 
let $ \theta \in \R^{ \fd } $. 
Then 
$
  \cG( \theta ) \in ( \bD \cL_{ \infty } )( \theta )
$
\cfload.
\end{prop}
\begin{cproof}{prop:loss:gradient:subdiff}
\Nobs that \cref{prop:risk:ae:diff} ensures that 
there exists 
$ E \in \cB( \R^{ \fd } ) $ 
which satisfies  
$
  \int_{ \R^{ \fd } \backslash E } 1 \, \d \vartheta = 0 
$, 
which satisfies for all $ \vartheta \in E $ 
that $\cL_\infty$ is differentiable at $\vartheta$, 
and which satisfies for all 
$\vartheta \in E$ that 
\begin{equation}
\label{eq:differentiable_property}
  ( \nabla \cL_{ \infty } )( \vartheta ) 
  = \cG( \vartheta )
  .
\end{equation}
\Nobs that \cref{eq:differentiable_property} 
and \cref{lem:subdifferential:c1} 
prove for all $ \vartheta \in E $ 
that 
\begin{equation}
\label{eq:generalized_gradient_subdifferential}
  \cG( \vartheta ) \in ( \cD \cL_{ \infty } )( \vartheta )
  .
\end{equation}
\Moreover the fact that 
$ \int_{ \R^{ \fd } \backslash E } 1 \, \d \vartheta = 0 $ 
and \cref{lem:gradient:approx:sequence} imply 
that there exists 
$
  \vartheta = ( \vartheta_n )_{ n \in \N } \colon \N \to E 
$ 
which satisfies 
\begin{equation}
\label{eq:proof_construction_of_vartheta_sequence}
\textstyle
  \limsup_{ n \to \infty } 
  \bigl(
    \| \vartheta_n - \theta \| 
    + 
    \| \cG( \vartheta_n ) - \cG( \theta ) \|
  \bigr)
  = 0
  .
\end{equation}
\Nobs that \cref{eq:proof_construction_of_vartheta_sequence,eq:generalized_gradient_subdifferential} 
demonstrate that 
$
  \cG( \theta ) 
  \in 
  ( \bD \cL_{ \infty } )( \theta )
$.
\end{cproof}

Finally, as a consequence of \cref{prop:loss:gradient:subdiff} we show in \cref{cor:cG_equal_to_gradient} that on every open set on which the risk function $\cL_\infty$ is continuously differentiable its gradient agrees with $\cG$.
This fact will be used in the convergence analysis of GD processes in \cref{sec:convergence_GD}.

\cfclear
\begin{cor}
\label{cor:cG_equal_to_gradient}
\cfadd{def:limit:subdiff}
Assume \cref{setting:dnn}. 
Then it holds for all 
$
  \theta \in 
  \cup_{ 
    U \subseteq \R^{ \fd }, 
    \, U \text{ is open}, \, 
    ( \cL_{ \infty } )|_U \in C^1( U, \R )
  }
  U
$
that 
$
  \cG( \theta ) = ( \nabla \cL_{ \infty } )( \theta )
$
\cfload.
\end{cor}
\begin{cproof}{cor:cG_equal_to_gradient}
\Nobs that 
\cref{item:subdiff_item_iv}
in 
\cref{lem:subdifferential:c1} 
(applied with 
$
  n \with \fd 
$,
$  
  f \with \cL_{ \infty }
$
in the notation of \cref{lem:subdifferential:c1}) 
and \cref{prop:loss:gradient:subdiff}
ensure that for all open $ U \subseteq \R^{ \fd } $ 
and all $ \theta \in U $
with 
$ ( \cL_{ \infty } )|_U \in C^1( U, \R) $
it holds that 
\begin{equation}
  \cG( \theta ) 
  \in 
  ( \bbD \cL_{ \infty } )( \theta )
  =
  \{
    ( \nabla \cL_{ \infty } )( \theta )
  \}  
  .
\end{equation}
\Hence for all open $ U \subseteq \R^{ \fd } $ 
and all $ \theta \in U $
with 
$ ( \cL_{ \infty } )|_U \in C^1( U, \R) $
that 
$
  \cG( \theta ) = 
  ( \nabla \cL_{ \infty } )( \theta )
$. 
\end{cproof}

\section{Suitable piecewise rational functions}
\label{section:piecewise:rational}

In this section we identify
in \cref{eq:def_function_amn}
in \cref{def:function:amn}
a suitable subclass of the class of semi-algebraic functions 
which is closed under integration 
(see \cref{prop:integrals:amn} 
in \cref{ssec:closedness_integration} below) 
and which contains the realization functions of deep ReLU ANNs 
(see \cref{prop:realization:bmn} 
in \cref{ssec:realization_DNNs} below). 
The fact that functions in this class of suitable piecewise rational functions are semi-algebraic is established in \cref{prop:amn:semialgebraic} below.
We also summarize in \cref{subsection:semialgebraic:def} some basic facts regarding semi-algebraic sets and functions.
The results from this section will be employed in \cref{sec:piecewise_polynomial} below to establish that the considered risk function in the training of deep ANNs with ReLU activation are semi-algebraic.

Closedness under integration is not a trivial issue due to the fact that, 
in general, 
the integral of a semi-algebraic function 
is not necessarily semi-algebraic 
(in fact, in general not even globally subanalytic, see Kaiser~\cite{Kaiser2013}). 
Our analysis of the integrals of the functions considered in \cref{def:function:amn} below crucially relies on the fact that they are piecewise rational on regions separated by hyperplanes in the $x$-component. This property is also satisfied by the realization functions of ANNs with ReLU activation.

The function class in 
\cref{def:function:amn}
and some of the results in this section are 
inspired by the findings in our previous article 
Eberle et al.~\cite[Section~4]{EberleJentzenRiekertWeiss2021}. 
In particular, 
\cref{def:function:amn}
extends 
\cite[Definition~4.6]{EberleJentzenRiekertWeiss2021}, 
\cref{prop:amn:semialgebraic} in \cref{ssec:amn:semialgebraic} below 
extends 
\cite[Lemma~4.7]{EberleJentzenRiekertWeiss2021}, 
and 
\cref{prop:integrals:amn}
in \cref{ssec:closedness_integration} below 
extends 
\cite[Proposition~4.8]{EberleJentzenRiekertWeiss2021}.

\subsection{Suitable piecewise rational functions}

\cfclear
\begin{definition}[Vector spaces of suitable piecewise rational functions]
\label{def:function:amn}
Let $ m, n \in \N_0 $, $ \delta \in (0,\infty] $. 
Then we denote by $ \scrF_{ m, n, \delta } $ 
the $ \R $-vector space given by
\begin{multline}
\label{eq:def_function_amn}
\cfadd{def:polynomial}\cfadd{def:rational:function}\cfadd{def:degree}
  \scrF_{ m, n, \delta } 
  = 
  \operatorname{span}_{ \R 
  }\!\Biggl( 
    \Biggl\{ 
      F \colon \R^m \times \R^n \to \R \colon 
      \Biggl[  
        \exists \, r \in \N, \,
        R \in \ratio_{ m, \delta } , \,
        Q \in \{ q \in \polyn_n \colon \deg(q) \leq \delta \} ,
\\
        P = ( P_{ i, j } )_{ (i,j) \in \cu{1, 2, \ldots, r } \times \cu{0, 1, \ldots, n } } 
        \subseteq \polyn_m 
        \colon 
        \biggl(
          \forall \, \theta \in \R^m, \, 
          x = (x_1, \ldots, x_n) \in \R^n \colon 
        \\
          f( \theta , x ) 
          = R( \theta ) Q( x ) 
          \biggl[
            \textstyle\prod\limits_{ i = 1 }^r 
            \indicator{ [0,\infty) }(
              P_{ i, 0 }( \theta ) 
              + 
              \smallsum_{ j = 1 }^n 
              P_{ i, j }( \theta ) x_j 
            )
          \biggr] 
        \biggr)
        \Biggr] 
    \Biggr\} 
  \Biggr)
\end{multline}
\cfload.
\end{definition}

In \cref{eq:def_function_amn} above we denote by $\operatorname{span}
_\R$ the linear span of the given functions $F \colon \R^m \times \R^n \to \R$
with coefficients in $\R$, i.e., the $\R$-vectorspace generated by these functions.

\Nobs that functions in $\scrF_{m,n,\delta}$ depend on two vectors $\theta \in \R^m$ and $x \in \R^n$. In the considered deep learning framework this will be applied in the situation where $\theta$ is the parameter vector of a suitable ANN, $x$ is the input vector of the ANN, and $f(\theta,x)$ is the output.

\subsection{Elementary properties of suitable piecewise rational functions}

\cfclear
\begin{lemma} 
\label{lem:bmn:amn}
Let $ m, n \in \N_0 $. 
Then 
\begin{enumerate}[label = (\roman*)]
\item 
\label{item:properties_F_i}
it holds for all $ \delta_1, \delta_2 \in (0,\infty] $ 
with $ \delta_1 \leq \delta_2 $ 
that 
$
  \scrR_{ n, \delta_1 } 
  \subseteq 
  \scrR_{ n, \delta_2 } 
$, 
\item 
\label{item:properties_F_ib}
it holds for all $ \delta \in (0,1] $ 
that 
$
  \scrR_{ n, \delta } 
  =
  \scrP_n
$, 
\item 
\label{item:properties_F_ii}
it holds for all $ \delta_1, \delta_2 \in (0,\infty] $ 
with $ \delta_1 \leq \delta_2 $ 
that 
$
  \scrF_{ m, n, \delta_1 } 
  \subseteq 
  \scrF_{ m, n, \delta_2 } 
$, 
\item 
\label{item:properties_F_iii}
it holds that 
$
  \scrF_{ m, n, 1 } \subseteq 
  \scrF_{ m, n, \infty }
$, 
\item 
\label{item:properties_F_iv}
it holds for all 
$ 
  f, g \in 
  \scrF_{ m, n, \infty }
$ 
that
\begin{equation}
\textstyle
  \bigl( 
    \R^m \times \R^n 
    \ni (\theta, x) \mapsto
    f(\theta,x) g(\theta,x) \in \R 
  \bigr) 
  \in
  \scrF_{ m, n, \infty }
  ,
\end{equation}
and 
\item 
\label{item:properties_F_v}
it holds for all $ \delta \in (0,\infty] $ 
that 
\begin{multline} 
  \scrF_{ m, n, \delta } 
  = 
  \operatorname{span}_{ \R 
  }\Bigl( 
    \Bigl\{ 
      F \colon \R^m \times \R^n \to \R \colon 
      \Bigl[  
        \exists \, r \in \N, \,
        A_1, A_2, \ldots, A_r \in \cu{ \cu{0}, [0, \infty ), (0, \infty )}, 
\\
        R \in \ratio_{ m, \delta } , \,
        Q \in \{ q \in \polyn_n \colon \deg(q) \leq \delta \} , 
        \,
        P = ( P_{ i, j } )_{ (i,j) \in \cu{1, 2, \ldots, r } \times \cu{0, 1, \ldots, n } } 
        \subseteq \polyn_m 
        \colon 
        \Bigl(
          \forall \, \theta \in \R^m \colon
        \\
          \forall \, x = (x_1, \ldots, x_n) \in \R^n \colon 
          f( \theta , x ) 
          = R( \theta ) Q( x ) 
          \br[\big]{ 
            \textstyle\prod_{ i = 1 }^r 
            \indicator{ A_i }(
              P_{ i, 0 }( \theta ) 
              + 
              \smallsum_{ j = 1 }^n 
              P_{ i, j }( \theta ) x_j 
            )
          } 
        \Bigr)
        \Bigr] 
    \Bigr\} 
  \Bigr)
\end{multline}
\cfadd{def:function:amn}
\end{enumerate}
\cfload.
\end{lemma}
\begin{cproof}{lem:bmn:amn}
\Nobs that 
\cref{eq:def_set_rational_functions}
and the fact 
that for all $ \delta_1, \delta_2 \in (0,\infty] $ 
with $ \delta_1 \leq \delta_2 $
it holds that 
\begin{equation}
  \{ q \in \scrP_n \colon \deg(q) < \delta_1 \}
  \subseteq
  \{ q \in \scrP_n \colon \deg(q) < \delta_2 \}
\end{equation}
establish  
\cref{item:properties_F_i}. 
\Nobs that 
\cref{eq:def_set_rational_functions} 
and the fact 
that for all $ \delta \in (0,1] $ 
it holds that 
\begin{equation}
  \{ q \in \scrP_n \colon \deg(q) < \delta \}
  =
  \{ q \in \scrP_n \colon \deg(q) = 0 \}
\end{equation}
prove
\cref{item:properties_F_ib}. 
\Nobs that 
\cref{eq:def_function_amn}, 
\cref{item:properties_F_i}, 
and the fact 
that for all $ \delta_1, \delta_2 \in (0,\infty] $ 
with $ \delta_1 \leq \delta_2 $
it holds that 
\begin{equation}
  \{ q \in \scrP_n \colon \deg(q) \leq \delta_1 \}
  \subseteq
  \{ q \in \scrP_n \colon \deg(q) \leq \delta_2 \}
\end{equation}
establish  
\cref{item:properties_F_ii}. 
\Nobs that 
\cref{item:properties_F_ii} 
proves  
\cref{item:properties_F_iii}. 
\Nobs that \cref{eq:def_function_amn} 
establishes 
\cref{item:properties_F_iv}. 
\Nobs that 
the fact that 
$ 
  \forall \, y \in \R \colon 
  \mathbbm{1}_{ \{ 0 \} }( y )
  =
  \mathbbm{1}_{ (-\infty,0] \cap [0,\infty) }( y )
  =
  \mathbbm{1}_{ (-\infty,0] }( y )
  \mathbbm{1}_{ [0,\infty) }( y )
  =
  \mathbbm{1}_{ [0,\infty) }( y )
  \mathbbm{1}_{ [0,\infty) }( -y )
$ 
shows that 
for all 
$ 
  P_0, P_1, \dots, P_n \in \polyn_m 
$
it holds that 
\begin{equation}
\label{eq:indicator_nur_0}
\textstyle
\begin{split}
&
\textstyle
  \mathbbm{1}_{
    \{ 0 \}
  }(
    P_0( \theta )
    +
    \sum_{ j = 1 }^n
    P_j( \theta ) x_j
  )
\\ & =
\textstyle
  \mathbbm{1}_{
    [0,\infty)
  }(
    P_0( \theta )
    +
    \sum_{ j = 1 }^n
    P_j( \theta ) x_j
  )
  \mathbbm{1}_{
    [0,\infty)
  }(
    - P_0( \theta )
    +
    \sum_{ j = 1 }^n
    ( - P_j( \theta ) ) x_j
  )
  .
\end{split}
\end{equation}
\Moreover 
the fact that 
$ 
  \forall \, y \in \R \colon 
  \mathbbm{1}_{ (0,\infty) }( y )
  =
  1 - \mathbbm{1}_{ (-\infty,0] }( y )
  =
  1 - \mathbbm{1}_{ [0,\infty) }( -y )
$ 
shows that 
for all 
$ 
  P_0, P_1, \dots, P_n \in \polyn_m 
$
it holds that 
\begin{equation}
\textstyle
% \begin{split}
% &
% \textstyle
  \mathbbm{1}_{
    (0,\infty)
  }(
    P_0( \theta )
    +
    \sum_{ j = 1 }^n
    P_j( \theta ) x_j
  )
% \\ & 
=
% \textstyle
  1
  -
  \mathbbm{1}_{
    [0,\infty)
  }(
    - P_0( \theta )
    +
    \sum_{ j = 1 }^n
    ( - P_j( \theta ) ) x_j
  )
  .
% \end{split}
\end{equation}
Combining \cref{eq:def_function_amn,eq:indicator_nur_0} 
\hence shows 
that for all 
$ \delta \in (0,\infty] $, 
$ r \in \N $,
$ A_1, A_2, \dots, A_r \in \{ \{ 0 \} , (0,\infty), [0,\infty) \} $, 
$ R \in \scrR_{ m, \delta } $,
$ Q \in \{ q \in \polyn_n \colon \deg( q ) \leq \delta \} $, 
$
  P = ( P_{ i, j } )_{
    \{ 1, 2, \dots, r \} \times \{ 0, 1, \dots, n \} 
  }
  \subseteq \polyn_m
$
it holds that 
\begin{equation}
\textstyle
  \bigl(
    \R^m \times \R^n
    \ni 
    ( \theta, x )
    \mapsto 
    R( \theta ) Q( x )
    \bigl[
      \prod_{ i = 1 }^r
      \mathbbm{1}_{ A_i }(
        P_0( \theta )
        +
        \sum_{ j = 1 }^n
        P_j( \theta ) x_j
      )
    \bigr]
    \in \R
  \bigr)
  \in \scrF_{ m, n, \delta }
  .
\end{equation}
This establishes \cref{item:properties_F_v}. 
\end{cproof}

\subsection{Semi-algebraic sets}
\label{subsection:semialgebraic:def}

In the following we gather several known definitions and elementary results 
regarding semi-algebraic sets and functions; 
cf., e.g., Bochnak et al.~\cite{BochnakCoste1998}, 
Coste~\cite{Coste2000},
Shiota~\cite{Shiota1997},
 and Van den Dries \& Miller~\cite{VanDries1996geometric}.

\begin{definition}[Set of polynomials]
	\label{def:polynomial}	
	Let $ n \in \N_0 $. 
	Then we denote by $ \polyn_n \subseteq C( \R^n, \R ) $ 
	the set\footnote{Note that $ \R^0 = \{ 0 \} $, 
		$ C( \R^0, \R ) = C( \{ 0 \}, \R ) $, and 
		$ \#( C( \R^0, \R ) ) = \#( C( \{ 0 \}, \R ) ) = \infty $. In particular, this shows 
		for all $ n \in \N_0 $ that $ \operatorname{dim}( \R^n ) = n $ 
		and $ \#( C( \R^n, \R ) ) = \infty $.} 
	of all polynomials from $ \R^n $ to $ \R $.
\end{definition}

\cfclear
\begin{definition}[Multidimensional semi-algebraic sets]
	\label{def:semialgebraic:set}
	Let $ n \in \N $ and let $ A \subseteq \R^n $ be a set. 
	Then we say that $ A $ is an $ n $-dimensional semi-algebraic set 
	if and only if there exist 
	$ M , N \in \N $ 
	and 
	$
	( 
	P_{ i, j, k } 
	)_{ 
		( i, j, k ) \in \cu{ 1, 2, \ldots, M } \times \cu{1, 2, \ldots, N } \times \cu{ 0, 1 } 
	} 
	\subseteq \polyn_n 
	$ 
	such that
	\begin{equation} 
	\label{eq:def_semi-algebraic_set}
	\cfadd{def:polynomial}
	\textstyle
	A = 
	\bigcup_{ i = 1 }^M 
	\bigl( 
	\bigcap_{ j = 1 }^N 
	\cu*{ 
		x \in \R^n \colon P_{ i, j, 0 }( x ) = 0 < P_{ i, j, 1 }( x ) 
	}
	\bigr)
	\end{equation}
	\cfload.
\end{definition}

Note that in \cref{eq:def_semi-algebraic_set} we have that $\cu{ 
	x \in \R^n \colon P_{ i, j, 0 }( x ) = 0 < P_{ i, j, 1 }( x ) 
} = \cu{ x \in \R^n \colon [ P_{i,j,0} (x ) = 0 \wedge P_{i,j,1} ( x ) > 0]} = \cu{x \in \R^n \colon P_{i,j,0} ( x ) = 0 } \cap \cu{x \in \R^n \colon P_{i,j,1} ( x ) > 0 }$.

The following properties of semi-algebraic sets are well-known and not hard to show from the definition; see, e.g., Shiota~\cite[(I.2.9)]{Shiota1997}.

\cfclear
\begin{prop}
\label{prop:semialgebraic:sets}
Let $m ,n \in \N$. Then 
\begin{enumerate} [label = (\roman*)]
	\cfadd{def:polynomial} \cfadd{def:semialgebraic:set}
	\item it holds for all $n$-dimensional semi-algebraic sets $A,B$ that $A \cup B$, $A \cap B$, and $\R^n \backslash A$ are $n$-dimensional semi-algebraic sets,
	\item it holds for every $n$-dimensional semi-algebraic set $A$ and every $m$-dimensional semi-algebraic set $B$ that $A \times B $ is an $(m+n)$-dimensional semi-algebraic set,
	\item it holds for every $P \in \scrP_n$ that $\cu{x \in \R^n \colon P ( x ) \ge 0 }$ is an $n$-dimensional semi-algebraic set,
	\item it holds for all $a \in \R^n$ that $\cu{a} \subseteq \R^n$ is an $n$-dimensional semi-algebraic set
\end{enumerate}
\cfload.
\end{prop}

\subsection{Semi-algebraic functions}

\cfclear
\begin{definition}[Semi-algebraic functions]
	\label{def:semialgebraic:function}
	Let $ m, n \in \N $ and let 
	$ f \colon \R^m \to \R^n $ be a function. 
	Then we say that $ f $ is a semi-algebraic function 
	(we say that $ f $ is semi-algebraic) 
	if and only if it holds that 
	$
	\operatorname{Graph}( f ) 
	$
	is an $ ( m + n ) $-dimensional 
	semi-algebraic set
	\cfadd{def:semialgebraic:set}\cfload.
\end{definition}

The next elementary result, \cref{lem:sum:semialgebraic},
is a direct consequence of, e.g., \cite[(I.2.9)]{Shiota1997} or \cite[Proposition 2.2.6]{BochnakCoste1998} (see, e.g., also Bierstone \& Milman~\cite[Section 1]{BierstoneMilman1988}).

\cfclear
\begin{lemma} \label{lem:sum:semialgebraic}
	Let $ n \in \N $ and let $ f \colon \R^n \to \R $ 
	and 
	$
	g \colon \R^n \to \R 
	$
	be semi-algebraic \cfadd{def:semialgebraic:function}\cfload.
	Then
	\begin{enumerate}[label = (\roman*)]
		\item \label{lem:sum:semialgebraic:item1} 
		it holds that $ \R^n \ni x \mapsto f(x) + g(x) \in \R $ is semi-algebraic and
		\item \label{lem:sum:semialgebraic:item2} 
		it holds that $ \R^n \ni x \mapsto f(x) g(x) \in \R $ is semi-algebraic. 
	\end{enumerate}
\end{lemma}
%\begin{cproof}{lem:sum:semialgebraic}
%	\Nobs[note] that, e.g., Coste~\cite[Corollary 2.9]{Coste2000} 
%	or \cite[Proposition 2.2.6]{BochnakCoste1998}
%	(see, e.g., also Bierstone \& Milman~\cite[Section 1]{BierstoneMilman1988}) 
%	establishes \cref{lem:sum:semialgebraic:item1,lem:sum:semialgebraic:item2}.
%\end{cproof}

\cfclear
\begin{lemma}
	\label{lem:indicator_semialgebraic}
	\cfadd{def:semialgebraic:set}
	Let $ n \in \N $ and let $ A \subseteq \R^n $ be an 
	$ n $-dimensional semi-algebraic set 
	\cfload. 
	\cfadd{def:semialgebraic:function}
	Then 
	$
	\R^n \ni x \mapsto \indicator{ A }( x ) \in \R
	$
	is semi-algebraic
	\cfload. 
\end{lemma}

\begin{cproof}{lem:indicator_semialgebraic}
	Throughout this proof let $ f \colon \R^n \to \R $ 
	satisfy for all $ x \in \R^n $ that
	\begin{equation}
	\label{eq:indicator_lemma:def_f}
	f(x) = \indicator{ A }( x )
	.
	\end{equation}
	\Nobs that \cref{eq:indicator_lemma:def_f} 
	shows that
	\begin{equation}
	\label{eq:graph_f_rep_proof}
	\operatorname{Graph}( f )
	% \\ & 
	=
	(
	A \times \{ 1 \}
	) 
	\cup 
	(
	( \R^n \backslash A ) \times \{ 0 \}
	) \subseteq \R^{n + 1 }
	.
	\end{equation}
	\Moreover \cref{prop:semialgebraic:sets} ensures that 
	$ \{ 0 \} $
	and 
	$ \{ 1 \} $ 
	are $ 1 $-dimensional semi-algebraic sets and that 
	$ \R^n \backslash A $ 
	is an $ n $-dimensional semi-algebraic set. 
	Combining this
	with \cref{prop:semialgebraic:sets}
	shows that 
	$
	A \times \{ 1 \}
	$
	and 
	$
	( \R^n \backslash A ) \times \{ 0 \}
	$
	are $ ( n + 1 ) $-dimensional semi-algebraic sets. 
	\cref{prop:semialgebraic:sets,eq:graph_f_rep_proof} \hence 
	show that $ \operatorname{Graph}( f ) $ is 
	an $ ( n + 1 ) $-dimensional semi-algebraic set. 
	This establishes that $ f $ is semi-algebraic. 
\end{cproof}

\subsection{Rational functions as semi-algebraic functions}

The next goal is to establish in \cref{prop:amn:semialgebraic} below that the functions in the classes $\scrF_{m,0 , \infty}$, $m \in \N$, are semi-algebraic.
As a preparation, we first recall in \cref{lem:rational:semialgebraic} below the simple fact that rational functions are semi-algebraic.

\cfclear
\begin{definition}[Degree\footnote{\Nobs that 
		$ \deg( \R^2 \ni (x_1,x_2) \mapsto x_1 x_2 \in \R ) = 2 $. 
		\Moreover for all $ P \in \polyn_0 $, $ x, y \in \R^0 = \{ 0 \} $ 
		it holds that $ P( x ) = P( y ) = P(0) $ and $ \deg(P) = 0 $.} 
	of a polynomial] 
	\label{def:degree}
	\cfadd{def:polynomial}
	Let $ n \in \N_0 $, $ P \in \polyn_n $ 
	\cfload. 
	Then we denote by 
	$
	\deg( P ) \in \N_0 
	$
	the degree of $ P $. 
\end{definition}

\cfclear
\begin{definition}[Sets of suitable rational functions]
	\label{def:rational:function}
	Let $ n \in \N_0 $, $ \delta \in (0,\infty] $. 
	Then we denote by $ \ratio_{ n, \delta } $ the set given by
	\begin{multline} 
	\label{eq:def_set_rational_functions}
	\cfadd{def:polynomial}
	\ratio_{ n, \delta } 
	= 
	\Biggl\{ 
	R \colon \R^n \to \R \colon 
	\Biggl(
	\exists \, P \in \polyn_n, \, Q \in \{ q \in \polyn_n \colon \deg( q ) < \delta \} 
	\colon 
	\\
	\Biggl[
	\forall \, x \in \R^n \colon  R(x) 
	= 
	\begin{cases}
	[ Q(x) ]^{ - 1 } P(x) 
	& 
	\colon Q( x ) \not= 0 
	\\[0.5ex]
	0 
	& 
	\colon Q ( x ) = 0 
	\end{cases} 
	\Biggr]
	\Biggr)
	\Biggr\}
	\end{multline}
	\cfload.
\end{definition}

\begin{lemma}
	\label{lem:rational:semialgebraic}
	Let $ n \in \N $, $ R \in \scrR_{ n, \infty } $. 
	Then $ R $ is semi-algebraic.
\end{lemma}
\begin{cproof}{lem:rational:semialgebraic}
	\Nobs that the assumption that $ R \in \scrR_{ n, \infty } $ 
	assures that there exist $ P, Q \in \polyn_n $ which satisfy for all 
	$ x \in \R^n $ that
	\begin{equation}
	\label{eq:property_rational_semi-algebraic}
	R(x) = 
	\begin{cases}
	\frac{P(x) }{ Q( x ) } & \colon Q( x ) \not= 0 
	\\
	0 & \colon Q ( x ) = 0
	.
	\end{cases}
	\end{equation}
	\Nobs that \cref{eq:property_rational_semi-algebraic} 
	ensures that 
	\begin{equation}
	\begin{split}
	&
	\operatorname{Graph}( R ) 
	 =
	\cu[\big]{
	(x,y) \in \R^n \times \R
	\colon
	\bigl(
	R(x) = y
	\bigr) }
	\\ &
	= \cu[\Big]{ 
	(x,y) \in \R^n \times \R
	\colon
	\br[\big]{ \rbr[\big]{
	R(x) = y }
	, \rbr[\big]{
	Q(x) = 0 } } }
	\\ &
	\quad
	\cup 
	\cu[\Big]{
	(x,y) \in \R^n \times \R
	\colon
	\br[\big]{\rbr[\big]{
	R(x) = y } , \rbr[\big]{
	Q(x) \neq 0  } } }
	\\ &
	=
	\cu[\Big]{
	(x,y) \in \R^n \times \R
	\colon
	\rbr[\big]{
	y = Q(x) = 0 } }
	\cup 
	\cu[\Big]{ 
	(x,y) \in \R^n \times \R
	\colon \br[\big]{ \rbr[\big]{
	P(x) = y Q(x) } , \rbr[\big]{
	Q(x) \neq 0 } } }
	.
	\end{split}
	\end{equation}
	\Hence 
	\begin{equation}
	\label{eq:graph_representation}
	\begin{split}
	\operatorname{Graph}( R ) 
	 & =
	\Bigl[
	\Bigl\{ 
	(x,y) \in \R^n \times \R
	\colon
	\bigl(
	y = 0
	\bigr) 
	\Bigr\}
	\cap
	\Bigl\{ 
	(x,y) \in \R^n \times \R
	\colon
	\bigl(
	Q(x) = 0
	\bigr) 
	\Bigr\}
	\Bigr]
	\\ & \quad 
	\cup 
	\Bigl[
	\Bigl\{ 
	(x,y) \in \R^n \times \R
	\colon
	\bigl(
	P(x) - y Q(x) = 0 
	< 
	[ Q(x) ]^2 
	\bigr) 
	\Bigr\}
	\Bigr]
	.
	\end{split}
	\end{equation}
	Combining this
	with \cref{eq:def_semi-algebraic_set} establishes that 
	$ R $ is semi-algebraic. 
\end{cproof}

\subsection{Suitable piecewise rational functions as semi-algebraic functions}
\label{ssec:amn:semialgebraic}

\cfclear
\begin{prop}  
\label{prop:amn:semialgebraic}
Let $ m \in \N $, $ f \in \scrF_{ m, 0, \infty } $ \cfadd{def:function:amn}\cfload. 
Then 
$
  \R^m \ni \theta \mapsto f(\theta,0) \in \R 
$ 
is semi-algebraic \cfadd{def:semialgebraic:function}\cfload.
\end{prop}
\begin{cproof}{prop:amn:semialgebraic}
\Nobs that \cref{eq:def_function_amn} 
and the assumption that $ f \in \scrF_{ m, 0, \infty } $ 
assure that there exist $ V \in \N $, $ r_1, r_2, \dots, r_V \in \N $, 
$ R_1, R_2, \dots, R_V \in \scrR_{ m, \infty} $, 
$ 
  P^1 = ( P^1_i )_{ i \in \{ 1, 2, \dots, r_1 \} } 
  \subseteq \scrP_m 
$,
$ 
  P^2 = ( P^2_i )_{ i \in \{ 1, 2, \dots, r_2 \} } 
  \subseteq \scrP_m 
$,
$ \dots $, 
$ 
  P^V = ( P^V_i )_{ i \in \{ 1, 2, \dots, r_V \} } 
  \subseteq \scrP_m 
$
which satisfy 
for all $ \theta \in \R^m $ that
\begin{equation}
\label{eq:f_semi_algebraic_property_rectified}
  f( \theta , 0 ) 
  = 
  \sum_{ v=1 }^V 
  \Bigg[ 
    R_v( \theta ) 
    \br*{
    \textstyle
      \prod\limits_{ i = 1 }^{ r_v } 
      \indicator{ [0, \infty ) }\big( 
        P^v_i( \theta ) 
      \big)
    }
  \Bigg]
  .
\end{equation}
\Nobs that \cref{eq:f_semi_algebraic_property_rectified} 
shows for all $ \theta \in \R^m $ that
\begin{equation}
\label{eq:f_semi_algebraic_property_rectified2}
  f( \theta , 0 ) 
  = 
  \sum_{ v=1 }^V 
  \Bigg[ 
    R_v( \theta ) 
    \br*{
    \textstyle
      \prod\limits_{ i = 1 }^{ r_v } 
      \indicator{ \{ \vartheta \in \R^m \colon P^v_i( \vartheta ) \geq 0 \} }( \theta )
    }
  \Bigg]
  .
\end{equation}
\Moreover \cref{prop:semialgebraic:sets} and \cref{lem:indicator_semialgebraic} 
prove that 
for all $ v \in \{ 1, 2, \dots, V \} $, 
$ i \in \{ 1, 2, \dots, r_v \} $ 
it holds that 
\begin{equation}
\label{eq:indicator_is_semi_algebraic}
  \R^m \ni \theta \mapsto 
  \indicator{ \{ \vartheta \in \R^m \colon P^v_i( \vartheta ) \geq 0 \} }( \theta )
  \in \R
\end{equation}
is semi-algebraic. 
\Moreover \cref{lem:rational:semialgebraic} assures 
that for all 
$
  v \in \{ 1, 2, \dots, V \}
$
it holds that 
$
  R_v
$
is semi-algebraic. 
Combining this and 
\cref{eq:indicator_is_semi_algebraic} 
with 
\cref{lem:sum:semialgebraic}
shows that for all $ v \in \{ 1, 2, \dots, V \} $ 
it holds that 
\begin{equation}
  \R^m 
  \ni 
  \theta 
  \mapsto 
  R_v( \theta ) 
  \br*{
  \textstyle
    \prod_{ i = 1 }^{ r_v } 
    \indicator{ \{ \vartheta \in \R^m \colon P^v_i( \vartheta ) \geq 0 \} }( \theta )
  }
  \in \R
\end{equation}
is semi-algebraic. 
\cref{lem:sum:semialgebraic} and \cref{eq:f_semi_algebraic_property_rectified2} 
\hence show that 
$
  \R^m \ni \theta \mapsto f( \theta, 0 ) \in \R
$
is semi-algebraic. 
\end{cproof}

\subsection{Closedness under parametric integration 
of suitable piecewise rational functions}
\label{ssec:closedness_integration}

\cfclear
\begin{prop} 
\label{prop:integrals:amn}
Let $ m, n \in \N $, $ a \in \R $, $ b \in ( a, \infty) $, 
$ f \in \scrF_{ m, n, \infty } $ \cfadd{def:function:amn}\cfload. 
Then
\begin{enumerate}[label = (\roman*)]
  \item 
  \label{prop:int:semialgebraic:item1} 
  it holds 
  for all $ \theta \in \R^m $, $ x_1, x_2, \dots, x_{ n - 1 } \in \R $ 
  that 
  $
    \int_a^b | f( \theta , x_1, x_2, \ldots, x_{ n - 1 }, x_n ) | \, \d x_n     
    < \infty 
  $
  and
  \item 
  \label{prop:int:semialgebraic:item2} 
  it holds that 
\begin{equation}
\textstyle
  \bigl( 
    \R^m \times \R^{n-1} \ni 
    ( \theta, x_1, \ldots, x_{ n - 1 } ) 
    \mapsto \int_a^b f ( \theta , x_1, x_2, \ldots, x_{ n - 1 }, x_n ) \, \d x_n 
    \in \R
  \bigr) 
  \in \scrF_{ m, n - 1, \infty } 
  .
\end{equation}
\end{enumerate}
\end{prop}

\begin{cproof}{prop:integrals:amn}
\Nobs that \cref{eq:def_function_amn} 
and the fact that 
$
  \{ q \in \polyn_n \colon \deg(q) \leq \infty \} 
  =
  \polyn_n
  \subseteq C( \R^n, \R ) 
$
prove that for all 
$ \theta \in \R^m $, $ r \in (0,\infty) $ it holds that 
\begin{equation}
  \sup\nolimits_{ 
    x \in [-r,r]^n
  }
  | f( \theta, x ) |
  < \infty 
  .
\end{equation}
This shows \cref{prop:int:semialgebraic:item1}. 
\Moreover \cite[Proposition 4.8]{EberleJentzenRiekertWeiss2021} 
and \cref{item:properties_F_v} in \cref{lem:bmn:amn} 
establish \cref{prop:int:semialgebraic:item2}. 
\end{cproof}

\subsection{Closedness under rectification of suitable piecewise rational functions}

The next result, \cref{lem:maximum:bmn}, establishes that the function classes $\scrF_{m,n,1}$ introduced in \cref{def:function:amn} above are closed under composition with the ReLU function.  This will be used to show in \cref{prop:realization:bmn} below
that these function classes contain the realization functions of DNNs with ReLU activation.

\cfclear
\begin{lemma} \label{lem:maximum:bmn}
Let $ m, n \in \N $, $ f \in \scrF_{ m, n, 1 } $ \cfadd{def:function:amn}\cfload.
Then 
\begin{equation}
\label{eq:closed_rectification}
  \big(
    \R^m \times \R^n 
    \ni 
    v
    \mapsto
    \max\cu{ f(v) , 0 } \in 
    \R 
  \big) 
  \in 
  \scrF_{ m, n, 1 } 
  .
\end{equation}
\end{lemma}
\begin{cproof}{lem:maximum:bmn}
\Nobs that \cref{eq:def_function_amn} 
and the assumption that 
$ f \in \scrF_{ m, n, 1 } $ 
ensure that there exist 
$ V \in \N $, $ r_1, r_2, \ldots, r_V \in \N $, 
$ R_1, R_2, \dots, R_V \in \scrR_{ m, 1 } $, 
$ Q_1, Q_2, \dots, Q_V \in \{ q \in \polyn_n \colon \deg( q ) \leq 1 \} $, 
$ 
  P^1 = ( P^1_{ i, j } )_{ (i,j) \in \{ 1, 2, \dots, r_1 \} \times \{ 0, 1, \dots, n \} } 
  \subseteq \scrP_m 
$,
$ 
  P^2 = ( P^2_{ i, j } )_{ (i,j) \in \{ 1, 2, \dots, r_2 \} \times \{ 0, 1, \dots, n \} } 
  \subseteq \scrP_m 
$,
$ \dots $, 
$ 
  P^V = ( P^V_{ i, j } )_{ (i,j) \in \{ 1, 2, \dots, r_V \} \times \{ 0, 1, \dots, n \} } 
  \subseteq \scrP_m 
$
which satisfy 
for all $ \theta \in \R^m $, $ x = (x_1, \dots, x_n) \in \R^n $ that
\begin{equation}
\label{eq:f_property_rectified}
  f( \theta , x ) 
  = 
  \sum_{ v=1 }^V 
  \Bigg[ 
    R_v( \theta ) Q_v( x ) 
    \br*{
    \textstyle
      \prod\limits_{ i = 1 }^{ r_v } 
      \indicator{ [0, \infty ) }\big( 
        P^v_{ i, 0 }( \theta ) 
        +
        \smallsum_{ j = 1 }^n P^v_{ i, j }( \theta ) x_j 
      \big)
    }
  \Bigg]
  .
\end{equation}
In the following let 
$
  p_v \colon \R^m \times \R^n \to \{ 0, 1 \} 
$, 
$ v \in \{ 1, 2, \dots, V \} $, 
satisfy for all 
$ v \in \{ 1, 2, \dots, V \} $, 
$ \theta \in \R^m $, 
$ x = (x_1, \dots, x_n) \in \R^n $
that 
\begin{equation}
\label{eq:def:p_v_rectifier}
\textstyle
  p_v( \theta, x ) 
  =
  \prod\limits_{ i = 1 }^{ r_v } 
  \indicator{ [0, \infty ) }\big( 
    P^v_{ i, 0 }( \theta ) 
    +
    \smallsum_{ j = 1 }^n P^v_{ i, j }( \theta ) x_j 
  \big) ,
\end{equation}
for every $ W \subseteq \{ 1, 2, \dots, V \} $ 
let 
$ \fp_W \in \{ 0, 1 \} $
satisfy 
\begin{equation}
\label{eq:rectified_cW_def}
\textstyle 
  \fp_W 
  =
  \biggl[ 
    \prod\limits_{ v \in W }
    p_v( \theta, x ) 
  \biggr] 
  \biggl[ 
    \prod\limits_{ v \in \{ 1, 2, \dots, V \} \backslash W }
    (
      1 - p_v( \theta, x )
    )
  \biggr] 
  ,
\end{equation}
and for every $ \theta \in \R^m $, 
$ x \in \R^n $ 
let 
$ 
  \scrV_{ \theta, x } \subseteq \N
$ 
satisfy 
\begin{equation}
\label{eq:def_script_V}
  \scrV_{ \theta, x } = \{ v \in \{ 1, 2, \dots, V \} \colon p_v( \theta, x ) = 1 \}
  .
\end{equation} 
\Nobs that
\cref{eq:rectified_cW_def,eq:def:p_v_rectifier,eq:def_script_V} 
assure that for all $ W \subseteq \{ 1, 2, \dots, V \} $ 
it holds that 
\begin{equation}
\label{eq:property_fp_W}
\textstyle
  \fp_W
  = 
  \begin{cases}
    1 
  &
    \colon W = \scrV_{ \theta, x }
  \\
    0
  &
    \colon W \neq \scrV_{ \theta, x } .
  \end{cases}
\end{equation}
Combining this with 
\cref{eq:f_property_rectified,eq:def:p_v_rectifier,eq:def_script_V} 
proves that for all $ \theta \in \R^m $, $ x \in \R^n $ 
it holds that 
\begin{equation}
\label{eq:f_property_rectified_2}
\begin{split}
  f( \theta, x ) 
& =
  \sum_{ v = 1 }^V
  R_v( \theta ) Q_v( x ) p_v( \theta, x )
  =
  \sum_{
    v \in \scrV_{ \theta, x }
  }
  R_v( \theta ) Q_v( x ) p_v( \theta, x ) 
\\ &
  =
  \sum_{
    v \in \scrV_{ \theta, x }
  }
  R_v( \theta ) Q_v( x ) 
  =
  \sum_{
    W \subseteq \{ 1, 2, \dots, V \} 
  }
  \left(
    \fp_W
    \textstyle
    \left[ 
      \sum\limits_{
        v \in W
      }
      R_v( \theta ) Q_v( x ) 
    \right]
  \right)
  .
\end{split}
\end{equation}
This and \cref{eq:property_fp_W} show that for all 
$ \theta \in \R^m $, $ x \in \R^n $ it holds that
\begin{equation}
\label{eq:f_property_rectified_3}
\begin{split}
  \max\{ f( \theta, x ) , 0 \}
& =
  \sum_{
    W \subseteq \{ 1, 2, \dots, V \} 
  }
  \left(
    \fp_W
    \textstyle
    \max\!\left\{ 
      \sum\limits_{
        v \in W
      }
      R_v( \theta ) Q_v( x ) 
      , 0
    \right\}
  \right)
  .
\end{split}
\end{equation}
The fact that for all $ r \in \R $ it holds that 
$
  \max\{ r, 0 \} = r \indicator{ [0,\infty) }( r )
$ 
\hence demonstrates that for all 
$ \theta \in \R^m $, $ x \in \R^n $ it holds that
\begin{equation}
\label{eq:f_property_rectified_4}
\begin{split}
  \max\{ f( \theta, x ) , 0 \}
& =
  \sum_{
    W \subseteq \{ 1, 2, \dots, V \} 
  }
  \left(
    \fp_W
    \textstyle
    \left[
      \sum\limits_{
        v \in W
      }
      R_v( \theta ) Q_v( x ) 
    \right]
    \indicator{ [0,\infty) }\!\left(
      \sum\limits_{
        v \in W
      }
      R_v( \theta ) Q_v( x ) 
    \right)
  \right)
\\ & =
  \sum_{
    W \subseteq \{ 1, 2, \dots, V \} 
  }
  \sum\limits_{
    w \in W
  }
  \left(
    R_w( \theta ) Q_w( x ) 
    \left[
      \fp_W
      \textstyle
      \indicator{ [0,\infty) }\!\left(
        \sum\limits_{
          v \in W
        }
        R_v( \theta ) Q_v( x ) 
      \right)
    \right]
  \right)
  .
\end{split}
\end{equation}
\Moreover 
\cref{eq:rectified_cW_def} 
shows that 
\begin{equation}
\begin{split}
  \fp_W 
  & =
\textstyle
  \biggl[ 
    \prod\limits_{ v \in W }
    p_v( \theta, x ) 
  \biggr] 
  \biggl[ 
    \sum\limits_{
      U \subseteq ( \{ 1, 2, \dots, V \} \backslash W )
    }
    ( - 1 )^{ \#( U ) }
    \bigl[
    \prod_{ v \in U }
    p_v( \theta, x )
    \bigr]
  \biggr] 
\\ & =
  \sum\limits_{
    U \subseteq ( \{ 1, 2, \dots, V \} \backslash W )
  }
  \left(
  ( - 1 )^{ \#( U ) }
\textstyle
    \biggl[ 
      \prod\limits_{ v \in W }
      p_v( \theta, x ) 
    \biggr] 
    \biggl[ 
      \prod\limits_{ v \in U }
      p_v( \theta, x )
    \biggr] 
  \right) .
\end{split}
\end{equation}
Combining this and \cref{eq:f_property_rectified_4} 
proves that 
for all $ \theta \in \R^m $, $ x \in \R^n $ 
it holds that 
\begin{multline}
\label{eq:f_property_rectified_5}
  \max\{ f( \theta, x ) , 0 \}
=
  \sum_{
    W \subseteq \{ 1, 2, \dots, V \} 
  }
  \sum\limits_{
    w \in W
  }
  \sum\limits_{
    U \subseteq ( \{ 1, 2, \dots, V \} \backslash W )
  }
  \Biggl(
    ( - 1 )^{ \#( U ) }
    R_w( \theta ) Q_w( x ) 
\\ 
    \cdot 
    \biggl[
      \textstyle
      \indicator{ [0,\infty) }\!\left(
        \sum\limits_{
          v \in W
        }
        R_v( \theta ) Q_v( x ) 
      \right)
    \biggr]
    \biggl[ 
      \prod\limits_{ v \in W }
      p_v( \theta, x ) 
    \biggr] 
    \biggl[ 
      \prod\limits_{ v \in U }
      p_v( \theta, x )
    \biggr] 
  \Biggr)
  .
\end{multline}
\Moreover 
\cref{item:properties_F_ib}
in 
\cref{lem:bmn:amn}, 
\cref{eq:def_function_amn}, 
\cref{eq:def:p_v_rectifier}, 
and 
the fact that 
$ 
  R_1, R_2, \dots, R_V \in \scrR_{ m, 1 }
$
show that
for all 
$ W \subseteq \{ 1, 2, \dots, V \} $, 
$
  U \subseteq ( \{ 1, 2, \dots, V \} \backslash W )
$
and all 
$ w \in W $
it holds that 
\begin{multline}
  \Biggl(
    \R^m \times \R^n \ni ( \theta, x ) 
    \mapsto 
    R_w( \theta ) Q_w( x ) 
    \biggl[
      \textstyle
      \indicator{ [0,\infty) }\!\left(
        \sum\limits_{
          v \in W
        }
        R_v( \theta ) Q_v( x ) 
      \right)
    \biggr]
\\
\textstyle
    \cdot 
    \biggl[ 
      \prod\limits_{ v \in W }
      p_v( \theta, x ) 
    \biggr] 
    \biggl[ 
      \prod\limits_{ v \in U }
      p_v( \theta, x )
    \biggr] 
    \in \R
  \Biggr)
  \in \scrF_{ m, n, 1 }
  .
\end{multline}
Combining this and \cref{eq:f_property_rectified_5} 
with \cref{eq:def_function_amn} 
establishes \cref{eq:closed_rectification}. 
\end{cproof}

\subsection{Realization functions of DNNs as suitable piecewise rational functions}
\label{ssec:realization_DNNs}

\cfclear 
\begin{prop} 
\label{prop:realization:bmn}
Assume \cref{setting:dnn}. 
Then it holds for all 
$ i \in \cu{1, 2, \ldots, \ell_L} $
that
\begin{equation} 
\label{eq:to_prove_DNN_realization_in_F}
\cfadd{def:function:amn}
  \bigl(
    \R^{ \fd } \times \R^{ \ell_0 } \ni( \theta, x ) 
    \mapsto \cN^{ L, \theta }_{ \infty, i }( x ) 
    \in \R 
  \bigr)
  \in \scrF_{ \fd, \ell_0, 1 }
\end{equation}
\cfload.
\end{prop}

\begin{cproof}{prop:realization:bmn}
\Nobs that \cref{eq:def_NN_realization} ensures  
for all $ k \in \N_0 $, 
$ \theta \in \R^{ \fd } $, 
$ i \in \{ 1, 2, \dots, \ell_{ k + 1 } \} $, 
$ x = ( x_1, \dots, x_{ \ell_0 } ) \in \R^{ \ell_0 } $ 
that
\begin{equation}
\label{eq:NN_realization_in_F}
  \cN_{ \infty, i }^{ k + 1, \theta }( x ) 
  =
  \begin{cases}
    \b{ k + 1, \theta }_i 
    + \smallsum_{ j = 1 }^{ \ell_k } 
    \w{ k + 1, \theta }_{ i, j } x_j .
  &
    \colon
    k = 0
  \\[1ex]
    \b{ k + 1, \theta }_i 
    + \smallsum_{ j = 1 }^{ \ell_k } 
    \w{ k + 1, \theta }_{ i, j } \max\{ \cN_{ \infty, j }^{ k, \theta }( x ) , 0 \}.
  &
    \colon
    k > 0
    .
  \end{cases}
\end{equation}
Next we claim that for all 
$ k \in \{ 1, 2, \dots, L \} $  
it holds that
\begin{equation}
\label{lem:realization:bmn:eq:inductclaim}
\textstyle
  \bigl(
    \bigcup_{ i = 1 }^{ \ell_k }
    \bigl\{
      \R^\fd \times \R^{ \ell_0 } \ni (\theta, x) 
      \mapsto \cN^{ k, \theta }_{ \infty, i }( x ) \in \R 
    \bigr\}
  \bigr)
  \subseteq 
  \scrF_{ \fd, \ell_0, 1 }
  .
\end{equation}
In the following we prove \cref{lem:realization:bmn:eq:inductclaim} by induction 
on $ k \in \{ 1, 2, \dots, L \} $. 
For the base case $ k = 1 $ 
\nobs that 
\cref{eq:NN_realization_in_F} assures that 
for all 
$ \theta \in \R^{ \fd } $, 
$ i \in \cu{ 1, 2, \ldots, \ell_1 } $, 
$ x = ( x_1, \dots, x_{ \ell_0 } ) \in \R^{ \ell_0 } $ 
it holds that
\begin{equation}
  \cN_{ \infty, i }^{ 1, \theta }( x ) 
  = 
  \b{ 1, \theta }_i 
  + \smallsum_{ j = 1 }^{ \ell_0 } 
  \w{ 1, \theta }_{ i, j } x_j .
\end{equation}
This establishes \cref{lem:realization:bmn:eq:inductclaim} in the case $ k = 1 $. 
For the induction step \nobs 
that \cref{lem:maximum:bmn} 
implies that 
for all $ k \in \N \cap (0,L) $, 
$ j \in \cu{ 1, 2, \ldots, \ell_k } $ 
with 
$
  (
    \bigcup_{ i = 1 }^{ \ell_k }
    \{
      \R^\fd \times \R^{ \ell_0 } \ni (\theta, x) 
      \mapsto \cN^{ k, \theta }_{ \infty, i }( x ) \in \R 
    \}
  )
  \subseteq 
  \scrF_{ \fd, \ell_0, 1 }
$
it holds that 
\begin{equation}
\label{eq:max_in_F_1_proof}
  \rbr[\big]{
    \R^{ \fd } \times \R^{ \ell_0 } 
    \ni ( \theta, x ) 
    \mapsto 
    \max\cu{ 
      \cN^{ k, \theta }_{ \infty, j }( x )
      , 0 
    } 
    \in \R 
  } 
  \in \scrF_{ \fd, \ell_0, 1 }.
\end{equation}
\Moreover \cref{eq:NN_realization_in_F} 
shows that 
for all $ k \in \N \cap (0,L) $, 
$ \theta \in \R^{ \fd } $, 
$ i \in \cu{ 1, 2, \ldots, \ell_{ k + 1 } } $, 
$ x \in \R^{ \ell_0 } $
we have that
\begin{equation}
\label{eq:cN_rep_proof_NN_realization_in_F}
  \cN_{ \infty, i }^{ k + 1, \theta }( x ) 
  = 
  \b{ k + 1, \theta }_i 
  + 
  \smallsum_{ j = 1 }^{ \ell_k } 
  \w{ k + 1, \theta }_{ i, j } 
  \max\cu{ 
    \cN^{ k, \theta }_{ \infty, j }( x ), 0 
  }
  .
\end{equation}
\Moreover \cref{eq:def_function_amn,eq:max_in_F_1_proof} 
demonstrate 
that for all 
$ k \in \N \cap (0,L) $, 
$ i \in \cu{ 1, 2, \ldots, \ell_{ k + 1 } } $, 
$ j \in \cu{ 1, 2, \ldots, \ell_k } $
with 
$
  (
    \bigcup_{ v = 1 }^{ \ell_k }
    \{
      \R^\fd \times \R^{ \ell_0 } \ni (\theta, x) 
      \mapsto \cN^{ k, \theta }_{ \infty, v }( x ) \in \R 
    \}
  )
  \subseteq 
  \scrF_{ \fd, \ell_0, 1 }
$ 
it holds that 
\begin{equation}
  \bigl(
    \R^{ \fd } \times \R^{ \ell_0 }
    \ni 
    ( \theta, x )
    \mapsto 
    \w{ k + 1, \theta }_{ i, j } 
    \max\cu{ 
      \cN^{ k, \theta }_{ \infty, j }( x ), 0 
    }
    \in \R
  \bigr)
  \in 
  \scrF_{ \fd, \ell_0, 1 }
  .
\end{equation}
The fact that 
$
  \scrF_{ \fd, \ell_0, 1 }
$
is an $ \R $-vector space 
and \cref{eq:def_function_amn} 
hence show that for all 
$ k \in \N \cap (0,L) $, 
$ i \in \cu{ 1, 2, \ldots, \ell_{ k + 1 } } $ 
with 
$
  (
    \bigcup_{ v = 1 }^{ \ell_k }
    \{
      \R^\fd \times \R^{ \ell_0 } \ni (\theta, x) 
      \mapsto \cN^{ k, \theta }_{ \infty, v }( x ) \in \R 
    \}
  )
  \subseteq 
  \scrF_{ \fd, \ell_0, 1 }
$ 
it holds that 
\begin{equation}
\textstyle
  \bigl(
    \R^{ \fd } \times \R^{ \ell_0 }
    \ni 
    ( \theta, x )
    \mapsto 
    \b{ k + 1, \theta }_i 
    +
    \sum_{ j = 1 }^{ \ell_k }
    \w{ k + 1, \theta }_{ i, j } 
    \max\cu{ 
      \cN^{ k, \theta }_{ \infty, j }( x ), 0 
    }
    \in \R
  \bigr)
  \in 
  \scrF_{ \fd, \ell_0, 1 }
  .
\end{equation}
This and \cref{eq:cN_rep_proof_NN_realization_in_F} 
assure that for all 
$ k \in \N \cap (0,L) $
with 
$
  (
    \bigcup_{ i = 1 }^{ \ell_k }
    \{
      \R^\fd \times \R^{ \ell_0 } \ni (\theta, x) 
      \mapsto \cN^{ k, \theta }_{ \infty, i }( x ) \in \R 
    \}
  )
  \subseteq 
  \scrF_{ \fd, \ell_0, 1 }
$ 
it holds that 
\begin{equation}
\textstyle
  \bigl(
    \bigcup_{ i = 1 }^{ \ell_{ k + 1 } }
    \bigl\{
      \R^{ \fd } \times \R^{ \ell_0 }
      \ni 
      ( \theta, x )
      \mapsto 
      \cN^{ k + 1, \theta }_{ \infty, i }( x ) 
      \in \R
    \bigr\} 
  \bigr)
  \subseteq 
  \scrF_{ \fd, \ell_0, 1 }
  .
\end{equation}
Induction thus proves \cref{lem:realization:bmn:eq:inductclaim}. 
\Nobs that \cref{lem:realization:bmn:eq:inductclaim} 
establishes \cref{eq:to_prove_DNN_realization_in_F}. 
\end{cproof}

\section{Piecewise polynomial functions}
\label{sec:piecewise_polynomial}

In this section we establish 
in \cref{cor:loss:semialgebraic}
in \cref{ssec:semi_algebraic_risk} below 
that in the set-up of \cref{setting:dnn} 
in \cref{subsection:dnn:framework} above
we have,  
under 
the assumption that 
the measure $ \mu $ is absolutely continuous 
with density 
% function 
$
  \dens \colon [a,b]^{ \ell_0 } \to \R 
$
and 
the assumption that 
the density function
$ \dens \colon [a,b]^{ \ell_0 } \to \R $
and every component 
of the target function 
$ 
  f = ( f_1, \dots, f_{ \ell_L } ) \colon 
  [a,b]^{ \ell_0 } \to \R^{ \ell_L } 
$
are piecewise polynomial 
in the sense of \cref{def:multidim:piece:polyn} 
in \cref{subsec:piecewise_polynomial} below, 
that the risk function 
$
  \cL_{ \infty } \colon [a,b]^{ \ell_0 } \to \R
$
is semi-algebraic. 
In \cref{sec:KL} below we will employ 
\cref{cor:loss:semialgebraic}
to conclude that for every 
ANN parameter vector 
$ \theta \in \R^{ \fd } $ 
we have that 
the risk function 
$
  \cL_{ \infty } \colon [a,b]^{ \ell_0 } \to \R
$
satisfies a generalized 
Kurdyka-\L ojasiewicz inequality 
on a neighbourhood 
of $ \theta $. 

Throughout this work we consider fully connected feedforward ANNs.
For different network architectures such as convolutional neural networks (CNNs) or residual neural networks
it might be possible to establish analogous results by suitably adapting our arguments.

Our proof of \cref{cor:loss:semialgebraic} mainly relies 
on \cref{prop:integrals:amn} in \cref{ssec:closedness_integration} above, 
on \cref{prop:realization:bmn} in \cref{ssec:realization_DNNs} above, 
as well as 
on the fact that for all $ m \in \N $
it holds that 
functions in 
$ \scrF_{ m, 0, \infty } $, 
are semi-algebraic 
according to 
\cref{prop:amn:semialgebraic} in \cref{ssec:amn:semialgebraic} above.

Some of the concepts and results in this section 
are inspired by our previous article 
Eberle et al.~\cite[Section~4]{EberleJentzenRiekertWeiss2021}. 
In particular, \cref{def:multidim:piece:polyn} 
is a slight extension of 
\cite[Definition~4.9]{EberleJentzenRiekertWeiss2021} 
and \cref{cor:loss:semialgebraic} 
extends \cite[Corollary~4.10]{EberleJentzenRiekertWeiss2021} 
from the situation of shallow ReLU ANNs with one hidden layer
to deep ReLU ANNs with an arbitrarily large number of 
hidden layers. 

It should also be noted that \cref{cor:loss:semialgebraic} is a novel contribution mainly due to the fact that we consider the true risk defined by the integral over the entire input data. If one considers the empirical risk (calculated from a finite set of input data data) an analogous result is already known, cf.~Davis et al.~\cite[Corollary 5.11]{Davis2018stochastic}.

\subsection{Piecewise polynomial functions}
\label{subsec:piecewise_polynomial}

\cfclear
\begin{definition}[Piecewise polynomial functions] 
\label{def:multidim:piece:polyn}
Let $ d \in \N $, 
let $ A \subseteq \R^d $ and $ B \subseteq \R $ be sets, 
and let $ f \colon A \to B $ be a function. 
Then we say that $ f $ is piecewise polynomial 
if and only if there exist $n \in \N$, 
$ \alpha_1, \alpha_2, \ldots, \alpha_n \in \R^{ n \times d } $, 
$ \beta_1, \beta_2, \ldots, \beta_n \in \R^n $,
$ P_1, P_2, \ldots, P_n \in \polyn_d $ 
such that for all $ x \in A $
it holds that
\begin{equation} 
\label{eq:def_piecewise_polynomial}
\cfadd{def:polynomial}
  f(x) = 
  \smallsum_{ i = 1 }^n 
%   \br*{ 
    P_i(x) \indicator{[0, \infty )^n} \rbr{ \alpha_i x + \beta_i } 
%   } 
\end{equation}
\cfload.
\end{definition}

\subsection{Characterization results for piecewise polynomial functions}

The following results, \cref{lem:piecewise_polynomial_charac} to \cref{prop:piecewise_polynomial_extensionB},
are elementary consequences of the definition of piecewise polynomial functions. They will be employed in the proof of \cref{cor:loss:semialgebraic} to show that the risk function is semi-algebraic if the density function and every component of the target function are piecewise polynomial.

\cfclear
\begin{lemma}
\label{lem:piecewise_polynomial_charac}
\cfadd{def:multidim:piece:polyn}
Let $ d \in \N $, 
let $ A \subseteq \R^d $ be a set, 
and let 
$
  f \colon A \to \R
$
be a function. 
Then the following three statements are equivalent:
\begin{enumerate}[label=(\roman*)]

\item
\label{item:piecewise_polynomial_charac:i}
It holds that $ f $ is piecewise polynomial 
\cfload. 

\item 
\label{item:piecewise_polynomial_charac:ii}
There exist 
$ m, n \in \N $, 
$ \alpha_1, \alpha_2, \dots, \alpha_m \in \R^{ n \times d } $, 
$ \beta_1, \beta_2, \dots, \beta_m \in \R^n $,
$ P_1, P_2, \dots, P_m \in \polyn_d $
such that for all $ x \in A $ it holds that 
\begin{equation}
\textstyle
  f(x) = \sum_{ i = 1 }^m P_i( x ) \mathbbm{1}_{ [0,\infty)^n }( \alpha_i x + \beta_i )
  .
\end{equation}

\item 
\label{item:piecewise_polynomial_charac:iii}
There exist 
$ n \in \N $, 
$ m_1, m_2, \dots, m_n \in \N $, 
$ \alpha_1 \in \R^{ m_1 \times d } $, 
$ \alpha_2 \in \R^{ m_2 \times d } $, $ \dots $, 
$ \alpha_n \in \R^{ m_n \times d } $, 
$ \beta_1 \in \R^{ m_1 } $, 
$ \beta_2 \in \R^{ m_2 } $, 
$ \dots $, 
$ \beta_m \in \R^{ m_n } $,
$ P_1, P_2, \dots, P_n \in \polyn_d $
such that for all $ x \in A $ it holds that 
\begin{equation}
\textstyle
  f(x) = 
  \sum_{ i = 1 }^n P_i( x ) 
  \mathbbm{1}_{ [0,\infty)^{ m_i } }( 
    \alpha_i x + \beta_i 
  )
  .
\end{equation}

\end{enumerate}
\end{lemma}
\begin{cproof}{lem:piecewise_polynomial_charac}
Throughout this proof for every $ m \in \N $ let 
$ {\bf e}^{ (m) }_1, {\bf e}^{ (m) }_2, \dots, {\bf e}^{ (m) }_m \in \R^n $ 
satisfy 
$
  {\bf e}^{ (m) }_1 = (1, 0, \dots, 0)
$,
$
  {\bf e}^{ (m) }_2 = (0, 1, 0, \dots, 0)
$,
$ \dots $,
$
  {\bf e}^{ (m) }_m = (0, \dots, 0, 1)
$
and for every $ m, n \in \N $ 
with $ m \geq n $ 
let $ A_{ m, n } \in \R^{ m \times n } $ 
satisfy for all 
$ x = ( x_1, \dots, x_n ) \in \R^n $ 
that 
\begin{equation}
\label{eq:def_A_m_n_matrix}
\textstyle
  A_{ m, n } x 
  = 
  \sum\limits_{ i = 1 }^n x_i {\bf e}_i^{ ( m ) }
  .
\end{equation}
\Nobs that 
\cref{eq:def_A_m_n_matrix} ensures that 
for all 
$ n \in \N $, 
$ m_1, m_2, \dots, m_n \in \N $, 
$ \alpha_1 \in \R^{ m_1 \times d } $, 
$ \alpha_2 \in \R^{ m_2 \times d } $, $ \dots $, 
$ \alpha_n \in \R^{ m_n \times d } $, 
$ \beta_1 \in \R^{ m_1 } $, 
$ \beta_2 \in \R^{ m_2 } $, 
$ \dots $, 
$ \beta_m \in \R^{ m_n } $,
$ P_1, P_2, \dots, P_n \in \polyn_d $, 
$ x \in A $
it holds that 
\begin{equation}
\begin{split}
\textstyle
&
  \sum_{ i = 1 }^n P_i( x ) 
  \mathbbm{1}_{ [0,\infty)^{ m_i } }( 
    \alpha_i x + \beta_i 
  )
\\ &
  =
  \sum_{ i = 1 }^n 
  \bigl[
    P_{ \min\{ i, n \} }( x ) 
  \bigr]
  \bigl[
    \mathbbm{1}_{ 
      [0,\infty)^{ m_{ \min\{ i, n \} } } 
    }\bigl( 
      \alpha_{ \min\{ i, n \} } x + \beta_{ \min\{ i, n \} } 
    \bigr)
  \bigr]
\\ &
  =
  \sum\limits_{ i = 1 }^{ m_1 + \ldots + m_n } 
  \bigl[
    P_{ \min\{ i, n \} }( x ) 
    \mathbbm{1}_{
      [1,n]
    }( i )
  \bigr]
  \bigl[
    \mathbbm{1}_{ 
      [0,\infty)^{ m_{ \min\{ i, n \} } } 
    }\bigl( 
      \alpha_{ \min\{ i, n \} } x + \beta_{ \min\{ i, n \} } 
    \bigr)
  \bigr]
\\ &
  =
  \sum\limits_{ i = 1 }^{ 
    m_1 + \ldots + m_n 
  } 
  \bigl[
    P_{ \min\{ i, n \} }( x ) 
    \mathbbm{1}_{
      [1,n]
    }( i )
  \bigr]
\\ & \cdot 
  \bigl[
    \mathbbm{1}_{ 
      [0,\infty)^{ m_1 + \ldots + m_n } 
    }\bigl( 
      \bigl[ 
        A_{ 
          m_1 + \ldots + m_n, m_{ \min\{ i, n \} } 
        } 
        \alpha_{ \min\{ i, n \} } 
      \bigr]
      x 
      + 
      \bigl[ 
        A_{ 
          m_1 + \ldots + m_n, m_{ \min\{ i, n \} } 
        } 
        \beta_{ \min\{ i, n \} } 
      \bigr]
    \bigr)
  \bigr]
  .
\end{split}
\end{equation}
Combining this with \cref{eq:def_piecewise_polynomial} 
establishes that 
for all 
$ n \in \N $, 
$ m_1, m_2, \dots, m_n \in \N $, 
$ \alpha_1 \in \R^{ m_1 \times d } $, 
$ \alpha_2 \in \R^{ m_2 \times d } $, $ \dots $, 
$ \alpha_n \in \R^{ m_n \times d } $, 
$ \beta_1 \in \R^{ m_1 } $, 
$ \beta_2 \in \R^{ m_2 } $, 
$ \dots $, 
$ \beta_m \in \R^{ m_n } $,
$ P_1, P_2, \dots, P_n \in \polyn_d $
with 
$  
  \forall \, x \in A \colon
  f(x) = \sum_{ i = 1 }^n P_i( x ) 
  \mathbbm{1}_{ [0,\infty)^{ m_i } }( \alpha_i x + \beta_i )
$
it holds that 
$ f $ is piecewise polynomial. 
\end{cproof}

\subsection{Sums and products of piecewise polynomial functions}

\cfclear
\begin{lemma}
\label{lem:piecewise_polynomial_product}
\cfadd{def:multidim:piece:polyn}
Let $ d \in \N $, 
let $ A \subseteq \R^d $ be a set, 
and let 
$
  f \colon A \to \R
$
and 
$
  g \colon A \to \R
$
be piecewise polynomial 
\cfload. 
Then 
\begin{enumerate}[label=(\roman*)]
\item 
\label{item:piecewise_polynomial_product_i}
it holds that 
$ A \ni x \mapsto f(x) + g(x) \in \R $
is piecewise polynomial and 
\item 
\label{item:piecewise_polynomial_product_ii}
it holds that 
$ A \ni x \mapsto f(x) g(x) \in \R $
is piecewise polynomial.
\end{enumerate}
\end{lemma}
\begin{cproof}{lem:piecewise_polynomial_product}
\Nobs that 
\cref{eq:def_piecewise_polynomial}, 
the assumption that $ f $ is piecewise polynomial, 
and the assumption that $ g $ is piecewise polynomial 
ensure that there exist $ n, m \in \N $, 
$ \alpha_1, \alpha_2, \ldots, \alpha_n \in \R^{ n \times d } $, 
$ \fa_1, \fa_2, \ldots, \fa_m \in \R^{ m \times d } $, 
$ \beta_1, \beta_2, \ldots, \beta_n \in \R^n $,
$ \fb_1, \fb_2, \ldots, \fb_m \in \R^m $,
$ 
  P_1, P_2, \ldots, P_n, \fP_1, \fP_2, \dots, \fP_m \in \polyn_d 
$
which satisfy for all $ x \in A $ that
\begin{equation} 
\label{eq:f_and_g_representation_piecewise_polynomial}
  f(x) = 
  \smallsum_{ i = 1 }^n 
    P_i(x) \indicator{[0, \infty )^n} \rbr{ \alpha_i x + \beta_i } 
\qqandqq
  g(x)
  =
  \smallsum_{ i = 1 }^m 
  \fP_i(x) \indicator{[0, \infty )^m} \rbr{ \fa_i x + \fb_i } 
  .
\end{equation}
\Nobs that \cref{eq:f_and_g_representation_piecewise_polynomial} 
assures for all $ x \in A $ that 
\begin{equation}
  f(x) + g(x)
  =
  \bigl[ 
    \smallsum_{ i = 1 }^n 
    P_i(x) \indicator{[0, \infty )^n} \rbr{ \alpha_i x + \beta_i } 
  \bigr]
  +
  \bigl[
    \smallsum_{ i = 1 }^m 
    \fP_i(x) \indicator{[0, \infty )^m} \rbr{ \fa_i x + \fb_i } 
  \bigr]
\end{equation}
and 
\begin{equation}
\begin{split}
  f(x) g(x)
& =
  \bigl[ 
    \smallsum_{ i = 1 }^n 
    P_i(x) \indicator{[0, \infty )^n} \rbr{ \alpha_i x + \beta_i } 
  \bigr]
  \bigl[
    \smallsum_{ i = 1 }^m 
    \fP_i(x) \indicator{[0, \infty )^m} \rbr{ \fa_i x + \fb_i } 
  \bigr]
\\ & =
  \smallsum_{ 
    (i,j) 
    \in 
    \{ 1, 2, \dots, n \} \times 
    \{ 1, 2, \dots, m \}
  } 
  \bigl[ 
    P_i(x) 
    \fP_j(x) 
  \bigr]
  \bigl[
    \indicator{ [0, \infty)^n } \rbr{ \alpha_i x + \beta_i } 
    \indicator{ [0, \infty)^m } \rbr{ \fa_j x + \fb_j } 
  \bigr]
  .
\end{split}
\end{equation}
Combining this with \cref{lem:piecewise_polynomial_charac} 
establishes 
\cref{item:piecewise_polynomial_product_i,item:piecewise_polynomial_product_ii}.
\end{cproof}

\subsection{Indicator functions as piecewise polynomial functions}

\cfclear
\begin{lemma}
\label{lem:piecewise_polynomial_indicator}
\cfadd{def:multidim:piece:polyn}
Let $ d \in \N $, 
$ a_1, a_2, \dots, a_d \in \R $, 
$ b_1 \in [a_1,\infty) $, 
$ b_2 \in [a_2,\infty) $,
$ \dots $,
$ b_d \in [a_d,\infty) $
and let 
$
  f \colon \R^d \to \R
$
satisfy for all $ x \in \R^d $ that 
$
  f(x) 
  =
  \mathbbm{1}_{ [a_1, b_1] \times \ldots \times [a_d,b_d] }( x )
%   = \prod_{ i = 1 }^d \mathbbm{1}_{ [a_i,b_i] }( x_i )
$. 
Then  
$ f $ is piecewise polynomial 
\cfload. 
\end{lemma}
\begin{cproof}{lem:piecewise_polynomial_indicator}
Throughout this proof 
let 
$ \alpha_1, \alpha_2, \dots, \alpha_{ 2 d } \in \R^{ ( 2 d ) \times d } $ 
satisfy for all 
$ i \in \{ 1, 2, \dots, 2 d \} $, 
$ x = ( x_1, \dots, x_d ) \in \R^d $
that 
\begin{equation}
\label{eq:piecewise_polynomial_def_alpha}
  \alpha_i x 
  =
  ( x_1, x_2, \dots, x_d, - x_1, - x_2, \dots, - x_d )
  ,
\end{equation}
let $ \beta_1, \beta_2, \dots, \beta_{ 2 d } \in \R^{ 2 d } $ 
satisfy for all $ i \in \{ 1, 2, \dots, 2 d \} $ that 
\begin{equation}
\label{eq:piecewise_polynomial_def_beta}
  \beta_i = ( - \alpha_1, - \alpha_2, \dots, - \alpha_d, \beta_1, \beta_2, \dots, \beta_d ) ,
\end{equation}
and let $ P_i \colon \R^d \to \R $, 
$ i \in \{ 1, 2, \dots, 2 d \} $, 
satisfy for all 
$ i \in \N \cap (1,2d] $, 
$ x \in \R^d $ that 
$
  P_1( x ) = 1
$
and 
$
  P_i( x ) = 0
$.
\Nobs that 
\cref{eq:piecewise_polynomial_def_alpha,eq:piecewise_polynomial_def_beta} 
ensure that 
\begin{equation}
\begin{split}
&
  \times_{ i = 1 }^d [a_i, b_i] 
\\ & 
  =
  \Bigl\{
    x = (x_1, \dots, x_d) \in \R^d 
    \colon
    \Bigl(
      \forall \, i \in \{ 1, 2, \dots, d \} \colon
      a_i \leq x_i \leq b_i
    \Bigr)
  \Bigr\}
\\ &
  =
  \Bigl\{
    x = (x_1, \dots, x_d) \in \R^d 
    \colon
    \Bigl(
      \forall \, i \in \{ 1, 2, \dots, d \} \colon
      \Bigl[
        \bigl( 
          x_i - a_i \in [0,\infty)
        \bigr) ,
        \bigl(
          - x_i + b_i \in [0,\infty)
        \bigr)
      \Bigr]
    \Bigr)
  \Bigr\}
\\ &
  =
  \bigl\{
    x \in \R^d 
    \colon
    \bigl(
      \alpha_1 x + \beta_1 \in [0,\infty)^{ 2 d }
    \bigr)
  \bigr\}
  .
\end{split}
\end{equation}
\Hence for all $ x \in \R^d $ that 
\begin{equation}
\begin{split}
  f(x)
&
  =
  \mathbbm{1}_{
    \{
      y \in \R^d 
      \colon
      \alpha_1 y + \beta_1 \in [0,\infty)^{ 2 d }
    \}
  }( x )
  =
  \mathbbm{1}_{
    [0,\infty)^{ 2 d }
  }( \alpha_1 x + \beta_1 )
  =
  P_1(x)
  \mathbbm{1}_{
    [0,\infty)^{ 2 d }
  }( \alpha_1 x + \beta_1 )
\\ &
=
  \sum_{ i = 1 }^{ 2 d }
  P_i(x)
  \mathbbm{1}_{
    [0,\infty)^{ 2 d }
  }( \alpha_i x + \beta_i )
  .
\end{split}
\end{equation}
Combining this with the fact that 
$
  P_1, P_2, \dots, P_{ 2 d } \in \polyn_d
$
establishes that 
$ f $ is piecewise polynomial. 
\end{cproof}

\subsection{Extensions of piecewise polynomial functions}

\cfclear
\begin{lemma}
\label{lem:piecewise_polynomial_extension}
\cfadd{def:multidim:piece:polyn}
Let $ d \in \N $, 
let $ A \subseteq \R^d $ be a set, 
and let 
$
  f \colon A \to \R
$
be piecewise polynomial 
\cfload. 
Then there exists a piecewise polynomial $ F \colon \R^d \to \R $ 
such that $ F|_A = f $. 
\end{lemma}
\begin{cproof}{lem:piecewise_polynomial_indicator}
\Nobs that 
\cref{eq:def_piecewise_polynomial} 
and the assumption that $ f $ is piecewise polynomial 
ensure that there exist $ n \in \N $, 
$ \alpha_1, \alpha_2, \ldots, \alpha_n \in \R^{ n \times d } $, 
$ \beta_1, \beta_2, \ldots, \beta_n \in \R^n $,
$ P_1, P_2, \ldots, P_n \in \polyn_d $ 
which satisfy for all $ x \in A $ that
\begin{equation} 
\label{eq:piecewise_polynomial_extension1}
  f(x) = 
  \smallsum_{ i = 1 }^n 
%   \br*{ 
    P_i(x) \indicator{[0, \infty )^n} \rbr{ \alpha_i x + \beta_i } 
%   } 
  .
\end{equation}
In the following let $ F \colon \R^d \to \R $ satisfy 
for all $ x \in \R^d $ that
\begin{equation}
\label{eq:piecewise_polynomial_extension2}
  F(x) =
  \smallsum_{ i = 1 }^n 
%   \br*{ 
    P_i(x) \indicator{[0, \infty )^n} \rbr{ \alpha_i x + \beta_i } 
%   } 
  .
\end{equation}
\Nobs that 
\cref{eq:def_piecewise_polynomial,eq:piecewise_polynomial_extension1,eq:piecewise_polynomial_extension2} 
assure that $ F $ is piecewise polynomial. 
\Moreover 
\cref{eq:piecewise_polynomial_extension1,eq:piecewise_polynomial_extension2} 
establish that $ F|_A = f $. 
\end{cproof}

\cfclear
\begin{prop}
\label{prop:piecewise_polynomial_extensionB}
\cfadd{def:multidim:piece:polyn}
Let $ d \in \N $, 
$ a_1, a_2, \dots, a_d \in \R $,
$ b_1 \in [a_1,\infty) $,
$ b_2 \in [a_2,\infty) $,
$ \dots $,
$ b_d \in [a_d,\infty) $, 
$ A = [a_1,b_1] \times \ldots \times [a_d, b_d] $, 
let 
$
  f \colon A \to \R
$ 
be piecewise polynomial, 
and let 
$ F \colon \R^d \to \R $ 
satisfy for all $ x \in A $, $ y \in \R^d \backslash A $ that 
$ 
  F(x) = f(x) 
$
and 
$
  F(y) = 0
$
\cfload. 
Then $ F $ is piecewise polynomial. 
\end{prop}
\begin{cproof}{prop:piecewise_polynomial_extensionB}
Throughout this proof 
let $ \fg \colon \R^d \to \R $ satisfy for all 
$ x \in \R^d $ that 
$ 
  \fg(x) = \mathbbm{1}_{ A }( x ) 
$.
\Nobs that \cref{lem:piecewise_polynomial_indicator} ensures that 
$ \fg $ is piecewise polynomial. 
\Moreover \cref{lem:piecewise_polynomial_extension} assures 
that there exists a piecewise polynomial $ \ff \colon \R^d \to \R $
which satisfies 
\begin{equation}
\label{eq:frak_f_piecewise_polynomial_extension}
  \ff|_A = f
  .
\end{equation}
\Nobs that 
\cref{eq:frak_f_piecewise_polynomial_extension}
shows for all $ x \in \R^d $ that 
\begin{equation}
\label{eq:product_ff_and_fg_is_F}
  \ff(x) \fg(x)
  =
  \ff(x) \mathbbm{1}_A(x)
  =
  F(x)
\end{equation}
\Moreover \cref{lem:piecewise_polynomial_product}, 
the fact that $ \ff $ is piecewise polynomial, 
and the fact that $ \fg $ is piecewise polynomial 
demonstrate that 
$ \R^d \ni x \mapsto \ff(x) \fg(x) \in \R $
is piecewise polynomial. 
Combining this with \cref{eq:product_ff_and_fg_is_F} 
establishes that $ F $ is piecewise polynomial. 
\end{cproof}

\subsection{Piecewise polynomial functions as suitable piecewise rational functions}

We next establish in \cref{prop:piecewise_polynomial} that every piecewise polynomial function in the sense of \cref{def:multidim:piece:polyn} is contained in a class of suitable piecewise rational functions introduced in \cref{def:function:amn}.

\cfclear
\begin{prop}
\label{prop:piecewise_polynomial}
\cfadd{def:multidim:piece:polyn}
Let $ m, n \in \N $ 
and let $ f \colon \R^n \to \R $ 
be piecewise polynomial \cfload. Then 
\begin{equation}
\cfadd{def:function:amn}
\label{eq:piecewise_polynomial_are_Fnminfty}
  \big(
    \R^m \times \R^n \ni ( \theta, x ) \mapsto f(x) \in \R
  \big)
  \in 
  \scrF_{ m, n, \infty }
\end{equation}
\cfload. 
\end{prop}
\begin{cproof}{prop:piecewise_polynomial}
\Nobs that the fact that for all 
$ N \in \N $, $ v = ( v_1, \dots, v_N ) \in \R^N $ 
it holds that 
$
  \indicator{ [0,\infty)^N }( v ) 
  =
  \prod_{ i = 1 }^N
  \indicator{ [0,\infty) }( v_i )
$
assures that 
for all 
$ N \in \N $, 
$
  \alpha = 
  ( 
    \alpha_{ i, j } 
  )_{ 
    (i,j) \in \{ 1, \dots, N \} \times \{ 1, \dots, n \}
  } 
  \in \R^{ N \times n }
$,
$
  \beta = ( \beta_1, \dots, \beta_N ) \in \R^N
$, 
$ 
  P \in \scrP_n
$,
$
  x = ( x_1, \dots, x_n ) \in \R^n
$
it holds that 
\begin{equation}
\textstyle 
  P(x) 
  \indicator{ [0,\infty)^N }( \alpha x + \beta )
  =
  P(x) 
  \bigl[
    \prod_{ 
      i = 1
    }^N
    \mathbbm{1}_{ [0,\infty) }( \beta_i + \sum_{ j = 1 }^n \alpha_{ i, j } x_j )
  \bigr]
  .
\end{equation}
Combining this with 
\cref{eq:def_function_amn} 
demonstrates that 
for all 
$ N \in \N $, 
$
  \alpha 
  \in \R^{ N \times n }
$,
$
  \beta 
  \in \R^N
$, 
$ 
  P \in \scrP_n
$
it holds that 
\begin{equation}
\textstyle 
  \big(
    \R^m \times \R^n \ni ( \theta, x ) 
    \mapsto 
      P(x) 
      \mathbbm{1}_{ [0,\infty)^N }( \alpha x + \beta )
    \in 
    \R
  \big)
  \in \scrF_{ m, n, \infty }
  .
\end{equation}
This and 
\cref{eq:def_function_amn} 
assure that for all 
$ N \in \N $, 
$
  \alpha_1, \alpha_2, \dots, \alpha_N 
  \in \R^{ N \times n }
$,
$
  \beta_1, \beta_2, \dots, \beta_N
  \in \R^N
$, 
$ 
  P_1, P_2, \dots, P_N \in \scrP_n
$
it holds that 
\begin{equation}
\textstyle 
  \big(
    \R^m \times \R^n \ni ( \theta, x ) 
    \mapsto 
      \sum_{ i = 1 }^N
      \br{
        P_i(x) 
        \mathbbm{1}_{ [0,\infty)^N }( \alpha_i x + \beta_i )
      }
    \in 
    \R
  \big)
  \in \scrF_{ m, n, \infty }
  .
\end{equation}
Combining this and 
\cref{eq:def_piecewise_polynomial}
establishes \cref{eq:piecewise_polynomial_are_Fnminfty}. 
\end{cproof}

\subsection{Semi-algebraic risk functions 
in the training of deep ANNs}
\label{ssec:semi_algebraic_risk}

Finally, we combine the previous results to establish the main result of this section.

\cfclear
\begin{cor} 
\label{cor:loss:semialgebraic}
Assume \cref{setting:dnn},
assume for all $ i \in \cu{ 1, 2, \ldots, \ell_L } $ that $ f_i $ is piecewise polynomial, 
let $ \dens \colon [ a, b ]^{ \ell_0 } \to \R $ be piecewise polynomial, 
and assume for all $ E \in \cB( [ a, b ]^{ \ell_0 } ) $ 
that $ \mu ( E ) = \int_E \dens ( x ) \, \d x $ \cfadd{def:multidim:piece:polyn}\cfload.
Then $\cL_\infty$ is semi-algebraic \cfadd{def:semialgebraic:function}\cfload.
\end{cor}
\begin{cproof}{cor:loss:semialgebraic}
Throughout this proof let 
$ 
  F = (F_1, \ldots, F_{ \ell_L } ) \colon \R^{ \ell_0 } \to \R^{ \ell_L } 
$
and 
$
  \fP \colon \R^{ \ell_0 } \to \R 
$ 
satisfy for all $ i \in \cu{ 1, 2, \ldots, \ell_L } $, 
$ x \in \R^{ \ell_0 } $ 
that
\begin{equation} 
\label{cor:loss:semialgebraic:eqtilde}
  F_i ( x ) = 
  \begin{cases} 
    f_i ( x ) 
  & 
    \colon x \in [a , b ]^{ \ell_0 } 
  \\
    0 
  & 
    \colon x \notin [a , b ]^{ \ell_0 }
  \end{cases}
\qqandqq
  \fP( x ) = 
  \begin{cases} 
    \dens ( x ) 
  & 
    \colon x \in [ a, b ]^{ \ell_0 } 
  \\
    0 
  & 
    \colon x \notin [ a, b ]^{ \ell_0 }
    .
  \end{cases}
\end{equation}
\Nobs that 
\cref{prop:piecewise_polynomial_extensionB,cor:loss:semialgebraic:eqtilde}
assure for all $ i \in \{ 1, 2, \dots, \ell_L \} $ 
that 
$ F_i $ and $ \fP $ are piecewise polynomial. 
\cref{prop:piecewise_polynomial} \hence ensures 
for all 
$ i \in \cu{ 1, 2, \ldots, \ell_L } $ 
that
\begin{equation} 
\label{eq:density_and_target_function_in_F}
\cfadd{def:function:amn}
  \bigl\{ 
    \bigl( 
      \R^{ \fd } \times \R^{ \ell_0 } 
      \ni ( \theta, x ) \mapsto F_i( x ) \in \R 
    \bigr)
    ,
    \bigl(
      \R^{ \fd } \times \R^{ \ell_0 } 
      \ni ( \theta, x ) \mapsto \fP( x ) \in \R 
    \bigr)
  \bigr\}
  \subseteq
  \scrF_{ \fd, \ell_0, \infty } 
\end{equation}
\cfload. 
\Moreover \cref{prop:realization:bmn} and 
\cref{item:properties_F_iii} in \cref{lem:bmn:amn} 
demonstrate for all $ i \in \{ 1, 2, \dots, \ell_L \} $ that 
$
  ( 
    \R^{ \fd } \times \R^{ \ell_0 } 
    \ni ( \theta, x ) 
    \mapsto \cN_{ \infty, i }^{ L, \theta }( x )
    \in \R
  )
  \in 
  \scrF_{ \fd, \ell_0, \infty }
$. 
Combining this, 
\cref{eq:def_function_amn}, 
and 
\cref{eq:density_and_target_function_in_F} 
with 
\cref{item:properties_F_iv} in \cref{lem:bmn:amn} 
establish that
\begin{equation}
  \bigl(
    \R^\fd \times \R^{\ell_0} 
    \ni 
    ( \theta , y ) 
%     = 
%     ( \theta, y_1, \dots, y_{ \ell_0 } )
    \mapsto 
%     \norm{\cN_\infty^{ L, \theta }( x ) - F( x ) }^2 \fP( x ) 
%     = 
    \smallsum_{ i = 1 }^{ \ell_L } 
    \bigl[
      (
        \cN_{ \infty, i }^{ L, \theta }( y ) - F_i( y ) 
      )^2
      \fP( y ) 
    \bigr]
    \in \R 
  \bigr) 
  \in 
  \scrF_{ \fd, \ell_0, \infty } .
\end{equation}
\cref{prop:integrals:amn} 
and induction \hence prove that
\begin{multline}
\label{eq:risk_in_F}
\textstyle
  \bigl(
    \R^{ \fd } \times \R^0 \ni 
    (\theta, x) 
    \mapsto 
%   \cL_\infty( \theta ) 
%   = 
    \int_a^b \int_a^b \cdots \int_a^b 
    \sum_{ i = 1 }^{ \ell_L }
    \bigl[
      ( \cN^{ L, \theta }_{ \infty, i }( y_1, \dots, y_{ \ell_0 } ) 
\\
      - F_i( y_1, \dots, y_{ \ell_0 } ) )^2 
      \fP( y_1, \dots, y_{ \ell_0 } ) 
    \bigr]
    \, \d y_{ \ell_0 } \cdots \, \d y_2 \, \d y_1 \in \R 
  \bigr)
  \in \scrF_{ \fd, 0, \infty } .
\end{multline}
\Moreover Fubini's theorem 
and \cref{eq:def_risk_function} show 
for all $ \theta \in \R^{ \fd } $
that 
\begin{equation}
\begin{split}
&
  \cL_{ \infty }( \theta )
% \\ & 
\textstyle 
  =
  \int_{ [a,b]^{ \ell_0 } }
  \| \cN^{ L, \theta }_{ \infty }( y ) - F( y ) \|^2 
  \fP( y ) 
  \, \d y
\\
&
\textstyle
  =
  \int_a^b \int_a^b \cdots \int_a^b 
  \| 
    \cN^{ L, \theta }_{ \infty }( y_1, \dots, y_{ \ell_0 } ) 
    - 
    F( y_1, \dots, y_{ \ell_0 } ) 
  \|^2 \fP( y_1, \dots, y_{ \ell_L } ) 
  \, \d y_d \cdots \, \d y_2 \, \d y_1 
\\ &
\textstyle
  =
  \int_a^b \int_a^b \cdots \int_a^b 
  \sum_{ i = 1 }^{ \ell_L }
  \bigl[
    ( 
      \cN^{ L, \theta }_{ \infty, i }( y_1, \dots, y_{ \ell_0 } ) 
      - 
      F_i( y_1, \dots, y_{ \ell_0 } ) 
    )^2 
    \fP( y_1, \dots, y_{ \ell_0 } )
  \bigr]
  \, \d y_{ \ell_0 } \cdots \, \d y_2 \, \d y_1 
  .
\end{split}
\end{equation}
This and \cref{eq:risk_in_F} disclose that
$
  (
    \R^{ \fd } \times \R^0 \ni 
    (\theta, x) 
    \mapsto 
    \cL_\infty( \theta ) 
    \in \R 
  )
  \in \scrF_{ \fd, 0, \infty } 
$.
\cref{prop:amn:semialgebraic} \hence yields that $ \cL_\infty $ is semi-algebraic.
\end{cproof}

\section{Generalized Kurdyka-\L ojasiewicz inequalities 
for the training of deep ANNs}
\label{sec:KL}

The main result of this section is \cref{prop:loss:lojasiewicz} below,
which reveals that under the assumption that the distribution of the input data has a piecewise polynomial density and that the target function is piecewise polynomial an appropriately generalized Kurdyka-\L ojasiewicz inequality for the risk function is satisfied. 
We prove \cref{prop:loss:lojasiewicz} by combining Bolte et al.~\cite[Theorem 3.1]{BolteDaniilidis2006}) with the fact that the considered risk function $ \cL_\infty$ is semi-algebraic (cf.~\cref{cor:loss:semialgebraic} above). Since \cite[Theorem 3.1]{BolteDaniilidis2006} is formulated for subanalytic functions (cf.~\cref{def:subanalytic_function} below), we state in
\cref{prop:semialgebraic_subanalytic} below the well-known fact that every semi-algebraic function is subanalytic.
We also formulate in \cref{lem:lower_semicontinuous} below the fact that the nonsmooth slope defined in \cref{eq:def_bfM_in_statement} below is lower semi-continuous, which is well-known in the literature (see \cite{BolteDaniilidis2006}). Only for completeness we include in this article a detailed proof of \cref{lem:lower_semicontinuous}.
As a simple consequence of \cref{lem:lower_semicontinuous} we show in \cref{cor:lower_semicontinuous} below that the Kurdyka-\L ojasiewicz inequality always holds around non-critical points, which is also known (cf.~Remark~3.2 in \cite{BolteDaniilidis2006}).

For ANNs with analytic activation functions the risk function was shown to be analytic in Dereich \& Kassing~\cite[Theorem 4.2]{DereichKassing2021} (for an arbitrary compactly supported input distribution). It therefore satisfies an analogous Kurdyka-\L ojasiewicz inequality.

\subsection{Semi-analytic and subanalytic sets}

\begin{definition}[Set of real analytic functions]
\label{def:analytic}	
Let $ n \in \N $ and let $ U \subseteq \R^n $ be an open set. 
Then we denote by $ \analytic_U \subseteq C^{ \infty }( U, \R ) $ 
the set of all real analytic functions from $ U $ to $ \R $.
\end{definition}

For the next notions see, e.g., Bolte et al.~\cite[Definition~2.1]{BolteDaniilidis2006} and Van den Dries \& Miller~\cite{VanDries1996geometric}.

\cfclear
\begin{definition}[Multidimensional semi-analytic sets]
\label{def:semi_analytic_set}
Let $ n \in \N $ and let $ A \subseteq \R^n $ be a set. 
Then we say that $ A $ is an $ n $-dimensional semi-analytic set 
if and only if for all $ v \in \R^n $ there exist 
% $ U \in \{ V \subseteq \R^n \colon V \text{ is open} \} $, 
$M , N \in \N $, 
an open $ U \subseteq \R^n $, 
and 
$
  ( 
    P_{ i, j, k } 
  )_{ 
    ( i, j, k ) \in \cu{ 1, 2, \ldots, M } \times \cu{1, 2, \ldots, N} \times \cu{ 0, 1 } 
  } 
  \subseteq \analytic_U 
$ 
such that 
$ v \in U $ and 
\begin{equation} 
\label{eq:def_semi_analytic_set}
\cfadd{def:analytic}
\textstyle
  A 
  \cap 
  U
  = 
  \bigcup_{ i = 1 }^M
  \bigl( 
    \bigcap_{ j = 1 }^N 
    \cu*{ 
      x \in U \colon P_{ i, j, 0 }( x ) = 0 < P_{ i, j, 1 }( x ) 
    }
  \bigr)
\end{equation}
\cfload.
\end{definition}

\cfclear
\begin{definition}[Multidimensional subanalytic sets]
\label{def:subanalytic_set}
Let $ n \in \N $ and let $ A \subseteq \R^n $ be a set. 
Then we say that $ A $ is an $ n $-dimensional subanalytic set 
if and only if for all $ v \in A $ there exist 
$ m \in \N $, 
an open $ U \subseteq \R^n $, 
and 
a bounded $ ( n + m ) $-dimensional semi-analytic set
$ B \subseteq \R^{ n + m } $
such that 
$ v \in U $ and 
\begin{equation} 
\label{eq:def_sub_analytic_set}
\cfadd{def:semi_analytic_set}
\textstyle
  A 
  \cap 
  U
  = 
  \{ 
    x \in \R^n \colon 
    (
      \exists \, y \in \R^m \colon
      (x,y) \in B
    )
  \} 
\end{equation}
\cfload.
\end{definition}

\subsection{Subanalytic functions}

\cfclear
\begin{definition}[Subanalytic functions]
\label{def:subanalytic_function}
Let $ m, n \in \N $ and let 
$ f \colon \R^m \to \R^n $ be a function. 
Then we say that $ f $ is a subanalytic function 
(we say that $ f $ is subanalytic) 
if and only if it holds that 
$
  \operatorname{Graph}( f ) 
$
is an $ ( m + n ) $-dimensional 
subanalytic set
\cfadd{def:subanalytic_set}\cfload.
\end{definition}

\cfclear 
\begin{prop}
\label{prop:semialgebraic_subanalytic}
\cfadd{def:semialgebraic:function}
Let $ m, n \in \N $ and let 
$ f \colon \R^m \to \R^n $ be semi-algebraic \cfload. 
Then $ f $ is subanalytic
\cfadd{def:subanalytic_function}\cfload. 
\end{prop}
\begin{cproof}{prop:semialgebraic_subanalytic}
\Nobs that the assumption that $ f $ is semi-algebraic 
demonstrates that $ \operatorname{Graph}( f ) $ is an 
$ ( m + n ) $-dimensional semi-algebraic set
\cfadd{def:semialgebraic:set}\cfload. 
Moreover,
it is well-known that every semi-algebraic set is subanalytic (cf., e.g., \cite[Section 2.5]{VanDries1996geometric}).
Hence, we obtain that
$ \operatorname{Graph}( f ) $ is an 
$ ( m + n ) $-dimensional subanalytic set.
\end{cproof}

\subsection{Lower semi-continuity of the nonsmooth slope}

\cfclear
\begin{lemma}
\label{lem:lower_semicontinuous}
Let $ \fd \in \N $, $ f \in C( \R^{ \fd }, \R ) $, 
let 
$
  \bfM \colon \R^{ \fd } \to [0,\infty] 
$
satisfy for all $ \theta \in \R^{ \fd } $ that 
\begin{equation}
\label{eq:def_bfM_in_statement}
\cfadd{def:limit:subdiff}
  \bfM( \theta ) = 
  \inf\bigl(
    \bigl\{
      r \in \R \colon 
      (
        \exists \,
        h \in ( \bbD f )( \theta ) 
        \colon
        r =
        \norm{h} 
      )
    \bigr\} 
    \cup \{ \infty \}
  \bigr) 
  ,
\end{equation}
and let 
$ \theta = ( \theta_n )_{ n \in \N_0 } \colon \N_0 \to \R^{ \fd } $
satisfy 
$
  \limsup_{ n \to \infty } \| \theta_n - \theta_0 \| = 0
$ 
\cfload. 
Then 
$
  \liminf_{ n \to \infty } \bfM( \theta_n ) \geq \bfM( \theta_0 ) 
$.
\end{lemma}
\begin{cproof}{lem:lower_semicontinuous}
Throughout this proof let 
$ \bfm \in [0,\infty] $ satisfy 
$ \bfm = \liminf_{ n \to \infty } \bfM( \theta_n ) $
and 
assume without loss of generality that 
\begin{equation}
\label{eq:bfM_limit_smaller_than_infty}
\textstyle
  \bfm < \infty 
  .
\end{equation}
\Nobs that \cref{eq:bfM_limit_smaller_than_infty} assures 
that there exists a strictly increasing 
$ N \colon \N \to \N $ which satisfies 
\begin{equation}
\textstyle 
\label{eq:construction_of_bfM_subsequence1}
  \limsup_{ n \to \infty } 
  |
    \bfM( \theta_{ N(n) } )
    -
    \bfm 
  |
  = 0 
\qqandqq
  \sup_{ n \in \N }
  \bfM( \theta_{ N(n) } ) < \infty 
  .
\end{equation}
\Nobs that \cref{eq:def_bfM_in_statement,eq:construction_of_bfM_subsequence1} 
prove that there exists 
$ h = ( h_n )_{ n \in \N } \colon \N \to \R^{ \fd } $ 
which satisfies for all $ n \in \N $ that 
\begin{equation}
\label{eq:construction_of_bfM_subsequence2}
  h_n \in ( \bbD f )( \theta_{ N(n) } )
\qqandqq 
  \| h_n \| \leq \bfM( \theta_{ N(n) } ) + n^{ - 1 }
  .
\end{equation}
\Nobs that 
\cref{eq:construction_of_bfM_subsequence1,eq:construction_of_bfM_subsequence2} 
demonstrate that there exist 
$ \bfh \in \R^{ \fd } $ 
and a strictly increasing 
$
  M \colon \N \to \N
$
which satisfy 
\begin{equation}
\label{eq:convergent_subsequence}
\textstyle 
  \limsup_{ n \to \infty }
  \|
    h_{ M(n) }
    - 
    \bfh 
  \|
  = 0 .
\end{equation}
\Nobs that 
\cref{eq:construction_of_bfM_subsequence2}, 
\cref{eq:convergent_subsequence}, 
and the assumption that 
$
  \limsup_{ n \to \infty } 
  \| \theta_n - \theta_0 \| = 0
$ 
demonstrate that 
$
  \limsup_{ k \to \infty } 
  (
    \| 
      h_{ M(k) }
      -
      \bfh 
    \|
    +
    \|
      \theta_{ N( M(k) ) }
      -
      \theta_0
    \|
  )
  = 0
$
and 
$
  \forall \, k \in \N \colon 
  h_{ M(k) }
  \in 
  ( \bbD \cL_{ \infty } )( 
    \theta_{ N( M( k ) ) }
  )
$. 
Combining this and 
\cref{lem:limiting_derivatives} 
(applied with 
$
  n \with \fd 
$, 
$
  f \with f 
$, 
$
  x_0 \with \theta_0 
$, 
$ 
  ( x_k )_{ k \in \N }
  \with 
  ( \theta_{ N( M( k ) ) } )_{ k \in \N }
$, 
$
  y_0 \with \bfh 
$, 
$
  ( y_k )_{ k \in \N }
  \with 
  ( h_{ M(k) } )_{ k \in \N }
$
in the notation of \cref{lem:limiting_derivatives}) 
proves that 
$
  \bfh \in ( \bbD f )( \theta_0 )
$. 
This, 
\cref{eq:def_bfM_in_statement}, 
\cref{eq:construction_of_bfM_subsequence1}, 
\cref{eq:construction_of_bfM_subsequence2}, 
and 
\cref{eq:convergent_subsequence}
show that 
\begin{equation}
\textstyle 
  \bfM( \theta_0 ) \leq \| \bfh \| =
  \limsup_{ n \to \infty } \| h_{ M(n) } \| 
  \leq 
  \limsup_{ n \to \infty } 
  (
    \bfM( \theta_{ N( M(n) ) } )
    +
    [ M(n) ]^{ - 1 }
  )
  = \bfm 
  .
\end{equation}
\end{cproof}

\cfclear
\begin{cor}
\label{cor:lower_semicontinuous}
Let $ \fd \in \N $, $ f \in C( \R^{ \fd }, \R ) $, 
let 
$
  \bfM \colon \R^{ \fd } \to [0,\infty] 
$
satisfy for all $ \theta \in \R^{ \fd } $ that 
\begin{equation}
\label{eq:def_bfM_in_statement2}
\cfadd{def:limit:subdiff}
  \bfM( \theta ) = 
  \inf\bigl(
    \bigl\{
      r \in \R \colon 
      (
        \exists \,
        h \in ( \bbD f )( \theta ) 
        \colon
        r =
        \norm{h} 
      )
    \bigr\} 
    \cup \{ \infty \}
  \bigr) 
  ,
\end{equation}
let 
$ \vartheta \in \R^{ \fd } $ 
satisfy 
$
  0 \notin ( \bbD f )( \vartheta )
$, 
and let $ \fa \in [0,1) $
\cfload. 
Then there exist $ \varepsilon, \const \in (0,\infty) $ such that 
for all 
$ 
  \theta \in 
  \{ \psi \in \R^{ \fd } \colon \| \psi - \vartheta \| < \varepsilon \} 
$
it holds that 
$
  |
    f( \theta ) - f( \vartheta )
  |^{ \fa }
  \leq \const \bfM( \theta )
$.
\end{cor}
\begin{cproof}{cor:lower_semicontinuous}
\Nobs that \cref{item:subdiff_item_v} in \cref{lem:subdifferential:c1}, \cref{eq:def_bfM_in_statement2}, 
and the assumption that $ 0 \notin ( \bbD f )( \vartheta ) $ 
prove that $ \bfM( \vartheta ) > 0 $. 
Combining this with \cref{lem:lower_semicontinuous} proves that there 
exists $ \varepsilon \in (0,\infty) $ such that 
for all 
$
  \theta \in 
  \{ 
    \psi \in \R^{ \fd } \colon 
    \| \psi - \vartheta \| < \varepsilon 
  \} 
$
it holds that 
$
  0 
  < 
  \frac{ \bfM( \vartheta ) }{ 2 } 
  \leq 
  \bfM( \theta )
$. 
\Moreover the assumption that $ f \in C( \R^{ \fd }, \R ) $ 
assures that there exists $ \varepsilon \in (0,\infty) $ 
such that for all 
$
  \theta \in 
    \{ 
      \psi \in \R^{ \fd } \colon 
      \| \psi - \vartheta \| < \varepsilon 
    \} 
$
it holds that 
$
  | 
    f( \theta ) - f( \vartheta )
  |^{ \fa }
  \leq 1
$. 
\end{cproof}

\subsection{Generalized Kurdyka-\L ojasiewicz inequalities 
for the training of deep ANNs}
\label{subsection:loja}

%Using \cref{cor:lower_semicontinuous} we are now able to prove the generalized Kurdyka-\L ojasiewicz inequality for the risk function around every point, not necessarily critical.

\cfclear
\begin{prop}[Generalized {\L}ojasiewicz inequalities] 
\label{prop:loss:lojasiewicz} 
Assume \cref{setting:dnn}, 
assume for all $ i \in \cu{ 1, 2, \ldots, \ell_L } $ that $ f_i $ is piecewise polynomial, 
let $ \dens \colon [ a, b ]^{ \ell_0 } \to \R $ be piecewise polynomial, 
assume for all $ E \in \cB( [ a, b ]^{ \ell_0 } ) $ 
that $ \mu ( E ) = \int_E \dens ( x ) \, \d x $,
and let $ \vartheta \in \R^{ \fd } $ \cfadd{def:multidim:piece:polyn}\cfload. 
Then there exist $ \varepsilon, \const \in (0, \infty) $, 
$ \fa \in [0, 1) $ 
such that for all 
$ 
  \theta \in 
  \{ 
    \psi \in \R^{ \fd } \colon 
    \| \psi - \vartheta \| < \varepsilon 
  \} 
$, 
$ \alpha \in [\fa, 1] $
it holds that
\begin{equation} 
\label{prop:loss:loja:eqclaim}
  \abs{ 
    \cL_{ \infty }( \theta ) - \cL_{ \infty }( \vartheta ) 
  }^{ \alpha } 
  \leq \const 
  \norm{ \cG( \theta ) } 
  .
\end{equation}
\end{prop}
\cfclear
\begin{cproof}{prop:loss:lojasiewicz}
Throughout this proof 
for every $ \varepsilon \in (0,\infty) $ 
let $ B_{ \varepsilon } \subseteq \R^{ \fd } $ 
satisfy 
\begin{equation}
\label{eq:def_B_eps}
  B_{ \varepsilon } =
  \{ 
    \theta \in \R^{ \fd } \colon 
    \| \theta - \vartheta \| < \varepsilon
  \} 
\end{equation}
and let 
$ \bfM \colon \R^{ \fd } \to [0, \infty] $ 
satisfy for all $ \theta \in \R^{ \fd } $ 
that
\begin{equation} 
\label{eq:def_function_M_in_proof}
\cfadd{def:limit:subdiff}
  \bfM( \theta ) 
  = 
  \inf\bigl(
    \bigl\{
      r \in \R \colon 
      (
        \exists \,
        h \in ( \bbD \cL_{ \infty } )( \theta ) 
        \colon
        r =
        \norm{h} 
      )
    \bigr\} 
    \cup \{ \infty \}
  \bigr)
\end{equation}
\cfload.
\Nobs that \cref{prop:loss:gradient:subdiff} implies 
for all $ \theta \in \R^{ \fd } $ that 
$
  \cG( \theta ) \in ( \bbD \cL_{ \infty } )( \theta )
$. 
Combining this with \cref{eq:def_function_M_in_proof} shows 
that for all $ \theta \in \R^{ \fd } $ 
it holds that
\begin{equation}
\label{eq:bfM_estimate_in_proof}
  \bfM( \theta ) \leq \norm{ \cG( \theta ) }
  .
\end{equation}
\Moreover \cref{lem:realization:lip} implies that 
\begin{equation}
\label{eq:risk_is_continuous_in_proof_Bolte}
  \cL_{ \infty } 
  \in 
  C( \R^{ \fd }, \R )
  .
\end{equation}
\Hence for all $ \varepsilon \in (0,\infty) $, $ r \in [0,1] $ that  
\begin{equation}
\label{eq:in_proof_sup_difference_estimate}
\textstyle 
    \sup\nolimits_{ \psi \in B_{ \varepsilon } } 
    \bigl(
      \abs{
        \cL_{ \infty }( \psi ) 
        - 
        \cL_{ \infty }( \vartheta ) 
      }^r 
    \bigr)
  \leq 
  \max\bigl\{ 
    1 ,
    \sup\nolimits_{ \psi \in B_{ \varepsilon } } 
      \abs{
        \cL_{ \infty }( \psi ) 
        - 
        \cL_{ \infty }( \vartheta ) 
      }
  \bigr\} 
  < \infty 
  .
\end{equation}
\Moreover \cref{cor:loss:semialgebraic} 
assures that 
$ \cL_{ \infty } $ is semi-algebraic. 
\cref{prop:semialgebraic_subanalytic} 
\hence proves that 
$ \cL_{ \infty } $ is subanalytic. 
Combining this, 
\cref{eq:def_B_eps}, 
\cref{eq:def_function_M_in_proof}, 
\cref{eq:risk_is_continuous_in_proof_Bolte}, 
\cref{cor:lower_semicontinuous} 
(applied with 
$ f \with \cL_{ \infty } $ 
in the notation of 
\cref{cor:lower_semicontinuous}), 
and 
Bolte et al.~\cite[Theorem 3.1 and (4)]{BolteDaniilidis2006} 
(applied with 
$ n \with \fd $, 
$
  f \with 
  ( 
    \R^{ \fd } \ni \theta \mapsto \cL_{ \infty }( \theta ) 
    \in \R \cup \{ \infty \} 
  )
$
in the notation of Bolte et al.~\cite[Theorem 3.1]{BolteDaniilidis2006})
ensures that there exist $ \varepsilon, \const \in (0, \infty) $, 
$ \fa \in [0, 1) $ 
which satisfy for all 
$ 
  \theta \in B_{ \varepsilon } 
%   \{ \psi \in \R^{ \fd } \colon \| \psi - \vartheta \| < \varepsilon \} 
$ 
that
\begin{equation} 
\label{prop:loss:loja:eq1}
  \abs{
    \cL_{ \infty }( \theta ) - \cL_{ \infty }( \vartheta ) 
  }^{ \fa } 
  \leq \const \bfM( \theta ) 
  .
\end{equation}
\Nobs that \cref{eq:bfM_estimate_in_proof,prop:loss:loja:eq1} 
assure for all $ \theta \in B_{ \varepsilon } $ that 
$
  \abs{
    \cL_{ \infty }( \theta ) - \cL_{ \infty }( \vartheta ) 
  }^{ \fa } 
  \leq \const \| \cG( \theta ) \| 
$.
Combining this with 
\cref{eq:in_proof_sup_difference_estimate} 
demonstrates that for all 
$ \theta \in B_{ \varepsilon } $, $ \alpha \in [ \fa, 1 ] $ 
it holds that
\begin{equation}
\begin{split}
&
  \abs{ 
    \cL_{ \infty }( \theta ) 
    - 
    \cL_{ \infty }( \vartheta ) 
  }^{ \alpha } 
\le 
  \abs{
    \cL_{ \infty }( \theta ) - 
    \cL_{ \infty }( \vartheta ) 
  }^{ \fa } 
  \bigl[ 
    \sup\nolimits_{ \psi \in B_{ \varepsilon } } 
    \bigl(
      \abs{
        \cL_{ \infty }( \psi ) 
        - 
        \cL_{ \infty }( \vartheta ) 
      }^{ \alpha - \fa } 
    \bigr)
  \bigr] 
\\
& 
  \leq
  \abs{
    \cL_{ \infty }( \theta ) - 
    \cL_{ \infty }( \vartheta ) 
  }^{ \fa } 
  \bigl[ 
    \max\bigl\{ 
      1 ,
      \sup\nolimits_{ \psi \in B_{ \varepsilon } } 
        \abs{ 
          \cL_{ \infty }( \psi ) - \cL_{ \infty }( \vartheta ) 
        }
    \bigr\}
  \bigr] 
\\
& 
  \leq 
  \const 
  \bigl[ 
    \max\bigl\{ 
      1 ,
      \sup\nolimits_{ \psi \in B_{ \varepsilon } } 
        \abs{ 
          \cL_{ \infty }( \psi ) - \cL_{ \infty }( \vartheta ) 
        }
    \bigr\}
  \bigr] 
  \| \cG( \theta ) \|  
  < \infty 
  .
\end{split}
\end{equation}
\end{cproof}

\section{Convergence analysis for solutions of GF differential equations}
\label{section:gf:loja}

In \cref{prop:gf:conv:local} below we establish an abstract local convergence result for GF processes under the assumption that a Kurdyka-\L ojasiewicz inequality is satisfied.
The arguments used in the proof of \cref{prop:gf:conv:local} are essentially well-known in the scientific literature; see, e.g.,
 Kurdyka et al.~\cite[Section 1]{KurdykaMostowski2000},
 Bolte et al.~\cite[Theorem 4.5]{BolteDaniilidis2006},
  Absil et al.~\cite[Theorem 2.2]{AbsilMahonyAndrews2005},
  or our previous article Eberle et al.~\cite{EberleJentzenRiekertWeiss2021} (see also \cite{DereichKassing2021SDE} for a version for SDEs). 
  
 The above mentioned works \cite{AbsilMahonyAndrews2005,BolteDaniilidis2006,DereichKassing2021SDE,KurdykaMostowski2000} assume that the objective function is $C^1$ or satisfies some other regularity conditions (in \cite{BolteDaniilidis2006} the objective function is required to be lower-$C^2$ or convex).
 Some works, e.g. \cite{AbsilMahonyAndrews2005}, also assume a certain weak decrease condition for the objective function.
  These assumptions are not necessary for our proof of \cref{prop:gf:conv:local}.
  In fact, we do not even assume that $\cG$ is a subgradient of the objective function $\cL$ at every point.
  The only regularity we need is the chain rule $\cL( \Theta_t ) = \cL( \Theta_0 ) - \int_0^t \norm{ \cG( \Theta_s ) }^2 \, \d s$ for all $t \in [0, \infty ) $. Therefore, our result is not implied by the mentioned previous works.
   %Moreover, we establish precise convergence rates with explicit constants.

  In \cref{cor:gf:conv:local} below we then prove a simplified version of \cref{prop:gf:conv:local}.
  Afterwards, in \cref{prop:gf:global:abstract} and \cref{cor:gf:global:abstract} below we derive
  global convergence of every non-divergent GF trajectory.
  Finally, in \cref{theo:gf:conv:simple} below we combine \cref{cor:gf:global:abstract} with the Kurdyka-\L ojasiewicz inequality for the risk function in \cref{prop:loss:lojasiewicz} and the fact that the generalized gradient is a limiting subdifferential of the risk function (cf.~\cref{prop:loss:gradient:subdiff}) to establish the convergence of GF in the considered deep ANN framework and, thereby, prove \cref{thm:main_thm} from the introduction.

\subsection{Abstract local convergence results for GF processes}
\label{subsection:gf:local:conv}

\begin{prop}
	\label{prop:gf:conv:local}
	Let 
	$ \fd \in \N $, 
	$ \vartheta \in \R^{ \fd } $, 
	$ \consttt \in \R $, 
	$ \const, \varepsilon \in (0, \infty) $, 
	$ \alpha \in (0,1) $, 
	$ \Theta \in C( [0,\infty), \R^{ \fd } ) $, 
	$ \cL \in C( \R^{ \fd }, \R ) $, 
	let $ \cG \colon \R^{ \fd } \to \R^{ \fd } $ 
	be measurable, 
	assume 
	for all $ t \in [0,\infty) $ 
	that
	$
	\cL( \Theta_t ) = \cL( \Theta_0 ) - \int_0^t \| \cG( \Theta_s ) \|^2 \, \d s
	$
	and
	$
  	\Theta_t = \Theta_0 - \int_0^t \cG( \Theta_s ) \, \d s
	$,
	and assume for all 
	$
	\theta \in \R^{ \fd } 
	$
	with 
	$
	\norm{ \theta - \vartheta } < \varepsilon
	$
	that
	\begin{equation} 
	\label{prop:gf:conv:local:eq:loj}
	\abs{ 
		\cL( \theta ) - \cL( \vartheta ) 
	}^\alpha 
	\leq 
	\const \norm{ \cG( \theta ) } ,
	\quad 
	\consttt = 
	|
	\cL( \Theta_0 ) - \cL( \vartheta )
	|,
	\quad 
	\const 
	( 1 - \alpha )^{ - 1 }
	\consttt^{ 1 - \alpha }
	+
	\norm{ \Theta_0 - \vartheta }
	< 
	\varepsilon
	,
	\end{equation}
	and 
	$
	\inf_{ 
		t \in 
		\cu{ 
			s \in  [0, \infty) \colon 
			\forall \, r \in [0,s] \colon 
			\norm{ \Theta_r - \vartheta } < 
			\varepsilon
		} 
	} 
	\cL( \Theta_t ) \geq \cL( \vartheta )
	$. 
	Then there exists 
	$ \psi \in \cL^{ - 1 }( \cu{ \cL( \vartheta ) } ) $ 
	such that 
	\begin{enumerate}[label = (\roman*)]
		\item 
		\label{prop:quantitative:item_i}
		it holds for all $ t \in [0, \infty) $ that 
		$
		\norm{ \Theta_t - \vartheta } < \varepsilon
		$,
		\item 
		\label{prop:quantitative:item_ii}
		it holds for all $ t \in [0,\infty) $ that
		$
		0 \leq 
		\cL( \Theta_t ) - \cL( \psi ) 
		\leq 
		\const^2
		\consttt^2
		(
		\indicator{
			\{ 0 \}
		}( \consttt )
		+
		\const^2 \consttt
		+ 
		\consttt^{
			2 \alpha 
		}
		t 
		)^{ - 1 }
		$,
		and 
		\item 
		\label{prop:quantitative:item_iii}
		it holds for all $ t \in [0,\infty) $ that 
		\begin{equation}
		\label{prop:gf:conv:local:eq:statement}
		\begin{split}
		\norm{
			\Theta_t - \psi 
		}
		& \leq
		\smallint_t^{ \infty }
		\norm{ 
			\cG( \Theta_s )
		}
		\,
		\d s
		\leq
		\const 
		\rbr{ 1 - \alpha }^{ - 1 } 
		[ 
		\cL( \Theta_t ) - \cL( \psi )
		]^{ 1 - \alpha } 
		\\ & \leq 
		\const^{ 3 - 2 \alpha } 
		\consttt^{ 2 - 2 \alpha }
		\rbr{ 1 - \alpha }^{ - 1 } 
		(
		\indicator{
			\{ 0 \}
		}( \consttt )
		+
		\const^2 \consttt
		+ 
		\consttt^{
			2 \alpha 
		}
		t 
		)^{ \alpha - 1 }
		.
		\end{split}
		\end{equation}
	\end{enumerate}
\end{prop}
The assumption that
$\inf_{ t \in \cu{ s \in  [0, \infty) \colon 
		\forall \, r \in [0,s] \colon 
		\norm{ \Theta_r - \vartheta } < 
		\varepsilon	} } \cL( \Theta_t ) \geq \cL( \vartheta )$
means that for all $t \in [0 , \infty )$ which satisfy that the trajectory of $\Theta$ remains within distance $\varepsilon$ of $\vartheta$ until time $t$, it holds that $\cL( \Theta_t ) \geq \cL( \vartheta )$.
This assumption is in particular satisfied if $\vartheta$ is a local minimum of $\cL$ with $\forall \, \theta \in \cu{\psi \in \R^\fd \colon \norm{ \psi - \vartheta } < \varepsilon} \colon \cL ( \theta ) \ge \cL ( \vartheta )$.
But the statement of \cref{prop:gf:conv:local} also covers more general cases, since we only assume this lower bound for the values $\cL ( \Theta_t )$ and not for all values of $\cL$ in a neighborhood of $\vartheta$.

\begin{cproof}{prop:gf:conv:local}
	Throughout this proof let $ \bL \colon [0, \infty) \to \R $ 
	satisfy for all $ t \in [0, \infty) $ that 
	\begin{equation} 
	\label{prop:gf:convergence:eq:defl}
	\bL(t) = \cL ( \Theta_t ) - \cL ( \vartheta ) ,
	\end{equation}
	let $ \bB \subseteq \R^{ \fd } $ satisfy 
	\begin{equation}
	\bB = \{ \theta \in \R^{ \fd } \colon \| \theta - \vartheta \| < \varepsilon \} ,
	\end{equation}
	let $ T \in [0, \infty] $ 
	satisfy
	\begin{equation} 
	\label{prop:gf:convergence:eq:tout}
	T = 
	\inf\rbr*{ 
		\cu*{ t \in [0, \infty ) \colon \Theta_t \notin \bB } 
		\cup \cu{\infty} 
	} ,
	\end{equation}
	let $ \tau \in [0,T] $ satisfy
	\begin{equation}
	\label{prop:gf:convergence:eq:tau}
	\tau = 
	\inf\rbr*{ 
		\cu*{ t \in [0,T) \colon \bL( t ) = 0 } \cup \cu{T } 
	} 
	,
	\end{equation}
	let 
	$ 
	\mathscr{g} = ( \mathscr{g}_t )_{ t \in [0,\infty) }
	\colon [0,\infty) \to [0,\infty] 
	$ 
	satisfy for all $ t \in [0,\infty) $ that 
	$
	\mathscr{g}_t = \int_t^{ \infty } \norm{ \cG( \Theta_s ) } \, \d s
	$, 
	and let 
	$
	\constt \in \R 
	$ satisfy 
	$
	\constt 
	=
	\const^2 
	\consttt^{ ( 2 - 2 \alpha ) } 
	$.
	In the first step of our proof 
	of 
	\cref{prop:quantitative:item_i,prop:quantitative:item_ii,prop:quantitative:item_iii}
	we show that for all $ t \in [0,\infty) $
	it holds that 
	\begin{equation}
	\label{eq:to_prove}
	\Theta_t \in \bB
	.
	\end{equation}
	For this we observe that 
	\cref{prop:gf:conv:local:eq:loj}, 
	the triangle inequality, 
	and 
	the assumption that 
	for all $ t \in [0,\infty) $ it holds that 
	$
	\Theta_t = \Theta_0 - \int_0^t \cG( \Theta_s ) \, \d s
	$
	ensure that 
	for all $ t \in [0,\infty) $ we have that
	\begin{equation} 
	\label{eq:theta_t}
	\begin{split}
	&
	\norm{ \Theta_t - \vartheta } 
	\leq  
	\norm{ \Theta_t - \Theta_0 } 
	+ 
	\norm{ \Theta_0 - \vartheta }  
	\leq
	\norm*{ \int_{0}^t \cG ( \Theta_s )  \, \d s} 
	+
	\norm{ \Theta_0 - \vartheta }  
	\\ & \leq 
	%   \left[ 
	\int_0^t \norm{\cG(\Theta_s ) } \, \d s 
	%   \right]
	+
	\norm{ \Theta_0 - \vartheta }  
	% \\ & 
	< 
	%   \left[ 
	\int_0^t \norm{\cG(\Theta_s ) } \, \d s 
	%   \right]
	-
	\const 
	( 1 - \alpha )^{ - 1 }
	\abs{
		\cL( \Theta_0 ) - \cL( \vartheta ) 
	}^{ 1 - \alpha } 
	+
	\varepsilon
	.
	\end{split}
	\end{equation} 
	To establish \cref{eq:to_prove}, it thus sufficient 
	to prove that 
	$
	\int_0^T \norm{\cG(\Theta_s ) } \, \d s 
	\leq 
	\const 
	( 1 - \alpha )^{ - 1 }
	\abs{
		\cL( \Theta_0 ) - \cL( \vartheta ) 
	}^{ 1 - \alpha } 
	$. 
	We will accomplish this 
	by employing an appropriate differential inequality 
	for a fractional power of the function $ \bL $ in 
	\cref{prop:gf:convergence:eq:defl} 
	(see \cref{prop:gf:conv:eq:derivative} below for details). 
	For this we need several technical preparations. 
	More formally, 
	\nobs that 
	\cref{prop:gf:convergence:eq:defl} 
	and the assumption that for all $ t \in [0,\infty) $ it holds that 
	\begin{equation}
	\label{eq:weak_chain_rule}
	\cL( \Theta_t ) = \cL( \Theta_0 ) - \int_0^t \| \cG( \Theta_s ) \|^2 \, \d s
	\end{equation} 
	imply that for almost all $ t \in [0,\infty) $ it holds 
	that $ \bL $ is differentiable at $ t $ and 
	satisfies 
	\begin{equation}
	\label{eq:derivative_bL}
	\bL'( t ) = \tfrac{ \d }{ \d t }( \cL( \Theta_t ) ) 
	= - \norm{ \cG( \Theta_t ) }^2
	.
	\end{equation}
	Moreover, \nobs that the assumption that 
	$
	\inf_{ 
		t \in 
		\cu{ 
			s \in  [0, \infty) \colon 
			\forall \, r \in [0,s] \colon 
			\norm{ \Theta_r - \vartheta } < 
			\varepsilon
		} 
	} 
	\cL( \Theta_t ) \geq \cL( \vartheta )
	$
	assures for all $ t \in [0, T) $ that 
	\begin{equation}
	\label{eq:Lgeq0}
	\bL(t) \geq 0 
	.
	\end{equation}
	Combining this with 
	\cref{prop:gf:conv:local:eq:loj}, 
	\cref{prop:gf:convergence:eq:defl}, 
	and 
	\cref{prop:gf:convergence:eq:tau}
	demonstrates for all $t \in [0, \tau)$ that
	\begin{equation}
	0 < [ \bL( t ) ]^{ \alpha } 
	= \abs{ \cL( \Theta_t ) - \cL( \vartheta ) }^{ \alpha } 
	\leq \const \norm{\cG ( \Theta_t ) } 
	.
	\end{equation}
	The chain rule and \cref{eq:derivative_bL} 
	hence prove that for almost all $ t \in [0, \tau) $ it holds that
	\begin{equation} 
	\label{prop:gf:conv:eq:derivative}
	\begin{split}
	\tfrac{ \d }{ \d t } 
	( [ \bL( t ) ]^{ 1 - \alpha } ) 
	&
	= ( 1 - \alpha ) [ \bL( t ) ]^{ - \alpha } \rbr{ - \norm{\cG(\Theta_t)}^2  } 
	\\
	&
	\leq 
	- ( 1 - \alpha ) 
	\const^{ - 1 } 
	\norm{ \cG( \Theta_t ) }^{ - 1 } 
	\norm{ \cG( \Theta_t ) }^2 
	= - \const^{ - 1 } ( 1 - \alpha ) 
	\norm{ \cG( \Theta_t ) } .
	\end{split}
	\end{equation}
	Next \nobs that \cref{eq:weak_chain_rule} ensures that 
	$ [0, \infty ) \ni t \mapsto \bL(t) \in \R $ is absolutely continuous. 
	This and the fact that for all $ r \in (0, \infty) $ it holds that 
	$ [ r , \infty ) \ni y \mapsto y^{ 1 - \alpha } \in \R $ is Lipschitz continuous 
	demonstrate that for all $ t \in [0, \tau) $ 
	it holds that $ [0, t] \ni s \mapsto [ \bL( s ) ]^{ 1 - \alpha } \in \R $ 
	is absolutely continuous. 
	Combining this with \cref{prop:gf:conv:eq:derivative} shows 
	that for all $ s, t \in [0, \tau ) $ with $ s \le t $ it holds that
	\begin{equation} 
	\label{eq:integral_inequality}
	\int_s^t 
	\norm{ \cG( \Theta_u ) } 
	\, \d u 
	\leq 
	- \const \rbr{ 1 - \alpha }^{ - 1 } 
	\rbr{  
		[ \bL( t ) ]^{ 1 - \alpha } - [ \bL ( s ) ]^{ 1 - \alpha } 
	} 
	\leq 
	\const 
	\rbr{ 1 - \alpha }^{ - 1 } 
	[ \bL( s ) ]^{ 1 - \alpha } .
	\end{equation}
	In the next step we \nobs that 
	\cref{eq:weak_chain_rule}
	ensures that 
	$
	[0, \infty) \ni t \mapsto \cL ( \Theta_t ) \in \R 
	$ 
	is non-increasing. 
	This and \cref{prop:gf:convergence:eq:defl} 
	prove that $ \bL $ is non-increasing. 
	Combining \cref{prop:gf:convergence:eq:tau} 
	and \cref{eq:Lgeq0} hence implies 
	that for all 
	$ t \in [\tau , T) $ 
	it holds that $ \bL(t) = 0 $. 
	Therefore, we obtain for all 
	$ t \in ( \tau, T ) $ 
	that 
	\begin{equation}
	\label{eq:derivative_zero}
	\bL'( t ) = 0
	.
	\end{equation}
	This and \cref{eq:derivative_bL} assure that 
	for almost all $ t \in (\tau, T) $ it holds that 
	\begin{equation}
	\label{eq:G_vanishes}
	\cG( \Theta_t ) = 0 
	.
	\end{equation}
	Combining this with \cref{eq:integral_inequality} demonstrates that 
	for all $ s, t \in [0, T) $ with $ s \le t $ 
	it holds that 
	\begin{equation} 
	\label{prop:gf_conv:eq:integrated}
	\int_s^t \norm{ \cG( \Theta_u ) } \, \d u 
	\leq 
	\const \rbr{ 1 - \alpha}^{-1} [ \bL( s ) ]^{ 1 - \alpha } 
	.
	\end{equation}
	Hence, we obtain 
	for all $ t \in [0, T) $ 
	that 
	\begin{equation}
	\label{prop:gf_conv:eq:integrated1b}
	\int_0^t \norm{ \cG( \Theta_u ) } \, \d u 
	\leq 
	\const \rbr{ 1 - \alpha}^{-1} [ \bL( 0 ) ]^{ 1 - \alpha } 
	.
	\end{equation}
	In addition, \nobs that \cref{prop:gf:conv:local:eq:loj} 
	assures that $ \Theta_0 \in \bB $. 
	Combining this with \cref{prop:gf:convergence:eq:tout} 
	proves that $ T > 0 $. 
	This, \cref{prop:gf_conv:eq:integrated1b}, 
	and \cref{prop:gf:conv:local:eq:loj} 
	demonstrate that 
	\begin{equation} 
	\label{prop:gf_conv:eq:integrated2}
	\int_0^T \norm{ \cG( \Theta_u ) } \, \d u 
	\leq 
	\const \rbr{ 1 - \alpha}^{-1} [ \bL( 0 ) ]^{ 1 - \alpha } 
	< 
	\varepsilon 
	< \infty 
	.
	\end{equation}
	Combining \cref{prop:gf:convergence:eq:tout} and \cref{eq:theta_t} 
	hence assures that 
	\begin{equation}
	\label{eq:T_infinity} 
	T = \infty 
	.
	\end{equation}
	This establishes \cref{eq:to_prove}. 
	In the next step of our proof of 
	\cref{prop:quantitative:item_i,prop:quantitative:item_ii,prop:quantitative:item_iii}
	we verify that 
	$ \Theta_t \in \R^{ \fd } $, $ t \in [0,\infty) $, 
	is convergent (see \cref{eq:convergent} below). 
	For this observe that 
	the assumption that 
	for all $ t \in [0,\infty) $ it holds that 
	$
	\Theta_t = \Theta_0 - \int_0^t \cG( \Theta_s ) \, \d s
	$
	demonstrates that for all 
	$ r, s, t \in [0,\infty) $ 
	with $ r \leq s \leq t $ 
	it holds that 
	\begin{equation}
	\label{eq:Cauchy_family}
	\norm{ \Theta_t - \Theta_s } 
	=
	\norm*{ 
		\int_s^t
		\cG( \Theta_u ) \, \d u
	}
	\leq 
	\int_s^t 
	\norm{ 
		\cG( \Theta_u )
	}
	\, \d u
	\leq 
	\int_r^{ \infty }
	\norm{
		\cG( \Theta_u )
	}
	\, \d u
	= \mathscr{g}_r 
	.
	\end{equation}
	Moreover, note that 
	\cref{prop:gf_conv:eq:integrated2}
	and 
	\cref{eq:T_infinity} 
	assure that 
	$
	\infty >
	\mathscr{g}_0
	\geq 
	\limsup_{ r \to \infty } 
	\mathscr{g}_r = 0
	$. 
	Combining this with 
	\cref{eq:Cauchy_family} proves that there exist 
	$ \psi \in \R^{ \fd } $ 
	which satisfies 
	\begin{equation}
	\label{eq:convergent}
	\limsup\nolimits_{ t \to \infty } \norm{ \Theta_t - \psi } = 0 .
	\end{equation}
	In the next step of our proof of 
	\cref{prop:quantitative:item_i,prop:quantitative:item_ii,prop:quantitative:item_iii}
	we show that $ \cL( \Theta_t ) $, $ t \in [0,\infty) $, 
	converges to $ \cL( \psi ) $ with convergence order $ 1 $. 
	We accomplish this bringing a suitable differential inequality 
	for the reciprocal of the function $ \bL $ 
	in \cref{prop:gf:convergence:eq:defl} into play 
	(see \cref{eq:d_dt_reciprocal_L} below for details). More specifically, 
	\nobs that 
	\cref{eq:derivative_bL}, 
	\cref{eq:T_infinity}, 
	\cref{prop:gf:convergence:eq:tout}, 
	and \cref{prop:gf:conv:local:eq:loj} 
	demonstrate that for almost all $ t \in [0, \infty) $ 
	it holds that
	\begin{equation}
	\bL'( t ) 
	= - \norm{ \cG( \Theta_t ) }^2 
	\leq - \const^{ - 2 } [ \bL( t ) ]^{ 2 \alpha } .
	\end{equation}
	Hence, we obtain that $ \bL $ is non-increasing. 
	This shows for all $ t \in [0,\infty) $ 
	that 
	$ \bL(t) \leq \bL( 0 ) $. 
	This and the fact that for all $ t \in [0, \tau) $ 
	it holds that $ \bL( t ) > 0 $ 
	show that for almost all $ t \in [0, \tau) $ 
	we have that
	\begin{equation}
	\bL'( t )
	\leq 
	- \const^{ - 2 } 
	[ \bL( t ) ]^{ ( 2 \alpha - 2 ) } 
	[ \bL( t ) ]^{ 2 } 
	\leq 
	- \const^{ - 2 } 
	[ \bL( 0 ) ]^{ ( 2 \alpha - 2 ) } 
	[ \bL( t ) ]^{ 2 } 
	=
	- \constt^{ - 1 }
	[ \bL( t ) ]^{ 2 } 
	.
	\end{equation}
	Therefore, we obtain that for almost all $ t \in [0, \tau) $ 
	it holds that
	\begin{equation}
	\label{eq:d_dt_reciprocal_L}
	\frac{ \d }{ \d t } 
	\rbr*{ 
		\frac{ 
			\constt  
		}{ 
			\bL( t ) 
		} 
	} 
	= 
	- 
	\rbr*{
		\frac{
			\constt \, 
			\bL'( t )
		}{
			[ \bL( t ) ]^2 
		} 
	}
	\geq 1 .
	\end{equation}
	Moreover, \nobs that the fact that 
	for all $ t \in [0, \tau) $ 
	it holds that 
	$ [0, t] \ni s \mapsto \bL( s ) \in (0, \infty) $ 
	is absolutely continuous 
	proves that for all $ t \in [0,\tau) $ 
	we have that 
	$ [0, t] \ni s \mapsto \constt [ \bL( s ) ]^{ - 1 } \in (0, \infty) $ 
	is absolutely continuous. 
	This and \cref{eq:d_dt_reciprocal_L} imply
	for all $ t \in [0,\tau) $ that 
	$
	\frac{ 
		\constt
	}{
		\bL(t)
	}
	-
	\frac{ 
		\constt
	}{
		\bL(0)
	}
	\geq t
	$.
	Hence, we obtain for all $ t \in [0,\tau) $ 
	that 
	$
	\frac{ \constt }{ \bL(t) } 
	\geq 
	\frac{ \constt }{ \bL(0) } 
	+ t 
	$. 
	Therefore, we get for all 
	$ t \in [0,\tau) $ that 
	$
	\constt 
	\,
	(
	\frac{ \constt }{ \bL(0) } 
	+ t 
	)^{ - 1 }
	\geq 
	\bL(t)
	$. 
	This implies for all $ t \in [0,\tau) $ 
	that 
	\begin{equation}
	\begin{split}
	\bL(t)
	& \leq
	\constt 
	\,
	(
	\constt 
	[ \bL(0) ]^{ - 1 }
	+ 
	t 
	)^{ - 1 }
	=
	\const^2
	\consttt^{
		2 - 2 \alpha 
	}
	(
	\const^2
	\consttt^{
		1 - 2 \alpha 
	}
	+ 
	t 
	)^{ - 1 }
	=
	\const^2
	\consttt^2
	(
	\const^2 \consttt
	+ 
	\consttt^{
		2 \alpha 
	}
	t 
	)^{ - 1 }
	.
	\end{split}
	\end{equation}
	The fact that for all 
	$ t \in [\tau, \infty) $ 
	it holds that $ \bL(t) = 0 $ 
	and \cref{prop:gf:convergence:eq:tau} 
	therefore prove 
	that for all $ t \in [0, \infty) $ 
	it holds that 
	\begin{equation}
	\label{prop:gf:convergence:eq:lossest}
	0 
	\leq
	\bL( t ) 
	\leq 
	\const^2
	\consttt^2
	(
	\indicator{
		\{ 0 \}
	}( \consttt )
	+
	\const^2 \consttt
	+ 
	\consttt^{
		2 \alpha 
	}
	t 
	)^{ - 1 }
	.
	\end{equation}
	Next \nobs that \cref{eq:convergent} and 
	the assumption that $ \cL \in C( \R^{ \fd }, \R ) $ 
	assure that 
	$
	\limsup_{ t \to \infty } 
	\abs{
		\cL( \Theta_t ) - \cL( \psi ) 
	} = 0
	$. 
	Combining this with 
	\cref{prop:gf:convergence:eq:lossest} 
	demonstrates that 
	$ \cL( \psi ) = \cL( \vartheta ) $. 
	This and \cref{prop:gf:convergence:eq:lossest} 
	ensure for all $ t \in [0,\infty) $ 
	that
	\begin{equation}
	\label{prop:gf:convergence:eq:claim2}
	0 
	\leq
	\cL( \Theta_t ) - \cL( \psi ) 
	\leq 
	\const^2
	\consttt^2
	(
	\indicator{
		\{ 0 \}
	}( \consttt )
	+
	\const^2 \consttt
	+ 
	\consttt^{
		2 \alpha 
	}
	t 
	)^{ - 1 }
	.
	\end{equation}
	In the final step of our proof 
	of 
	\cref{prop:quantitative:item_i,prop:quantitative:item_ii,prop:quantitative:item_iii}
	we establish convergence rates for the 
	real numbers  
	$ \norm{ \Theta_t - \psi } $, 
	$ t \in [0, \infty) $. 
	\Nobs that \cref{eq:convergent}, 
	\cref{eq:Cauchy_family}, 
	and \cref{prop:gf_conv:eq:integrated} 
	assure
	for all $ t \in [0,\infty) $ 
	that 
	\begin{equation}
	\norm{
		\Theta_t - \psi 
	}
	=
	\norm*{
		\Theta_t 
		- 
		\left[ 
		\lim\nolimits_{ s \to \infty } 
		\Theta_s 
		\right]
	}
	=
	\lim\nolimits_{ s \to \infty }
	\norm{
		\Theta_t - \Theta_s 
	}
	\leq 
	\mathscr{g}_t
	\leq 
	\const 
	\rbr{ 1 - \alpha }^{ - 1 } [ \bL( t ) ]^{ 1 - \alpha } 
	.
	\end{equation}
	This and \cref{prop:gf:convergence:eq:claim2} 
	ensure for all $ t \in [0,\infty) $ 
	that 
	\begin{equation}
	\begin{split}
	\norm{
		\Theta_t - \psi 
	}
	& \leq 
	\mathscr{g}_t
	\leq
	\const 
	\rbr{ 1 - \alpha }^{ - 1 } 
	[ 
	\cL( \Theta_t ) - \cL( \psi )
	]^{ 1 - \alpha } 
	\\ & \leq 
	\const 
	\rbr{ 1 - \alpha }^{ - 1 } 
	\left[ 
	\const^2
	\consttt^2
	(
	\indicator{
		\{ 0 \}
	}( \consttt )
	+
	\const^2 \consttt
	+ 
	\consttt^{
		2 \alpha 
	}
	t 
	)^{ - 1 }
	\right]^{ 1 - \alpha } 
	\\ & =
	\const^{ 3 - 2 \alpha } 
	\consttt^{ 2 - 2 \alpha }
	\rbr{ 1 - \alpha }^{ - 1 } 
	(
	\indicator{
		\{ 0 \}
	}( \consttt )
	+
	\const^2 \consttt
	+ 
	\consttt^{
		2 \alpha 
	}
	t 
	)^{ \alpha - 1 }
	.
	\end{split}
	\end{equation}
	Combining this with \cref{eq:to_prove,prop:gf:convergence:eq:claim2} 
	establishes 
	\cref{prop:quantitative:item_i,prop:quantitative:item_ii,prop:quantitative:item_iii}. 
\end{cproof}

\begin{cor}
\label{cor:gf:conv:local}
Let 
$ \fd \in \N $, 
$ \vartheta \in \R^{ \fd } $, 
$ \consttt \in [0,1] $, 
$ \const, \varepsilon \in (0, \infty) $, 
$ \alpha \in (0,1) $, 
$ \Theta \in C( [0,\infty), \R^{ \fd } ) $, 
$ \cL \in C( \R^{ \fd }, \R ) $, 
let $ \cG \colon \R^{ \fd } \to \R^{ \fd } $ 
be measurable, 
assume 
for all $ t \in [0,\infty) $ 
that
$
  \cL( \Theta_t ) = \cL( \Theta_0 ) - \int_0^t \| \cG( \Theta_s ) \|^2 \, \d s
$
and 
$
  \Theta_t = \Theta_0 - \int_0^t \cG( \Theta_s ) \, \d s
$, 
and assume for all 
$
  \theta \in \R^{ \fd } 
$
with 
$
  \norm{ \theta - \vartheta } < \varepsilon
$
that
\begin{equation} \label{cor:gf:conv:local:eq:ass}
  \abs{ 
    \cL( \theta ) - \cL( \vartheta ) 
  }^\alpha 
  \leq 
  \const \norm{ \cG( \theta ) } ,
\quad 
  \consttt = 
  |
    \cL( \Theta_0 ) - \cL( \vartheta )
  |,
\quad 
  \const 
  ( 1 - \alpha )^{ - 1 }
  \consttt^{ 1 - \alpha }
  +
  \norm{ \Theta_0 - \vartheta }
  < 
  \varepsilon
  ,
\end{equation}
and 
$
  \inf_{ 
    t \in 
    \cu{ 
      s \in  [0, \infty) \colon 
      \forall \, r \in [0,s] \colon 
      \norm{ \Theta_r - \vartheta } < 
      \varepsilon
    } 
  } 
  \cL( \Theta_t ) \geq \cL( \vartheta )
$. 
	Then there 
	exists 
	$ \psi \in \cL^{ - 1 }( \cu{ \cL( \vartheta ) } ) $ 
	such that for all $ t \in [0,\infty) $ 
	it holds that 
	$
	\norm{ \Theta_t - \vartheta } < \varepsilon
	$,
	$
	0 \leq 
	\cL( \Theta_t ) - \cL( \psi ) 
	\leq 
	(
	1
	+ 
	\const^{
		- 2
	}
	t 
	)^{ - 1 }
	$,
	and 
	\begin{equation}
	\label{cor:eq:statement}
	\norm{
		\Theta_t - \psi 
	}
	\leq
	\smallint\nolimits_t^{ \infty }
	\norm{ 
		\cG( \Theta_s )
	}
	\,
	\d s
	\leq 
	\const
	\rbr{ 1 - \alpha }^{ - 1 } 
	(
	1
	+ 
	\const^{
		- 2
	}
	t 
	)^{ \alpha - 1 }
	.
	\end{equation}
\end{cor}

\begin{cproof}{cor:gf:conv:local}
	Observe that \cref{prop:gf:conv:local} ensures that 
	exists 
	$ \psi \in \cL^{ - 1 }( \cu{ \cL( \vartheta ) } ) $ 
	which satisfies that
	\begin{enumerate}[label = (\roman*)]
		\item 
		\label{cor:quantitative:proof:item_i}
		it holds for all $ t \in [0, \infty) $ that 
		$
		\norm{ \Theta_t - \vartheta } < \varepsilon
		$,
		\item 
		\label{cor:quantitative:proof:item_ii}
		it holds for all $ t \in [0,\infty) $ that
		$
		0 \leq 
		\cL( \Theta_t ) - \cL( \psi ) 
		\leq 
		\const^2
		\consttt^2
		(
		\indicator{
			\{ 0 \}
		}( \consttt )
		+
		\const^2 \consttt
		+ 
		\consttt^{
			2 \alpha 
		}
		t 
		)^{ - 1 }
		$,
		and 
		\item 
		\label{cor:quantitative:proof:item_iii}
		it holds for all $ t \in [0,\infty) $ that 
		\begin{equation}
		\begin{split}
		\norm{
			\Theta_t - \psi 
		}
		& \leq
		\smallint_t^{ \infty }
		\norm{ 
			\cG( \Theta_s )
		}
		\,
		\d s
		\leq
		\const 
		\rbr{ 1 - \alpha }^{ - 1 } 
		[ 
		\cL( \Theta_t ) - \cL( \psi )
		]^{ 1 - \alpha } 
		\\ & \leq 
		\const^{ 3 - 2 \alpha } 
		\consttt^{ 2 - 2 \alpha }
		\rbr{ 1 - \alpha }^{ - 1 } 
		(
		\indicator{
			\{ 0 \}
		}( \consttt )
		+
		\const^2 \consttt
		+ 
		\consttt^{
			2 \alpha 
		}
		t 
		)^{ \alpha - 1 }
		.
		\end{split}
		\end{equation}
	\end{enumerate}
	\Nobs that \cref{cor:quantitative:proof:item_ii} 
	and the assumption that $ \consttt \leq 1 $ imply 
	that 
	for all $ t \in [0,\infty) $ it holds that
	\begin{equation}
	\label{eq:weak_estimate:proof:cor}
	0 
	\leq 
	\cL( \Theta_t ) - \cL( \psi ) 
	\leq 
	\consttt^2
	(
	\const^{ - 2 }
	\indicator{
		\{ 0 \}
	}( \consttt )
	+
	\consttt
	+ 
	\const^{ - 2 }
	\consttt^{
		2 \alpha 
	}
	t 
	)^{ - 1 }
	\leq 
	(
	1
	+ 
	\const^{ - 2 }
	t 
	)^{ - 1 }
	.
	\end{equation} 
	This and \cref{cor:quantitative:proof:item_iii} 
	ensure that 
	for all $ t \in [0,\infty) $ 
	it holds that 
	\begin{equation}
	\begin{split}
	\norm{
		\Theta_t - \psi 
	}
	& \leq
	\smallint_t^{ \infty }
	\norm{ 
		\cG( \Theta_s )
	}
	\,
	\d s
	\leq
	\const 
	\rbr{ 1 - \alpha }^{ - 1 } 
	[ 
	\cL( \Theta_t ) - \cL( \psi )
	]^{ 1 - \alpha } 
	\\ & \leq 
	\const 
	\rbr{ 1 - \alpha }^{ - 1 } 
	(
	1
	+ 
	\const^{ - 2 }
	t 
	)^{ \alpha - 1 }
	. 
	\end{split}
	\end{equation} 
	Combining this with 
	\cref{cor:quantitative:proof:item_i,eq:weak_estimate:proof:cor}
	establishes \cref{cor:eq:statement}. 
\end{cproof}

\subsection{Abstract global convergence results for GF processes}
\label{subsection:gf:global:conv}

We next employ \cref{cor:gf:conv:local} to establish under a Kurdyka-\L ojasiewicz assumption the convergence of every non-divergent GF trajectory. To prove \cref{prop:gf:global:abstract} we note that the trajectory must have a convergent subsequence with limit $\vartheta \in \R^\fd$. Hence, for a sufficiently large time the GF reaches a neighborhood of $\vartheta$ where the conditions of \cref{cor:gf:conv:local} in \cref{cor:gf:conv:local:eq:ass} are satisfied, and thus we get convergence of the entire trajectory.

\cfclear 
\begin{prop} 
\label{prop:gf:global:abstract}
Let 
$ \fd \in \N $, $ \Theta \in C( [0,\infty), \R^{ \fd } ) $, 
$ \cL \in C( \R^{ \fd }, \R ) $, 
let $ \cG \colon \R^{ \fd } \to \R^{ \fd } $ 
be measurable, 
% and locally bounded, 
assume that for all $ \vartheta \in \R^{ \fd } $ 
there exist 
$ \varepsilon, \const \in (0, \infty) $, 
$ \alpha \in (0, 1) $ such that 
for all $ \theta \in \R^{ \fd } $ 
with $ \norm{ \theta - \vartheta } < \varepsilon $ 
it holds that
$
  \abs{ \cL( \theta ) - \cL( \vartheta ) }^{ \alpha } 
  \leq \const \norm{ \cG( \theta ) } 
$,  
and assume for all $ t \in [0,\infty) $ 
that
\begin{equation} 
\label{prop:gf:global:eq:ass}
\textstyle
  \liminf_{ s \to \infty } \norm{ \Theta_s } < \infty,
\quad
  \cL( \Theta_t ) = 
  \cL( \Theta_0 ) - \int_0^t \| \cG( \Theta_s ) \|^2 \, \d s ,
\qandq
  \Theta_t = \Theta_0 - \int_0^t \cG( \Theta_s ) \, \d s
  .
\end{equation}	
Then there exist 
$ \vartheta \in \R^{ \fd } $, $ \fC, \tau, \beta \in (0, \infty) $ 
such that for all $ t \in [ \tau, \infty) $ 
it holds that
\begin{equation} 
\label{prop:gf:global:eq:claim}
  \norm{
    \Theta_t - \vartheta
  } 
  \leq 
%   \fC 
  \bigl( 
    1 + \fC^{ - 1 } ( t - \tau ) 
  \bigr)^{ - \beta } 
\qqandqq
  0
  \leq
  \cL( \Theta_t ) - \cL( \vartheta ) 
  \leq 
  \bigl( 
    1 + \fC^{ - 1 } ( t - \tau ) 
  \bigr)^{ - 1 } .
\end{equation}
\end{prop}
\begin{cproof}{prop:gf:global:abstract}
\Nobs that \cref{prop:gf:global:eq:ass} implies that 
$
  [0, \infty) \ni t \mapsto \cL( \Theta_t ) \in \R 
$ 
is non-increasing. 
\Hence that there exists $ \bfm \in [-\infty, \infty) $ 
which satisfies
\begin{equation} 
\label{proof:gf:global:eq:defm}
  \bfm = \limsup\nolimits_{ t \to \infty } \cL( \Theta_t ) 
  = 
  \liminf\nolimits_{ t \to \infty } 
  \cL( \Theta_t ) 
  = 
  \inf\nolimits_{ t \in [0, \infty) } 
  \cL( \Theta_t )
  .
\end{equation}
\Moreover that the assumption that 
$ 
  \liminf_{ t \to  \infty } \norm{ \Theta_t } < \infty 
$ 
ensures that there exist 
$
  \vartheta \in \R^{ \fd } 
$ 
and 
$
  \delta = ( \delta_n )_{ n \in \N } \colon \N \to [0, \infty) 
$ 
which satisfy 
\begin{equation} 
\label{proof:gf:global:eq:taun}
\textstyle
  \liminf_{ n \to \infty } \delta_n = \infty 
\qqandqq
  \limsup\nolimits_{ n \to \infty } \norm{ \Theta_{ \delta_n } - \vartheta } = 0 
  .
\end{equation} 
\Nobs that 
\cref{proof:gf:global:eq:defm}, 
\cref{proof:gf:global:eq:taun}, 
and the fact that $ \cL $ is continuous show that 
\begin{equation} 
\label{proof:gf:global:eq:bfm}
  \cL( \vartheta ) 
  = \bfm \in \R
\qqandqq 
  \forall \, t \in [0, \infty) \colon 
  \cL( \Theta_t ) \geq \cL( \vartheta ) 
  .
\end{equation}
Next let
$
  \varepsilon, \const \in (0, \infty) 
$, 
$
  \alpha \in (0, 1)
$
satisfy for all 
$
  \theta \in \R^{ \fd } 
$ 
with 
$
  \norm{ \theta - \vartheta } < \varepsilon 
$ 
that
\begin{equation}
\label{eq:loja_in_proof}
  \abs{ \cL( \theta ) - \cL( \vartheta ) }^{ \alpha } 
  \leq \const \norm{ \cG( \theta ) } 
  .
\end{equation}
\Nobs that \cref{proof:gf:global:eq:taun} and 
the fact that $ \cL $ is continuous demonstrate 
that there exist $ n \in \N $, $ \consttt \in [0,1] $ 
which satisfy
\begin{equation} 
\label{proof:gf:global:eq:consttt}
  \consttt = 
  \abs{
    \cL( \Theta_{ \delta_n } ) - \cL( \vartheta ) 
  } 
\qqandqq 
  \const 
  ( 1 - \alpha )^{ - 1 }
  \consttt^{ 1 - \alpha }
  +
  \norm{ \Theta_{ \delta_n } - \vartheta }
  < 
  \varepsilon
  .
\end{equation}
Next let $ \Phi \colon [0, \infty ) \to \R^{ \fd } $ 
satisfy for all $ t \in [0, \infty) $ that
\begin{equation} 
\label{proof:gf:global:eq:defphi}
  \Phi_t = \Theta_{ \delta_n + t }
  .
\end{equation}
\Nobs that \cref{prop:gf:global:eq:ass,proof:gf:global:eq:bfm,proof:gf:global:eq:defphi}
assure that for all $ t \in [0, \infty) $
it holds that
\begin{equation}
\textstyle
  \cL( \Phi_t ) = \cL( \Phi_0 ) - \int_0^t \norm{ \cG( \Phi_s ) }^2 \, \d s ,
\qquad
  \Phi_t = \Phi_0 - \int_0^t \cG ( \Phi_s ) \, \d s ,
\qqandqq
  \cL( \Phi_t ) 
  \ge \cL( \vartheta )
  .
\end{equation}
Combining this with 
\cref{eq:loja_in_proof}, 
\cref{proof:gf:global:eq:consttt}, 
\cref{proof:gf:global:eq:defphi}, 
and \cref{cor:gf:conv:local} 
(applied with 
$
  \Theta \with \Phi
$
in the notation of \cref{cor:gf:conv:local}) establishes 
that there exists 
$
  \psi \in \cL^{ - 1 }( \cu{ \cL( \vartheta ) } ) 
$
which satisfies 
for all $ t \in [0, \infty) $ that
\begin{equation} 
\label{proof:gf:global:eq:phiest}
  0 \leq 
  \cL( \Phi_t ) - \cL( \psi ) 
  \leq 
  ( 1 + \const^{- 2 } t  )^{ - 1 }
\qqandqq 
  \norm{ \Phi_t - \psi }
  \leq
  \const \rbr{ 1 - \alpha }^{ - 1 } ( 1 + \const^{- 2 } t )^{ \alpha - 1 } .
\end{equation} 
\Nobs that \cref{proof:gf:global:eq:defphi,proof:gf:global:eq:phiest} 
assure for all $ t \in [ 0, \infty) $ that 
$
  0 \leq 
  \cL( \Theta_{ \delta_n + t } ) - \cL( \psi ) 
  \leq 
  ( 1 + \const^{- 2 } t  )^{ - 1 }
$
and 
$
  \norm{ \Theta_{ \delta_n + t } - \psi }
  \leq
  \const \rbr{ 1 - \alpha }^{ - 1 } ( 1 + \const^{- 2 } t )^{ \alpha - 1 } 
$.
\Hence for all 
$ \tau \in [ \delta_n, \infty ) $, 
$ t \in [ \tau, \infty ) $ that 
\begin{equation}
\begin{split}
  0 
& \leq 
  \cL( \Theta_t ) - \cL( \psi ) 
  \leq 
  ( 1 + \const^{ - 2 } ( t - \delta_n ) )^{ - 1 }
  =
  ( 1 + \const^{ - 2 } ( t - \tau ) + \const^{ - 2 } ( \tau - \delta_n ) )^{ - 1 }
\\ &
\leq 
  ( 1 + \const^{ - 2 } ( t - \tau ) )^{ - 1 }
\end{split}
\end{equation}
and 
\begin{equation}
\label{eq:norm_Theta_t_psi_estimate}
\begin{split}
&
  \norm{ \Theta_t - \psi }
\leq
  \const \rbr{ 1 - \alpha }^{ - 1 } ( 1 + \const^{- 2 } ( t - \delta_n ) )^{ \alpha - 1 } 
= 
  \Bigl[
    \bigl[
      \const \rbr{ 1 - \alpha }^{ - 1 } 
    \bigr]^{ \frac{ 1 }{ \alpha - 1 } } 
    ( 1 + \const^{- 2 } ( t - \delta_n ) )
  \Bigr]^{ \alpha - 1 } 
\\ & =
  \biggl[
    \bigl[
      \const \rbr{ 1 - \alpha }^{ - 1 } 
    \bigr]^{ \frac{ 1 }{ \alpha - 1 } } 
    \bigl[
      1 
      +
      \const^{ - 2 } ( \tau - \delta_n ) 
    \bigr]
    + 
    \Bigl[
      \bigl[
        \const ( 1 - \alpha )^{ - 1 } 
      \bigr]^{ \frac{ 1 }{ 1 - \alpha } } 
      \const^2
    \Bigr]^{ - 1 }
    ( t - \tau ) 
  \biggr]^{ \alpha - 1 } 
  .
\end{split}
\end{equation}
Next let 
$ \scrC, \tau \in (0, \infty) $ 
satisfy 
\begin{equation}
\label{eq:def_tau_in_proof}
  \scrC 
  = 
  \max\bigl\{ 
    \fC^2 , 
    \bigl[
      \const ( 1 - \alpha )^{ - 1 } 
    \bigr]^{ \frac{ 1 }{ 1 - \alpha } } 
    \const^2
  \bigr\} 
\qqandqq
  \tau 
  =
  \delta_n 
  + 
  \const^2
  \bigl[
    \const \rbr{ 1 - \alpha }^{ - 1 } 
  \bigr]^{ \frac{ 1 }{ 1 - \alpha } } 
  .
\end{equation}
\Nobs that 
\cref{eq:norm_Theta_t_psi_estimate,eq:def_tau_in_proof} 
demonstrate for all $ t \in [ \tau, \infty ) $ that
\begin{equation}
\begin{split}
  0 
  \leq
  \cL( \Theta_t ) - \cL( \psi ) 
& \leq 
  ( 1 + \const^{ - 2 } ( t - \tau ) )^{ - 1 }
  \leq   
  ( 1 + \scrC^{ - 1 } ( t - \tau ) )^{ - 1 }
% 
%   =
%   ( 1 + \scrC^{ - 1 } ( \tau - \delta_n ) + \scrC^{ - 1 } ( t - \tau ) )^{ - 1 }
% \\ & \leq 
%   ( 1 + \scrC^{ - 1 } ( t - \tau ) )^{ - 1 }
\end{split}
\end{equation}
and 
\begin{equation}
\begin{split}
&
  \norm{ \Theta_t - \psi }
\leq
  \Bigl[
    \bigl[
      \const \rbr{ 1 - \alpha }^{ - 1 } 
    \bigr]^{ \frac{ 1 }{ \alpha - 1 } } 
    \bigl[
      1 
      +
      \const^{ - 2 } ( \tau - \delta_n ) 
    \bigr]
    + 
    \scrC^{ - 1 }
    ( t - \tau ) 
  \Bigr]^{ \alpha - 1 } 
\\ & =
  \Bigl[
    \bigl[
      \const \rbr{ 1 - \alpha }^{ - 1 } 
    \bigr]^{ \frac{ 1 }{ \alpha - 1 } } 
    \bigl[
      1 
      +
      \bigl[
        \const \rbr{ 1 - \alpha }^{ - 1 } 
      \bigr]^{ \frac{ 1 }{ 1 - \alpha } } 
    \bigr]
    + 
    \scrC^{ - 1 }
    ( t - \tau ) 
  \Bigr]^{ \alpha - 1 } 
% \\ & 
\leq 
  \bigl[
    1
    + 
    \scrC^{ - 1 }
    ( t - \tau ) 
  \bigr]^{ \alpha - 1 } 
  .
\end{split}
\end{equation}
\end{cproof}

The next result, \cref{cor:gf:global:abstract},
is a simplified version of \cref{prop:gf:global:abstract} where the sufficiently large finite time $\tau \in [0 , \infty )$ is incorporated in the constant $\scrC$.

\cfclear 
\begin{cor} \label{cor:gf:global:abstract}
	Let 
	$ \fd \in \N $, 
	$ \Theta \in C( [0,\infty), \R^{ \fd } ) $, 
	$ \cL \in C( \R^{ \fd }, \R ) $, 
	let $ \cG \colon \R^{ \fd } \to \R^{ \fd } $ 
	be measurable, 
	assume that for all $\vartheta \in \R^\fd$ there exist
	$\varepsilon, \const \in (0, \infty)$, $\alpha \in ( 0 , 1 )$ such that for all $\theta \in \R^\fd$ with $\norm{\theta - \vartheta } < \varepsilon$ it holds that
	$\abs{\cL ( \theta ) - \cL ( \vartheta ) } ^\alpha \leq \const \norm {\cG ( \theta ) }$,
	and assume for all $ t \in [0,\infty) $ 
	that
\begin{equation} 
	\label{cor:gf:global:eq:ass}
	\textstyle
	\liminf_{ s \to \infty } \norm{ \Theta_s } < \infty ,
	\quad
	\cL( \Theta_t ) = \cL( \Theta_0 ) - \int_0^t \| \cG( \Theta_s ) \|^2 \, \d s ,
	\qandq
	\Theta_t = \Theta_0 - \int_0^t \cG( \Theta_s ) \, \d s .
\end{equation}
Then there exist 
$ \vartheta \in \R^{ \fd } $, $ \scrC, \beta \in (0, \infty) $ 
which satisfy for all $ t \in [0, \infty) $ 
that
\begin{equation} 
\label{cor:gf:global:eq:claim}
  \norm{ \Theta_t - \vartheta } 
  \leq \scrC ( 1 + t )^{ - \beta } 
\qqandqq
  0 \leq \cL( \Theta_t ) - \cL( \vartheta ) \leq \scrC ( 1 + t )^{ - 1 } .
\end{equation}
\end{cor}
\begin{cproof}{cor:gf:global:abstract}
	\Nobs that
	\cref{prop:gf:global:abstract} demonstrates that
	there exist $\vartheta \in \R^\fd$,
	$ \fC, \tau , \beta \in (0, \infty)$ which satisfy for all $t \in [ \tau , \infty )$ that
	\begin{equation} \label{proof:cor:gf:convergence:eq1}
	\norm{\Theta_t - \vartheta} \leq \rbr[\big]{ 1 + \fC^{-1} ( t - \tau ) }^{- \beta} 
% 	\end{equation}
% 	and
% 	\begin{equation} \label{proof:cor:gf:convergence:eq2}
\qqandqq
	0 \leq \cL( \Theta_t ) - \cL( \vartheta ) 
	\leq \rbr[\big]{ 1 + \fC^{ - 1 } ( t - \tau ) }^{ - 1 } .
	\end{equation}
	In the following let $\scrC \in (0, \infty)$ satisfy
	\begin{equation} 
	\label{proof:cor:gf:conv:eq3}
	\scrC =
	\max\bigl\{
	  1 + \tau ,
	  ( 1 + \tau )^{ \beta } , 
	  \fC , 
	  \fC^{ \beta } , 
	  ( 1 + \tau )^{ \beta } 
	  \bigl( 
	    \sup\nolimits_{ s \in [0, \tau ] } \norm{ \Theta_s - \vartheta } 
	  \bigr),
	  ( 1 + \tau ) \abs{ \cL( \Theta_0 ) - \cL ( \vartheta ) }  
        \bigr\} .
	\end{equation}
	\Nobs that
	\cref{proof:cor:gf:convergence:eq1}, 
	\cref{proof:cor:gf:conv:eq3},
	and the fact that $ [0, \infty) \ni t \mapsto \cL( \Theta_t ) \in \R $ 
	is non-increasing show for all $ t \in [0, \tau] $ that
	\begin{equation}
	  \norm{ \Theta_t - \vartheta } 
	  \le 
	  \sup\nolimits_{ s \in [0, \tau] } 
	  \norm{ \Theta_s - \vartheta } 
	  \le 
	  \scrC ( 1 + \tau )^{ - \beta } 
	  \le 
	  \scrC ( 1 + t )^{ - \beta }
	\end{equation}
	and
	\begin{equation}
	  0 \leq \cL( \Theta_t ) - \cL( \vartheta ) 
	  \le \cL( \Theta_0 ) - \cL ( \vartheta ) 
	  \le \scrC ( 1 + \tau )^{ - 1 } 
	  \le \scrC ( 1 + t )^{ - 1 } . 
	\end{equation}
	\Moreover \cref{proof:cor:gf:convergence:eq1,proof:cor:gf:conv:eq3} imply for all $t \in [ \tau , \infty ) $ that
	\begin{equation}
	\begin{split}
	\norm{\Theta_t - \vartheta } 
	&\leq 
	\bigl( 1 + \fC^{ - 1 } ( t - \tau ) \bigr)^{ - \beta } 
	= \scrC 
	\bigl( 
	  \scrC^{ \nicefrac{ 1 }{ \beta } } 
	  + 
          \scrC^{ 
            \nicefrac{ 1 }{ \beta } 
          } 
          \fC^{ - 1 } 
          ( t - \tau ) 
        \bigr)^{ - \beta } 
\\
& 
\le  
  \scrC 
  \bigl(
    \scrC^{ 
      \nicefrac{ 1 }{ \beta } 
    } 
    + 
    t - \tau 
  \bigr)^{ - \beta }
\le 
  \scrC ( 1 + t )^{ - \beta } .
\end{split}
\end{equation}
\Moreover \cref{proof:cor:gf:convergence:eq1,proof:cor:gf:conv:eq3} demonstrate 
for all $ t \in [ \tau, \infty) $ that 
\begin{equation}
  0 \leq \cL( \Theta_t ) - \cL( \vartheta )  
  \le \scrC \rbr[\big]{ \scrC + \fC^{-1} \scrC ( t - \tau ) }^{ - 1 }
  \le \scrC \rbr[\big]{ \scrC - \tau + t }^{-1} \le \scrC ( 1 + t )^{ - 1 } .
\end{equation}
\end{cproof}

\subsection{Convergence of GF processes in the training of deep ANNs}
\label{sec:convergence_GF_processes}

Due to the Kurdyka-\L ojasiewicz inequality for the risk function from \cref{prop:loss:lojasiewicz} we are now able to apply \cref{cor:gf:global:abstract} to the risk function $\cL_\infty$ from \cref{setting:dnn}.

\cfclear
\begin{theorem} \label{theo:gf:conv:simple}
	Assume \cref{setting:dnn},
		assume for all $ i \in \cu{ 1, 2, \ldots, \ell_L } $ that $ f_i $ is piecewise polynomial, 
	let $ \dens \colon [ a, b ]^{ \ell_0 } \to \R $ be piecewise polynomial, 
	 assume for all $ E \in \cB( [ a, b ]^{ \ell_0 } ) $ 
	that $ \mu ( E ) = \int_E \dens ( x ) \, \d x $,
	and let $\Theta \in C ( [ 0 , \infty ) , \R^\fd )$ satisfy $\liminf_{t \to  \infty } \norm{\Theta_t } < \infty$ and $\forall \, t \in [0, \infty ) \colon \Theta_t = \Theta_0 - \int_0^t \cG ( \Theta_s ) \, \d s$ \cfadd{def:multidim:piece:polyn}\cfload.
	Then there exist $\vartheta \in \R^\fd$,
	$\scrC,  \beta \in (0, \infty)$ 
	with  $0 \in (\bbD \cL _\infty ) ( \vartheta )$\cfadd{def:limit:subdiff}
	such that for all $t \in [ 0 , \infty )$ 
	it holds that
	\begin{equation} \label{cor:gf:convergence:eq1}
	\norm{ \Theta_t - \vartheta } 
	\leq 
	\scrC ( 1 + t )^{ - \beta }
	\qqandqq
	0 \leq \cL_{ \infty }( \Theta_t ) - \cL_{ \infty }( \vartheta ) 
	\leq  \scrC ( 1 + t )^{ - 1 } 
	\end{equation}
	\cfload.
\end{theorem}
\begin{cproof}{theo:gf:conv:simple}
\Nobs that \cref{prop:G} shows that 
$ \cG $ is measurable. 
% \Moreover \cref{lem:approx:gradient:bounded} 
% establishes that $ \cG $ is locally bounded. 
\Moreover \cite[Lemma 3.7]{DNNReLUarXiv} ensures 
that for all $ t \in [0, \infty) $ it holds that
\begin{equation} 
\label{proof:gf:convergence:eq:lossintegral}
  \cL_{ \infty }( \Theta_t ) = 
  \cL_{ \infty }( \Theta_0 ) - \int_0^t \norm{ \cG( \Theta_s ) }^2 \, \d s .
\end{equation}
\Moreover \cref{lem:realization:lip} assures that 
\begin{equation}
\label{eq:risk_is_continuous}
  \cL_{ \infty } \in C( \R^{ \fd }, \R ) 
  .
\end{equation}
\Moreover \cref{prop:loss:lojasiewicz} shows that 
for all $ \vartheta \in \R^{ \fd } $ there exist
$ \varepsilon, \const \in (0, \infty) $, 
$ \alpha \in (0, 1) $ 
such that for all 
$ \theta \in \R^{ \fd } $ 
with $ \norm{ \theta - \vartheta } < \varepsilon $ 
it holds that
$
  \abs{ 
    \cL_{ \infty }( \theta ) - \cL_{ \infty }( \vartheta ) 
  }^{ \alpha } 
  \leq \const \norm{ \cG( \theta ) } 
$.
\cref{cor:gf:global:abstract}, 
the fact that 
$
  \cG
$
is measurable, 
\cref{proof:gf:convergence:eq:lossintegral}, 
and 
\cref{eq:risk_is_continuous} 
\hence demonstrate that
there exist 
$ \vartheta \in \R^{ \fd } $, 
$ \scrC, \beta \in (0, \infty) $ 
which satisfy for all $ t \in [0, \infty) $ that
\begin{equation} 
\label{eq:norm_estimate_in_final_theorem_GF}
  \norm{ \Theta_t - \vartheta } 
  \leq 
  \scrC 
  ( 1 + t )^{ - \beta } 
\qqandqq
  0 \leq 
  \cL_{ \infty }( \Theta_t ) - \cL_{ \infty }( \vartheta ) 
  \leq 
  \scrC ( 1 + t )^{ - 1 } .
\end{equation}
\Moreover \cref{proof:gf:convergence:eq:lossintegral} 
demonstrates that 
$ \int_0^{ \infty } 
  \norm{ 
    \cG( \Theta_s )
  }^2
  \, \d s < \infty $. 
\Hence 
$ \liminf_{ s \to \infty } \norm{ \cG( \Theta_s ) } = 0 $. 
This implies that there exists a strictly increasing 
$ \tau = ( \tau_n )_{ n \in \N } \colon \N \to [0,\infty) $ 
which satisfies 
\begin{equation}
\label{eq:construction_of_tau_n_in_final_theorem_GF}
\textstyle
  \liminf_{ n \to \infty } \tau_n = \infty 
\qqandqq
  \limsup_{ n \to \infty } 
  \norm{ \cG( \Theta_{ \tau_n } ) } = 0 .
\end{equation}
\Moreover \cref{eq:norm_estimate_in_final_theorem_GF} 
assures that 
$ \limsup_{ t \to \infty } \| \Theta_t - \vartheta \| = 0 $. 
Combining this with \cref{eq:construction_of_tau_n_in_final_theorem_GF} 
shows that 
\begin{equation}
\label{eq:construction_of_tau_n_in_final_theorem_GF2}
\textstyle
  \limsup_{ n \to \infty } 
  \norm{ \cG( \Theta_{ \tau_n } ) } 
  =
  \limsup_{ n \to \infty } 
  \norm{ \Theta_{ \tau_n } - \vartheta } 
  = 0 .
\end{equation}
\Moreover \cref{prop:loss:gradient:subdiff} assures that 
for all $ \theta \in \R^{ \fd } $ it holds that 
$
  \cG( \theta ) \in ( \bbD \cL_{ \infty } )( \theta )
$. 
\Hence for all $ n \in \N $ that 
\begin{equation}
  \cG( \Theta_{ \tau_n } ) \in ( \bbD \cL_{ \infty } )( \Theta_{ \tau_n } )
  .
\end{equation}
Combining this and \cref{eq:construction_of_tau_n_in_final_theorem_GF2} 
with \cref{lem:limiting_derivatives} 
demonstrates that $ 0 \in ( \bbD \cL_{ \infty } )( \vartheta ) $. 
Combining this with \cref{eq:norm_estimate_in_final_theorem_GF} 
establishes \cref{cor:gf:convergence:eq1}. 
\end{cproof}

\section{Convergence analysis for GD processes} 
\label{sec:convergence_GD}

In this section we establish in \cref{prop:gd:loja} below an abstract local convergence result for GD under a Kurdyka-\L ojasiewicz assumption. 
In the scientific literature related abstract convergence results for GD type processes under a \L ojasiewicz assumption can be found, e.g., in
Absil et al.~\cite{AbsilMahonyAndrews2005}, Attouch \& Bolte~\cite{AttouchBolte2009}, and
Dereich \& Kassing \cite{DereichKassing2021}.
Similar arguments have recently been employed in the analysis of optimization algorithms for tensor decomposition \cite{XuYin2013},
%minimizing fractional functions involving $\ell_1$ and $\ell_2$ norms \cite{ZengYuPong2021},
deep neural networks \cite{Davis2018stochastic},
and residual neural networks \cite{ZengLauLinYao2018}.
The latter two works consider the empirical risk, which is measured with respect to a finite set of training data, while we focus on the true risk defined as the expectation over the entire input distribution.

Except for \cite{XuYin2013} the above mentioned works prove convergence of GD type processes to a critical point, but do not show explicit convergence rates.
 On the other hand, the authors of \cite{XuYin2013} consider block coordinate descent under the assumption that the objective function is convex with respect to each block.
This property is in general not satisfied for objective functions that arise in the training of DNNs with ReLU activation.
The novel contribution of \cref{prop:gd:loja} is to establish a precise convergence rate with explicit constants
and without such convexity assumptions.

 To prove \cref{prop:gd:loja} we transfer the ideas from the continuous-time setting in \cref{section:gf:loja} to the discrete-time setting.
In addition, we require the descent statement in \cref{lem:descent} below.
\cref{lem:descent} below is well-known, see, e.g., Lei et al.~\cite[Lemma 1]{LeiHuLiTang2020},
Attouch et al.~\cite[Lemma 3.1]{AttouchBolteSvaiter2013},
 or Karimi et al.~\cite{Karimi2020linear}.
The elementary proof is only included for completeness.

 In \cref{cor:gd:loja,cor:gd:local:simple} below we establish two simplified versions of \cref{prop:gd:loja}, and as a consequence we obtain in \cref{cor:gd:multi:random:init} below the convergence of GD with random initializations in an abstract setting.
Afterwards, in \cref{prop:gd:random:dnn} below we derive the convergence of GD with random initializations in the considered deep ANN framework in \cref{setting:dnn} under the assumption that there exists a global minimum of the risk function around which suitable regularity assumptions are satisfied.
Our proof of \cref{prop:gd:random:dnn} also uses the ANN approximation result in \cref{prop:dnn:approx} below which, in turn, relies on the universal approximation theorem; cf., e.g., Leshno et al.~\cite{LeshnoLinPinkusSchocken1993}, Cybenko~\cite{Cybenko1989},
Hornik~\cite{Hornik1991}, Lu et al.~\cite{LuShenYangZhang2021}, and Shen et al.~\cite{ShenYangZhang2020}.
As a consequence of \cref{prop:gd:random:dnn} we obtain \cref{theo:gd:random:dnn} below and, thereby, prove \cref{theo:intro:convergence} from the introduction.
Finally, in \cref{prop:gd:random:snn} below we combine \cref{prop:gd:random:dnn} with the existence result for global minima in \cref{cor:existence:regular} to establish the convergence of GD with random initializations in the case of shallow ANNs. As a consequence of \cref{prop:gd:random:snn}
we derive \cref{cor:gd:random:shallow} below and, thereby, prove \cref{theo:intro:random:init} from the introduction.

\subsection{One-step descent property for GD processes}
\begin{lemma} \label{lem:descent}
	Let $\fd \in \N$,
	$L \in \R$,
	let $U \subseteq \R^\fd$ be open and convex,
	let $f \in C^1 (U, \R)$,
	and assume for all $x,y \in U$ that $\norm{ ( \nabla f ) ( x ) - ( \nabla f ) ( y )} \le L \norm{x-y}$.
	Then it holds for all $x,y \in U$ that
	\begin{equation}
	f(y) \le  f(x) + \spro{ (\nabla f ) ( x ) , y - x } + \tfrac{L}{2} \norm{x-y}^2.
	\end{equation}
\end{lemma}
\begin{cproof}{lem:descent}
	\Nobs that the fundamental theorem of calculus,
	the Cauchy-Schwarz inequality,
	and the assumption that for all $x,y \in U$ it holds that $\norm{ ( \nabla f ) ( x ) - ( \nabla f ) ( y )} \le L \norm{x-y}$ 
	assure that for all $x,y \in U$ we have that
	\begin{equation}
	\begin{split}
	f(y) - f(x) 
	&= \br[\big]{ f(x+r(y-x)) }_{r=0}^{r=1} = \int_0^1 \spro{ (\nabla f ) ( x + r ( y - x ) ) , y - x } \, \d r \\
	&= \spro{ (\nabla f ) ( x ) , y - x ) } + \int_0^1 \spro{ (\nabla f ) ( x + r ( y - x ) ) - ( \nabla f ) ( x ), y - x } \, \d r \\
	& \le \spro{ (\nabla f ) ( x ) , y - x ) } + \int_0^1 \abs*{ \spro{ (\nabla f ) ( x + r ( y - x ) ) - ( \nabla f ) ( x ), y - x } } \, \d r \\
	& \le  \spro{ (\nabla f ) ( x ) , y - x ) } + \br*{\int_0^1 \norm{ (\nabla f ) ( x + r ( y - x ) ) - ( \nabla f ) ( x )} \, \d r } \norm{y-x} \\
	& \le \spro{ (\nabla f ) ( x ) , y - x ) } + L \norm{y-x} \br*{ \int_0^1 \norm{r(y-x)} \, \d r } \\
	&= \spro{ (\nabla f ) ( x ) , y - x ) } + \tfrac{L}{2} \norm{x-y}^2.
	\end{split}
	\end{equation}
\end{cproof}

\begin{cor} \label{cor:descent:general}
	Let $\fd \in \N$,
	$L , \gamma  \in \R$,
	let $U \subseteq \R^\fd$ be open and convex,
	let $f \in C^1 (U, \R)$,
	and assume for all $x,y \in U$ that $\norm{ ( \nabla f ) ( x ) - ( \nabla f ) ( y )} \le L \norm{x-y}$.
	Then it holds for all $x \in U $ with $x - \gamma ( \nabla f ) ( x ) \in U$ that
	\begin{equation} \label{cor:descent:general:eq:claim}
	f(x - \gamma ( \nabla f ) ( x )) 
	\le f ( x ) + \gamma \big( \tfrac{L \gamma}{2} - 1 \big) \norm{( \nabla f ) ( x ) }^2 .
	\end{equation}
\end{cor}
\begin{cproof}{cor:descent:general}
	\Nobs that \cref{lem:descent} ensures for all $x \in U $ with $x - \gamma ( \nabla f ) ( x ) \in U$ that
	\begin{equation}
	\begin{split}
	f(x - \gamma ( \nabla f ) ( x )) 
	&\le f(x) + \spro{ ( \nabla f ) ( x ) , - \gamma ( \nabla f ) (x)} + \tfrac{L}{2} \norm{\gamma (\nabla f ) (x )}^2 \\
	&= f(x) - \gamma \norm{ ( \nabla f ) ( x ) }^2 + \tfrac{L \gamma^2}{2} \norm{ ( \nabla f ) ( x ) }^2.
	\end{split}
	\end{equation}
	This establishes \cref{cor:descent:general:eq:claim}.
\end{cproof}

\begin{cor} \label{cor:descent}
	Let $\fd \in \N$,
	$L \in (0 , \infty )$,
	$\gamma \in [0, L^{-1} ]$,
	let $U \subseteq \R^\fd$ be open and convex,
	let $f \in C^1 (U, \R)$,
	and assume for all $x,y \in U$ that $\norm{ ( \nabla f ) ( x ) - ( \nabla f ) ( y )} \le L \norm{x-y}$.
	Then it holds for all $x \in U $ with $x - \gamma ( \nabla f ) ( x ) \in U$ that
	\begin{equation} \label{cor:descent:eq:claim}
	f(x - \gamma ( \nabla f ) ( x )) \le f(x) - \tfrac{\gamma}{2} \norm{ ( \nabla f ) ( x ) }^2 \le f(x).
	\end{equation}
\end{cor}
\begin{cproof}{cor:descent}
	\Nobs that \cref{cor:descent:general},
	the fact that $\gamma \ge 0$,
	and the fact that $\frac{L \gamma}{2} - 1 \le - \frac{1}{2}$
	establish
	\cref{cor:descent:eq:claim}.
\end{cproof}

\subsection{Abstract local convergence results for GD processes}
\begin{prop} 
\label{prop:gd:loja}
Let $ \fd \in \N $, $ \consttt \in \R $,
$ \varepsilon, L, \const \in (0, \infty) $,
	$\alpha \in (0,1)$,
	$\gamma \in (0, L^{-1}]$,
	$\vartheta \in \R^\fd$,
	let $\bB \subseteq \R^\fd$ satisfy $\bB = \cu{ \theta \in \R^\fd \colon \norm{\theta - \vartheta } < \varepsilon}$,
	let $\cL \in C ( \R^\fd , \R)$ satisfy $\cL |_{\bB} \in C^1 ( \bB , \R)$,
	let $\cG \colon \R^\fd \to \R^\fd$ satisfy for all $\theta \in \bB$ that $\cG (\theta ) = ( \nabla \cL ) ( \theta )$,
	assume $\cG ( \vartheta ) = 0$,
	assume for all $\theta_1, \theta_2 \in \bB$ that $\norm{ \cG ( \theta_1 ) - \cG ( \theta_2 )} \le L \norm{\theta_1 - \theta_2}$,
	let $\Theta \colon \N_0 \to \R^\fd$ satisfy for all $n \in \N_0$ that $\Theta_{n+1} = \Theta_n - \gamma \cG (\Theta_n )$,
	and assume for all $\theta \in \bB$ that
	\begin{equation} \label{prop:gd:eq:ass:loja}
	\abs{ 
		\cL( \theta ) - \cL( \vartheta ) 
	}^\alpha 
	\leq 
	\const \norm{ \cG( \theta ) } ,
	\quad 
	\consttt = 
	\abs{\cL( \Theta_0 ) - \cL( \vartheta )},
	\quad 
	2 \const 
	( 1 - \alpha )^{ - 1 }
	\consttt^{ 1 - \alpha }
	+
	\norm{ \Theta_0 - \vartheta }
	< 
	\tfrac{\varepsilon}{\gamma L + 1 },
	\end{equation}
	and $\inf_{n \in \cu{m \in \N_0 \colon \forall \, k \in \N_0 \cap [0, m ] \colon \Theta_k \in \bB}} \cL (\Theta_n  ) \ge \cL ( \vartheta )$.
	Then there exists $\psi \in \cL^{-1} ( \cu{\cL ( \vartheta )}) \cap \cG^{-1} ( \cu{0})$ such that 
	\begin{enumerate} [label = (\roman*)]
		\item \label{prop:gd:item1} it holds for all $n \in \N_0$ that $\Theta_n \in \bB$,
		\item \label{prop:gd:item2} it holds for all $n \in \N_0$ that
		$
		0 \le \cL ( \Theta_n ) - \cL ( \psi )  
		\le 
		2 \const^2 \consttt^2 
		\rbr{ 
			\indicator{ \cu{0}} ( \consttt ) + \consttt^{2 \alpha} n \gamma + 2 \const^2 \consttt 
		}^{ - 1 } ,
		$
		and
		\item \label{prop:gd:item3} it holds for all $n \in \N_0$ that
		\begin{equation}
		\begin{split}
		\norm{\Theta_n - \psi } 
		&\le \smallsum_{k=n}^\infty \norm{\Theta_{k+1} - \Theta_k } 
		\le 2 \const \rbr{1 - \alpha}^{-1}
		\abs{\cL ( \Theta_n ) - \cL ( \psi ) }^{1 - \alpha} \\
		& \le
		2^{2 - \alpha} \const^{3 - 2 \alpha} \consttt^{2 - 2 \alpha }(1 - \alpha)^{-1}
		\rbr{ 
			\indicator{ \cu{0}} ( \consttt ) + \consttt^{2 \alpha} n \gamma + 2 \const^2 \consttt 
		}^{ \alpha - 1 } .
		\end{split}
		\end{equation}
	\end{enumerate}
\end{prop}
Observe that the assumption that $\inf_{n \in \cu{m \in \N_0 \colon \forall \, k \in \N_0 \cap [0, m ] \colon \Theta_k \in \bB}} \cL (\Theta_n  ) \ge \cL ( \vartheta )$ is in particular satisfied
 if $\vartheta$ is a local minimum of $\cL$ with $\forall \, \theta \in  \bB \colon \cL ( \theta ) \ge \cL ( \vartheta )$. Hence \cref{prop:gd:loja} implies as a consequence a local convergence result of GD to a local minimum.
But our assumption also covers more general cases, since we only require an estimate on the values of $\cL ( \Theta_n)$ and not for all values $\cL ( \theta)$ with $\theta \in \bB$.

\begin{cproof}{prop:gd:loja}
	Throughout this proof 
	let $T \in \N_0 \cup \cu{\infty}$ satisfy
	\begin{equation}
	T = \inf \rbr*{ \cu{n \in \N_0 \colon \Theta_n \notin \bB } \cup \cu{\infty}} 
	,
	\end{equation} 
	let $ \bL \colon \N_0 \to \R $ 
	satisfy for all $ n \in \N_0 $ that 
	$ 
	\bL( n ) = \cL( \Theta_n ) - \cL( \vartheta ) 
	$, 
	and let $ \tau \in \N_0 \cup \cu{ \infty } $ satisfy
	\begin{equation} \label{prop:gd:eq:deftau}
	\tau = \inf \rbr*{ \cu{n \in \N_0 \cap [0, T ) \colon \bL(n) = 0} \cup \cu{T}}.
	\end{equation}
	In the first step of our proof we verify that $T=\infty$,
	i.e., that the Gd iterates remain inside the neighborhood $\bB$ at all times.
	\Nobs that the 
	assumption that $\cG ( \vartheta ) = 0$
	implies for all $\theta \in \bB$ that
	\begin{equation} \label{prop:gd:eq:stepsize}
	\gamma \norm{\cG ( \theta ) } = \gamma  \norm{\cG ( \theta ) - \cG ( \vartheta ) } \le  \gamma L \norm{\theta - \vartheta } .
	\end{equation}
	This,
	the fact that $\norm{\Theta_0 - \vartheta }  < \varepsilon$,
	and the fact that 
	\begin{equation}
	\norm{\Theta_1 - \vartheta } \le \norm{\Theta_1 - \Theta_0 } + \norm{\Theta_0 - \vartheta } = \gamma \norm{\cG ( \Theta_0 ) } + \norm{\Theta_0 - \vartheta} \le (\gamma L + 1) \norm{\Theta_0 - \vartheta}  < \varepsilon
	\end{equation}
	ensure that $T \ge 2$.
	Next \nobs that the assumption that $\inf_{n \in \cu{m \in \N_0 \colon \forall \, k \in \N_0 \cap [0, m ] \colon \Theta_k \in \bB}} \cL (\Theta_n  ) \ge \cL ( \vartheta )$ proves for all $n \in \N_0 \cap [0, T)$ that $\bL ( n ) \ge 0$.
	%We intend to prove that $\tau = \infty$. We thus assume for the sake of contradiction that $\tau < \infty$.
	In addition, \nobs that the fact that $\bB \subseteq \R^\fd$ is open and convex, \cref{cor:descent},
	and \cref{prop:gd:eq:ass:loja} demonstrate for all $n \in \N_0 \cap [0, T - 1 )$ that
	\begin{equation} \label{prop:gd:eq:descentest}
	\begin{split}
	\bL ( n+1 ) - \bL (n) 
	&= \cL ( \Theta_{n+1} ) - \cL ( \Theta_n) \le - \tfrac{\gamma}{2} \norm{\cG ( \Theta_n ) } ^2 
	= - \tfrac{1}{2} \norm{\cG ( \Theta_n ) } \norm{\gamma \cG ( \Theta_n ) } \\
	&= - \tfrac{1}{2} \norm{\cG ( \Theta_n ) } \norm{\Theta_{n+1} - \Theta_n } 
	\le - (2 \const)^{-1} \abs{\cL ( \Theta_n ) - \cL ( \vartheta ) }^\alpha  \norm{\Theta_{n+1} - \Theta_n } \\
	& = - (2 \const)^{-1} [ \bL ( n ) ] ^{\alpha}  \norm{\Theta_{n+1} - \Theta_n } \le 0 .
	\end{split}
	\end{equation}
	Therefore, we obtain that $\N_0 \cap [0, T ) \ni n \mapsto \bL(n) \in [0, \infty)$ is non-increasing.
	Combining this with \cref{prop:gd:eq:deftau} shows for all $n \in \N_0 \cap [\tau , T )$ that $\bL(n) = 0$.
	This and \cref{prop:gd:eq:descentest}
	demonstrate for all $n \in \N_0 \cap [\tau , T - 1 )$ that $0 = \bL ( n+1 ) - \bL ( n ) \le - \frac{\gamma}{2} \norm{\cG ( \Theta_n ) } ^2 \le 0$.
	The fact that $\gamma > 0$ therefore assures 
	for all $n \in \N_0 \cap [\tau , T - 1 )$  that $\cG ( \Theta_n ) = 0$.
	Hence, we obtain for all $n \in \N_0 \cap[\tau , T)$ that 
	\begin{equation} \label{prop:gd:eq:theta:stop}
	\Theta_n  = \Theta_{\tau}.
	\end{equation}
	In addition, \nobs that \cref{prop:gd:eq:descentest,prop:gd:eq:deftau} ensure for all $n \in \N_0 \cap [0, \tau) \cap [0, T - 1 )$ that
	\begin{equation}
	\begin{split}
	\norm{\Theta_{n+1} - \Theta_n } 
	&\le \frac{2 \const ( \bL (n) - \bL (n+1))}{[\bL (n)]^\alpha} = 2 \const \int_{\bL (n+1)}^{\bL (n)} [\bL (n)]^{-\alpha} \, \d u \\
	& \le 2 \const \int_{\bL (n+1)}^{\bL (n)} u^{-\alpha} \, \d u = \frac{2 \const \rbr*{[\bL (n)]^{1 - \alpha} - [ \bL ( n+1) ]^{1 - \alpha}}}{1 - \alpha}.
	\end{split}
	\end{equation}
	This and \cref{prop:gd:eq:theta:stop} show for all $n  \in \N_0 \cap [0, T - 1 )$ that
	\begin{equation}
	\norm{\Theta_{n+1} - \Theta_n } \le \frac{2 \const \rbr*{[\bL (n)]^{1 - \alpha} - [ \bL ( n+1) ]^{1 - \alpha}}}{1 - \alpha}.
	\end{equation}
	Combining this with the triangle inequality
	proves for all $m , n \in \N_0 \cap [0, T )$ with $m \le n$ that
	\begin{equation} \label{prop:gd:eq:triangleest}
	\begin{split}
	\norm{\Theta_n - \Theta_m} 
	& \le \sum_{k=m}^{n-1} \norm{\Theta_{k+1} - \Theta_k } \le \frac{2 \const}{1 - \alpha} \br*{ \sum_{k=m}^{n-1} \rbr*{[\bL ( k ) ]^{1 - \alpha} - [ \bL ( k + 1) ]^{1 - \alpha}} } \\
	& 
	= \frac{2 \const \rbr*{[ \bL ( m ) ]^{1 - \alpha} - [ \bL ( n ) ]^{1 - \alpha} } }{1 - \alpha}
	\le \frac{2 \const [ \bL ( m ) ]^{1 - \alpha}}{1 - \alpha}.
	\end{split}
	\end{equation}
	This and \cref{prop:gd:eq:ass:loja} demonstrate for all $n \in \N_0 \cap [0, T )$ that
	\begin{equation} \label{prop:gd:eq:dist:theta0}
	\norm{\Theta_n - \Theta_0 } \le \frac{2 \const [ \bL ( 0 )]^{1 - \alpha}}{1 - \alpha} = \frac{2 \const \abs{\cL ( \Theta_0 ) - \cL ( \vartheta )}^{1 - \alpha}}{1 - \alpha} = 2 \const 
	( 1 - \alpha )^{ - 1 }
	\consttt^{ 1 - \alpha } .
	\end{equation}
	Combining this with \cref{prop:gd:eq:stepsize}, \cref{prop:gd:eq:ass:loja},
	and the triangle inequality shows for all $n \in \N_0 \cap [0, T )$ that
	\begin{equation} \label{prop:gd:eq:induct:dist}
	\begin{split}
	\norm{\Theta_{n+1} - \vartheta }
	& \le \norm{\Theta_{n+1} - \Theta_n } + \norm{\Theta_n - \vartheta} 
	= \gamma \norm{ \cG ( \Theta_n ) } + \norm{\Theta_n - \vartheta } \\
	& \le ( \gamma L + 1 ) \norm{\Theta_n - \vartheta } \le( \gamma L + 1 ) (\norm{\Theta_n - \Theta_0} + \norm{\Theta_0 - \vartheta }) \\
	& \le ( \gamma L + 1 ) ( 2 \const 
	( 1 - \alpha )^{ - 1 }
	\consttt^{ 1 - \alpha } + \norm{\Theta_0 - \vartheta}) <  \varepsilon.
	\end{split}
	\end{equation}
	Hence, we obtain that 
	\begin{equation} \label{prop:gd:eq:tau:infty}
	T = \infty .
	\end{equation}
	Combining this with \cref{prop:gd:eq:ass:loja} 
	and \cref{prop:gd:eq:triangleest} proves that
	\begin{equation} \label{prop:gd:eq:traj:length}
	\sum_{k=0}^\infty \norm{\Theta_{k+1} - \Theta_k } = \lim_{n \to \infty} \br*{	\sum_{k=0}^n \norm{\Theta_{k+1} - \Theta_k } }
	\le \frac{2 \const [ \bL ( 0 )]^{1 - \alpha}}{1 - \alpha} = \frac{2 \const \consttt^{1 - \alpha}}{1 - \alpha} < \varepsilon < \infty. 
	\end{equation}
	Therefore, we obtain that there exists $\psi \in \R^\fd$ which satisfies
	\begin{equation} \label{prop:gd:eq:limitpsi}
	\limsup\nolimits_{n \to \infty} \norm{\Theta_n - \psi } = 0.
	\end{equation}
	This establishes convergence of the GD process. We next deduce explicit convergence rates.
	\Nobs that \cref{prop:gd:eq:induct:dist,prop:gd:eq:tau:infty,prop:gd:eq:limitpsi} imply that $\norm{\psi - \vartheta } \le ( \gamma L + 1 ) ( 2 \const 
	( 1 - \alpha )^{ - 1 }
	\consttt^{ 1 - \alpha } + \norm{\Theta_0 - \vartheta}) <  \varepsilon $.
	Therefore, we obtain that $\psi \in \bB$.
	Next \nobs that \cref{prop:gd:eq:descentest},
	\cref{prop:gd:eq:ass:loja},
	and the fact that for all $n \in \N_0$ 
	it holds that $ \bL( n ) \le \bL ( 0 ) = \consttt $ ensure 
	that for all $ n \in \N_0 \cap [0, \tau) $ we have that
	\begin{equation} \label{prop:gd:eq:ln:descent}
	- \bL ( n ) 
	\le 
	\bL (n+1) - \bL (n) 
	\le - \tfrac{\gamma}{2} \norm{\cG ( \Theta_n)}^2 \le - \tfrac{\gamma}{2 \const^2} [\bL ( n ) ]^{2 \alpha} \le -  \tfrac{\gamma}{2 \const^2 \consttt^{2 - 2 \alpha}} [ \bL ( n ) ]^{2 }.
	\end{equation}
	This assures for all $ n \in \N_0 \cap [0, \tau) $ that 
	$ 
	0 < \bL (n) \le \frac{ 2 \const^2 \consttt^{ 2 - 2 \alpha } }{ \gamma } 
	$. 
	Combining this and \cref{prop:gd:eq:ln:descent} demonstrates 
	for all $ n \in \N_0 \cap [0, \tau - 1 ) $ that
	\begin{equation}
	\begin{split}
	\frac{ 1 }{ \bL(n) } - \frac{ 1 }{ \bL(n+1) } 
	& 
	\le 
	\frac{ 1 }{ 
		\bL(n) 
	} 
	- 
	\frac{ 1 }{ 
		\bL(n) 
		( 
		1 - 
		\tfrac{ 
			\gamma 
		}{ 2 \const^2 \consttt^{ 2 - 2 \alpha } 
		} 
		\bL( n ) 
		) 
	} 
	=
	\frac{
		\big( 
		1 
		- 
		\frac{ \gamma }{ 2 \const^2 \consttt^{ 2 - 2 \alpha } } 
		\bL( n )
		\big) 
		- 1
	}{
		\bL(n) 
		\big( 
		1 
		- 
		\frac{ \gamma }{ 2 \const^2 \consttt^{ 2 - 2 \alpha } } 
		\bL( n )
		\big)
	}
	\\ &   =
	\frac{
		- 
		\frac{ \gamma }{ 2 \const^2 \consttt^{ 2 - 2 \alpha } } 
	}{
		\big( 
		1 
		- 
		\frac{ \gamma }{ 2 \const^2 \consttt^{ 2 - 2 \alpha } } 
		\bL( n )
		\big)
	}
	=
	- 
	\frac{ 
		1 
	}{ 
		(
		\tfrac{ 2 \const^2 \consttt^{ 2 - 2 \alpha } 
		}{ \gamma } 
		- \bL(n)
		)
	} 
	< 
	- \frac{ \gamma }{ 2 \const^2 \consttt^{ 2 - 2 \alpha } } 
	.
	\end{split}
	\end{equation}
	Therefore,
	we get for all $n \in \N_0 \cap [0, \tau )$ that
	\begin{equation}
	\frac{1}{ \bL ( n ) } 
	= \frac{1}{\bL ( 0 ) } + \sum_{k=0}^{n-1} \br*{ \frac{1}{\bL ( k+1 )} - \frac{1}{\bL ( k ) } }
	\ge \frac{1}{ \bL (0)} + \frac{n \gamma}{2 \const^2 \consttt^{2 - 2 \alpha} } = \frac{1}{\consttt} +  \frac{n \gamma}{2 \const^2 \consttt^{2 - 2 \alpha} }.
	\end{equation}
	Hence, we obtain for all $n \in \N_0 \cap [0, \tau )$ that $ \bL (n) \le \frac{2 \const^2 \consttt^{2 - 2 \alpha} }{n \gamma + 2 \const^2 \consttt^{1 - 2 \alpha}}$. Combining this with 
	the fact that for all $n \in \N_0 \cap  [ \tau , \infty )$ it holds that $\bL (n) = 0$ shows that for all $n \in \N_0$ we have that
	\begin{equation} \label{prop:gd:eq:ln:finalest}
	\bL (n) \le \frac{2 \const^2 \consttt^{2} }{\indicator{ \cu{0 } } ( \consttt) + \consttt ^{2 \alpha} n \gamma + 2 \const^2 \consttt} .
	\end{equation}
	This, \cref{prop:gd:eq:limitpsi}, and the assumption that $\cL$ is continuous prove that 
	\begin{equation} \label{prop:gd:eq:limit:risk}
	\cL( \psi ) 
	= 
	\lim\nolimits_{n \to \infty} \cL( \Theta_n ) 
	= \cL ( \vartheta ) .
	\end{equation} 
	Combining this with \cref{prop:gd:eq:ln:finalest}
	assures for all $n \in \N_0$ that
	\begin{equation} \label{prop:gd:eq:risk:rate}
	0 \le 
	\cL ( \Theta_n ) - \cL ( \psi ) 
	\le \frac{2 \const^2 \consttt^{2} }{\indicator{ \cu{0} } ( \consttt) + \consttt ^{2 \alpha} n \gamma + 2 \const^2 \consttt} .
	\end{equation}
	Furthermore, \nobs that 
	the fact that $ \bB \ni \theta \mapsto \cG ( \theta ) \in \R^\fd $ is continuous,
	% 	\cref{prop:gd:eq:tau:infty},
	the fact that $ \psi \in \bB$,
	and \cref{prop:gd:eq:limitpsi} imply that 
	\begin{equation} \label{prop:gd:eq:limit:gradient}
	\cG( \psi ) 
	=
	\lim\nolimits_{ n \to \infty } \cG( \Theta_n )
	= 
	\lim\nolimits_{ n \to \infty } 
	( \gamma^{ - 1 } ( \Theta_{ n } - \Theta_{ n + 1 } ) ) 
	= 0 .
	\end{equation}
	Next \nobs that \cref{prop:gd:eq:ln:finalest} and \cref{prop:gd:eq:triangleest} 
	ensure for all $n \in \N_0$ that
	\begin{equation}
	\begin{split}
	\norm{\Theta_n - \psi } 
	&= \lim_{m \to \infty} \norm{\Theta_n - \Theta_m } \le \sum_{k=n}^\infty \norm{\Theta_{k+1} - \Theta_k }
	\le \frac{2 \const [ \bL ( n ) ]^{1 - \alpha}}{1 - \alpha} \\
	& \le 
	\frac{ 2^{2 - \alpha} \const^{3 - 2 \alpha} \consttt^{2  - 2 \alpha } }{(1 - \alpha)
		\rbr{ \indicator{ \cu{0}} ( \consttt ) + \consttt^{2 \alpha} n \gamma + 2 \const^2 \consttt } ^{1 - \alpha } }.
	\end{split}
	\end{equation}
	Combining this with \cref{prop:gd:eq:limit:risk}, 
	\cref{prop:gd:eq:tau:infty}, \cref{prop:gd:eq:limit:gradient}, 
	and \cref{prop:gd:eq:risk:rate} establishes \cref{prop:gd:item1,prop:gd:item2,prop:gd:item3}.
\end{cproof}

The next result, \cref{cor:gd:loja}, specializes \cref{prop:gd:loja} to the case where $\vartheta \in \R^\fd$ is a local minimum of $\cL$ in the sense that for all $\theta \in \bB$ we have that $\cL ( \theta ) \ge \cL ( \vartheta )$, where $\bB$ is a suitable neighborhood of $\vartheta$.

\begin{cor} \label{cor:gd:loja}
	Let $\fd \in \N$,
	$\consttt \in [0 , 1 ]$,
	$\varepsilon, L, \const \in (0, \infty)$,
	$\alpha \in (0,1)$,
	$\gamma \in (0, L^{-1}]$,
	$\vartheta \in \R^\fd$,
	let $ \bB \subseteq \R^{ \fd } $ satisfy 
	$
	  \bB = \cu{ \theta \in \R^\fd \colon \norm{\theta - \vartheta } < \varepsilon } 
	$,
	let $ \cL \in C( \R^{ \fd }, \R) $ satisfy $ \cL|_{\bB} \in C^1( \bB , \R) $,
	let $\cG \colon \R^\fd \to \R^\fd$ satisfy for all $\theta \in \bB$ that 
	$ \cG( \theta ) = ( \nabla \cL ) ( \theta )$,
	assume for all $\theta_1, \theta_2 \in \bB$ that $\norm{ \cG ( \theta_1 ) - \cG ( \theta_2 )} \le L \norm{\theta_1 - \theta_2}$,
	let $ \Theta = ( \Theta_n )_{ n \in \N_0 } \colon \N_0 \to \R^{ \fd } $ satisfy for all $n \in \N_0$ that $\Theta_{n+1} = \Theta_n - \gamma \cG (\Theta_n )$,
	and assume for all  $\theta \in \bB$ that
	\begin{equation} \label{cor:gd:eq:ass:loja}
	\abs{ 
		\cL( \theta ) - \cL( \vartheta ) 
	}^\alpha 
	\leq 
	\const \norm{ \cG( \theta ) } ,
	\quad 
	\consttt = 
	\abs{\cL( \Theta_0 ) - \cL( \vartheta )},
	\quad 
	2 \const 
	( 1 - \alpha )^{ - 1 }
	\consttt^{ 1 - \alpha }
	+
	\norm{ \Theta_0 - \vartheta }
	< 
	\tfrac{\varepsilon}{\gamma L + 1 },
	\end{equation}
	and $ \cL( \theta ) \geq \cL( \vartheta ) $.
	Then there exists $\psi \in \cL^{-1} ( \cu{\cL ( \vartheta )}) \cap \cG^{-1} ( \cu{0})$ such that for all $n \in \N_0$ it holds that $\Theta_n \in \bB$, $0 \le \cL ( \Theta_n ) - \cL ( \psi ) \le 2 ( 2 + \const^{-2} \gamma n )^{-1}$, and
	\begin{equation}
	\norm{\Theta_n - \psi } \le \smallsum_{k=n}^\infty \norm{\Theta_{k+1} - \Theta_k } \le 
	2^{2 - \alpha} \const (1 - \alpha)^{-1} ( 2 + \const^{-2} \gamma n )^{\alpha - 1} .
	\end{equation}
\end{cor}
\begin{cproof}{cor:gd:loja}
	\Nobs that
	the fact that $\cL (\vartheta ) = \inf_{\theta \in \bB} \cL ( \theta )$
	ensures that $\cG ( \vartheta ) = ( \nabla \cL ) ( \vartheta ) = 0$
	and
	$\inf_{n \in \cu{m \in \N_0 \colon \forall \, k \in \N_0 \cap [0, m ] \colon \Theta_k \in \bB}} \cL (\Theta_n  ) \ge \cL ( \vartheta )$.
	Combining this with \cref{prop:gd:loja} implies that there exists
	$\psi \in \cL^{-1} ( \cu{\cL ( \vartheta )}) \cap \cG^{-1} ( \cu{0})$ such that 
	\begin{enumerate}[label = (\Roman*)]
		\item \label{cor:gd:item1} it holds for all $n \in \N_0$ that $\Theta_n \in \bB$,
		\item \label{cor:gd:item2} it holds for all $n \in \N_0$ that
		$0 \le \cL ( \Theta_n ) - \cL ( \psi )  \le \frac{ 2 \const^2 \consttt^{2} }{\indicator{ \cu{0}} ( \consttt ) + \consttt^{2 \alpha} n \gamma + 2 \const^2 \consttt } $,
		and
		\item \label{cor:gd:item3} it holds for all $n \in \N_0$ that
		\begin{equation}
		\begin{split}
		\norm{\Theta_n - \psi } 
		&\le \sum_{k=n}^\infty \norm{\Theta_{k+1} - \Theta_k } 
		\le \frac{2 \const \abs{\cL ( \Theta_n ) - \cL ( \psi ) }^{1 - \alpha}}{1 - \alpha} \\
		& \le
		\frac{ 2^{2 - \alpha} \const^{3 - 2 \alpha} \consttt^{2  - 2 \alpha } }{(1 - \alpha)
			\rbr{ \indicator{ \cu{0}} ( \consttt ) + \consttt^{2 \alpha} n \gamma + 2 \const^2 \consttt } ^{1 - \alpha } }.
		\end{split}
		\end{equation}
	\end{enumerate}
	\Nobs that \cref{cor:gd:item2} and the assumption that $\consttt \le 1$ show for all $n \in \N_0$ that
	\begin{equation}
	0 \le \cL ( \Theta_n ) - \cL ( \psi ) \le  2  \consttt^{2} \rbr*{\const^{-2} \indicator{ \cu{0}} ( \consttt ) + \const^{-2} \consttt^{2 \alpha} n \gamma + 2 \consttt }^{-1} \le 2 ( 2 + \const^{-2} \gamma n )^{-1}.
	\end{equation}
	This and \cref{cor:gd:item3} demonstrate for all $n \in \N_0$ that
	\begin{equation}
	\begin{split}
	\norm{\Theta_n - \psi } 
	\le \sum_{k=n}^\infty \norm{\Theta_{k+1} - \Theta_k } 
	\le \frac{2 \const \abs{\cL ( \Theta_n ) - \cL ( \psi ) }^{1 - \alpha}}{1 - \alpha} 
	\le 
	\left[ 
	\frac{ 2^{2 - \alpha} \const }{1 - \alpha } 
	\right] 
	( 2 + \const^{-2} \gamma n )^{\alpha - 1} .
	\end{split}
	\end{equation}
\end{cproof}

\subsection{Abstract global convergence results for GD processes}

In \cref{cor:gd:local:simple} we reformulate \cref{cor:gd:loja}
to show that around every local minimum point which admits a Kurdyka-\L ojasiewicz inequality and a certain regularity condition there exists an open neighborhood such that the risk of every GD sequence started in this neighborhood converges with rate $1$ to the risk of the local minimum.

\begin{cor} \label{cor:gd:local:simple}
	Let $\fd \in \N$,
	$\varepsilon, L, \const \in (0, \infty)$,
	$\alpha \in (0,1)$,
	$\vartheta \in \R^\fd$,
	let $\bB \subseteq \R^\fd$ satisfy $\bB = \cu{ \theta \in \R^\fd \colon \norm{\theta - \vartheta } < \varepsilon}$,
	let $\cL \in C ( \R^\fd , \R)$ satisfy $\cL |_{\bB} \in C^1 ( \bB , \R)$,
	let $\cG \colon \R^\fd \to \R^\fd$ satisfy for all $\theta \in \bB$ that $\cG (\theta ) = ( \nabla \cL ) ( \theta )$,
	assume for all $\theta_1, \theta_2 \in \bB$ that $\norm{ \cG ( \theta_1 ) - \cG ( \theta_2 )} \le L \norm{\theta_1 - \theta_2}$,
	for every $\theta \in \R^\fd$, $\gamma \in \R$ 
	let 
	$ 
	  \Theta^{ \gamma, \theta } = 
	  ( \Theta^{ \gamma, \theta }_n )_{ n \in \N_0 } \colon \N_0 \to \R^\fd
	$ 
	satisfy for all $ n \in \N_0 $ that 
	$
	  \Theta_0^{ \gamma, \theta } = \theta
	$
	and
	$\Theta^{\gamma , \theta }_{n+1} = \Theta^{\gamma , \theta }_n - \gamma \cG (\Theta^{\gamma , \theta }_n )$,
	and assume for all  $\theta \in \bB$ that
	\begin{equation} \label{cor:gd:simple:eq:ass:loja}
	\abs{ 
		\cL( \theta ) - \cL( \vartheta ) 
	}^\alpha 
	\leq 
	\const \norm{ \cG( \theta ) }
	\qqandqq \cL( \theta ) \geq \cL( \vartheta ) .
	\end{equation}
	Then there exist
	$\delta , \scrC \in (0, \infty)$
	such that for all
	$\theta \in \cu{\psi \in \R^\fd \colon \norm{\psi - \vartheta } < \delta }$,
	$\gamma \in (0, L^{-1}]$,
	$n \in \N_0$
	 it holds that $0 \le \cL ( \Theta_n^{\gamma , \theta } ) - \cL ( \vartheta ) \le \scrC(1 + \gamma n ) ^{-1}$. 
\end{cor}

\begin{cproof}{cor:gd:local:simple}
	\Nobs that the fact that $\cL$ is continuous ensures that there exist
	$\consttt \in [0,1]$,
	 $\delta \in (0, \varepsilon )$ which satisfy for all $\theta \in \cu{\psi \in \R^\fd \colon \norm{\psi - \vartheta } < \delta }$, $\gamma \in (0 , L^{-1}]$ that 
	\begin{equation}
	\label{eq:global_convergence_GD_construction_of_delta}
	\consttt = 
	\abs{\cL( \theta ) - \cL( \vartheta )}
	\qqandqq
	2 \const 
	( 1 - \alpha )^{ - 1 }
	\consttt^{ 1 - \alpha }
	+
	\norm{ \theta - \vartheta }
	< \tfrac{\varepsilon}{2} \le
	\tfrac{\varepsilon}{\gamma L + 1 } .
	\end{equation}
	\Nobs that \cref{eq:global_convergence_GD_construction_of_delta} 
	and \cref{cor:gd:loja} 
	(applied for every 
	$ \theta \in \{ \psi \in \R^{ \fd } \colon \norm{ \psi - \vartheta } < \delta \} $, 
	$ \gamma \in ( 0, L^{ - 1 } ] $ 
	with 
	$ \varepsilon \with \delta $, 
	$ \gamma \with \gamma $, 
	$ \bB \with \{ \psi \in \R^{ \fd } \colon \norm{ \psi - \vartheta } < \delta \} $,
	$ \Theta \with \Theta^{ \gamma, \theta } $
	in the notation of \cref{cor:gd:loja})
	demonstrate that for all
	$ \theta \in \{ \psi \in \R^{ \fd } \colon \norm{ \psi - \vartheta } < \delta \} $,
	$\gamma \in (0, L^{-1}]$
	there exists $\psi \in \cL^{-1} ( \cu{\cL ( \vartheta )})$
	such that for all
	$n \in \N_0$
	it holds that
	\begin{equation}
	  0 \le 
	  \cL( \Theta_n^{ \gamma, \theta } ) - \cL( \psi ) 
	  \le 2 ( 2 + \fC^{ - 2 } \gamma n )^{ - 1 } .
	\end{equation}
	\Hence for all
	$\theta \in \cu{\psi \in \R^\fd \colon \norm{\psi - \vartheta } < \delta }$,
	$\gamma \in (0, L^{-1}]$,
	$n \in \N_0$
	that
	\begin{equation}
	0 
	\le  \cL( \Theta_n^{\gamma , \theta } ) - \cL( \vartheta ) 
	\le 
	2 ( 2 + \fC^{ - 2 } \gamma_n )^{ - 1 }
	\le
	2 ( \min \cu{ 2, \fC^{ - 2 } } ( 1 + \gamma n ) )^{ - 1 }
	= \max \cu{ 1 , 2 \fC ^2 } ( 1 + \gamma n )^{-1}.
	\end{equation}
\end{cproof}

\subsection{Abstract convergence result for GD with random initializations}

The next result, \cref{cor:gd:multi:random:init}, establishes convergence in probability of the GD method with multiple random initalizations under a \L ojasiewicz type assumption.
The proof relies on \cref{cor:gd:local:simple} and the fact that for a sufficiently high number of initilizations at least one of the GD trajectories will start in a suitable open domain of attraction with high probability.

\begin{cor} \label{cor:gd:multi:random:init}
Let $ \fd \in \N $,
$ \varepsilon, L, \const, \gamma \in (0, \infty) $,
$ \alpha \in (0,1) $,
$ \vartheta \in \R^{ \fd } $
satisfy
$ \gamma L \le 1 $,
let $\bB \subseteq \R^\fd$ satisfy $\bB = \cu{ \theta \in \R^\fd \colon \norm{\theta - \vartheta } < \varepsilon}$,
let $\cL \in C ( \R^\fd , \R)$ satisfy $\cL |_{\bB} \in C^1 ( \bB , \R)$,
let $\cG \colon \R^\fd \to \R^\fd$ satisfy for all $\theta \in \bB$ that $\cG (\theta ) = ( \nabla \cL ) ( \theta )$,
assume for all $\theta_1, \theta_2 \in \bB$ that 
$
  \norm{ \cG( \theta_1 ) - \cG( \theta_2 ) } 
  \le L \norm{ \theta_1 - \theta_2 } 
$,
assume for all $ \theta \in \bB $ that
\begin{equation} 
  \abs{ 
    \cL( \theta ) - \cL( \vartheta ) 
  }^{ \alpha }
  \leq 
  \const \norm{ \cG( \theta ) }
\qqandqq 
  \cL( \theta ) \geq \cL( \vartheta ) ,
\end{equation}
let $ ( \Omega, \cF, \P) $ be a probability space, 
for every $ K, n \in \N_0 $
let $ \Theta^K_n \colon  \Omega \to \R^{ \fd } $
% , $ K, n \in \N_0 $, 
and 
$
  \bfk^K_n \colon \Omega \to \N 
$
% ,$
%   K, n \in \N_0 
% $, 
be random variables,
assume that $ \Theta_0^K $, $ K \in \N $, are i.i.d.,
assume for all $ \delta \in (0, \infty) $ that
$
  \P( \norm{ \Theta_0^1 - \vartheta } < \delta ) > 0 
$,
and assume for all $ K \in \N $, $ n \in \N_0 $,
$ \omega \in \Omega $ that
\begin{equation} 
\label{cor:gd:random:eq:defk}
  \Theta_{ n + 1 }^K( \omega ) = 
  \Theta_n^K( \omega ) 
  - \gamma \cG( \Theta_n^K( \omega ) ) 
\qqandqq 
  \bfk^K_n( \omega) 
  \in \arg\min\nolimits_{ \kappa \in \{ 1, 2, \ldots, K \} } 
  \cL( \Theta_n^{ \kappa }( \omega ) ) .
\end{equation}
Then
\begin{equation}
\label{cor:gd:multi:random:init:eq}
  \liminf\nolimits_{ K \to \infty } 
  \P\bigl(
    \limsup\nolimits_{ n \to \infty } 
    \cL(
      \Theta^{ \bfk^K_n }_n 
    )
    \le \cL( \vartheta )  
  \bigr) = 1 .
\end{equation}
\end{cor}
\begin{cproof}{cor:gd:multi:random:init}
\Nobs that \cref{cor:gd:random:eq:defk} shows for all $ K \in \N $ that
\begin{equation} 
\label{cor:gd:random:proof:eq1}
\begin{split}
& 
  \P\bigl(
    \limsup\nolimits_{ n \to \infty } \cL\rbr[\big]{ \Theta^{  \bfk^{ K }_n}_n } 
    \le \cL ( \vartheta) 
  \bigr)
\\
& 
  \geq 
  \P\rbr*{ 
    \exists \, \kappa \in \{ 1, 2, \ldots, K \} \colon 
%     \br[\big]{ 
      \limsup\nolimits_{ n \to \infty } 
      \cL( \Theta_n^{ \kappa } ) \le \cL( \vartheta ) 
%     } 
  }
  .
\end{split}
\end{equation}
\Moreover \cref{cor:gd:local:simple} demonstrates 
that there exist $ \delta, \scrC \in (0, \infty) $
which satisfy for all 
$ \kappa \in \N $, 
$ \omega \in \cu{ w \in \Omega \colon \norm{ \Theta^{ \kappa }_0( w ) - \vartheta } < \delta } $,
$ n \in \N_0 $
that 
$ 
  0 \le \cL( \Theta_n^{ \kappa }( \omega ) ) - \cL( \vartheta ) 
  \le \scrC ( 1 + \gamma n )^{ - 1 }
$.
\Hence for all $ \kappa \in \N $, 
$ 
  \omega \in 
  \cu{ w \in \Omega \colon \norm{ \Theta^{ \kappa }_0( w ) - \vartheta } < \delta } 
$
that 
\begin{equation}
\textstyle
  \limsup_{ n \to \infty } \cL( \Theta_n^{ \kappa }( \omega ) )
  \leq \cL( \vartheta )
  .
\end{equation}
This shows for all $ \kappa \in \N $ that 
$
  \{ 
    \omega \in \Omega \colon 
    \norm{ \Theta^{ \kappa }_0( \omega ) - \vartheta } < \delta 
  \} 
  \subseteq
  \{ 
    \omega \in \Omega \colon
    \limsup_{ n \to \infty } \cL( \Theta_n^{ \kappa }( \omega ) )
    \leq \cL( \vartheta )
  \}
$. 
\Hence for all $ K \in \N $ that
\begin{equation} 
\label{cor:gd:random:proof:eq2}
\begin{split}
& 
  \P\bigl(
    \exists \, \kappa \in \{ 1, 2, \ldots, K \} \colon 
    \norm{ \Theta_0^{ \kappa } - \vartheta } < \delta 
  \bigr)
  =
  \P\bigl(
    \cup_{ \kappa = 1 }^K
    \{ 
      \norm{ \Theta_0^{ \kappa } - \vartheta } < \delta 
    \}
  \bigr)
\\ &
  \le 
  \P\bigl(
    \cup_{ \kappa = 1 }^K 
    \{ 
      \limsup\nolimits_{ n \to \infty } 
      \cL( \Theta_n^{ \kappa } ) \le \cL ( \vartheta ) 
    \}
  \bigr)
\\ &
  =
  \P\bigl(
    \exists \, \kappa \in \{ 1, 2, \ldots, K \} \colon 
%     \br[\big]{ 
      \limsup\nolimits_{ n \to \infty } 
      \cL( \Theta_n^{ \kappa } ) \le \cL ( \vartheta ) 
%     } 
  \bigr)
  .
\end{split}
\end{equation}
\Moreover the assumption that $ \Theta_0^{ \kappa } $, $ \kappa \in \N $, 
are i.i.d.\ proves that for all $ K \in \N $ we have that
\begin{equation} 
\label{cor:gd:random:proof:eq3}
\begin{split}
  \P\bigl(
    \exists \, \kappa \in \{ 1, 2, \ldots, K \} \colon 
    \norm{ \Theta_0^{ \kappa } - \vartheta } < \delta 
  \bigr) 
& 
  = 1 - 
  \P\bigl( 
    \forall \, \kappa \in \{ 1, 2, \ldots, K \} \colon 
    \norm{ \Theta_0^{ \kappa } - \vartheta } \ge \delta 
  \bigr) 
\\
&
  = 1 - 
  \bigl[ 
    \P\bigl( 
      \norm{ \Theta_0^1 - \vartheta } \ge \delta 
    \bigr) 
  \bigr]^K
  .
\end{split}
\end{equation}
The fact that 
$
  \P( 
    \norm{ \Theta_0^1 - \vartheta } \ge \delta 
  ) 
  = 
  1 - 
  \P( 
    \norm{ \Theta_0^1 - \vartheta } < \delta 
  )
  < 1 
$
\hence implies that 
\begin{equation} 
  \liminf\nolimits_{ K \to \infty } 
  \P\bigl( 
    \exists \, \kappa \in \{ 1, 2, \ldots, K \} \colon 
    \norm{ \Theta_0^{ \kappa } - \vartheta } < \delta 
  \bigr) = 1
  .
\end{equation}
Combining this 
with \cref{cor:gd:random:proof:eq1,cor:gd:random:proof:eq2} 
establishes \cref{cor:gd:multi:random:init:eq}.
\end{cproof}

\subsection{Approximation results for deep ANNs}

We next show an $L^2$-universal approximation result for shallow ANNs. In \cref{lem:dnn:shallow:approx} the target function is not necessarily continuous and takes values in a multidimensional space  $\R^\delta$.
We establish \cref{lem:dnn:shallow:approx} by employing the universal approximation theorem for $\R$-valued functions in Leshno et al.~\cite[Proposition~1 in Section~4]{LeshnoLinPinkusSchocken1993}. The proof is only included for completeness.

\begin{lemma} 
\label{lem:dnn:shallow:approx}
Let $ d, \delta \in \N $, $ a \in \R $, $ b \in (a, \infty) $, 
$ \varepsilon \in (0,\infty) $ 
and let 
$ f = ( f_1, \ldots, f_{ \delta } ) \colon [a,b]^d \to \R^{ \delta } $
and $ \dens \colon [a,b]^d \to [0, \infty) $ be bounded and measurable. 
Then there exist
$ n \in \N $,
$ \fw_1, \fw_2, \ldots, \fw_n \in \R^{ 1 \times d } $,
$ \fb_1, \fb_2, \ldots, \fb_n \in \R $,
$ \fv_1, \fv_2, \ldots, \fv_n \in \R^{ \delta } $ 
such that
\begin{equation}
\textstyle
  \int_{ [a,b]^d } 
  \norm{
    f( x )
    -
    \smallsum_{ i = 1 }^n 
    \fv_i 
    \max\cu{ 
      \fw_i x + \fb_i , 0 
    }  
  }^2 
  \,
  \dens( x ) \, \d x 
  < \varepsilon .
\end{equation}
\end{lemma}
\begin{cproof}{lem:dnn:shallow:approx}
Throughout this proof let 
$ 
  \mu \colon \cB( \R^d ) \to [0, \infty] 
$
satisfy for all 
$ E \in \cB( \R^d ) $
that
\begin{equation}
\label{eq:def_mu_measure_in_proof}
\textstyle
  \mu( E ) = \int_{ [a,b]^d \cap E } \dens( x ) \, \d x 
\end{equation}
and 
for every $ i \in \{ 1, 2, \dots, \delta \} $
let 
$ 
  F_i \colon \R^d \to \R
$
satisfy for all 
$ x \in [a,b]^d $, 
$ y \in \R^d \backslash [a,b]^d $
that
\begin{equation}
\label{eq:def_F_in_proof}
  F_i( x ) = f_i( x )
\qqandqq
  F_i( y ) = 0 .
\end{equation}
\Nobs that \cref{eq:def_mu_measure_in_proof} 
and the assumption that $ \dens $ 
is bounded and measurable 
ensure that $ \mu $ is a 
finite, absolutely continuous, and compactly supported measure. 
The assumption that $ f $ is bounded \hence implies 
that for all $ i \in \cu{ 1, 2, \ldots, \delta } $ 
it holds that 
$
  \int_{ \R^d } \abs{ F_i( x ) }^2 \, \mu( \d x ) < \infty 
$.
Combining this, the universal approximation theorem 
(cf.\ Leshno et al.~\cite[Proposition~1 in Section~4]{LeshnoLinPinkusSchocken1993}
(applied with 
$ \sigma \with ( \R \ni x \mapsto \max \cu{ x, 0 } \in \R ) $,
$ \mu \with \mu $, $ p \with 2 $ 
in the notation of \cite[Proposition~1 in Section~4]{LeshnoLinPinkusSchocken1993})), 
\cref{eq:def_mu_measure_in_proof}, 
\cref{eq:def_F_in_proof}, 
and the fact that 
$ \mu $ is 
a finite, absolutely continuous, and compactly supported measure 
proves that 
for every $ i \in \cu{ 1, 2, \ldots, \delta } $
there exist $ \bfn^{ (i) } \in \N $,
$ 
  \bfw_1^{ (i) }, \bfw_2^{ (i) }, \ldots, \bfw_{ \bfn^{ (i) } }^{ (i) } \in \R^{ 1 \times d } 
$,
$
  \bfb_1^{ (i) }, \bfb_2^{ (i) }, \ldots, 
  \bfb_{ \bfn^{ (i) } }^{ (i) }, \bfv_1^{ (i) }, \bfv_2^{ (i) }, 
  \ldots, \bfv_{ \bfn^{(i)} }^{ (i) } \in \R 
$ 
which satisfy
\begin{equation} 
\label{lem:shallow:approx:eq1}
\begin{split}
&
\textstyle
  \int_{ [a,b]^d } 
  \abs{ 
    f_i( x ) 
    -
    \smallsum_{ k = 1 }^{ \bfn^{ (i) } } 
    \bfv_k^{ (i) } 
    \max\cu{ 
      \bfw_k^{ (i) } x + \bfb_k^{ (i) }, 0 
    } 
  }^2 
  \, \dens ( x ) \, \d x 
\\
&
\textstyle
  = 
  \int_{ \R^d } 
  | 
    F_i( x ) 
    -
    \smallsum_{ k = 1 }^{ \bfn^{ (i) } } 
    \bfv_k^{ (i) } 
    \max\cu{ 
      \bfw_k^{ (i) } x + \bfb_k^{ (i) }, 0 
    } 
  |^2 \, \mu( \d x )
  < 
  \frac{ \varepsilon }{ \delta } 
  .
\end{split}
\end{equation}
In the following 
let $ e_1, e_2, \dots, e_{ \delta } \in \R^{ \delta } $
satisfy 
$
  e_1 = ( 1, 0, 0, \dots, 0 )
$, 
$
  e_2 = ( 0, 1, 0, \dots, 0 )
$,
$ \dots $,
$ 
  e_{ \delta } = ( 0, \dots, 0, 1 ) 
$,
let $ n \in \N $ satisfy 
$ n = \sum_{ i = 1 }^\delta \bfn^{ (i) } $, 
and let 
$ \fw_1, \fw_2, \ldots, \fw_n \in \R^{ 1 \times d } $,
$ \fb_1, \fb_2, \ldots, \fb_n \in \R $,
$ \fv_1, \fv_2, \ldots, \fv_n \in \R^{ \delta } $
satisfy for all
$ j \in \cu{ 1, 2, \ldots, \delta } $, 
$ k \in \{ 1, 2, \dots, \bfn^{ (j) } \} $
that
\begin{equation}
\label{eq:def_v_in_proof_universal_approximation}
  \fw_{ k + \sum_{ i = 1 }^{ j - 1 } \bfn^{ (i) } } 
  = \bfw^{ (j) }_k ,
\qquad
  \fb_{ k + \sum_{ i = 1 }^{ j - 1 } \bfn^{ (i) } } 
  = \bfb^{ (j) }_k ,
\qqandqq
  \fv_{ k + \sum_{ i = 1 }^{ j - 1 } \bfn^{ (i) } } 
  = \bfv^{ (j) }_k e_j
%   ( \underbrace{ 0, \ldots, 0 }_{ j - 1 }, \bfv^{ (j) }_k, 0, \ldots, 0 )
  .
\end{equation}
\Nobs that \cref{eq:def_v_in_proof_universal_approximation} 
assures that for all $ c_1, c_2, \ldots, c_n \in \R $,
$ x \in [a,b]^d $ 
it holds that
\begin{equation}
  \norm[\big]{
    f(x)
    -
    \smallsum_{ i = 1 }^n 
    c_i \fv_i 
  }^2 
  = 
  \smallsum_{ j = 1 }^{ \delta } 
  \abs[\big]{
    f_j( x ) 
    -
    \smallsum_{ k = 1 }^{ \bfn^{ (j) } } 
    \bfv_k^{ (j) } 
    c_{ 
      k 
      +
      \sum_{ i = 1 }^{ j - 1 } 
      \bfn^{ (i) } 
    } 
  }^2.
\end{equation}
Combining this with \cref{lem:shallow:approx:eq1} 
and \cref{eq:def_v_in_proof_universal_approximation} demonstrates that
\begin{equation}
\begin{split}
& 
\textstyle
  \int_{ [a,b]^d } 
  \norm[\big]{ 
    f( x ) 
    -
    \smallsum_{ i = 1 }^n 
    \bigl[
      \max\cu{ 
        \fw_i x + \fb_i , 0 
      } 
    \bigr]
    \fv_i 
  }^2 
  \dens( x ) \, \d x 
\\
&
\textstyle
  = 
  \sum\limits_{ j = 1 }^{ \delta } 
  \int\limits_{ [a,b]^d } 
    \abs[\big]{ 
      f_j( x ) 
      -
      \smallsum_{ k = 1 }^{ \bfn^{ (j) } } 
      \bfv_k^{ (j) } 
      \max\cu{ 
        \fw_{ k + \sum_{ i = 1 }^{ j - 1 } \bfn^{ (i) } } 
        x 
        + 
        \fb_{ k + \sum_{ i = 1 }^{ j - 1 } \bfn^{ (i) } } 
        , 0 
      } 
    }^2 
  \dens ( x ) \, \d x 
\\
&
\textstyle
  = 
  \sum\limits_{ j = 1 }^{ \delta } 
  \int\limits_{ [a,b]^d } 
    \abs[\big]{ 
      f_j( x ) 
      -
      \smallsum_{ k = 1 }^{ \bfn^{ (j) } } 
      \bfv_k^{ (j) } 
      \max\cu{ 
        \bfw_k^{ (j) } x 
        + 
        \bfb_k^{ (j) } , 0 
      } 
    }^2 
  \dens ( x ) \, \d x 
  < \delta \bigl[ \frac{ \varepsilon }{ \delta } \bigr] 
  = \varepsilon
  .
\end{split}
\end{equation}
\end{cproof}

As a consequence of \cref{lem:dnn:shallow:approx}
we show in \cref{prop:dnn:approx} a universal approximation result for deep ANNs as the width increases to infinity.

\begin{prop} \label{prop:dnn:approx}
Let 
$ d, \delta \in \N $, 
$ a \in \R $, $ b \in [a,\infty) $, 
$ ( \rho_{ \fa } )_{ \fa \in \N } \subseteq ( \N \cap (1,\infty) ) $,
let 
$
  \ell^{ \fa } = ( \ell^{ \fa }_0, \ell^{ \fa }_1, \dots, \ell^{ \fa }_{ \rho_{ \fa } } ) 
  \allowbreak 
  \in 
  \{ d \} \times \N^{ \rho_{ \fa } - 1 } \times \{ \delta \}$, 
$ \fa \in \N $, 
satisfy 
\begin{equation} 
\label{prop:dnn:approx:eq:architecture}
\textstyle
  \liminf_{ \fa \to \infty }
  \min\cu{
    \ell^{ \fa }_1, \ell^{ \fa }_2, \ldots , \ell^{ \fa }_{ \rho_{ \fa } - 1 } 
  }
  = \infty ,
\end{equation}
for every $ \fa \in \N $
let 
$
  \fd_{ \fa } = 
  \sum_{ k = 1 }^{ \rho_{ \fa } } 
  \ell^{ \fa }_k ( \ell^{ \fa }_{ k - 1 } + 1 ) 
$, 
let
$ f \colon [a,b]^d \to \R^{ \delta } $ 
and 
$ \dens \colon [ a, b ]^d \to [0, \infty) $
be bounded and measurable,
for every $ \fa \in \N $, $ k \in \{ 1, 2, \dots, \rho_{ \fa } \} $, 
$ \theta = ( \theta_1, \dots, \theta_{ \fd_{ \fa } } ) \in \R^{ \fd_{ \fa } } $
let 
$ 
  \fw^{ k, \theta }_{ \fa } 
  = 
  ( 
    \fw^{ k, \theta }_{ \fa, i, j } 
  )_{ 
    (i,j) \in 
    \{ 1, \ldots, \ell_k^{ \fa } \} \times 
    \{ 1, \ldots, \ell_{ k - 1 }^{ \fa } \} } \in \R^{ \ell_k^{ \fa } \times 
    \ell_{ k - 1 }^{ \fa } 
  } 
$
and 
$ 
  \fb^{ k, \theta }_{ \fa } = 
  ( 
    \fb^{ k, \theta }_{ \fa, 1 }, \dots, \fb^{ k, \theta }_{ \fa, \ell_k^{ \fa } } 
  )
  \in \R^{ \ell_k^{ \fa } }
$
satisfy for all 
$ i \in \{ 1, 2, \ldots, \ell_k^{ \fa } \} $, 
$ j \in \{ 1, 2, \ldots, \ell_{ k - 1 }^{ \fa } \} $ 
that
\begin{equation}
  \fw^{ k, \theta }_{\fa , i, j } 
  =
  \theta_{ 
    (i-1) \ell^\fa_{k-1} + j 
    + \sum_{h=1}^{k-1} \ell^\fa_h (\ell^\fa_{h-1} + 1)
  }
\qqandqq
  \fb^{ k, \theta }_{\fa , i} 
  =
  \theta_{ 
    \ell^{ \fa }_k \ell^{ \fa }_{ k - 1 } + i + 
    \sum_{ h = 1 }^{ k - 1 } \ell^{ \fa }_h ( \ell^{ \fa }_{ h - 1 } + 1 ) 
  } 
  ,
\end{equation}
let 
$ 
  \mathfrak{M} \colon ( \cup_{ n \in \N } \R^n )
  \to ( \cup_{ n \in \N } \R^n ) 
$
satisfy for all 
$ n \in \N $,
$ x = (x_1, \ldots, x_n ) \in \R^n $
that 
$ 
  \mathfrak{M}( x ) = ( \max\{ x_1, 0 \}, \ldots, \max\{ x_n, 0 \} ) 
$, 
for every 
$ \fa \in \N $, 
$ \theta \in \R^{ \fd_{ \fa } } $ 
let 
$ 
  \mathcal{N}^{ k, \theta }_{ \fa } 
  \colon \R^d \to \R^{ \ell^{ \fa }_k } 
$, 
$ k \in \N \cap [1, \rho_{ \fa } ]
$, 
satisfy for all 
$ k \in \N \cap [1, \rho_{ \fa } ) $,
$ x \in \R^d $
that 
\begin{equation}
  \mathcal{N}^{ 1, \theta }_{ \fa }( x ) 
  =
  \fb^{ 1, \theta }_{ \fa } 
  +
  \fw^{ 1, \theta }_{ \fa } x 
\qqandqq
  \mathcal{N}^{ k+1, \theta }_{ \fa }( x ) 
  =
  \fb^{ k + 1, \theta }_{ \fa } 
  +
  \fw^{ k + 1, \theta }_{ \fa }\big(
    \mathfrak{M}( \mathcal{N}^{ k, \theta }_{ \fa }( x ) ) 
  \big)
  ,
\end{equation} 
and for every $ \fa \in \N $
let 
$ \cL_{ \fa } \colon \R^{ \fd_{ \fa } } \to \R $
satisfy for all
$ \theta \in \R^{ \fd_{ \fa } } $ 
that 
$ 
  \cL_{ \fa }( \theta ) 
  =
  \int_{ [a,b]^{ d } } 
  \norm{ \mathcal{N}_{ \fa }^{ \rho_\fa, \theta }( x ) - f(x) }^2 \dens( x )
  \, \d x 
$.
Then
\begin{equation}
  \limsup\nolimits_{ \fa \to \infty } 
  \inf\nolimits_{ \theta \in \R^{ \fd_{ \fa } } } 
  \cL_{ \fa }( \theta ) = 0 .
\end{equation}
\end{prop}
\begin{cproof2}{prop:dnn:approx}
Throughout this proof let $ \varepsilon \in (0, \infty) $.
\Nobs that \cref{lem:dnn:shallow:approx}
proves that there exist $ n \in \N $, $ \bfw \in \R^{ n \times d } $,
$ \bfb \in \R^n $, $ \bfv \in \R^{ \delta \times n } $
which satisfy
\begin{equation} 
\label{prop:dnn:approx:eq1}
\textstyle
  \int_{ [a,b]^d } 
  \norm{
    \bfv \fM( \bfw x + \bfb ) - f ( x ) 
  }^2 \dens( x ) \, \d x 
  < \varepsilon .
\end{equation}
\Moreover \cref{prop:dnn:approx:eq:architecture}
assures that there exists $ \bfA \in \N $
which satisfies for all $ \fa \in \N \cap [ \bfA , \infty) $, 
$ i \in \N \cap ( 1, \rho_{ \fa } ) $
that
$ 
  \ell_1^{ \fa } \ge n 
$
and
$
  \ell_i^{ \fa } \ge 2 \delta 
$.
Combining this with Beck et al.~\cite[Lemma 2.10]{BeckJentzenKuckuck2022} 
(applied for every $ \fa \in \N \cap [ \bfA, \infty ) $ with
$ L \with 2 $, 
$ (l_0, l_1, l_2 ) \with ( d, n, \delta ) $, 
$ d \with n ( d + 1 ) + \delta ( n + 1 ) $, 
$ \fL \with \rho_\fa $, 
$ ( \fl_0, \fl_1, \ldots, \fl_{ \fL } ) \with \ell^{ \fa } $,
$ \fd \with \fd_{ \fa } $ 
in the notation of \cite[Lemma 2.10]{BeckJentzenKuckuck2022}) 
shows for every 
$ \fa \in \N \cap [ \bfA, \infty) $
that there exists $ \theta_{ \fa } \in \R^{ \fd_{ \fa } } $
which satisfies for all $ x \in \R^d $
that
\begin{equation}
\label{eq:construction_of_identity_DNN}
  \cN_{ \fa }^{ \rho_{ \fa }, \theta_{ \fa } }( x ) 
  = \bfv \fM( \bfw x + \bfb ) 
  .
\end{equation}
\Nobs that \cref{prop:dnn:approx:eq1,eq:construction_of_identity_DNN} 
ensure for all $ \fa \in \N \cap [\bfA , \infty ) $
that
\begin{equation}
\textstyle
  \inf_{ \vartheta \in \R^{ \fd_{ \fa } } } 
  \cL_{ \fa }( \vartheta ) 
  \le 
  \cL_{ \fa }( \theta_{ \fa } )
  =
  \int_{ [a,b]^d } 
    \norm{ \cN_{ \fa }^{ \rho_{ \fa }, \theta_{ \fa } }( x ) - f( x ) 
  }^2 \dens( x ) \, \d x 
  < \varepsilon .
\end{equation}
\end{cproof2}

\subsection{Convergence of GD with random initializations in the training of deep ANNs}
\label{sec:GD_random_init}

We next combine the Kurdyka-\L ojasiewicz inequality from \cref{prop:loss:lojasiewicz}
with the abstract convergence result for GD with random initializations from \cref{cor:gd:multi:random:init}
to prove convergence in probability of GD with random initializations for deep ANNs with a fixed architecture. In \cref{prop:gd:random:dnn} the parameter vector $\vartheta \in \R^\fd$ is assumed to be a local minimum of the risk function $\cL_\infty$ in a neighborhood of which the regularity assumptions in \cref{cor:gd:multi:random:init} are satisfied. The convergence holds for every sufficiently small positive learning rate $\gamma \in (0 , \fg ]$.

\cfclear
\begin{prop}
\label{prop:gd:random:dnn}
Assume \cref{setting:dnn},
assume for all $ i \in \cu{ 1, 2, \ldots, \ell_L } $ that $ f_i $ is piecewise polynomial, 
let $ \dens \colon [ a, b ]^{ \ell_0 } \to \R $ be piecewise polynomial\cfadd{def:multidim:piece:polyn}, 
assume for all $ E \in \cB( [ a, b ]^{ \ell_0 } ) $ 
that $ \mu( E ) = \int_E \dens ( x ) \, \d x $,
let $ U \subseteq \R^{ \fd } $ be open, 
assume $ ( \cL_{ \infty } )|_U \in C^1( U, \R) $, 
let $ \fG \colon \R^{ \fd } \to \R^{ \fd } $ 
satisfy for all 
$ 
  \theta \in U
$
that 
$ 
  \fG( \theta ) = ( \nabla \cL_{ \infty } )( \theta ) 
$,
assume that $ \fG|_U $ is locally Lipschitz continuous,
let $ \vartheta \in U $ satisfy 
$ 
  \cL_{ \infty }( \vartheta ) 
%   = \inf_{ \theta \in \R^{ \fd } } \cL_{ \infty }( \theta ) 
  = \inf_{ \theta \in U } \cL_{ \infty }( \theta ) 
$,
let $ ( \Omega, \cF, \P ) $ be a probability space, 
for every $ K, n \in \N_0 $, $ \gamma \in \R $ 
let 
$ \Theta^{ K, \gamma }_n \colon \Omega \to \R^{ \fd } $
and 
$
  \bfk^{ K, \gamma }_n \colon \Omega \to \N 
$
be random variables,
assume for all $ \gamma \in \R $ that $ \Theta_0^{ K, \gamma } $, $ K \in \N $, 
are i.i.d., assume for all
$ \gamma, \delta \in (0, 1) $
that
$
  \P( 
    \norm{ \Theta_0^{ 1, \gamma } - \vartheta } < \delta 
  ) > 0 
$,
and assume for all $ K \in \N $, $ n \in \N_0 $, $ \gamma \in \R $, $ \omega \in \Omega $ 
that
\begin{equation} 
\label{prop:gd:random:eq:defk}
  \Theta_{ n + 1 }^{ K, \gamma }( \omega ) 
  = 
  \Theta_n^{ K, \gamma }( \omega ) - \gamma \fG( \Theta_n^{ K, \gamma }( \omega ) ) 
\qandq 
  \bfk^{ K, \gamma }_n( \omega) 
  \in 
  \arg\min\nolimits_{
    \kappa \in \{ 1, 2, \ldots, K \} 
  } 
  \cL_{ \infty }( 
    \Theta_n^{ \kappa, \gamma }( \omega ) 
  )
\end{equation}
\cfload.
Then there exists $ \fg \in (0, \infty) $ 
such that for all $ \gamma \in (0, \fg] $ 
it holds that
\begin{equation}
  \liminf\nolimits_{ K \to \infty } 
  \P\Bigl( 
    \limsup\nolimits_{ n \to \infty } 
    \cL_{ \infty }( 
      \Theta^{ \bfk^{ K, \gamma }_n , \gamma }_n 
    ) 
    \le \inf\nolimits_{ \theta \in U } \cL_{ \infty }( \theta ) 
  \Bigr)
  = 1 .
\end{equation}
\end{prop}
\begin{cproof}{prop:gd:random:dnn}
\Nobs that \cref{cor:cG_equal_to_gradient} 
assures for all open $ V \subseteq \R^{ \fd } $ 
and all $ \theta \in V $
with 
$ ( \cL_{ \infty } )|_V \in C^1( V, \R) $
that 
$
  \cG( \theta ) = 
  ( \nabla \cL_{ \infty } )( \theta )
$. 
The assumption that $ U \subseteq \R^{ \fd } $ is open, 
the assumption that $ ( \cL_{ \infty } )|_U \in C^1( U, \R ) $, 
and the assumption that for all $ \theta \in U $ it holds that 
$ \fG( \theta ) = ( \nabla \cL_{ \infty } )( \theta ) $ 
\hence demonstrates that for all 
$ \theta \in U $ it holds that 
\begin{equation}
\label{eq:gradient_smooth_on_U_in_proof}
  \cG( \theta ) = 
  ( \nabla \cL_{ \infty } )( \theta )
  = 
  \fG( \theta )
  .
\end{equation}
\Moreover \cref{prop:loss:lojasiewicz}, 
the assumption that $ U \subseteq \R^{ \fd } $ is open,
and the assumption that $ \fG|_U $ is locally Lipschitz continuous 
assure that there exist 
$ \bL, \varepsilon, \fC \in (0, \infty) $, $ \alpha \in (0, 1) $ 
which satisfy for all 
$
  v, w 
  \in 
  \cu{ 
    \psi \in \R^{ \fd } \colon \norm{ \psi - \vartheta } < \varepsilon 
  } 
$ 
that
\begin{equation}
\label{eq:property_eps_Loja_in_proof}
  v \in U , 
\qquad
  \abs{
    \cL_{ \infty }( v ) - \cL_{ \infty }( \vartheta ) 
  }^{ \alpha } 
  \le 
  \fC \norm{ \cG( v ) } ,
\qqandqq
  \norm{
    \fG( v ) - \fG( w ) 
  } 
  \le \bL \norm{ v - w } 
  .
\end{equation}
\Moreover \cref{lem:realization:lip} shows that 
$
  \cL_{ \infty } \in C( \R^{ \fd }, \R )
$. 
Combining this, 
\cref{eq:gradient_smooth_on_U_in_proof}, 
\cref{eq:property_eps_Loja_in_proof}, 
\cref{cor:gd:multi:random:init} 
(applied for every $ \gamma \in (0, \bL^{ - 1 } ] \cap (0,1) $ 
with 
$
  \fd \with \fd 
$, 
$
  \varepsilon \with \varepsilon
$, 
$
  L \with \bL
$, 
$ 
  \fC \with \fC
$, 
$
  \gamma \with \gamma 
$,
$
  \alpha \with \alpha 
$, 
$
  \vartheta \with \vartheta
$, 
$
  \cL \with \cL_{ \infty }
$, 
$
  \cG \with \fG
$
in the notation of \cref{cor:gd:multi:random:init}), 
and the assumption that 
$ 
  \cL_{ \infty }( \vartheta ) = \inf_{ \theta \in U } \cL_{ \infty }( \theta )
$
assures that for all $ \gamma \in (0, \bL^{ - 1 } ] \cap (0,1) $ 
it holds that
\begin{equation}
\begin{split}
&
  \liminf\nolimits_{ K \to \infty } 
  \P\Bigl(
    \limsup\nolimits_{ n \to \infty } 
    \cL_{ \infty }\rbr[\big]{ 
      \Theta^{ \bfk^{ K, \gamma }_n , \gamma }_n 
    } 
    \le 
    \inf\nolimits_{ \theta \in U } \cL_{ \infty }( \theta ) 
  \Bigr) 
\\
& 
  =
  \liminf\nolimits_{ K \to \infty } 
  \P\Bigl(
    \limsup\nolimits_{ n \to \infty } 
    \cL_{ \infty }\rbr[\big]{ 
      \Theta^{ \bfk^{ K, \gamma }_n , \gamma }_n 
    } 
    \le 
    \cL_{ \infty }( \vartheta ) 
  \Bigr) 
  = 1 .
\end{split}
\end{equation}
\end{cproof}

As a consequence of \cref{prop:gd:random:dnn}
and the universal approximation result from \cref{prop:dnn:approx}
we verify in \cref{theo:gd:random:dnn}
that the risk of the GD method with 
random initializations converges in probability 
to $ 0 $ as the number of GD steps,
the number of random initializations,
and the width of the ANNs increase to $ \infty $
and as the step size of the GD method decreases to $ 0 $.
In \cref{cor:gd:random:dnn:item1} we establish convergence in probability, and as a consequence we obtain in \cref{cor:gd:random:dnn:item2} convergence with respect to the metric $ \E [ \min\{ | X - Y |, 1 \} ] $ on the space of random variables.

\cfclear
\begin{theorem}
\label{theo:gd:random:dnn}
Let $ d, \delta \in \N $, $ a \in \R $, $ b \in [a,\infty) $, 
$ ( \rho_{ \fa } )_{ \fa \in \N } \subseteq ( \N \cap (1,\infty) ) $,
let 
$ 
  \ell^{ \fa } = ( \ell^{ \fa }_0, \ell^{ \fa }_1, \dots, \ell^{ \fa }_{ \rho_{ \fa } } ) \in 
  \{ d \} \times \N^{ \rho_{ \fa } - 1 } \times \{ \delta \}
$, 
$ \fa \in \N $, 
satisfy 
\begin{equation}
\textstyle
  \liminf_{ \fa \to \infty }
  \min \cu{\ell^\fa_1, \ell^\fa_2, \ldots , \ell^{\fa}_{\rho_\fa - 1 } }
  = \infty,
\end{equation}
for every $ \fa \in \N $
let 
$
  \fd_{ \fa } = 
  \sum_{ k = 1 }^{ \rho_{ \fa } } 
  \ell^{ \fa }_k ( \ell^{ \fa }_{ k - 1 } + 1 )
$, 
let 
$ f = (f_1, \ldots, f_\delta) \colon [a,b]^d \to \R^{ \delta } $ 
and 
$
  \dens \colon [a,b]^d \to [0,\infty)
$
be functions, 
assume for all 
$ i \in \cu{ 1, 2, \dots, \delta } $ 
that $ f_i $ and $ \dens $ are piecewise polynomial\cfadd{def:multidim:piece:polyn},
for every $ \fa \in \N $, 
$ k \in \{ 1, 2, \dots, \rho_{ \fa } \} $, 
$ 
  \theta = ( \theta_1, \dots, \theta_{ \fd_{ \fa } } ) \in \R^{ \fd_{ \fa } } 
$
let 
$ 
  \fw^{ k, \theta }_{ \fa }
  = 
  ( 
    \fw^{ k, \theta }_{ \fa, i, j } 
  )_{ 
    (i,j) \in 
    \{ 1, \ldots, \ell_k^{ \fa } \} \times 
    \{ 1, \ldots, \ell_{ k - 1 }^{ \fa } \}
  }
  \in \R^{ \ell_k^{ \fa } \times \ell_{ k - 1 }^{ \fa } } 
$
and 
$ 
  \fb^{ k, \theta }_{ \fa } 
  = 
  ( 
    \fb^{ k, \theta }_{ \fa, 1 }, \dots, \fb^{ k, \theta }_{ \fa, \ell_k^{ \fa } } 
  )
  \in \R^{ \ell_k^{ \fa } }
$
satisfy for all 
$ i \in \{ 1, 2, \ldots, \ell_k^{ \fa } \} $, 
$ j \in \{ 1, 2, \ldots, \ell_{ k - 1 }^{ \fa } \} $ 
that
\begin{equation}
	\fw^{ k, \theta }_{\fa , i, j } 
	=
	\theta_{ 
		(i-1) \ell^\fa_{k-1} + j 
		+ \sum_{h=1}^{k-1} \ell^\fa_h (\ell^\fa_{h-1} + 1)
	}
	\qqandqq
	\fb^{ k, \theta }_{\fa , i } 
	=
	\theta_{ 
		\ell^\fa_k \ell^\fa_{ k - 1 } + i + \sum_{ h = 1 }^{ k - 1 } \ell^\fa_h (\ell^\fa_{ h - 1 } + 1 ) 
	} ,
	\end{equation}
	let 
	$ 
	\mathfrak{M} \colon ( \cup_{ n \in \N } \R^n )
	\to ( \cup_{ n \in \N } \R^n ) 
	$
	satisfy for all 
	$ n \in \N $,
	$ x = (x_1, \ldots, x_n ) \in \R^n $
	that 
	$ 
	\mathfrak{M}( x ) = ( \max\{ x_1, 0 \}, \ldots, \max\{ x_n, 0 \} ) 
	$, 
	for every 
	$ \fa \in \N $, 
	$ \theta \in \R^{ \fd_{ \fa } } $ 
	let 
	$ 
	\mathcal{N}^{ k, \theta }_{ \fa } 
	\colon \R^d \to \R^{ \ell^{ \fa }_k } 
	$, 
	$ k \in \N \cap [1, \rho_{ \fa } ]
	$, 
	satisfy for all 
	$ k \in \N \cap [1, \rho_{ \fa } ) $,
	$x \in \R^d$
	that 
	\begin{equation}
	\mathcal{N}^{ 1, \theta }_{ \fa }( x ) 
	=
	\fb^{ 1, \theta }_{ \fa } 
	+
	\fw^{ 1, \theta }_{ \fa } x 
	\qquad
	\text{and}
	\qquad
	\mathcal{N}^{ k+1, \theta }_{ \fa }( x ) 
	=
	\fb^{ k + 1, \theta }_{ \fa } 
	+
	\fw^{ k + 1, \theta }_{ \fa }\big(
	\mathfrak{M}( \mathcal{N}^{ k, \theta }_{ \fa }( x ) ) 
	\big)
	,
	\end{equation} 
	for every $ \fa \in \N $
	let 
	$ \cL_{ \fa } \colon \R^{ \fd_{ \fa } } \to \R $
	satisfy for all
	$ 
	\theta \in \R^{ \fd_{ \fa } } 
	$ 
	that 
	$ 
	\cL_{ \fa }( \theta ) 
	=
	\int_{ [a,b]^{ d } } 
	\norm{ \mathcal{N}_{ \fa }^{ \rho_\fa, \theta }( x ) - f(x) }^2 
	\, 
	\dens( x )
	\,
	\d x
	$, 
	for every $ \fa \in \N $ 
	let 
	$ 
	\vartheta_{ \fa } \in 
	( \cL_{ \fa } )^{ - 1 }( 
	\{ 
	\inf_{ \theta \in \R^{ \fd_{ \fa } } } 
	\cL_{ \fa }( \theta )
	\} 
	)
	$, 
	$ \varepsilon_{ \fa } \in (0,1) $ 
	satisfy that
	$
	\cL_{ \fa }|_{ 
		\{ 
		\theta \in \R^{ \fd_{ \fa } } \colon 
		\norm{ \theta - \vartheta_{ \fa } } 
		< \varepsilon_{ \fa } 
		\}
	} 
	$
	has a Lipschitz continuous derivative, 
	for every $ \fa \in \N $ 
	let 
	$ \cG_{ \fa } \colon \R^{ \fd_{ \fa } } \to \R^{ \fd_{ \fa } } $
	satisfy for all 
$ 
  \theta \in 
  \cup_{ 
    U \subseteq \R^{ \fd_{ \fa } } , 
    \,  U \text{ is open}, 
    \, \cL_{ \fa }|_U \in C^1( U , \R ) 
  } 
  U
$
that 
$ 
  \cG_{ \fa }( \theta ) = ( \nabla \cL_{ \fa } )( \theta ) 
$,
let $ ( \Omega, \cF, \P ) $ be a probability space, 
for every 
$ n, \fa, K \in \N_0 $, $ \gamma \in \R $
let 
$ 
  \Theta^{ K, \gamma }_{ \fa, n } \colon \Omega \to  \R^{ \fd_{ \fa } } 
$
and 
$ 
  \bfk^{ K, \gamma }_{ \fa, n } \colon \Omega \to \N 
$ 
be random variables,
	assume for all $ \fa \in \N $, $ \gamma \in \R $ 
	that $ \Theta_{ \fa, 0 }^{ K, \gamma } $, $ K \in \N $, are i.i.d., 
	assume for all $ \fa \in \N $, $ \gamma, r \in (0,1) $, 
	$ \theta \in \R^{ \fd_{ \fa } } $ that
	$
	\P( 
	\norm{ \Theta^{ 1, \gamma }_{ \fa, 0 } - \theta } < r
	)
	> 0
	$, 
	and assume for all 
	$ n \in \N_0 $, $ \fa, K \in \N $, $ \gamma \in \R $, $ \omega \in \Omega $ 
	that
	\begin{equation} 
	\Theta_{ \fa, n + 1 }^{ K, \gamma }( \omega ) 
	= \Theta_{ \fa, n }^{ K, \gamma }( \omega ) 
	- \gamma \cG_{ \fa }( \Theta_{ \fa, n }^{ K, \gamma }( \omega ) ) 
	\quad\text{and}\quad
	\bfk^{ K, \gamma }_{ \fa, n }( \omega) 
	\in 
	\argmin\nolimits_{ \kappa \in \{ 1, 2, \ldots, K \} } 
	\cL_{ \fa }( \Theta_{ \fa, n }^{ \kappa, \gamma }( \omega ) ) 
	\end{equation}
	\cfload.
Then
\begin{enumerate} [label = (\roman*)]
\item 
\label{cor:gd:random:dnn:item1}
there exist
$
  \bfA 
  \colon \R \to \R 
$
and 
$ 
  \fg 
  \colon \N \to (0, \infty) 
$ 
such that 
\begin{equation}  
\textstyle
  \inf_{ \varepsilon \in (0,\infty) }
  \inf_{
    \fa \in \N \cap [ \bfA( \varepsilon ), \infty)
  }
  \inf_{ \gamma \in (0, \fg( \fa ) ] }
  \liminf\nolimits_{ K \to \infty } 
  \P\bigl(
    \limsup\nolimits_{ n \to \infty } 
    \cL_{ \fa }\bigl(
      \Theta^{ 
        \bfk^{ K, \gamma }_{ \fa, n }, \gamma
      }_{ \fa, n } 
    \bigr) 
    \le \varepsilon  
  \bigr) = 1 
\end{equation}
and
\item 
\label{cor:gd:random:dnn:item2}
it holds that
\begin{equation}
  \limsup\nolimits_{ \fa \to \infty }
  \limsup\nolimits_{ \gamma \searrow 0 }
  \limsup\nolimits_{ K \to \infty }
  \E\bigl[ 
    \limsup\nolimits_{ n \to \infty } 
    \min\bigl\{ 
      \cL_{ \fa }\rbr[\big]{ 
      \Theta^{ \bfk^{ K, \gamma }_{ \fa, n }, \gamma }_{ \fa, n } } 
      , 1
    \bigr\}
  \bigr]= 0 .
\end{equation}
\end{enumerate} 
\end{theorem}
\begin{cproof}{theo:gd:random:dnn}
Throughout this proof for every $ \fa \in \N $ let 
$ U_{ \fa } \subseteq \R^{ \fd_{ \fa } } $ satisfy 
\begin{equation}
  U_{ \fa } = 
  \{ 
    \theta \in \R^{ \fd_{ \fa } }
    \colon 
    \| \theta - \vartheta_{ \fa } \| < \varepsilon_{ \fa }
  \}
  .
\end{equation}
\Nobs the assumption that $ \dens $ is piecewise polynomial 
and the assumption that for all $ i \in \cu{ 1, 2, \ldots, \delta } $ 
it holds that $ f_i $ is piecewise polynomial 
imply that $ f $ and $ \dens $ are bounded and measurable. 
Combining this, 
\cite[Lemma~2.4]{DNNReLUarXiv}, 
\cref{prop:gd:random:dnn} 
(applied for every $ \fa \in \N $ with 
$
  a \with a
$, 
$
  b \with b
$,
$
  \scrA \with \nicefrac{ 1 }{ 2 }
$, 
$
  \scrB \with 1 
$, 
$
  (
    \N_0 \ni k 
    \mapsto 
    \ell_k \in \N
  )
  \with 
  (
    \N_0 \ni k
    \mapsto 
    \ell_{ \min\{ k, \rho_{ \fa } \}^{ \fa } } 
    \in \N
  )
$, 
$ 
  L \with \rho_{ \fa } 
$, 
$
  \fd \with \fd_{ \fa }
$, 
$ 
  f \with f
$, 
$
  \mu \with 
  ( \cB( [a,b]^d ) \ni E \mapsto \int_E \dens(x) \, dx \in [0,\infty] )
$,
$
  U \with U_{ \fa }
$, 
$
  \fG \with \cG_{ \fa }
$, 
$
  \vartheta \with \vartheta_{ \fa }
$
in the notation of \cref{prop:gd:random:dnn}), 
and the fact that for all $ \fa \in \N $ it holds that 
\begin{equation}
\textstyle
  \cL_{ \fa }( \vartheta_{ \fa } ) 
  = 
  \inf_{ \theta \in \R^{ \fd_{ \fa } } } \cL_{ \fa }( \theta )
  = 
  \inf_{ \theta \in U_{ \fa } } \cL_{ \fa }( \theta )
\end{equation}
proves that there
there exists $ \fg \colon \N \to (0, \infty) $ 
which satisfies for all $ \fa \in \N $, $ \gamma \in (0, \fg( \fa ) ] $ 
that
\begin{equation}
\label{eq:limit_to_global_minimum_in_proof}
  \liminf\nolimits_{ K \to \infty } 
  \P\Big(
    \limsup\nolimits_{ n \to \infty } 
    \cL_{ \fa }\bigl( 
      \Theta^{ \bfk^{ K, \gamma }_{ \fa, n }, \gamma }_{ \fa, n } 
    \bigr) 
    \le 
    \inf\nolimits_{ \theta \in \R^{ \fd_{ \fa } } } 
    \cL_{ \fa }( \theta ) 
  \Bigr) = 1 .
\end{equation}
\cref{prop:dnn:approx} \hence 
establishes that there exists 
$ \bfA \colon \R \to \R $ 
which satisfies for all 
$ \varepsilon \in (0,\infty) $, $ \fa \in \N \cap [ \bfA( \varepsilon ) , \infty ) $
that
\begin{equation}
\label{eq:global_minimum_small_in_proof}
  \inf\nolimits_{ \theta \in \R^{ \fd_{ \fa } } } 
  \cL_{ \fa }( \theta ) 
  \le \varepsilon 
  .
\end{equation}
\Nobs that \cref{eq:limit_to_global_minimum_in_proof} 
and \cref{eq:global_minimum_small_in_proof} 
assure that for all 
$ \varepsilon \in (0,\infty) $, 
$ \fa \in \N \cap [ \bfA( \varepsilon ), \infty ) $, 
$ \gamma \in ( 0, \fg( \fa ) ] $
it holds that
\begin{equation} 
\label{cor:gd:random:dnn:eq1}
\begin{split}
&
  \liminf\nolimits_{ K \to \infty } 
  \P\Big(
    \limsup\nolimits_{ n \to \infty } 
    \cL_{ \fa }\rbr[\big]{ 
      \Theta^{ \bfk^{ K, \gamma }_{ \fa, n }, \gamma }_{ \fa, n } 
    } \le \varepsilon 
  \Big) 
\\ &
  \geq
  \liminf\nolimits_{ K \to \infty } 
  \P\Big(
    \limsup\nolimits_{ n \to \infty } 
    \cL_{ \fa }\rbr[\big]{ 
      \Theta^{ \bfk^{ K, \gamma }_{ \fa, n }, \gamma }_{ \fa, n } 
    } \le 
    \inf\nolimits_{ \theta \in \R^{ \fd_{ \fa } } } 
    \cL_{ \fa }( \theta ) 
  \Big) 
  = 1 .
\end{split}
\end{equation}
This establishes \cref{cor:gd:random:dnn:item1}. 
\Nobs that for all 
$ \varepsilon \in (0,\infty) $
and all 
random variables $ Z_n \colon \Omega \to [0, \infty) $, $ n \in \N $,
it holds that
\begin{equation}
\label{cor:in_proof_general_RV}
\begin{split}
%   \limsup\nolimits_{n \to \infty}
%   \E \br*{ \min \cu*{Z_n , 1 } }
% &
%   \le 
&
  \E\br*{ 
    \limsup\nolimits_{ n \to \infty }
    \min\cu*{ Z_n, 1 } 
  } 
\leq
  \E\br*{ 
    \min\cu*{ 
      \limsup\nolimits_{ n \to \infty } Z_n , 1 
    } 
  } 
\\
&
  \leq 
  \E\bigl[
    \min\cu*{ 
      \limsup\nolimits_{ n \to \infty } Z_n , 1 
    } 
    \mathbbm{1}_{
      \{ 
        \limsup\nolimits_{ n \to \infty } Z_n > \varepsilon
      \}
    }
  \bigr] 
\\ &
  +
  \E\bigl[
    \min\cu*{ 
      \limsup\nolimits_{ n \to \infty } Z_n , 1 
    } 
    \mathbbm{1}_{
      \{ 
        \limsup\nolimits_{ n \to \infty } Z_n \leq \varepsilon
      \}
    }
  \bigr] 
\\ &
  \leq
  \E\bigl[
    \mathbbm{1}_{
      \{ 
        \limsup\nolimits_{ n \to \infty } Z_n > \varepsilon
      \}
    }
  \bigr] 
  +
  \E\bigl[
    \min\cu*{ 
      \varepsilon , 1 
    } 
    \mathbbm{1}_{
      \{ 
        \limsup\nolimits_{ n \to \infty } Z_n \leq \varepsilon
      \}
    }
  \bigr] 
\\ &
  \leq
  \P\bigl(
    \limsup\nolimits_{ n \to \infty } Z_n > \varepsilon
  \bigr)
  +
  \min\cu*{ 
    \varepsilon , 1 
  } 
\le 
  \P\rbr*{
    \limsup\nolimits_{ n \to \infty } Z_n > \varepsilon 
  } 
  + \varepsilon 
  .	
\end{split}
\end{equation}
\Moreover \cref{cor:gd:random:dnn:eq1} assures that for all 
$ \varepsilon \in (0,\infty) $, 
$ \fa \in \N \cap [ \bfA( \varepsilon ) , \infty ) $, 
$ \gamma \in (0, \fg( \fa ) ] $
it holds that
\begin{equation} 
\begin{split}
&
  \limsup\nolimits_{ K \to \infty } 
  \P\Bigl(
    \limsup\nolimits_{ n \to \infty } 
    \cL_{ \fa }\rbr[\big]{ 
      \Theta^{ \bfk^{ K, \gamma }_{ \fa, n }, \gamma }_{ \fa, n } 
    } > \varepsilon 
  \Bigr) 
\\ &
  =
  \limsup\nolimits_{ K \to \infty } 
  \Bigl[
    1 
    -
    \P\Bigl(
      \limsup\nolimits_{ n \to \infty } 
      \cL_{ \fa }\rbr[\big]{ 
        \Theta^{ \bfk^{ K, \gamma }_{ \fa, n }, \gamma }_{ \fa, n } 
      } 
      \leq \varepsilon 
    \Bigr) 
  \Bigr]
  = 0 .
\end{split}
\end{equation}
Combining this with \cref{cor:in_proof_general_RV}
ensures that 
for all $ \varepsilon \in (0,\infty) $, 
$ \fa \in \N \cap [ \bfA( \varepsilon ) , \infty ) $, 
$ \gamma \in (0, \fg( \fa ) ] $
it holds that
\begin{equation}
  \limsup\nolimits_{ K \to \infty }
  \E\bigl[
    \limsup\nolimits_{ n \to \infty }
    \min\bigl\{
      \cL_{ \fa }\rbr[\big]{ 
        \Theta^{ \bfk^{ K, \gamma }_{ \fa, n }, \gamma }_{ \fa, n } 
      } , 1 
    \bigr\}
  \bigr]
  \le \varepsilon .
\end{equation}
This establishes \cref{cor:gd:random:dnn:item2}.
\end{cproof}

\subsection{Convergence of GD with random initializations in the training of shallow ANNs}
\label{subsec:convergence_GD_shallow_ANNs}

In this section we employ the general convergence results for deep ANNs from \cref{sec:GD_random_init}
to establish convergence of the risk of the GD method for shallow ANNs.
This time the regularity assumptions can be omitted, since they follow from the existence result for regular global minima for shallow ANNs in \cref{cor:existence:regular}.

\cfclear
\begin{prop}
\label{prop:gd:random:snn}
Assume \cref{setting:dnn},
assume $ L = 2 $, 
% $ \ell_0 = \ell_2 = 1 $, 
assume $ f \in C( [a,b], \R ) $,
assume that $ f $ is piecewise polynomial, 
let $ \dens \colon [ a, b ] \to \R $ 
be piecewise polynomial\cfadd{def:multidim:piece:polyn}, 
assume for all $ E \in \cB( [ a, b ] ) $ 
that $ \mu ( E ) = \int_E \dens ( x ) \, \d x $,
let $ \fG \colon \R^{ \fd } \to \R^{ \fd } $ 
satisfy for all 
$ 
  \theta \in 
  \{ 
    \vartheta \in \R^{ \fd } 
    \colon 
    \cL_{ \infty }\ \text{is}\ \allowbreak 
    \text{differentiable}\ \allowbreak \text{at}\ \vartheta
  \}
$
that 
$ 
  \fG( \theta ) = ( \nabla \cL_{ \infty } )( \theta ) 
$,
let $ ( \Omega, \cF, \P ) $ be a probability space, 
for every 
$ K, n \in \N_0 $, $ \gamma \in \R $
let 
$
  \Theta^{ K, \gamma }_n \colon \Omega \to \R^{ \fd } 
$
and 
$
  \bfk^{ K, \gamma }_n \colon \Omega \to \N 
$
be random variables,
assume for all $ \gamma \in \R $ that 
$ 
  \Theta_0^{ K, \gamma } 
$, $ K \in \N $, 
are i.i.d., 
assume for all 
$ \gamma, \delta \in (0, 1) $
that
$
  \P( 
    \| \Theta_0^{ 1, \gamma } - \vartheta \| < \delta 
  ) > 0 
$,
and assume for all $ K \in \N $, $ n \in \N_0 $,
$ \gamma \in \R $, $ \omega \in \Omega $ that
\begin{equation} 
\label{prop:gd:random:snn:eq:defk}
  \Theta_{ n + 1 }^{ K, \gamma }( \omega ) 
  = 
  \Theta_n^{ K, \gamma }( \omega ) - \gamma \fG( \Theta_n^{ K, \gamma }( \omega ) ) 
\qandq 
  \bfk^{ K, \gamma }_n( \omega ) 
  \in \arg\min\nolimits_{ \kappa \in \{ 1, 2, \ldots, K \} } 
  \cL_{ \infty }( \Theta_n^{ \kappa, \gamma }( \omega ) )
\end{equation}
\cfload.
Then there exists $ \fg \in (0, \infty) $ 
such that for all $ \gamma \in (0, \fg] $ 
it holds that
\begin{equation}
  \liminf\nolimits_{ K \to \infty } 
  \P\Bigl( 
    \limsup\nolimits_{ n \to \infty } 
    \cL_{ \infty }\bigl( 
      \Theta^{ \bfk^{ K, \gamma }_n , \gamma }_n 
    \bigr) 
    = 
    \inf\nolimits_{ \theta \in \R^{ \fd } } \cL_{ \infty }( \theta ) 
  \Bigr) = 1 .
\end{equation}
\end{prop}
\begin{cproof}{prop:gd:random:snn}
\Nobs that 
the assumption that $ f \in C( [a,b], \R ) $ 
and 
the assumption that $ f $ is piecewise polynomial imply 
that $ f $ is Lipschitz continuous. 
\Moreover 
\cref{lim_R} 
and the assumption that 
$
  \sup_{ r \in [1,\infty) }
  \sup_{ x \in \R } 
  | ( \Rect_r )'( x ) | < \infty 
$
assure that for all $ x \in \R $ 
it holds that 
$ 
  ( \cup_{ r \in \N } \{ \Rect_r \} ) \subseteq C^1( \R, \R ) 
$,
$ 
  \Rect_{ \infty }( x ) = \max\{ x, 0 \} 
$,
$
  \sup_{ r \in \N }
  \sup_{ y \in [ - |x|, |x| ] } 
  | (\Rect_r)'( y ) | < \infty 
$, 
and
\begin{equation}
  \limsup\nolimits_{ r \to \infty 
  }( 
    | 
      \Rect_r(x) - \Rect_\infty( x ) 
    |
    +
    | 
      ( \Rect_r )'( x ) - \indicator{ (0,\infty) }( x ) 
    | 
  ) = 0 
  .
\end{equation}
Combining this, \cref{cor:existence:regular}, 
and \cref{G3} in \cref{prop:G} with the fact that $ f $ is Lipschitz continuous 
shows that there exist 
$ \vartheta \in \R^{ \fd } $
and an open $ U \subseteq \R^{ \fd } $
which satisfy that 
\begin{enumerate}[label = (\roman*)]
\item 
\label{item:i_in_proof_shallow_ANNs}
it holds that 
$
  ( \cL_{ \infty } )|_U \in C^1( U, \R ) 
$, 
\item 
\label{item:ii_in_proof_shallow_ANNs}
it holds for all $ \theta \in U $ that 
$ \cG( \theta ) = ( \nabla \cL_{ \infty } )( \theta ) $, 
\item 
\label{item:iii_in_proof_shallow_ANNs}
it holds that 
$ \cG|_U $ 
is locally Lipschitz continuous, 
\item 
\label{item:iv_in_proof_shallow_ANNs}
it holds that $ \vartheta \in U $, 
and
\item 
\label{item:v_in_proof_shallow_ANNs}
it holds that 
$ 
  \cL_{ \infty }( \vartheta ) = 
  \inf_{ \theta \in \R^{ \fd } } 
  \cL_{ \infty }( \theta ) 
$.
\end{enumerate}
\Nobs that 
\cref{item:iv_in_proof_shallow_ANNs} and 
\cref{item:v_in_proof_shallow_ANNs} 
ensure that 
\begin{equation}
\label{eq:global_minimum_representation_in_proof}
\textstyle
  \cL_{ \infty }( \vartheta ) 
  = 
  \inf_{ \theta \in \R^{ \fd } } 
  \cL_{ \infty }( \theta ) 
  =
  \inf_{ \theta \in U } 
  \cL_{ \infty }( \theta ) 
  .
\end{equation}
\Moreover  
\cref{item:i_in_proof_shallow_ANNs}, 
\cref{item:ii_in_proof_shallow_ANNs}, 
and the assumption that 
for all 
$ 
  \theta \in 
  \{ 
    \vartheta \in \R^{ 3 \width + 1 } 
    \colon 
    \cL_{ \infty }\ \text{is}\ \allowbreak 
    \text{differentiable}\ \allowbreak \text{at}\ \vartheta
  \}
$
it holds that 
$ 
  \fG( \theta ) = ( \nabla \cL_{ \infty } )( \theta ) 
$
assure that for all $ \theta \in U $ it holds that 
\begin{equation}
\label{eq:representation_of_fG_in_proof}
  \cG( \theta ) = ( \nabla \cL_{ \infty } )( \theta ) = \fG( \theta  ) 
  .
\end{equation}
\Hence that $ \fG|_U = \cG|_U $. 
This and \cref{item:iii_in_proof_shallow_ANNs} ensure that 
$ \fG|_U $ is locally Lipschitz continuous. 
Combining this, 
\cref{item:i_in_proof_shallow_ANNs}, 
\cref{eq:global_minimum_representation_in_proof}, 
\cref{eq:representation_of_fG_in_proof}, 
\cref{prop:gd:random:dnn}, 
and 
the fact that $ U \subseteq \R^{ \fd } $ is open 
proves that there exists $ \fg \in (0, \infty) $ 
such that for all $ \gamma \in (0, \fg] $ 
it holds that
\begin{equation}
\begin{split}
&
  \liminf\nolimits_{ K \to \infty } 
  \P\Bigl( 
    \limsup\nolimits_{ n \to \infty } 
    \cL_{ \infty }( 
      \Theta^{ \bfk^{ K, \gamma }_n , \gamma }_n 
    ) 
    = \inf\nolimits_{ \theta \in \R^{ \fd } } \cL_{ \infty }( \theta ) 
  \Bigr)
\\ &
  =
  \liminf\nolimits_{ K \to \infty } 
  \P\Bigl( 
    \limsup\nolimits_{ n \to \infty } 
    \cL_{ \infty }( 
      \Theta^{ \bfk^{ K, \gamma }_n , \gamma }_n 
    ) 
    \le \inf\nolimits_{ \theta \in \R^{ \fd } } \cL_{ \infty }( \theta ) 
  \Bigr)
\\ &
=
  \liminf\nolimits_{ K \to \infty } 
  \P\Bigl( 
    \limsup\nolimits_{ n \to \infty } 
    \cL_{ \infty }( 
      \Theta^{ \bfk^{ K, \gamma }_n , \gamma }_n 
    ) 
    \le \inf\nolimits_{ \theta \in U } \cL_{ \infty }( \theta ) 
  \Bigr)
  = 1 .
\end{split}
\end{equation}
\end{cproof}

\begin{cor}
\label{cor:gd:random:shallow}
Let $ N \in \N $, 
$ \fx_0, \fx_1, \ldots, \fx_N, a, b \in \R $
satisfy 
$ a = \fx_0 < \fx_1 < \ldots < \fx_N = b $, 
let 
$ 
  f \in C( [a,b], \R )
$, 
let 
$ \dens \colon [a,b] \to [0, \infty) $ be a function, 
assume for all $ i \in \{ 1, 2, \ldots, N \} $ that 
$ f|_{ ( \fx_{i-1}, \fx_i ) } $ 
and 
$ \dens|_{ ( \fx_{i-1}, \fx_i ) } $ 
are polynomials, 
for every $ \width \in \N $ 
		let 
		$ \cL_\width \colon \R^{3 \width + 1} \to \R $
		satisfy for all 
		$ \theta = ( \theta_1, \ldots, \theta_{ 3 \width + 1 } ) \in \R^{ 3 \width + 1 } $ 
		that
		\begin{equation} 
		\label{theo:gd:random:init:eq2}
		\cL_\width( \theta ) = 
		\textstyle
		\int_a^b \rbr[\big]{ 
			f ( x ) 
			- 
			\theta_{\fd} - 
			\smallsum_{j=1}^\width 
			\theta_{2 \width + j} 
			\max\{ \theta_{j} x + \theta_{\width + j } , 0 \} }^2 
		\dens ( x ) \, \d x,
		\end{equation} 
		for every $ \width \in \N $
		let 
		$ \cG_\width \colon \R^{ 3 \width + 1 } \to \R^{ 3 \width + 1 } $ 
		satisfy for all 
		$ 
		\theta \in 
		\{ 
		\vartheta \in \R^{3 \width + 1} 
		\colon 
		\cL_{ \width }
		\ \text{is}\ \allowbreak \text{differentiable}\ \allowbreak \text{at}\ \vartheta
		\}
		$
		that 
		$ 
		\cG_{ \width }( \theta ) = ( \nabla \cL_{ \width } )( \theta ) 
		$,
		let $ (\Omega, \cF, \P) $ be a probability space, 
		for every 
		$ n, \width, K \in \N_0 $, $ \gamma \in \R $
		let 
		$ \Theta^{ K, \gamma }_{ \width, n } \colon \Omega \to \R^{ 3 \width + 1 } $
		and 
		$ \bfk^{ K, \gamma }_{ \width, n } \colon \Omega \to \N $ 
		be random variables,
		assume for all $ \width \in \N $, $ \gamma \in \R $ that 
		$ \Theta_{ \width, 0 }^{ K, \gamma } $, $ K \in \N $, are i.i.d., 
		assume for all $ \width \in \N $, $ \gamma, r \in (0,1) $, 
		$ \theta \in \R^{ 3 \width + 1 } $ that
		$
		\P( 
		\| \Theta^{ 1 , \gamma }_{ \width, 0 } - \theta \| < r
		)
		> 0
		$, 
		and assume for all 
		$ n, \width \in \N_0 $, $ K \in \N $, $ \gamma \in \R $, $ \omega \in \Omega $ 
		that
		\begin{equation} 
		\label{theo:gd:random:init:eq3}
		\Theta_{ \width, n+1 }^{ K, \gamma }( \omega ) 
		= \Theta_{ \width, n }^{ K, \gamma }( \omega ) 
		- \gamma \cG_\width( \Theta_{ \width, n }^{ K, \gamma }( \omega ) ) 
		\quad\text{and}\quad 
		\bfk^{ K, \gamma }_{ \width, n }( \omega) 
		\in 
		\argmin\nolimits_{ \kappa \in \{ 1, 2, \ldots, K \} } 
		\cL_{ \width }( \Theta_{ \width, n }^{ \kappa, \gamma }( \omega ) ) 
		  .
		\end{equation}
Then
\begin{enumerate} [label = (\roman*)]
\item \label{theo:gd:random:init:item1}
there exist 
$ \bfH \colon \R \to \R $ and 
$ \fg \colon \N \to (0,\infty) $ such that
\begin{equation}  
\textstyle
  \inf_{ \varepsilon \in (0,\infty) }
  \inf_{ \width \in \N \cap [ \bfH( \varepsilon ) , \infty) } 
  \inf_{ \gamma \in (0, \fg( \width ) ] }
  \liminf\nolimits_{ K \to \infty } 
  \P\bigl( 
    \limsup\nolimits_{ n \to \infty } 
    \cL_{ \width }\bigl(
      \Theta^{ \bfk^{ K, \gamma }_{ \width, n }, \gamma }_{ \width, n }
    \bigr) 
    \le \varepsilon
  \bigr)
  = 1 
\end{equation}
and
\item 
\label{theo:gd:random:init:item2}
it holds that
\begin{equation}  
  \limsup\nolimits_{ \width \to \infty } 
  \limsup\nolimits_{ \gamma \searrow 0 }
  \limsup\nolimits_{ K \to \infty } 
  \mathbb{E}\bigl[
    \limsup\nolimits_{ n \to \infty }
    \min\bigl\{
      \cL_{ \width }\bigl(
        \Theta^{ \bfk^{ K, \gamma }_{ \width, n }, \gamma }_{ \width, n }
      \bigr) 
      , 1 
    \bigr\}
  \bigr] 
  = 0 .
\end{equation}
\end{enumerate}
\end{cor}
\begin{cproof}{cor:gd:random:shallow}
\Nobs that \cref{cor:existence:regular} demonstrates that 
for every $ \width \in \N $ there exist 
$ \vartheta_{ \width } \in \R^{ 3 \width + 1 } $, 
$ \bL_{ \width } \in \R $, 
and an open 
$ V_{ \width } \subseteq \R^{ 3 \width + 1 } $
which satisfy that 
\begin{enumerate}[label = (\Roman*)]
\item
\label{item:in_proof_item_I}
it holds that 
$ \vartheta_{ \width } \in V_{ \width } $, 
\item 
\label{item:in_proof_item_II}
it holds that 
$ 
  \cL_{ \width }( \vartheta_{ \width } ) = \inf_{ \psi \in \R^{ 3 \width + 1 } } \cL_{ \width }( \psi )
$, 
\item 
\label{item:in_proof_item_III}
it holds that 
$
  \cL_{ \width }|_{ V_{ \width } } \in C^1( V_{ \width }, \R )
$, 
and 
\item 
\label{item:in_proof_item_IV}
it holds for all $ \theta_1, \theta_2 \in V_{ \width } $ that 
$
  \| 
    ( \nabla \cL_{ \width } )( \theta_1 ) 
    - 
    ( \nabla \cL_{ \width } )( \theta_2 ) 
  \|
  \le \bL_{ \width } 
  \| \theta_1 - \theta_2 \| 
$.
\end{enumerate}
\Moreover the fact that for all $ \width \in \N $,
$ 
  \theta \in 
  \{ 
    \vartheta \in \R^{ 3 \width + 1 } 
    \colon 
    \cL_{ \width }
    \ \text{is}\ \allowbreak \text{differentiable}\ \allowbreak \text{at}\ \vartheta
  \}
$
it holds that 
\begin{equation}
  \cG_{ \width }( \theta ) = ( \nabla \cL_{ \width } )( \theta ) 
\end{equation}
proves that for all $ \width \in \N $, 
$ 
  \theta \in 
  \cup_{ 
    U \subseteq \R^{ 3 \width + 1 }, \, 
    U \text{ is open}, \, \cL_{ \width }|_U \in C^1( U, \R ) 
  } 
  V
$
it holds that 
\begin{equation}
  \cG_{ \width }( \theta ) = ( \nabla \cL_{ \width } )( \theta ) 
  .
\end{equation}
Combining this, 
\cref{item:in_proof_item_I}, 
\cref{item:in_proof_item_II}, 
\cref{item:in_proof_item_III}, 
\cref{item:in_proof_item_IV}, 
and 
\cref{cor:gd:random:dnn:item1} in \cref{theo:gd:random:dnn} 
(applied with 
$ d \with 1 $, 
$ \delta \with 1 $, 
$ a \with a $,
$ b \with b $, 
$ 
  ( \N \ni \fa \mapsto \rho_{ \fa } \in \N \cap (1,\infty) ) 
  \with 
  ( \N \ni \fa \mapsto 2 \in \N \cap (1,\infty) ) 
$,
$
  ( 
    \N \ni \fa \mapsto \ell^{ \fa } 
    \in 
    \N^3
  )
  \with 
  ( 
    \N \ni \fa \mapsto (d, \fa, \delta)
    \in 
    \N^3
  )
$,
$
  ( \N \ni \fa \mapsto \fd_{ \fa } \in \N )
  \with 
  ( \N \ni \fa \mapsto ( 3 \fa + 1 ) \in \N )
$, 
$
  ( \N \ni \fa \mapsto \vartheta_{ \fa } \in ( \cup_{ k \in \N } \R^k ) )
  \with 
  (
    \N \ni \fa \mapsto \vartheta_{ \fa } \in ( \cup_{ k \in \N } \R^k )
  )
$
in the notation of \cref{theo:gd:random:dnn})
establishes \cref{theo:gd:random:init:item1,theo:gd:random:init:item2}.
\end{cproof}

\subsubsection*{Acknowledgments}
This work has been funded by the Deutsche Forschungsgemeinschaft 
(DFG, German Research Foundation) under Germany's Excellence Strategy 
EXC 2044-390685587, Mathematics M\"{u}nster: Dynamics-Geometry-Structure.
This project has been partially supported by the startup fund project 
of Shenzhen Research Institute of Big Data under grant 
No.\ T00120220001.

%\bibliographystyle{plainurl}
%\bibliography{loja_bib}

\end{document}